\documentclass{amsbook}
\usepackage{graphicx}
\usepackage{amssymb}
\usepackage{amsfonts}
\usepackage{epsfig}
\swapnumbers
\newtheorem{theorem}{Theorem}[section]
\newtheorem{lemma}[theorem]{Lemma}

\newtheorem{conjecture}[theorem]{Conjecture}
\newtheorem{corollary}[theorem]{Corollary}
\newtheorem{proposition}[theorem]{Proposition}
\theoremstyle{definition}

\theoremstyle{assumption}
\newtheorem{assumption}[theorem]{Assumption}
\newtheorem{remark}[theorem]{Remark}

 \theoremstyle{plain}    
 
 \numberwithin{section}{part}
 \numberwithin{equation}{part} 
 \numberwithin{figure}{section} 
 \theoremstyle{plain}    
 \theoremstyle{plain}    
 \theoremstyle{remark}    
 \newtheorem*{acknowledgement*}{Acknowledgement} 
\newcommand{\cC}{{\mathcal C}}
\newcommand{\cD}{{\mathcal D}}
\newcommand{\cE}{{\mathcal E}}
\newcommand{\cF}{{\mathcal F}}

\newcommand{\cI}{{\mathcal I}}

\newcommand{\cK}{{\mathcal K}}
\newcommand{\cL}{{\mathcal L}}
\newcommand{\cM}{{\mathcal M}}

\newcommand{\cO}{{\mathcal O}}
\newcommand{\cP}{{\mathcal P}}

\newcommand{\cT}{{\mathcal T}}
\newcommand{\cU}{{\mathcal U}}

\newcommand{\cW}{{\mathcal W}}
\newcommand{\cX}{{\mathcal X}}
\newcommand{\te}{{\theta}}
\newcommand{\Te}{{\Theta}}

\newcommand{\vrho}{{\varrho}}
\newcommand{\vsig}{{\varsigma}}
\newcommand{\Om}{{\Omega}}
\newcommand{\om}{{\omega}}
\newcommand{\ve}{{\varepsilon}}
\newcommand{\del}{{\delta}}
\newcommand{\Del}{{\Delta}}
\newcommand{\gam}{{\gamma}}
\newcommand{\Gam}{{\Gamma}}
\newcommand{\vf}{{\varphi}}
\newcommand{\vp}{{\varpi}}

\newcommand{\io}{{\iota}}

\newcommand{\sig}{{\sigma}}
\newcommand{\al}{{\alpha}}
\newcommand{\be}{{\beta}}
\newcommand{\ka}{{\kappa}}
\newcommand{\la}{{\lambda}}
\newcommand{\La}{{\Lambda}}

\newcommand{\bbN}{{\mathbb N}}

\newcommand{\bbR}{{\mathbb R}}
\newcommand{\bbT}{{\mathbb T}}

\newcommand{\bbI}{{\mathbb I}}

\newcommand{\bfM}{{\bf M}}

\newcommand{\bfr}{{\bf r}}

\newcommand{\brX}{{\bar X}}
\newcommand{\brZ}{{\bar Z}}
\newcommand{\brB}{{\bar B}}
\newcommand{\brU}{{\bar U}}

\newcommand{\brcW}{{\bar\cW}}

\newcommand{\SRBx}{{\mu_x^{\mbox{\tiny{SRB}}}}}
\newcommand{\SRBz}{{\mu_z^{\mbox{\tiny{SRB}}}}}
\newcommand{\SRBwx}{{\mu_{w,x}^{\mbox{\tiny{SRB}}}}}
\newcommand{\SRBwoi}{{\mu_{w,\cO_i}^{\mbox{\tiny{SRB}}}}}
\makeindex
\begin{document}
\frontmatter
\title{Large Deviations and Adiabatic Transitions\\
 for Dynamical Systems and Markov Processes\\
  in Fully Coupled Averaging.}%
 \vskip 0.1cm 
 \author{ Yuri Kifer}
\address{
 Institute of Mathematics, The Hebrew University, Jerusalem 91904, Israel}%
\email{ kifer@math.huji.ac.il}%
\thanks{The author was partially supported by US--Israel BSF}

\date{\today}
\begin{abstract}\noindent
The work treats dynamical systems given by ordinary differential equations
in the form $\frac {dX^\ve(t)}{dt}=\ve B(X^\ve(t),Y^\ve(t))$ where fast motions
$Y^\ve$ depend on the slow motion $X^\ve$ (coupled with it) and they are either
given by another differential equation $\frac {dY^\ve(t)}{dt}=b(X^\ve(t),
Y^\ve(t))$ or perturbations of an appropriate parametric family of Markov
processes with freezed slow variables.
In the first case we assume that the fast motions are hyperbolic for each 
freezed slow variable and in the second case we deal with Markov processes
such as random evolutions which are combinations of diffusions and continuous
time Markov chains. First, we study large deviations of the slow motion 
$X^\ve$ from its averaged (in fast variables $Y^\ve$) approximation $\brX^\ve.$
The upper large deviation bound justifies the averaging approximation on the 
time scale of order $1/\ve$, called the averaging principle, in the sense of 
convergence in measure (in the first case) or in probability (in the second
case) but our real goal is to obtain both the upper and the lower large 
deviations bounds which together with some Markov property type arguments 
(in the first case) or with the real Markov property (in the second case) 
enable us to study (adiabatic) behavior of the slow motion on the much longer 
exponential in 
$1/\ve$ time scale, in particular, to describe its fluctuations in a vicinity
 of an attractor of the averaged motion and its rare (adiabatic) transitions 
between neighborhoods of such attractors. When the fast motion $Y^\ve$ does
 not depend on the slow one we arrive at a simpler averaging setup studied in
 numerous papers but the above fully coupled case, which better describes real
phenomena, leads to much more complicated problems.
\end{abstract}
\maketitle
\setcounter{page}{4}
\tableofcontents
\include{preface}

\markboth{Y.Kifer}{Averaging} 
\renewcommand{\theequation}{\arabic{part}.\arabic{section}.\arabic{equation}}
\pagenumbering{arabic}

\mainmatter

\part[Hyperbolic Fast Motions]{Hyperbolic Fast Motions}\label{part1}
\setcounter{equation}{0}

\section{Introduction}\label{sec1.1}\setcounter{equation}{0}

Many real systems can be viewed as a combination of slow and fast motions
which leads to complicated double scale equations. Already in the 19th 
century in applications to celestial mechanics it was well understood
(though without rigorous justification) that a good approximation of the
slow motion can be obtained by averaging its parameters in fast variables.
Later, averaging methods were applied in signal processing and, rather
recently, to model climate--weather interactions \index{climate--weather
system} (see \cite{Ha}, \cite{CMP},
\cite{Ha1} and \cite{Ki4}). The classical setup of averaging justified 
rigorously in \cite{BM} presumes that the fast motion \index{fast motion}
does not depend on the slow one \index{slow motion}
and most of the work on averaging treats this case only. On the other
hand, in real systems both slow and fast motions depend on each other which 
leads to the more difficult fully coupled case which we study here.
This setup emerges, in particular, in perturbations of Hamiltonian systems
\index{Hamiltonian system}
which leads to fast motions on manifolds of constant energy and slow
motions across them.

In this work we consider a system of differential equations for 
$X^\ve=X^\ve_{x,y}$ and $Y^\ve=Y^\ve_{x,y},$
\begin{equation}\label{1.1.1}
\frac {dX^\ve(t)}{dt}=\ve B(X^\ve(t),Y^\ve(t)),\quad
\frac {dY^\ve(t)}{dt}= b(X^\ve(t),Y^\ve(t))
\end{equation}
with initial conditions $X^\ve(0)=x,\, Y^\ve(0)=y$ on the product 
$\bbR^d\times\bfM$ where $\bfM$ is a compact $n_\bfM$-dimensional $C^2$ 
Riemannian manifold and $B(x,y),\, b(x,y)$ are smooth in $x,\, y$ families of
 bounded vector fields on $\bbR^d$ and on $\bfM,$ respectively, so that $y$
serves as a parameter for $B$ and $x$ for $b$. The solutions
of (\ref{1.1.1}) determine the flow of diffeomorphisms $\Phi_\ve^t$ on 
$\bbR^d\times\bfM$ acting by $\Phi_\ve^t(x,y)=(X^\ve_{x,y}(t),Y^\ve_{x,y}(t))$.
Taking $\ve=0$ we arrive at the flow $\Phi^t=\Phi_0^t$ acting by $\Phi^t(x,y)
=(x,F_x^ty)$ where $F^t_x$ is another family of flows given by $F_x^ty=
Y_{x,y}(t)$ with $Y=Y_{x,y}=Y^0_{x,y}$ being the solution of
\begin{equation}\label{1.1.2}
\frac {dY(t)}{dt}=b(x,Y(t)),\quad Y(0)=y.
\end{equation}
 It is natural to view the flow $\Phi^t$ as describing 
an idealised physical system where parameters $x=(x_1,...,x_d)$ are assumed
to be constants (integrals) of motion \index{constants of motion}
while the perturbed flow $\Phi_\ve^t$
is regarded as describing a real system where evolution of these parameters 
is also taken into consideration. Essentially, the proofs of this paper work 
also in the slightly more general case when $B$ and $b$ in (\ref{1.1.1}) 
together with their derivatives depend Lipschitz continuously on $\ve$ 
(cf. \cite{Ki7}) but in order to simplify notations and estimates we do not
consider this generalisation here.

 Assume that the limit
\begin{equation}\label{1.1.3}
\bar B(x)=\bar B_y(x)=\lim_{T\to\infty}T^{-1}\int_0^TB(x,F^t_xy)dt
\end{equation}
exists and it is the same for "many" $y'$s. For instance, suppose that
$\mu_x$ is an ergodic invariant measure of the flow $F^t_x$ then the limit 
(\ref{1.1.3}) exists for $\mu_x-$almost all $y$ and is equal to
\begin{equation*}
\bar B(x)=\bar B_{\mu_x}(x)=\int B(x,y)d\mu_x(y).
\end{equation*}
If $b(x,y)$ does not, in fact, depend on $x$ then $F^t_x=F^t$ and $\mu_x=\mu$
are also independent of $x$ and we arrive at the classical uncoupled setup. 
In this case Lipschitz continuity of $B$ implies already that $\bar B(x)$ is 
also Lipshitz continuous  in $x$, and so there exists a unique solution
$\brX=\brX_{x}$ of the averaged equation \index{averaged equation}
\begin{equation}\label{1.1.4} 
\frac {d\brX^\ve(t)}{dt}=\ve\brB(\brX^\ve(t)),\quad \brX^\ve(0)=x.
\end{equation}
In this case the standard averaging principle \index{averaging principle}
says (see \cite{SV}) that for $\mu$-almost all $y$,
\begin{equation}\label{1.1.5}
\lim_{\ve\to 0}\sup_{0\leq t\leq T/\ve}|X^\ve_{x,y}(t)-\brX_x^\ve(t)|=0.
\end{equation}

As the main motivation for the study of averaging is the setup of
perturbations described above we have to deal in real problems with the fully
coupled system (\ref{1.1.1}) \index{fully coupled system}
which only in very special situations can be
reduced by some change of variables to a much easier uncoupled case where 
the fast motion does not depend on the slow one. Observe that
in the general case (\ref{1.1.1}) the averaged vector field $\brB(x)$
in (\ref{1.1.3}) may even not be continuous in $x$, let alone Lipschitz,
and so (\ref{1.1.4}) may have many solutions or none at all. Moreover,
there may exist no natural well dependent on $x\in\bbR^d$ family of 
invariant measures $\mu_x$ since dynamical systems $F^t_x$ may have
rather different properties for different $x$'s. Even when all measures
$\mu_x$ are the same the averaging principle often does not hold true
 in the form (\ref{1.1.5}), for instance, in the presence of resonances
 (see \cite{LM} and \cite{Ki8}). Thus even basic results on approximation 
 of the slow motion by the averaged one in the fully coupled case cannot
 be taken for granted and they should be formulated in a different 
 way requiring usually stronger and more specific assumptions.

If convergence in (\ref{1.1.3}) is uniform in $x$ and $y$ then (see, for
instance, \cite{Ki4}) any limit point $\brZ(t)=\brZ_{x}(t)$ as $\ve\to 0$
of $Z^\ve_{x,y}(t)=X^\ve_{x,y}(t/\ve)$ is a solution of the averaged
equation
\begin{equation}\label{1.1.6}
\frac {d\brZ(t)}{dt}=\brB(\brZ(t)),\quad\brZ(0)=x.
\end{equation}
It is known that the limit in (\ref{1.1.3}) is uniform in $y$ 
if and only if the flow $F^t_x$ on $\bfM$ is uniquely ergodic, i.e. it 
possesses
a unique invariant measure, which occurs rather rarely. Thus, the uniform
convergence in (\ref{1.1.3}) assumption is too restrictive and excludes many 
interesting cases. Probably, the first relatively general result on fully
coupled averaging is due to Anosov \cite{An} (see also \cite{LM} and 
\cite{Ki4}). Relying
on the Liouville theorem he showed that if each flow $F^t_x$ preserves a
 probability measure $\mu_x$ on $\bfM$ having a $C^1$ dependent on $x$
 density with respect to the Riemannian volume $m$ on $\bfM$ and $\mu_x$
 is ergodic for Lebesgue almost all $x$ then for any
  $\del>0,$ 
\begin{equation}\label{1.1.7}
\mbox{mes}\{ (x,y):\,\,\sup_{0\leq t\leq T/\ve}|X^\ve_{x,y}(t)-\brX^\ve_x(t)|>
\del\}\to 0\,\,\mbox{as}\,\,\ve\to 0,
\end{equation}
where mes is the product of $m$ and the Lebesgue measure in a relatively 
compact domain $\cX\subset\bbR^d$. An example in Appendix to \cite{Ki8}
shows that, in general, this convergence in measure cannot be strengthened
to the convergence for almost all initial conditions and, moreover, in this
example the convergence (\ref{1.1.5}) does not hold true for any initial
condition from a large open domain. Such examples exist due to the presence
of resonances \index{resonance}, more specifically to the "capture into
resonance" phenomenon, which is rather well understood in 
perturbations of integrable Hamiltonian systems. Resonances lead there to the 
wealth of ergodic invariant measures and to different time and space averaging. 
It turns out (see \cite{BK1}) that wealth of ergodic invariant
measures with nice properties (such as Gibbs measures) for Axiom A and 
expanding dynamical systems also yields in the fully coupled averaging
setup with the latter fast motions examples of nonconvergence as $\ve\to 0$
for large sets of initial conditions (see Remark \ref{rem1.2.12}).

In Hamiltonian systems, which are a classical object for applications of
averaging methods, the whole space is fibered into manifolds of constant
energy. For some mechanical systems these manifolds have negative curvature 
with respect to the natural metric and their motion is described by geodesic
flows there. Hyperbolic Hamiltonian systems were discussed, for instance, in
\cite{LMM} and a specific example of a particle in a magnetic field leading
to such systems was considered recently in \cite{Tai}. Of course, these
lead to Hamiltonian systems which are far from integrable.
Such situations fall in our framework and they are among main
motivations for this work. This suggests to consider the equation (\ref{1.1.1})
 on a (locally trivial) fiber bundle $\cM=\{(x,y):\, x\in U,\, y\in M_x\}$ with
a base $U$ being an open subset in a Riemannian manifold $N$ and fibers $M_x$
 being diffeomorphic compact Riemannian manifolds (see \cite{Ta}). On the other
  hand, $\cM$ has a local product structure and if $\| B\|$ is bounded then
  the slow motion stays in one chart during time intervals of order $\Del/\ve$
  with $\Del$ small enough. Hence, studying behavior of solutions of 
  (\ref{1.1.1}) on each such time interval separately we come back to the
  product space $\bbR^d\times\bfM$ setup and will only have to piece
  results together to see the picture on a larger time interval of length
  $T/\ve.$
 
 We assume in the first part of this work that $b(x,y)$ is $C^2$ in $x$ and $y$
  and that for
 each $x$ in a closure of a relatively compact domain $\cX$ the flow $F^t_x$
 is Anosov or, more generally, Axiom A in a neighborhood of an attractor
 $\La_x.$ Let $\SRBx$ be the Sinai-Ruelle-Bowen (SRB) invariant
 measure \index{SRB measure}
 of $F^t_x$ on $\La_x$ and set $\brB(x)=\int B(x,y)d\SRBx(y).$ It is 
 known (see \cite{Co}) that the vector field $\brB(x)$ is Lipschitz 
 continuous in $x,$ and so the averaged equations (\ref{1.1.4}) and 
 (\ref{1.1.6}) have unique solutions $\brX^\ve(t)$
 and $\brZ(t)=\brX^\ve(t/\ve).$ Still, in general, the measures $\SRBx$
 are singular with respect to the Riemannian volume on $\bfM$, and so the
 method of \cite{An} cannot be applied here. We proved in \cite{Ki7} that,
 nevertheless, (\ref{1.1.7}) still holds true in this case, as well, and, 
 moreover, the measure in (\ref{1.1.7}) can be estimated by $e^{-c/\ve}$ with
 some $c=c(\del)>0$. The convergence (\ref{1.1.7}) itself without an exponential
  estimate can be proved by another method (see \cite{Ki9})
 which can be applied also to some partially hyperbolic fast motions
 \index{partially hyperbolic}. An extension of the averaging principle in
 the sense of convergence of Young measures \index{Young measures} is discussed
 in Section \ref{sec1.11}.
 
 Once the convergence of 
 $Z^\ve_{x,y}(t)=X^\ve_{x,y}(t/\ve)$ to $\brZ_{x}(t)=\brX^\ve_{x}(t/\ve)$
  as $\ve\to 0$ is established it is interesting to study the asymptotic
  behavior of the normalized error
 \begin{equation}\label{1.1.8}
 V_{x,y}^{\ve,\te}(t)=\ve^{\te-1}(Z^\ve_{x,y}(t)-\brZ_x(t)),\quad
 \te\in [\frac 12,1].
 \end{equation}
 Namely, in our situation it is natural to study the distributions
  $m\{ y:\, V_{x,y}^{\ve,\te}(\cdot)\in A\}$ as $\ve\to 0$ where $m$ is the
   normalized Riemannian volume on $\bfM$ and $A$ is a Borel subset in the
   space $C_{0T}$ of continuous paths $\vf(t),\,t\in[0,T]$ on $\bbR^d.$ We 
   will obtain in this work large deviations bounds for
 $V_{x,y}^{\ve}=V_{x,y}^{\ve,1}$ which will give, in particular, the result
 from \cite{Ki7} saying that
 \begin{equation}\label{1.1.9}
 m\{ y:\,\|V_{x,y}^{\ve}\|_{0,T}>\del\}\to 0\quad\mbox{as}\quad\ve\to 0
 \end{equation}
 exponentially fast in $1/\ve$ where $\|\cdot\|_{0,T}$ is the uniform norm
 on $C_{0T}.$ However, the main goal of this work is not to provide another
 derivation of (\ref{1.1.9}) but to obtain precise upper and lower large 
 deviations bounds which not only estimate measure of sets of initial
 conditions for which the slow motion $Z^\ve$ exhibits substantially
 different behavior than the averaged one $\bar Z$ \index{averaged motion} 
 but also enable us to go further and to investigate much longer exponential
 in $1/\ve$ time behavior of 
 $Z^\ve$. Namely, we will be able to study exits of the slow motion from a 
 neighborhood of an attractor of the averaged one and transitions of $Z^\ve$ 
 between basins of attractors of $\bar Z$. Such evolution, which
 becomes visible only on much longer than $1/\ve$ time scales, is usually
 called adiabatic in the framework of averaging. In the simpler case when
 the fast motion does not depend on the slow one such results were 
 discussed in \cite{Ki2}. Still, even in this uncoupled situation descriptions
 of transitions of the slow motion between attractors of the averaged one
 were not justified rigorously both in the Markov processes case of
 \cite{Fre} and in the dynamical systems case of \cite{Ki2}. Extending
 these technique to three scale equations may exhibit stochastic resonance
 \index{stochastic resonance}
 type phenomena producing a nearly periodic motion of the slowest motion 
 which is described in Section \ref{sec1.10} below.
 These problems seem to be important in the study of
  climate--weather interactions and they were discussed in \cite{CMP} and 
  \cite{Ha1} in the framework of a model describing transitions between
  steady climatic states with weather evolving as a fast chaotic system
  and climate playing the role of the slow motion.
   Such "very long" time description of the slow motion is usually
  impossible in the traditional averaging setup which deals with perturbations
  of integrable Hamiltonian systems.
  In the fully coupled situation we cannot work just with one 
  hyperbolic flow but have to consider continuously changing fast motions which
  requires a special technique. In particular, the full flow $\Phi_\ve^t$ on 
  $\bbR^d\times\bfM$ defined above and viewed as a small perturbation of the 
  partially hyperbolic system $\Phi^t$ plays an important role in our 
  considerations. It is somewhat surprising that the "very long time" behavior
  of the slow motion which requires certain "Markov property type" arguments
   still can be described in the fully coupled setup which involves
   continuously changing fast hyperbolic motions. It turns out that the 
   perturbed system still possesses semi-invariant expanding cones 
   \index{expanding cones} and foliations and a certain volume lemma 
   \index{volume lemma} type result on expanding leaves \index{expanding leaves}
   plays an important role in our argument for transition from small time 
   were perturbation techniques still works to the long and "very long" time
    estimates.
 
 It is plausible that moderate deviations \index{moderate deviations}
 type results can be
 proved for $V_{x,y}^{\ve,\te}$ when $1/2<\te<1$ and that the distribution
 of $V_{x,y}^{\ve,1/2}(\cdot)$ in $y$ converges to the distribution of a 
 Gaussian diffusion process in $\bbR^d$. Still, this requires somewhat
  different methods and it will not be discussed here. In this regard we
  can mention limit theorems obtained in \cite{CD} for a system of two
  heavy and light particles which leads to an averaging setup for a billiard
  flow. For the simpler case when $b$ does not depend on $x,$ i.e. when all 
  flows $F_x^t$ are the same, the moderate deviations and Gaussian 
  approximations results were obtained previously in 
   \cite{Ki3}. Related results in this uncoupled situation concerning
   Hasselmann's nonlinear (strong) diffusion approximation of the slow motion
   $X^\ve$ were obtained in \cite{Ki8}. 
  
  We consider also the discrete time case where (\ref{1.1.1}) is replaced by 
  difference equations \index{difference equations}
  for sequences $X^\ve(n)=X^\ve_{x,y}(n)$ and
  $Y^\ve(n)=Y^\ve_{x,y}(n),\, n=0,1,...$ so that
  \begin{eqnarray}\label{1.1.10}
  &X^\ve(n+1)-X^\ve(n)=\ve B(X^\ve(n),Y^\ve(n)),\\
  &Y^\ve(n+1)= F_{X^\ve(n)}Y^\ve(n),\,\, X^\ve(0)=x,Y^\ve(0)=y
  \nonumber
  \end{eqnarray}
  where $B :\cX\times \bfM\to\bbR^d$ is Lipschitz in both variables
  and the maps $F_x:\bfM\to\bfM$ are smooth and depend 
  smoothly on the parameter $x\in\bbR^d$. Introducing the map 
  \[
  \Phi_\ve(x,y)=(X^\ve_{x,y}(1),Y^\ve_{x,y}(1))=(x+\ve B(x,y),F_xy)
  \]
  we can also view this setup as a perturbation of the map $\Phi(x,y)=
  (x,F_xy)$ describing an ideal system where parameters $x\in\bbR^d$ do
  not change. Assuming that $F_x,\, x\in\bbR^d$ are $C^2$ depending on
  $x$ families of either $C^2$ expanding transformations or $C^2$ Axiom A 
  diffeomorphisms in a neighborhood of an attractor $\La_x$ we will derive
  large deviations estimates for the difference $X^\ve_{x,y}(n)-
  \brX^\ve_x(n)$ where $\brX^\ve=\brX^\ve_x$ solves the equation
  \begin{equation}\label{1.1.11}
  \frac {d\brX^\ve(t)}{dt}=\ve\brB(\brX^\ve(t)),\,\,\,\brX^\ve(0)=x
  \end{equation}
  where $\brB(x)=\int B(x,y)d\SRBx(y)$ and $\SRBx$ is the 
  corresponding SRB invariant measure of $F_x$ on $\La_x$. The discrete
  time results are obtained, essentially, by simplifications of the
  corresponding arguments in the continuous time case which enable us
  to describe "very long" time behavior of the slow motion also in the
  discrete time case. Since our methods work not only for fast motions being
  Axiom A diffeomorphisms but also when they are expanding transformations
  \index{expanding transformations}
  we can construct simple examples satisfying conditions of our theorems and
  exhibiting corresponding effects. In particular, we produce in Section 
 \ref{sec1.9} computational examples which demonstrate transitions of the slow
  motion between neighborhoods of attractors of the averaged system.  
  
  A series of related results for the case when ordinary differential
  equations in (\ref{1.1.1}) are replaced by fully coupled stochastic
  differential equations appeared in \cite{Kh2}, \cite{Ve1}--\cite{Ve3}, 
  \cite{PV}, and \cite{Ba2}. Hasselmann's nonlinear (strong) diffusion 
  approximation of the slow motion in the fully coupled stochastic differential
  equations setup was justified in \cite{BK}. When the fast process does not 
  depend on the slow one such results were obtained in \cite{Kh1}, \cite{Fre}, 
  and \cite{Ki6}. Especially relevant for our results here is \cite{Ve2}
  and we employ some elements of the probabilistic strategy from this
  paper. Still, the methods there are quite different from ours and they are
  based heavily, first, on the Markov property of processes emerging there and, 
  secondly, on uniformity and nondegeneracy of the fast diffusion term 
assumptions which cannot be satisfied in our circumstances as our deterministic
  fast motions are very degenerate from this point of view. Note that the proof
  in \cite{Ve2} contains a vicious cycle and substantial gaps which recently 
  were essentially fixed in \cite{Ve3}. Some of the dynamical systems technique
  here resembles 
 \cite{Ki2} but the dependence of the fast motion on the slow one complicates
 the analysis substantially and requires additional machinery. 
 A series of results on Cramer's type  asymptotics for fully coupled averaging 
 with Axiom A diffeomorphisms as fast motions appeared recently in 
 \cite{Ba1}--\cite{Ba4}. Observe that the methods there do not work
  for continuous time Axiom A dynamical systems considered here, they cannot
  lead, in principle, to the standard large deviations estimates of our
  work and they deal with deviations of $X^\ve$ from the averaged motion only 
  at the last moment and not of its whole path. Various limit theorems for
  the difference equations setup (\ref{1.1.10}) with partially hyperbolic
  fast motions were obtained recently in \cite{Dol0} and \cite{Dol}.
  
  The study of deviations from the averaged motion in the fully coupled case
  seems to be quite important for applications, especially, from 
  phenomenological point of view. In addition to perturbations of Hamiltonian 
  systems mentioned above there are many non Hamiltonian systems which are 
  naturally to consider from the beginning as a combination of fast and slow 
  motions. For instance, Hasselmann \cite{Ha} based his model of 
  weather--climate interaction on the assumption that weather is a fast 
  chaotic motion depending on climate as a slow motion which differs from the 
  corresponding averaged motion mainly by a diffusion term. Though, as shown
  in \cite{Ki6}, \cite{BK} and \cite{Ki8}, this diffusion error term
  does not help in the study of large deviations which are responsible for
  rare transitions of the slow motion between attractors of the averaged one,
  the latter phenomenon can be described in our framework and it seems to be 
  important in certain models of climate fluctuations (see \cite{CMP} and
  \cite{Ha1}). Very slow nearly periodic motions appearing in the stochastic
  resonance framework considered in Section \ref{sec1.10} may also fit into
 this subject in the discussion on "ice ages". Of course, it is hard to believe
   that real world chaotic systems can be described precisely by an Anosov or
   Axiom A flow but one may take comfort in the Chaotic Hypothesis \cite{Ga}:
   " A chaotic mechanical system can be regarded for practical purposes as a
   topologically mixing Anosov system".

 \section{Main results}\label{sec1.2}\setcounter{equation}{0}

Let $F^t$ be a $C^2$ flow on a compact Riemannian manifold $\bfM$ given by
a differential equation
\begin{equation}\label{1.2.1}
\frac {dF^ty}{dt}=b(F^ty),\,\, F^0y=y.
\end{equation}
A compact $F^t-$invariant set
$\La\subset \bfM$ is called hyperbolic \index{hyperbolic set}
if there exists $\ka>0$ and the 
splitting $T_\La \bfM=\Gam^s\oplus\Gam^0\oplus\Gam^u$ into the continuous 
subbundles $\Gam^s,\Gam^0,\Gam^u$ of the tangent bundle $T\bfM$ restricted
to $\La,$ the splitting is invariant with respect to the differential $DF^t$
of $F^t,$ $\Gam^0$ is the one dimensional subbundle generated by the vector
field $b$, and there is $t_0>0$ such that for all
$\xi\in\Gam^s,\,\eta\in\Gam^u,$ and $t\geq t_0,$
\begin{equation}\label{1.2.2}
\| DF^t\xi\|\leq e^{-\ka t}\|\xi\|\quad\mbox{and}\quad\|DF^{-t}\eta\|\leq 
e^{-\ka t}\|\eta\|.
\end{equation}
A hyperbolic set $\La$ is said to be basic hyperbolic \index{basic hyperbolic
set} if the periodic orbits 
of $F^t|_\La$ are dense in $\La,$ $F^t|_\La$ is topologically transitive,
and there exists an open set $U\supset\La$ with $\La=\cap_{-\infty<t<\infty}
F^tU.$ Such a $\La$ is called a basic hyperbolic attractor \index{hyperbolic
 attractor} if for some open
 set $U$ and $t_0>0,$
\begin{equation*}
F^{t_0}\brU\subset U\quad\mbox{and}\quad\cap_{t>0}F^tU=\La
\end{equation*}
where $\brU$ denotes the closure of $U.$ If $\La=\bfM$ then $F^t$ is called an 
Anosov flow.
\begin{assumption}\label{ass1.2.1} The family $b(x,\cdot)$ in (\ref{1.1.2}) 
consists of $C^2$ vector fields on a compact $n_\bfM$-dimensional Riemannian 
manifold $\bfM$ with uniform $C^2$ dependence on the parameter $x$ belonging to 
a neighborhood of the closure $\bar\cX$ of a relatively compact open connected 
set $\cX\subset\bbR^d$. Each flow $F_x^t,\, x\in\bar\cX$ on $\bfM$ given by 
\begin{equation}\label{1.2.3}
\frac {dF_x^ty}{dt}=b(x,F_x^ty),\quad F_x^0y=y
\end{equation}
possesses a basic hyperbolic attractor $\La_x$ with a splitting 
$T_{\La_x}\bfM=\Gam^s_x\oplus\Gam^0_x\oplus\Gam^u_x$ satisfying (\ref{1.2.2})
with the same $\ka>0$ and there exists an
open set $\cW\subset \bfM$ and $t_0>0$ such that
\begin{equation}\label{1.2.4}
\La_x\subset \cW,\,\, F^t_x\brcW\subset \cW\,\forall t\geq t_0,\,\,\mbox{and}\,\,
\cap_{t>0}F^t_x\cW=\La_x\,\,\forall x\in\bar\cX.
\end{equation}
\end{assumption}

Let $J_x^u(t,y)$ be the absolute value of the
Jacobian of the linear map $DF_x^t(y):\Gam_x^u(y)\to
\Gam^u_x(F^t_xy)$ with respect to the Riemannian inner products and set
\begin{equation}\label{1.2.5}
 \vf_x^u(y)=-\frac {dJ_x^u(t,y)}{dt}\big\vert _{t=0}.
 \end{equation}
 The function $\vf_x^u(y)$ is known to be H\" older continuous in $y,$ since
 the subbundles $\Gam^u_x$ are H\" older continuous (see \cite{BR} and
  \cite{KH}), and $\vf_x^u(y)$ is $C^1$ in $x$ (see \cite{Co}).
  
  Let $\cW$ satisfy (\ref{1.2.4}) and set $\cW^t_x=\{ y\in \cW:\, F^s_xy\in
  \bar\cW\,\,\,
  \forall s\in [0,t]\}.$ A set $E\subset \cW^t_x$ is called 
  $(\del,t)-$separated \index{$(\cdot,\cdot)$-separated set}
  for the flow $F_x$ if $y,z\in E$, $y\ne z$ imply $d(F^s_xy, F_x^sz)>\del$
  for some $s\in[0,t]$, where $d(\cdot,\cdot)$ is the distance function on
  $\bfM.$ For each continuous function $\psi$ on $\cW$ set
  $P_x(\psi,\del,t)=\sup\{\sum_{y\in E}\exp\int_0^t\psi(F^s_xy)ds:\, E\subset
  \cW_x^t\,\,\mbox{is}\,\,(\del,t)-\mbox{separated for}\,\, F_x\},$
  $P_x(\psi,\del,t)=0$ if $\cW_x^t=\emptyset,$ and
  $$
  P_x(\psi,\del)=\limsup_{t\to\infty}\frac 1t\log P_x(\psi,\del,t).
  $$
  The latter is monotone in $\del,$ and so the limit
  \begin{equation*}
  P_x(\psi)=\lim_{\del\to 0}P_x(\psi,\del)
  \end{equation*}
  exists and it is called the topological pressure \index{topological pressure}
  of $\psi$ for the flow 
  $F^t_x.$ Let $\cM_x$ denotes the space of $F^t_x-$invariant probability
  measures \index{invariant measure}
  on $\La_x$ then (see, for instance, \cite{KH}) the following
  variational principle \index{variational principle}
  \begin{equation}\label{1.2.6}
  P_x(\psi)=\sup_{\mu\in\cM_x}(\int\psi d\mu+h_\mu(F^1_x))
  \end{equation}
  holds true where $h_\mu(F^1_x)$ is the Kolmogorov--Sinai entropy 
  \index{entropy} of the
  time-one map $F^1_x$ with respect to $\mu.$ If $q$ is a H\" older 
  continuous function on $\La_x$ then there exists a unique $F^t_x-$invariant
  measure $\mu^q_x$ on $\La_x,$ called the equilibrium state \index{equilibrium
  state} for $\vf_x^u+q,$
  such that 
  \begin{equation}\label{1.2.7}
  P_x(\vf_x^u+q)=\int (\vf_x^u+q)d\mu_x^q+h_{\mu_x^q}(F^1_x).
  \end{equation}
  We denote $\mu_x^0$ by $\SRBx$ since it is usually called the Sinai--Ruelle--
  Bowen (SRB) measure \index{SRB measure}. Since $\La_x$ are attractors we 
  have that $P_x(\vf_x^u)=0$ (see \cite{BR}).
  
  For any probability measure $\nu$ on $\brcW$ define \index{$I$-functional}
  \begin{equation}\label{1.2.8}
   I_x(\nu)=\left\{\begin{array}{ll}
  -\int\vf_x^ud\nu-h_\nu(F^1_x) &\mbox{if $\nu\in\cM_x$}\\
  \infty &\mbox{otherwise.}
  \end{array}\right.
  \end{equation}
 Then
 \begin{equation*}
 P_x(\vf^u_x+q)=\sup_{\nu}(\int qd\nu-I_x(\nu)).
 \end{equation*}
 Observe that by the Ruelle inequality (see, for instance, \cite{KH},
 Theorem S.2.13), $I_x(\nu)\geq 0$, and so in view of Assumption \ref{ass1.2.1}
 for any $\nu\in\cM_x$,
 \begin{equation}\label{1.2.9}
 I_x(\nu)\leq \sup_{y\in\La_x}|\vf^u_x(y)|\leq\sup_{x\in\bar\cX,y\in\La_x}|
 \vf^u_x(y)|<\infty .
 \end{equation}
  It is known that $h_\nu(F^1_x)$ is upper semicontinuous in $\nu$ since
  hyperbolic flows are entropy expansive ( see \cite{Bo1}). Thus $I_x(\nu)$
  is a lower semicontinuous functional in $\nu$ and it is also convex (and
  affine on $\cM_x$) since entropy $h_\nu$ is affine in $\nu$ (see, for 
  instance, \cite{Wa}). Hence, by the duality theorem \index{convex duality
  theorem}(see \cite{AE}, p.201),
  \begin{equation*}
  I_x(\nu)=\sup_{q\in\cC(\bfM)}(\int qd\nu-P_x(\vf^u_x+q)).
  \end{equation*}
  Observe that this formula can be proved more directly. Namely, if we
  define $I_x(\nu)$ by it in place of (\ref{1.2.8}) then (\ref{1.2.8})
  follows for $\nu\in\cM_x$ from Theorem 9.12 in \cite{Wa} and
  it is easy to show directly that $I_x(\nu)$ defined in this way equals
  $\infty$ for any finite signed measure $\nu$ which is not $F_x$-invariant.
  
  Since we assume that the vector field $B$ is $C^1$ in both arguments 
  (here only continuity in $y$ is needed) then
 for any $x,x'\in\cX$ and $\al,\be\in\bbR^d$ we can define $H(x,x',\be)=
 P_x(<\be, B(x',\cdot)>+\vf^u_x)$ and $H(x,\be)=H(x,x,\be)$. Then
  \begin{eqnarray}\label{1.2.10}
  &H(x,x',\be)=\sup_{\nu}\big(\int <\be,B(x',y)>d\nu(y)-I_x(\nu)\big)\\
  &=\sup_{\al\in\bbR^d}\big(<\al,\be>-L(x,x',\al)\big)\nonumber
  \end{eqnarray}
  where
  \begin{equation}\label{1.2.11}
  L(x,x',\al)=\inf\{ I_x(\nu):\,\int B(x',y)d\nu(y)=\al\}
  \end{equation}
  if $\nu\in\cM_x$ satisfying the condition in brackets exists and 
  $L(x,x',\al)=\infty$, otherwise. Since, $H(x,x',\be)$ is convex and 
  continuous the duality theorem (see \cite{AE}, p.201) yields that
  \begin{equation}\label{1.2.12}
  L(x,x',\al)=\sup_{\be\in\bbR^d}\big(<\al,\be>-H(x,x',\be)\big)
  \end{equation}
  provided there exists a probability measure $\nu\in\cM_x$ such that 
  $\int B(x',y)d\nu(y)=\al$ and $L(x,x',\al)=\infty$, otherwise. Clearly,
  $L(x,x',\al)$ is convex and lower semicontinuous in all arguments and,
  in particular, it is measurable. We set also $L(x,\al)= L(x,x,\al)$.
  
  Denote by $C_{0T}$ the space of continuous curves $\gam_t=\gam(t),\, t\in 
  [0,T]$ in $\cX$ which is the space of continuous maps of $[0,T]$ into $\cX.$
  For each absolutely continuous $\gam\in C_{0T}$ its velocity $\dot\gam_t$ 
  can be obtained as the almost everywhere limit of continuous functions
$n(\gam_{t+n^{-1}}-\gam_t)$ when $n\to\infty$. Hence $\dot\gam_t$ is 
measurable in $t$, and so we can set \index{$S$-functional}
  \begin{equation}\label{1.2.13}
  S_{0T}(\gam)=\int_0^TL(\gam_t,\dot{\gam}_t)dt=\int_0^T
  \inf\{ I_{\gam_t}(\nu):\,\dot{\gam}_t=\brB_\nu (\gam_t),\,
  \nu\in\cM_{\gam_t}\}dt,
  \end{equation}
  where $\brB_\nu(x)=\int B(x,y)d\nu(y),$ provided for Lebesgue almost all
  $t\in [0,T]$ there exists $\nu_t\in\cM_{\gam_t}$ for which $\dot{\gam}_t=
  \brB_{\nu_t}(\gam_t),$ and $S_{0T}(\gam)=\infty$ otherwise. It follows
  from \cite{BR} and \cite{Co} that 
  \[ 
  S_{0T}(\gam)\geq S_{0T}(\gam^u)=-\int_0^TP_{\gam_t^u}(\vf^u_{\gam_t^u})dt=0
  \]
   where $\gam_t^u$ is the unique solution of the equation
   \begin{equation}\label{1.2.14}
   \dot{\gam}_t^u=\brB(\gam_t^u),\quad\gam_0^u=x,
   \end{equation}
   where $\brB(z)=\brB_{\SRBz}(z),$ and the equality $S_{0T}(\gam)=0$ holds
   true if and only if $\gam=\gam^u.$
   
   Define the uniform metric on $C_{0T}$ by 
   \[\bfr_{0T}(\gam,\eta)=\sup_{0\leq t\leq T}|\gam_t-\eta_t|
   \]
   for any $\gam,\eta\in C_{0T}.$ Set 
   \[
   \Psi^a_{0T}(x)=\{\gam\in C_{0T}:\,\gam_0=x,\, S_{0T}(\gam)\leq a\}.
   \]
   Since $L(x,\al)$ is lower semicontinuous
   and convex in $\al$ and, in addition, $L(x,\al)=\infty$ if $|\al|>
   \sup_{y\in\bfM}|B(x,y)|$ we conclude that the conditions of Theorem 3 in 
   Ch.9 of \cite{IT} are satisfied as we can choose a fast growing minorant
   of $L(x,\al)$ required there to be zero in a sufficiently large ball and
   to be equal, say, $|\al|^2$ outside of it. As a result, it follows that
   $S_{0T}$ is lower semicontinuous functional on $C_{0T}$ with respect to
   the metric $\bfr_{0T}$, and so $\Psi^a_{0T}(x)$ is a closed set which plays 
   a crucial role in the large deviations \index{large deviations}
   arguments below.
   
   We suppose that the coefficients of (\ref{1.1.1}) satisfy the following
   \begin{assumption}\label{ass1.2.2}
   There exists $K>0$ such that
   \begin{equation}\label{1.2.15}
   \| B(x,y)\|_{C^1(\cX\times\bfM)}+\| b(x,y)\|_{C^2(\cX\times\bfM)}\leq K
   \end{equation}
   where $\|\cdot\|_{C^i(\cX\times\bfM)}$ is the $C^i$ norm of the corresponding
   vector fields on $\cX\times\bfM.$
   \end{assumption}
   Set $\cX_t=\{ x\in\cX:\, X^\ve_{x,y}(s)\in\cX$ and $\bar X_x^\ve(s)\in\cX$
   for all $y\in\bar\cW,\, s\in[0,t/\ve],\,\ve>0\}.$ Clearly, $\cX_t\supset
   \{ x\in\cX:\,\inf_{z\in\partial\cX}|x-z|\geq 2Kt\}.$
   The following is one of the main results of this paper.
   
   \begin{theorem}\label{thm1.2.3} Suppose that $x\in\cX_T$ and $X^\ve_{x,y}$, 
   $Y^\ve_{x,y}$ are solutions of (\ref{1.1.1}) with coefficients satisfying 
   Assumptions \ref{ass1.2.1} and \ref{ass1.2.2}. Set
    $Z^\ve_{x,y}(t)=X^\ve_{x,y}(t/\ve)$ then for any $a,\del,\la>0$ and every
    $\gam\in C_{0T},\,\gam_0=x$ there exists $\ve_0=\ve_0(x,\gam,a,\del,\la)
    >0$ such that for $\ve<\ve_0,$
    \begin{equation}\label{1.2.16}
    m\left\{ y\in \cW:\,\bfr_{0T}(Z^\ve_{x,y},\gam)<\del\right\}\geq
    \exp\left\{-\frac 1\ve(S_{0T}(\gam)+\la)\right\}
    \end{equation}
    and
    \begin{equation}\label{1.2.17}
    m\left\{ y\in \cW:\,\bfr_{0T}(Z^\ve_{x,y},\Psi^a_{0T}(x))\geq\del\right\}
    \leq\exp\left\{-\frac 1\ve(a-\la)\right\}
    \end{equation}
    where, recall, $m$ is the normalized Riemannian volume on $\bfM.$ The
    functional $S_{0T}(\gam)$ for $\gam\in C_{0T}$ is finite if and only if
    $\dot{\gam}_t=\brB_{\nu_t}(\gam_t)$ for $\nu_t\in\cM_{\gam_t}$ and
    Lebesgue almost all $t\in [0,T].$ Furthermore, $S_{0T}(\gam)$ achieves
    its minimum 0 only on $\gam^u$ satisfying (\ref{1.2.14})
     for all $t\in[0,T].$ Finally, for any $\del>0$ there
    exist $c(\del)>0$ and $\ve_0>0$ such that for all $\ve<\ve_0,$
    \begin{equation}\label{1.2.18}
    m\left\{ y\in \cW:\,\bfr_{0T}(Z^\ve_{x,y},\brZ_x)\geq\del\right\}\leq
    \exp\left(-\frac {c(\del)}\ve\right)
    \end{equation}
    where $\brZ_x=\gam^u$ is the unique solution of (\ref{1.2.14}).
    \end{theorem}
    
    Observe that (\ref{1.2.18}) (which was proved already in \cite{Ki7} by a
    less precise large deviations argument) follows from (\ref{1.2.17}) and the
    lower semicontinuity of the functional $S_{0T}$ and it
    says, in particular, that $Z^\ve_{x,\cdot}$
    converges to $\bar Z_x$ in measure on the space $(\cW,m)$ with respect
    to the metric $\bfr_{0T}$. It is naturally to ask whether
    we have here also the convergence for $m$-almost all $y\in\cW$. An example
    due to A.Neishtadt discussed in \cite{Ki8} shows that in the classical 
    situation of perturbations of integrable Hamiltonian systems, 
    \index{Hamiltonian system} in general,
    the averaging principle holds true only in the sense of convergence 
    in measure on the space of intitial conditions but not in the sense of
    the almost everywhere convergence. This example concerns the simple 
    system $\dot {I}=\ve(4+8\sin\vf-I),\,\dot {\vf}=I$ with the one dimensional
    slow motion $I$ and the fast motion $\vf$ evolving on the circle while
    the corresponding averaged motion $J=\bar I$ satisfies the equation 
    $\dot {J}=\ve(4-J)$. The resonance \index{resonance}
    occurs here only when $I=0$ but it 
    suffices already to create troubles in the averaging principle. Namely,
    it turns out that for any initial condition $(I_0,\vf_0)$ with 
    $-2<I_0<-1$ there exists a sequence $\ve_n\to 0$ such that
     $I^{\ve_n}_{I_0,\vf_0}(1/\ve_n)<J^{\ve_n}_{I_0}(1/\ve_n)-3/2$
    though, of course, convergence in measure holds true here (see
    \cite{LM}). Recently (see \cite{BK1} and Remark \ref{rem1.2.12}), 
    such nonconvergence examples were constructed for the difference
    equations averaging setup (\ref{1.1.10}) with expanding fast motions
    and there is no doubt that such examples exist also in the continuous
    time setup (\ref{1.1.1}) when fast motions are Axiom A flows as in this
    paper. Observe also that (\ref{1.2.16}) and (\ref{1.2.17}) remain true
    (with the same proof) if we replace $m$ there by $\SRBx$ but as an
    example in \cite{BK1} shows we cannot, in general, replace $m$ there by
    an arbitrary Gibbs measure $\mu_x$ of $F^t_x$.
    
     Next, let $V\subset\cX$ be a connected open set and put 
     $\tau^\ve_{x,y}(V)=\inf\{ t\geq 0:\, Z^\ve_{x,y}(t)
    \notin V\}$ where we take $\tau^\ve_{x,y}(V)=\infty$ if $X^\ve_{x,y}(t)\in
    V$ for all $t\geq 0.$ The following result follows directly from
    Theorem \ref{thm1.2.3}.
    \begin{corollary}\label{cor1.2.4} Under the conditions of Theorem 
    \ref{thm1.2.3} for any $T>0$ and $x\in V,$
    \begin{eqnarray*}
    &\lim_{\ve\to 0}\ve\log m\left\{ y\in \cW:\,\tau^\ve_{x,y}(V)<T\right\}\\
    &=-\inf\left\{ S_{0t}(\gam):\,\gam\in C_{0T},\, t\in[0,T],\,\gam_0=x,\,
    \gam_t\not\in V\right\}\nonumber.
    \end{eqnarray*}
    \end{corollary}
     
     Precise large deviations bounds such as (\ref{1.2.16}) and (\ref{1.2.17})
     of Theorem \ref{thm1.2.3} (which will be needed uniformly on certain
     unstable discs) are crucial in our study in Sections \ref{sec1.7}
     and \ref{sec1.8} of the "very long", i.e. exponential in 
    $1/\ve$, time "adiabatic" behaviour of the slow motion which cannot be 
    described usually in the traditional theory of averaging where only 
    perturbations of integrable Hamiltonian systems are considered. Namely, we
    will describe such long time behavior of $Z^\ve$ in terms of the function 
     \[
     R(x,z)=\inf_{t\geq 0,\gam\in C_{0t}}\{ S_{0t}(\gam):\,\gam_0=x,\,
     \gam_t=z\}
     \]
     under various assumptions on the averaged motion $\brZ.$
     Observe that $R$ satisfies the triangle inequality $R(x_1,x_2)+
     R(x_2,x_3)\geq R(x_1,x_3)$ for any $x_1,x_2,x_3\in\cX$ and it determines
     a semi metric on $\cX$ which measures "the difficulty'" for the slow
     motion to move from point to point in terms of the functional $S$. 
     
     Introduce the averaged flow $\Pi^t$ on $\cX_t$ by
     \begin{equation}\label{1.2.19}
     \frac {d\Pi^tx}{dt}=\bar B(\Pi^tx),\,\, x\in\cX_t
     \end{equation}
     where, recall, $\bar B(z)=\bar B_\SRBz(z)$ and $\bar B_\nu(z)=
     \int B(z,y)d\nu(y)$ for any probability measure $\nu$ on $\bfM$.
     Call a $\Pi^t$-invariant compact set $\cO\subset\cX$ an $S$-compact if 
     for any $\eta>0$ there exist $T_\eta\geq 0$ and  
      an open set $U_\eta\supset\cO$ such that whenever $x\in\cO$ and 
     $z\in U_\eta$ we can pick up $t\in[0,T_\eta]$ and $\gam\in C_{0t}$ 
     satisfying
     \[
     \gam_0=x,\,\, \gam_t=z\,\,\mbox{and}\,\, S_{0t}(\gam)\leq\eta.
     \]
     It is clear from this definition that $R(x,z)=0$ for any pair points
     $x,z$ of an $S$-compact $\cO$ and by the above triangle inequality
     for $R$ we see that $R(x,z)$ takes on the same value when $z\in\cX$
     is fixed and $x$ runs over $\cO$.
     We say that the vector field $B$ on $\cX\times\bfM$ is complete 
     \index{complete} at 
     $x\in\cX$ if the convex set of vectors $\{\be\bar B_\nu(x):\,\be\in[0,1],
     \,\nu\in\cM_x\}$ contains an open neigborhood of the origin in $\bbR^d$.
     In Lemma \ref{lem1.6.2} we will show that if $\cO\subset\cX$ is a compact
     $\Pi^t$-invariant set such that $B$ is complete at each $x\in\cO$ and
      either $\cO$ contains a dense orbit of the flow $\Pi^t$ (i.e. $\Pi^t$ is
     topologically transitive \index{topologically transitive}
     on $\cO$) or $R(x,z)=0$ for any $x,z\in\cO$
      then $\cO$ is an $S$-compact. Moreover, to ensure that $\cO$ is an 
      $S$-compact it suffices to assume that $B$ is complete only at some 
      point of $\cO$ and the flow $\Pi^t$ on $\cO$ is minimal \index{minimal},
     i.e. the $\Pi^t$-orbits of all points are dense in $\cO$ or, equivalently,
       for any $\eta>0$ there exists
       $T(\eta)>0$ such that the orbit $\{\Pi^tx,\, t\in[0,T(\eta)]\}$ of
       length $T(\eta)$ of each point $x\in\cO$ forms an $\eta$-net in $\cO$.
       The latter condition obviously holds true when $\cO$ is a fixed point
       or a periodic orbit of $\Pi^t$ but also, more generally, when $\Pi^t$ 
       on $\cO$ is uniquely ergodic \index{uniquely ergodic}
       (see \cite{KH}, \cite{Ma}, \cite{Wa}).
       Among well known examples of uniquely ergodic flows we can mention 
       irrational translations of tori and horocycle flows on surfaces of 
       negative curvature.
       
      A compact $\Pi^t$-invariant set $\cO\subset\cX$ is called an attractor
      \index{attractor}
     (for the flow $\Pi^t$) if there is an open set $U\supset\cO$ and $t_U>0$
     such that
     \[
     \Pi^{t_U}\bar U\subset U\,\,\mbox{and}\,\,\lim_{t\to\infty}\,
     \mbox{dist}(\Pi^tz,\cO)=0\,\,\mbox{for all}\,\, z\in U.
     \]
     For an attractor $\cO$ the set $V=\{ z\in\cX:\,\lim_{t\to\infty}\,
     \mbox{dist}(\Pi^tz,\cO)=0\}$, which is clearly open, is called the
     basin \index{basin}(domain of attraction) of $\cO$. An attractor which
      is also
     an $S$-compact will be called an $S$-attractor \index{$S$-attractor}.
     
     In what follows we will speak about connected open sets $V$ with 
     piecewise smooth boundaries $\partial V$. The latter can be introduced in 
     various ways but it will be convenient here to adopt the definition from 
     \cite{Cow} saying that $\partial V$ is the closure of a finite union
     of disjoint, connected, codimension one, extendible $C^1$ (open or
     closed) submanifolds of $\bbR^d$ which are called faces of the boundary.
     The extendibility condition means that the closure of each face is a part 
     of a larger submanifold of the same dimension which coincides with the
     face itself if the latter is a compact submanifold. This enables us to 
     extend fields of normal vectors to the boundary of faces and to speak 
     about minimal angles between adjacent faces which we assume to be 
     uniformly bounded away from zero or, in other words, angles between 
     exterior normals to adjacent faces at a point of intersection of their 
     closures are uniformly bounded away from $\pi$ and $-\pi$. The following
     result which will be proved in Section \ref{sec1.7} describes exits of the
     slow motion from neighborhoods of attractors of the averaged motion.

     \begin{theorem}\label{thm1.2.5} Let $\cO\subset\cX$ be an $S$-attractor
      of the flow $\Pi^t$ whose basin contains the closure $\bar V$ of a 
     connected open set $V$ with a piecewise smooth boundary $\partial V$ 
     such that $\bar V\subset\cX$ and assume that for each $z\in\partial V$
     there exists $\vp=\vp(z)>0$ and an $F^t_z-$invariant probability
      measure $\nu=\nu_z$ on $\La_z$ such that
     \begin{equation}\label{1.2.20}
     z+s\bar B(z)\in V\,\,\,\mbox{but}\,\,\, z+s\bar B_\nu(z)\in
     \bbR^d\setminus\bar V\,\,\mbox{for all}\,\, s\in(0,\vp],
     \end{equation}
     i.e. $\bar B(z)\ne 0,\,\bar B_\nu(z)\ne 0$ and the former vector
     points out into the interior while the latter into the exterior 
     of $V$. Set $R_{\partial}(z)=\inf\{ R(z,\tilde z):\,\tilde z\in
     \partial V\}$ and $\partial_{\min}(z)=\{\tilde z\in\partial V:\,
      R(z,\tilde z)=R_\partial(z)\}$. Then $R_\partial(z)$ takes on the 
      same value $R_\partial$ and $\partial_{\min}(z)$ coincides with
      the same compact nonempty set $\partial_{\min}$ for all $z\in\cO$
      while $R_\partial(x)\leq R_\partial$ for all $x\in V$. Furthermore,
      for any $x\in V$,
     \begin{equation}\label{1.2.21}
     \lim_{\ve\to 0}\ve\log\int_{\cW}\tau^\ve_{x,y}(V)dm(y)=R_{\partial}>0
     \end{equation}
     and for each $\al>0$ there exists $\la(\al)=\la(x,\al)>0$ such that for 
     all small $\ve>0$,
     \begin{equation}\label{1.2.22}
      m\big\{ y\in\cW:\, e^{(R_\partial-\al)/\ve}>\tau^\ve_{x,y}(V)\,\,
     \mbox{or}\,\,\tau^\ve_{x,y}(V)>e^{(R_\partial+\al)/\ve}\big\}
     \leq e^{-\la(\al)/\ve}.
     \end{equation}
     Next, set
     \[
     \Te^\ve_v(t)=\Te_v^{\ve,\del}(t)=\int_0^t
     \bbI_{V\setminus U_\del(\cO)}(Z^\ve_v(s))ds
     \]
     where $U_\del(\cO)=\{ z\in\cX:$ dist$(z,\cO)<\del\}$ and $\bbI_\Gam(z)
     =1$ if $z\in\Gam$ and $=0$, otherwise. Then for any $x\in V$ and 
     $\del>0$ there exists $\la(\del)=\la(x,\del)>0$ such that for all small
     $\ve>0$,
     \begin{equation}\label{1.2.23}
     m\big\{ y\in\cW:\,\Te^\ve_{x,y}(\tau^\ve_{x,y}(V))\geq e^{-\la(\del)/\ve}
     \tau^\ve_{x,y}(V)\big\}\leq  e^{-\la(\del)/\ve}.
     \end{equation}
     Finally, for every $x\in V$ and $\del>0$,
     \begin{equation}\label{1.2.24}
     \lim_{\ve\to 0}m\big\{ y\in\cW:\,\mbox{dist}\big(Z^\ve_{x,y}
     (\tau^\ve_{x,y}(V)),\partial_{\min}\big)\geq\del\big\}=0
     \end{equation}
     provided $R_\partial <\infty$ and the latter holds true if and only if
      for some $T>0$ there exists 
     $\gam\in C_{0T},\,\gam_0\in\cO,\,\gam_T\in\partial V$ such that 
     $\dot {\gam}_t=\bar B_{\nu_t}(\gam_t)$ for Lebesgue almost all 
     $t\in[0,T]$ with $\nu_t\in\cM_{\gam_t}$ then $R_\partial <\infty$.
     \end{theorem}
     
     Theorem \ref{thm1.2.5} asserts, in particular, that typically the slow 
     motion $Z^\ve$ performs rare (adiabatic) fluctuations \index{adiabatic
     fluctuations} in the vicinity of 
     an $S$-attractor $\cO$ since it exists from any domain $U\supset\cO$ with 
     $\bar U\subset V$ for the time much smaller than $\tau^\ve(V)$ (as the 
     corresponding number $R_\partial=R_{\partial U}$ will be smaller) and by 
     (\ref{1.2.23}) it can spend in $V\setminus U_\del(\cO)$ only small 
     proportion of time which implies that $Z^\ve$ exits from $U$ and returns 
     to $U_\del(\cO)$ (exponentially in $1/\ve$) many times before it finally 
     exits $V$. We observe that in the much simpler uncoupled setup 
     corresponding results in the case of $\cO$ being an attracting point were 
     obtained for a continuous time Markov chain and an Axiom A flow as fast 
     motions in \cite{Fre} and \cite{Ki2}, respectively, but the proofs there 
     rely on the lower semicontinuity of the function $R$ which does not hold
     true in general, and so extra conditions like $S$-compactness of $\cO$
     or, more specifically, the completness of $B$ at $\cO$ should be assumed
     there, as well. It is important to observe that the intuition based on
     diffusion type small random perturbations of dynamical systems should be
     applied with caution to problems of large deviations in averaging since
     the $S$-functional of Theorem \ref{thm1.2.3} describing them is more
     complex and have rather different properties than the corresponding
     functional emerging in diffusion type random perturbations of dynamical
     systems (see \cite{FW}). The reason for this is the deterministic nature
     of the slow motion $Z^\ve$ which unlike a diffusion can move only with
     a bounded speed, and moreover, even in order to ensure its "diffusive
     like" local behaviour (i.e. to let it go in many directions) some extra
     nondegeneracy type conditions on the vector field $B$ are required.
     
     Our next result describes rare (adiabatic) transitions of the slow motion
     $Z^\ve$ between basins of attractors of the averaged flow $\Pi^t$ which 
     we consider now in the whole $\bbR^d$ and impose certain conditions on the
     structure of its $\om$-limit set \index{$\om$-limit set}.
     \begin{assumption}\label{ass1.2.6} Assumptions \ref{ass1.2.1} and
     \ref{ass1.2.2} hold true for $\cX=\bbR^d$, the family $\{ F^t_x,\, t\leq 1,
     \, x\in\bbR^d\}$ is a compact set of diffeomorphisms in the $C^2$
     topology, 
     \begin{equation}\label{1.2.25}
     \| B(x,y)\|_{C^2(\bbR^d\times\bfM)}\leq K
     \end{equation}
     for some $K>0$ independent of $x,y$ and there exists $r_0>0$ such that
     \begin{equation}\label{1.2.26}
     \big(x,B(x,y)\big)\leq -K^{-1}\,\,\mbox{for any}\,\, y\in\cW
     \,\,\mbox{and}\,\, |x|\geq r_0.
     \end{equation}
     \end{assumption}
     
     The condition (\ref{1.2.26}) means that outside of some ball all vectors
     $B(x,y)$ have a bounded away from zero projection on the radial direction
     which points out to the origin. This condition can be weakened, for
     instance, it suffices that
     \[
     \lim_{d\to\infty}\inf\{ R(x,z):\,\mbox{dist}(x,z)\geq d\}=\infty
     \]
     but, anyway, we have to make some assumption which ensure that the slow
     motion stays (at least, for "most" initial points $y\in\cW$) in a compact
     region where really interesting dynamics takes place.
     
     Next, suppose that the $\om$-limit set of the averaged flow $\Pi^t$ is 
     compact and it consists of two parts, so that the first part is a finite 
     number of $S$-attractors $\cO_1,...,\cO_\ell$ whose basins
     $V_1,...,V_\ell$ have piecewise smooth boundaries $\partial V_1,...,
     \partial V_\ell$ and the remaining part of the $\om$-limit set is
     contained in $\cup_{1\leq j\leq\ell}\partial V_j$. We assume also
     that for any $z\in\cap_{1\leq i\leq k}\partial V_{j_i},\, k\leq\ell$
     there exist $\vp=\vp(z)>0$ and an $F^t_z$-invariant measures 
     $\nu_1,...,\nu_k$ such that
     \begin{equation}\label{1.2.27}
     z+s\bar B_{\nu_i}(z)\in V_{j_i}\,\,\mbox{for all}\,\, s\in(0,\vp]\,\,
     \mbox{and}\,\, i=1,...,k,
     \end{equation}
     i.e. $\bar B_{\nu_i}(z)\ne 0$ and it points out into the interior of
     $V_{j_i}$ which means that from any boundary point it is possible to
     go to any adjacent basin along a curve with an arbitrarily small
     $S$-functional. Let $\del>0$ be so small that the $\del$-neighborhood 
     $U_\del(\cO_i)=\{ z\in\cX:\,\mbox{dist}(z,\cO_i)<\del\}$
      of each $\cO_i$ is contained with its closure in the corresponding
      basin $V_i$. For any $x\in V_i$ set
     \[
     \tau^\ve_{x,y}(i)=\inf\big\{ t\geq 0:\, Z^\ve_{x,y}(t)\in
     \cup_{j\ne i}U_\del(\cO_j)\big\}.
     \]
     
      In Section \ref{sec1.8} we will derive the following result.
      
     \begin{theorem}\label{thm1.2.7} The function $R_{ij}(x)=\inf_{z\in V_j}
     R(x,z)$ takes on the same value $R_{ij}$ for all $x\in\cO_i,\, i\ne j$.
     Let $R_i=\min_{j\ne i,j\leq\ell}R_{ij}$. Then for any $x\in V_i$,
      \begin{equation}\label{1.2.28}
       \lim_{\ve\to 0}\ve\log\int_{\cW}\tau^\ve_{x,y}(i)dm(y)=R_i>0
     \end{equation}
      and for any $\al>0$ there exists $\la(\al)=\la(x,\al)>0$ such that for 
      all small $\ve>0$,
       \begin{equation}\label{1.2.29}
       m\big\{ y\in\cW:\,e^{(R_i-\al)/\ve}>
       \tau^\ve_{x,y}(i)\,\,
       \mbox{or}\,\,\tau^\ve_{x,y}(i)>e^{(R_i+\al)/\ve}
       \big\}\leq e^{-\la(\al)/\ve}.
      \end{equation}
       Next, set
     \[
     \Te^{\ve,i}_v(t)=\Te_v^{\ve,i,\del}(t)=\int_0^t
     \bbI_{V_i\setminus U_\del(\cO_i)}(Z^\ve_v(s))ds.
     \]
      Then for any $x\in V_i$ and $\del>0$ there exists $\la(\del)=
      \la(x,\del)>0$ such that for all small $\ve>0$,
     \begin{equation}\label{1.2.30}
     m\big\{ y\in\cW:\,\Te^{\ve,i}_{x,y}(\tau^\ve_{x,y}(i))\geq 
     e^{-\la(\del)/\ve}\tau^\ve_{x,y}(i)\big\}\leq 
     e^{-\la(\del)/\ve}.
     \end{equation}
     Now, suppose that the vector field $B$ is complete on 
     $\partial V_{i}$ for some $i\leq\ell$ (which strengthens (\ref{1.2.27})
     there) and the restriction of the $\om$-limit set of $\Pi^t$ to 
     $\partial V_{i}$ consists of a finite number of $S$-compacts. Assume
     also that there is a unique index $\io(i)\leq\ell,\, \io(i)\ne i$ such 
     that $R_{i}=R_{i\io(i)}$. Then for some $\la=\la(x)>0$ and all small 
     $\ve>0$,
      \begin{equation}\label{1.2.31}
      m\big\{ y\in\cW:\, Z^\ve_{x,y}(\tau^\ve_{x,y}(i))\not\in V_{\io(i)}
      \big\}\leq e^{-\la/\ve}.
      \end{equation}
      Finally, suppose that the above conditions hold true for all 
      $i=1,...,\ell$. Define $\io_0(i)=i$, $\tau^\ve_v(i,1)=\tau^\ve_v(i)$
      and recursively, 
      \[
      \io_k(i)=\io(\io_{k-1}(i))\,\,\mbox{and}\,\,\tau^\ve_v(i,k)=
      \tau^\ve_v(i,k-1)+\tau^\ve_{v_\ve(k-1)}
      \big(j(v_\ve(k-1))\big),
      \]
      where $v_\ve(k)=\Phi_\ve^{\ve^{-1}\tau^\ve_v(i,k)}v$,
      $j((x,y))=j$ if $x\in V_j$, and set $\Sigma^\ve_i(k,a)=\sum_{l=1}^k
      \exp\big((R_{\io_{l-1}(i),\io_l(i)}+a)/\ve\big)$. 
      Then for any $x\in V_i$ and $\al>0$ there 
      exists $\la(\al)=\la(x,\al)>0$ such that for all $n\in\bbN$ and 
      sufficiently small $\ve>0$,
       \begin{eqnarray}\label{1.2.32}
       &m\big\{ y\in\cW:\,\Sigma^\ve_i(k,-\al)>
       \tau^\ve_{x,y}(i,k)\,\,\mbox{or}\,\,\\
       &\tau^\ve_{x,y}(i,k)>\Sigma^\ve_i(k,\al)\,\,\mbox{for some}\,\,
       k\leq n\big\}\leq ne^{-\la(\al)/\ve}\nonumber
      \end{eqnarray}
      and for some $\la=\la(x)>0$,
       \begin{equation}\label{1.2.33}
      m\big\{ y\in\cW:\, Z^\ve_{x,y}(\tau^\ve_{x,y}(i,k))\not\in V_{\io_k(i)}
      \,\,\mbox{for some}\,\, k\leq n\big\}\leq ne^{-\la/\ve}.
      \end{equation}
      \end{theorem}
      
      Generically there exists only one index $\io(i)$ such that $R_i=
      R_{i\io(i)}$ and in this case Theorem \ref{thm1.2.7} asserts that
      $Z^\ve_{x,y},\, x\in V_i$ arrives (for "most" $y\in\cW$) at 
      $V_{\io(i)}$ after it leaves $V_i$. If $\cI(i)=\{ j:\, R_i=R_{ij}\}$
      contains more than one index then the method of the proof of Theorem
      \ref{thm1.2.7} enables us to conclude that in this case $Z^\ve_{x,y},\,
      x\in V_i$ arrives (for "most" $y\in\cW$) at $\cup_{j\in\cI(i)}V_j$
      after leaving $V_i$ but now we cannot specify the unique basin of
      attraction of one of $\cO_j$'s where $Z^\ve_{x,y}$ exits from $V_i$.
      If the succession function $\io$ is uniquely defined then it 
      determines an order of transitions of the slow motion $Z^\ve$
      between basins of attractors of $\bar Z$ and because of their finite
      number $Z^\ve$ passes them in certain cyclic order going around such
      cycle exponentially many in $1/\ve$ times while spending the total time
      in a basin $V_i$ which is approximately proportional to $e^{R_i/\ve}$.
      If there exist several cycles of indices $i_0,i_1,...,i_{k-1},i_k=i_0$ 
      where $i_j\leq\ell$ and $i_{j+1}=\io(i_j)$ then transitions between 
      different cycles may also be possible. In the uncoupled case with fast 
      motions being continuous time Markov chains a description of such 
      transitions via certain hierarchy of cycles appeared in
      \cite{Fre} and \cite{FW} without detailed proofs but relying on some
       heuristic arguments. In our fully coupled deterministic setup a 
      rigorous justification of the corresponding description seems to be 
      difficult in a more or less general situation though for some specific
      simple examples (as, for instance, those which are considered in Section
      \ref{sec1.9}) this looks feasible while it is not clear whether it is
      possible to describe in our situation a limiting as $t\to\infty$ 
      behaviour of the slow motion $Z^\ve(t)$ when $\ve$ is small but fixed.
      
      The proof of Theorems \ref{thm1.2.5}, \ref{thm1.2.7} and to certain extent
      also of Theorem \ref{thm1.2.3} rely, in particular, on certain "Markov
      property type" arguments which enable us to extend estimates on 
      relatively short time intervals to very long time intervals by,
      essentially, iterating them where the crucial role is played by a
      volume lemma type result of Section \ref{sec1.3} together with the
      technique of $(t,\del)$-separated sets and Bowen's $(t,\del)$-balls
      on unstable leaves of the perturbed flow $\Phi_\ve^t$. Moreover,
      the proof of (\ref{1.2.32}) and (\ref{1.2.33}) require certain rough
      strong Markov property type arguments which enable us to study the
      slow motion at subsequent hitting times $\tau^\ve_{x,y}(i,n)$ of
      small neighborhoods of attractors of the averaged motion.
      
      In order to produce a wide class of systems satisfying the conditions
      of Theorem \ref{thm1.2.7} we can choose, for instance, a vector field
      $\tilde B(x)$ on $\bbR^d$ whose $\om$-limit set satisfies the
      conditions stated above for the averaged system together with a
      family of vector fields $\hat B(x,y)$ on $\bbR^d$ (parametrized by
      $y\in\bfM$) such that $\int_{\bfM}\hat B(x,y)d\SRBx(y)\equiv 0$ and
      then set $B(x,y)=\tilde B(x)+\hat B(x,y)$. As a specific example
      we can take the flows $F_x^t,\, x\in\bbR^1_{-}=(-\infty,0)$ to be
      geodesic flows on the manifold $\bfM$ with (changing) constant
      negative curvature $x$, $\tilde B$ to be a one dimensional vector
      field on $\bbR^1$ and $\hat B(x,y)$ can be just a function $\hat B(y)$
      on $\bfM$ with zero integral with respect to the Lebesgue measure there.
      
      In Section \ref{sec1.9} we will derive similar results for the discrete
      time case where differential equations (\ref{1.1.1}) are replaced by
      difference equations (\ref{1.1.10}). Namely, recall that a compact subset
$\La$ of a compact Riemannian manifold $\bfM$ is called hyperbolic if it is
$F$-invariant and there exists $\ka>0$ and the splitting $T_\La \bfM=
\Gam^s\oplus\Gam^u$ into the continuous subbundles $\Gam^s,\Gam^u$ of the 
tangent bundle $T\bfM$ restricted to $\La,$ the splitting is invariant with 
respect to the differential $DF$ of $F,$ and there is $n_0>0$ such that for 
all $\xi\in\Gam^s,\,\eta\in\Gam^u,$ and $n\geq n_0$ the inequalities 
(\ref{1.2.2}) with $t$ replaced by $n$ hold true.
A hyperbolic set $\La$ is said to be basic hyperbolic if the periodic orbits 
of $F|_\La$ are dense in $\La,$ $F|_\La$ is topologically transitive,
and there exists an open set $U\supset\La$ with $\La=\cap_{-\infty<n<\infty}
F^nU.$ Such a $\La$ is called a basic hyperbolic attractor if for some open
 set $U$ and $n_0>0,$
\[
F^{n_0}\brU\subset U\quad\mbox{and}\quad\cap_{n>0}F^tU=\La
\]
where $\brU$ denotes the closure of $U.$ If $\La=\bfM$ then $F^t$ is called an 
Anosov flow. If $F$ is a $C^2$ endomorphism of $\bfM$ and there exists
$\ka>0$ such that $\|DF\xi\|\geq e^\ka\|\xi\|$ for all $\xi\in T\bfM$
then $F$ is called an expanding map (or transformation) of $\bfM$. It will
be convenient for our exposition to use the notation of the expanding
subbundle $\Gam^u$ also in the case of expanding maps where, of course,
$\Gam^u=T\bfM$. We replace now Assumption \ref{ass1.2.1} by the following one.

\begin{assumption}\label{ass1.2.8} The family $F_x=\Phi(x,\cdot)$ in (\ref{1.1.10}) 
consists of $C^2$-diffeomorphisms or endomorphisms of a compact 
$n_\bfM$-dimensional
Riemannian manifold $\bfM$ with uniform $C^2$ dependence on the parameter 
$x$ belonging to a neighborhood of the closure $\bar\cX$ of
 a relatively compact open connected set $\cX\subset\bbR^d$. All $F_x,\,
x\in\bar\cX$ are either expanding maps of $\bfM$ or diffeomorphisms  
possessesing basic hyperbolic attractors $\La_x$ with hyperbolic splittings  
satisfying (\ref{1.2.2}) with the same $\ka>0$ and there exists an
open set $\cW\subset \bfM$ and $n_0>0$ satisfying (\ref{1.2.4}) with $n$ in
place of $t$.
\end{assumption}

Let $J_x^u(y)$ be the absolute value of the
Jacobian of the linear map $DF_x(y):\Gam_x^u(y)\to
\Gam^u_x(F_xy)$ with respect to the Riemannian inner products and set
\begin{equation}\label{1.2.34}
 \vf_x^u(y)=-\log J_x^u(y).
 \end{equation}
 The function $\vf_x^u(y)$ is known to be H\" older continuous in $y,$ since
 the subbundles $\Gam^u_x$ are H\" older continuous (see \cite{KH}),
  and $\vf_x^u(y)$ is $C^1$ in $x$ (see \cite{Co}). The topological pressure
  $P_x(\psi)$ of a function $\psi$ for $F$ is defined similarly to the 
 continuous time (flow) case above but now time should run only over integers 
 and the integral $\int_0^t\psi(F^s_xy)ds$ should be replaced by the sum 
  $\sum_{k=0}^{n-1}\psi(F^ky)$ (see \cite{KH}). Again the variational 
  principle (\ref{1.2.6}) holds true and if $q$ is a H\" older 
  continuous function on $\La_x$ there exists a unique $F_x-$invariant
  measure $\mu^q_x$ on $\La_x,$ called the equilibrium state for $\vf_x^u+q$
  which satisfies (\ref{1.2.7}). In particular, $\mu^0_x=\SRBx$ is usually 
  called the Sinai--Ruelle--Bowen (SRB) measure. Since $\La_x$ are attractors 
  we have that $P_x(\vf_x^u)=0$ (see \cite{BR}) and the same holds true in the
  expanding case, as well. Next, we define $I_x(\nu)$,
  $H(x,x',\be)$, $H(x,\be)$, $L(x,x',\al)$, $L(x,\al)$, $S_{0T}(\gam)$,
  and $\gam^u$ as in (\ref{1.2.8}) and (\ref{1.2.10})--(\ref{1.2.14}). 
  In place of Assumption \ref{ass1.2.2} we will rely now on the similar one 
  concerning the equation (\ref{1.1.10}).
  \begin{assumption}\label{ass1.2.9}
   There exists $K>0$ such that
   \begin{equation}\label{1.2.35}
   \| B(x,y)\|_{C^1(\cX\times\bfM)}+\| \Phi(x,y)\|_{C^2(\cX\times\bfM)}
   \leq K
   \end{equation}
   where the first $\|\cdot\|_{C^1(\cX\times\bfM)}$ is the $C^1$ norm of the 
   corresponding vector fields on $\cX\times\bfM$ and the second expression
   is the $C^2$ norm (with respect to the corresponding Riemannian metrics) 
   of the map $\Phi:\cX\times\bfM\to\cX\times\bfM$ acting by $\Phi(x,y)=
   (x,F_xy)$.
   \end{assumption}
  \begin{theorem}\label{thm1.2.10} Assume that Assumptions \ref{ass1.2.8} and
  \ref{ass1.2.9} are satisfied and that $X^\ve(n)=X_{x,y}^\ve(n),\, n=0,1,2,...$ 
  is obtained by (\ref{1.1.10}). For $t\in[n,n+1]$ define $X^\ve(t)=(t-n)X^{\ve}
  (n+1)+(n+1-t)X^{\ve}(n)$ and set $Z^\ve_{x,y}(t)=X^\ve_{x,y}(t/\ve)$. Then
  Theorem \ref{thm1.2.3} and Corollary \ref{cor1.2.4} hold true with the 
  corresponding functionals $S_{0t}$. Theorems \ref{thm1.2.5} and \ref{thm1.2.7}
  hold true, as well, under the corresponding assumptions about the family
  $\{ F_x,\, x\in\cX\}$ (with $\cX=\bbR^d$ in the case of Theorem \ref{thm1.2.7})
  and about the averaged system (\ref{1.1.11}) (in particular, about its 
  attractors) in place of the system (\ref{1.1.6}).
  \end{theorem}
  In Section \ref{sec1.9} we exhibit computations which demonstrate the
  phenomenon of Theorem \ref{thm1.2.7} in the discrete time case for two
  simple examples where $F_xy$ are one dimensional maps $y\to 3y+x$ (mod 1)
  and the averaged equation has three attracting fixed points. 
  
  In the last Section \ref{sec1.10} we discuss a stochastic resonance type
  phenomenon which can be exhibited in three scale systems where fast 
  motions are hyperbolic flows (hyperbolic diffeomorphisms, expanding
  transformations) as above depending on the intermediate and slow motions
  while the intermediate motion performs rare transitions between attracting
  fixed points of corresponding averaged systems which under certain
  conditions creates a nearly periodic motion of the slow one dimensional
  motion.
  
  \begin{remark}\label{rem1.2.11} Computation or even estimates of functionals
  $S_{0T}(\gam)$ seem to be quite difficult already for simple discrete 
  (and, of course, more for continuous) time examples since this leads to
  complicated nonclassical variational problems. This is crucial in order
  to estimate numbers $R_{ij}$ which according to Theorem \ref{thm1.2.7}
  are responsible for transitions of the slow motion between basins of
  attractors of the averaged system.
  \end{remark}
   \begin{remark}\label{rem1.2.12} The estimate (\ref{1.2.18}) shows that
  $Z^\ve_{x,y}$ tends as $\ve\to 0$ to $\bar Z_x$  uniformly on $[0,T]$ 
  in the sense of convergence in measure $m$ considered on the space of
  initial conditions $y\in\bfM$. A natural question to ask is
  whether the convergence for almost all (fixed) initial conditions also takes 
  place in our circumstances. In \cite{BK1} we give a negative answer to 
  this question, in
  paricular, for the following simple discrete time example
  \begin{eqnarray*}
  &\big(X^\ve_{x,y}(n+1),Y^\ve_{x,y}(n+1)\big)\\
  &=\big(X^\ve_{x,y}(n)+\ve\sin\big(2\pi Y^\ve_{x,y}(n)\big),\,
   2Y^\ve_{x,y}(n)+X^\ve_{x,y}(n)\,\mbox{(mod 1)}\big).
  \end{eqnarray*}
  Identifying 0 and 1 we view $y$ variable as belonging to the circle in
  order to fit into our setup where the fast motion runs on a compact manifold.
  The averaged equation (\ref{1.1.11})
  has here zero in the right hand side so the averaged
  motion stays forever at the initial point. The discrete time version of
  (\ref{1.2.18}) asserted by Theorem \ref{thm1.2.10} implies that
  \begin{equation}\label{1.2.36} 
  \max_{0\leq n\leq 1/\ve}|X^\ve_{x,y}(n)-x|\to 0\,\,\mbox{as}\,\,\ve\to 0
  \end{equation}
  in the sense of convergence in (the Lebesgue) measure on the circle but 
  we show in \cite{BK1} that for each $x$ there is a set $\Gam_x$ of full
  Lebesgue measure on the circle such that if $y\in\Gam_x$ then $\limsup$
  as $\ve\to 0$ of the left hand side in (\ref{1.2.36}) is positive, i.e.
  there is no convergence for Lebesgue almost all $y$ there. Namely, it turns
  out that for almost all initial conditions there exists a sequence 
  $\ve_i\to 0$ such that the fast motion $Y^{\ve_i}_{x,y}(n)$ stays for a time
  of order $1/\ve_i$ close to an orbit $\{ 2^nv$ (mod $2\pi)\},\, n\geq 0$ of
  the doubling map with $v$ being a generic point with respect to a Gibbs
  invariant measure $\mu$ of this map satisfying $\int_0^1\sin(2\pi v)d\mu(v)
  \ne 0$ which prevents (\ref{1.2.36}).
  \end{remark}

  \section{Dynamics of $\Phi_\ve^t$}\label{sec1.3}\setcounter{equation}{0}

 For readers convenience we exhibit, first, in this section the setup and 
 necessary technical results from \cite{Ki7} and though their proofs can
 can be found in \cite{Ki7} we provide for completness and readers'
 convenience their slightly modified and corrected version also here. 

Any vector $\xi\in T(\bbR^d\times\bfM)=\bbR^d\oplus T\bfM$ can be uniquely
written as $\xi=\xi^\cX+\xi^\cW$ where $\xi^\cX\in T\bbR^d$ and $\xi^\cW\in
T\bfM$ and it has the Riemannian norm $|\|\xi|\|=|\xi^\cX|+\|\xi^\cW\|$ 
where $|\cdot|$ is the usual Euclidean norm on $\bbR^d$ and $\|\cdot\|$ is
 the Riemannian norm on $T\bfM$. The corresponding metrics on $\bfM$ and
 on $\bbR^d\times\bfM$ will be denoted by $d_\bfM$ and $dist$, respectively, so
 that if $z_1=(x_1,w_1),\, z_2=(x_2,w_2)\in\bbR^d\times\bfM$ then
 $dist(z_1,z_2)=|x_1-x_2|+d_\bfM(w_1,w_2).$
It is known (see \cite{Ro}) that the hyperbolic splitting 
$T_{\La_x} \bfM=\Gam_x^s\oplus\Gam_x^0
\oplus\Gam_x^u$ over $\La_x$ can be continuously extended to the splitting
$T_\cW\bfM=\Gam_x^s\oplus\Gam_x^0\oplus\Gam_x^u$ over $\cW$ which is forward
invariant with respect to $DF^s_x$
and satisfies exponential estimates with a uniform in $x\in\cX$ positive
exponent which we denote again by $\ka>0,$ i.e. we assume now that
\begin{equation}\label{1.3.1}
\| DF^t_x\xi\|\leq e^{-\ka t}\|\xi\|\quad\mbox{and}\quad\|DF_x^{-t}\eta\|\leq 
e^{-\ka t}\|\eta\|
\end{equation}
provided $\xi\in\Gam^s_x(w)$, $\eta\in\Gam^u_x(F^t_xw)$, $t\geq t_0$, and 
$w\in \cW.$ Moreover, by \cite{Co} (see also \cite{Ru}) we can choose these 
extensions so that $\Gam_x^s(w)$ and $\Gam_x^u(w)$ will be H\" older
continuous in $w$ and $C^1$ in $x$ in the 
corresponding Grassmann bundle. Actually, since $\cW$ is contained in the
basin of each attractor $\La_x$, any point $w\in\cW$ belongs to the stable 
manifold $W_x^s(v)$ of some point $v\in\La_x$ (see \cite{BR}), and so we 
 choose naturally
$\Gam_x^s(w)$ to be the tangent space to $W^s_x(v)$ at $w.$
Now each vector
$\xi\in T_{x,w}(\cX\times \cW)=T_x\cX\oplus T_w\cW$ can be represented uniquely
in the form $\xi=\xi^\cX +\xi^s+\xi^0+\xi^u$ with $\xi^\cX\in T_x\cX$,
$\xi^s\in\Gam^s_x(w)$, $\xi^0\in\Gam^0_x(w)$ and $\xi^u\in\Gam^u_x(w)$.
We denote also $\xi^{0s}=\xi^s+\xi^0$ and $\xi^{0u}=\xi^u+\xi^0$.
For each small $\ve,\al>0$ set $\cC^u(\ve,\al)=\{\xi\in T(\cX\times \cW):\,
\|\xi^{0s}\|\leq\ve\al^{-2}\|\xi^u\|\,\mbox{and}\, \|\xi^\cX\|\leq\ve\al^{-1}
\|\xi^u\|\}$ and $\cC^u_{x,w}(\ve,\al)=\cC^u(\ve,\al)\cap T_{x,w}(\cX\times \cW)$
which are unstable cones \index{unstable cones}
around $\Gam^u$ and $\Gam^u_x(w),$ respectively. 
Similarly, we define $\cC^s(\ve,\al)=\{\xi\in T(\cX\times \cW):\,
\|\xi^{0u}\|\leq\ve\al^{-2}\|\xi^s\|\,\mbox{and}\, \|\xi^\cX\|\leq\ve\al^{-1}
\|\xi^s\|\}$ and $\cC^s_{x,w}(\ve,\al)=\cC^s(\ve,\al)\cap T_{x,w}(\cX\times \cW)$
which are stable cones \index{stable cones}
around $\Gam^s$ and $\Gam^s_x(w),$ respectively. Put
$(\cX\times \cW)_t=\{ (x,w):\,\Phi_\ve^u(x,w)\in(\cX\times \cW)\,\,\forall 
u\in[0,t]\}$, where, recall, $\Phi_\ve^t$ is the flow determined by the
equations (\ref{1.1.1}).
 
 \begin{lemma}\label{lem1.3.1}.
  There exist $\al_0,\,\ve(\al),\, t_1>0$ such that
 if $z=(x,y)\in(\cX\times \cW)_t$ and
 $t\geq t_1,\,\al\leq\al_0,\,\ve\leq\ve(\al)$  then
 \begin{equation}\label{1.3.2}
 D_z\Phi_\ve^t\cC^u_z(\ve,\al)\subset\cC^u_{\Phi_\ve^tz}(\ve,\al),\quad
 \cC^s_z(\ve,\al)\supset D_z\Phi_\ve^{-t}\cC^s_{\Phi_\ve^tz}(\ve,\al),
 \end{equation}
 and for any $\xi\in\cC^u_z(\ve,\al),$ $\eta\in\cC^s_{\Phi_\ve^tz}(\ve,\al),$
 \begin{equation}\label{1.3.3}
  |\| D_z\Phi^t_\ve\xi|\|\geq e^{\frac 12\ka t}|\|\xi|\|,\quad
  |\| D_z\Phi^{-t}_\ve\eta|\|\geq e^{\frac 12\ka t}|\|\eta|\|.
  \end{equation}
  \end{lemma}
  \begin{proof} Let $\xi=\xi^\cX+\xi^u+\xi^{0s}\in T_z(\cX\times\bfM),$
  $D_z\Phi^t\xi^\cX=\zeta=\zeta^\cX+\zeta^u+\zeta^{0s}\in T_{\Phi^tz}(\cX
  \times\bfM),$ $z=(x,y),\, D_yF^t_x\xi^u=\eta^u,$ and $D_yF^t_x\xi^{0s}
  =\eta^{0s}.$ Then $D_z\Phi^t\xi=\zeta^\cX+(\zeta^u+\eta^u)+(\zeta^{0s}+
  \eta^{0s})$ and $\|\xi^\cX\|=\|\zeta^\cX\|,$ $\|\zeta^u\|\leq Ce^{Ct}
  \|\xi^\cX\|,$ $\|\zeta^{0s}\|\leq Ce^{Ct}\|\xi^\cX\|,$ $\|\eta^{0s}\|\leq
  C\|\xi^{0s}\|$ for some $C\geq 1$ independent of $\xi$ and $\|\eta^u\|\geq
  e^{\ka t}\|\xi^u\|$ if $t\geq t_0.$ Hence, for $t\geq t_0,$
  \begin{equation*}
  \|\zeta^u+\eta^u\|\geq\|\eta^u\|-\|\zeta^u\|\geq e^{\ka t}\|\xi^u\|-Ce^{Ct}
  \|\xi^\cX\|
  \end{equation*}
  and
  \begin{equation*}
  \|\zeta^{0s}+\eta^{0s}\|\leq\|\zeta^{0s}\|+\|\eta^{0s}\|\leq Ce^{Ct}
  \|\xi^\cX\|+C\|\xi^{0s}\|.
  \end{equation*}
  If $\xi\in\cC^u_z(\ve,\al)$ then
   $\|\xi^u\|\geq\al\ve^{-1}\|\xi^\cX\|$ and $\|\xi^u\|\geq\al^2
  \ve^{-1}\|\xi^{0s}\|.$ Hence, by the above,
  \begin{equation*}
  \|\zeta^u+\eta^u\|\geq\al\ve^{-1}(\frac 12e^{\ka t}-\ve\al^{-1}Ce^{Ct})
  \|\xi^\cX\|+\frac 12e^{\ka t}\al^2\ve^{-1}\|\xi^{0s}\|.
  \end{equation*}
  Set $t_1=\ka^{-1}\ln 6,$ choose $\al\leq 6^{-2C/\ka}$ and $\ve=\ve(\al)\leq
  \al^2/4C$. Then we obtain that $D_z\Phi^t\xi\in\cC^u_{\Phi^tz}(\ve,2\al)$ 
   for all $t\in[t_1,t_2]$, and so by continuity of the splitting
  $\Gam_x^s\oplus\Gam_x^0\oplus\Gam_x^u$ and by perturbation arguments it
  follows that $D_z\Phi^t_\ve\xi\in\cC^u_{\Phi^t_\ve z}(\ve,\al)$ for all $t\in
  [t_1,2t_1]$ provided $\ve$ is small enough. Repeating this argument
  for $t\in[it_1,(i+1)t_1]$, $i=2,3,..$ we conclude the proof of the first
  part of (\ref{1.3.2}) and its second part follows in the same way.
  
   Next, for $\xi\in\cC^u_z(\ve,\al)$ and $t\geq t_0,$
   \begin{eqnarray*}
  \quad\,\,\,\,\,|\| D_z\Phi^t\xi|\|\geq\|\eta^u\|-\|\zeta^\cX\|-\|\zeta^u\|- 
  \|\zeta^{0s}\|-\|\eta^{0s}\|\geq e^{\ka t}\|\xi^u\|&\\
  -(1+2Ce^{Ct})\|\xi^\cX\|-C\|\xi^{0s}\|\geq(e^{\ka t}-\al^{-1}\ve(1+2Ce^{Ct})
  -\al^{-2}\ve C)\|\xi^u\|&\nonumber\\
  \geq(e^{\ka t}-\al^2(1+2Ce^{Ct})-\al^{-2}\ve C)
  (1+\ve\al^{-1}+\ve\al^{-2})^{-1}|\|\xi|\|.&
   \nonumber\end{eqnarray*}
 Choose  $\al_0,\,\ve(\al)$ so small (for instance, $\ve(\al)=\al^3$) that for
  all $\al\leq\al_0$ and $\ve\leq\ve(\al)$,
 $$
 e^{\ka t}-\ve\al^{-1}(1+2Ce^{Ct})-\ve\al^{-2}C\geq (1+\ve\al^{-1}+\ve\al^{-2})
 e^{\frac 23\ka t}
 $$
 for all $t\in[t_1,2t_1].$ Then, $|\| D_z\Phi^t\xi|\|\geq e^{\frac 23\ka t}
 |\|\xi|\|$ for all such $t$, and so if $\ve$ small enough we have also
 $\| D_z\Phi^t_\ve\xi\|\geq e^{\frac 12\ka t}\|\xi\|$. Using (\ref{1.3.2}) and 
 repeating this argument for $D_z\Phi_\ve^{it_1}\xi,\, i=1,2,...$ in place of
  $\xi$ we derive (ii) for all $t\geq t_1$. The proof for stable cones
  $\cC^s_\ve(\ve,\al)$ is similar.
   \end{proof}
  
  For any linear subspace $\Xi$ of $T_z(\bbR^d\times\bfM)$ denote by 
  $J^\Xi_\ve(t,z)$ absolute value of the Jacobian of the linear map
   $D_z\Phi_\ve^t:\,\Xi\to D_z\Phi_\ve^t\Xi$ with respect to inner products
   induced by the Riemannian metric. For each $z=(x,y)\in\bbR^d\times\bfM$ 
   set also
   \begin{equation}\label{1.3.4}
   J^u_\ve(t,z)=\exp\big(-\int_0^t\vf^u_{X^\ve_{x,y}(s)}(Y^\ve_{x,y}(s))ds\big).
   \end{equation}
   Let $n_u$ be the dimension of $\Gam^u_x(w)$ which does not depend on $x$ and 
$w$ by continuity considerations. If $\Xi$ is an $n^u-$dimensional subspace
of $T_z(\cX\times \cW)$, $z=(x,y),$ and $\Xi\subset\cC^u_{x,y}(\ve,\al)$ then
it follows easily from Assumption \ref{ass1.2.2} and Lemma \ref{lem1.3.1} that
there exists a constant $C_1>0$ independent of $z\in\cX\times \cW$ and of a 
small $\ve$ such that for any $t\geq 0,$
\begin{equation}\label{1.3.5}
(1-C_1\ve)^t\leq J^\Xi_\ve(t,z)(J^u_\ve(t,z))^{-1}\leq(1+C_1\ve)^t.
\end{equation}

 Recall, that an embedded $C^k,\, k=1,2$ $l-$dimensional disc $D$ in
 $\bbR^d\times\bfM$, $l\leq d+n_\bfM$ is the image of an $l$-dimensional disc
 (ball) $K$ in $\bbR^{d+n_\bfM}$ centered at 0 under a $C^k$ diffeomorphism of
 a neighborhood of 0 in $\bbR^{d+n_\bfM}$ into $\bbR^d\times\bfM$ and we define
 the boundary $\partial D$ of $D$ as the image of the boundary $\partial K$ 
 of $K$ considered in the corresponding $l$-dimensional Euclidean subspace
 of $\bbR^{d+n_\bfM}$. Denote by $U(z,\rho)$ the ball in $\bbR^d
\times\bfM$ centered at $z$ and let $\cD_\ve^u(z,\al,\rho,C)$, $C\geq 1$ be the
set of all $C^1$ embedded $n_u-$dimensional closed discs $D\subset\cX\times \cW$ 
such that $z\in D$, $TD\subset\cC^u(\ve,\al)$ and if $v\in\partial D$ then 
$C\rho\leq d_D(v,z)\leq C^2\rho$ where $TD$ is the tangent bundle over $D$ and 
$d_D$ is the interior metric on $D.$ Each disc $D\in\cD^u_\ve(z,\al,\rho,C)$ 
will be called unstable or expanding \index{unstable disc} \index{expanding
disc} and, clearly, 
$D\subset U(z,C^2\rho)$ and if $\ve/\al^2$ and $\rho$ are small enough and
$C>1$ then 
dist$(v,z)\geq\rho$ for any $v\in\partial D.$
  Let $D\in\cD^u_\ve((x,y),\al,\rho,C)$ and $z\in D\subset\cX\times \cW$. Set
   $U^\ve_{D}(t,z,L)=\{\tilde z\in D:\, d_{\Phi_\ve^sD}
   (\Phi_\ve^sz,\Phi_\ve^s\tilde z) \leq L\quad\forall s\in[0,t]\}$ and
   let $\pi_1:\,\cX\times\cW\to\cX$ and $\pi_2:\,\cX\times\cW\to\cW$ be
   natural projections on the first and second factors, respectively. 
  
  \begin{lemma}\label{lem1.3.2} 
   Let $\ve,\,\al,\, t$, $(x,y)$ be as in Lemma \ref{lem1.3.1} and
   $T>0$. There exist $\rho_0,c,c_{\rho,T},C\geq 1$ such that if
    $\rho\leq\rho_0,\, D\in\cD_\ve^u((x,y),\al,\rho,C)$, $z\in D$,
     $V_{s,t}(z)=\Phi_\ve^sU^\ve_{D}(t,z,C\rho)$, $V_t(z)=V_{t,t}(z)$,
      $V=V_{0,t}(z)\subset D$ and $t\geq 0$ then 
   
   (i) $d_{V_{s,t}(z)}(\Phi_\ve^sv,\Phi_\ve^sz)\leq c^{-1}e^{-\frac 12\ka (t-s)}
   d_{V_{t}(z)}(\Phi_\ve^tv,\Phi_\ve^tz)\leq c^{-1}C\rho e^{-\frac 12\ka 
   (t-s)}$ for any $v\in V$ and $s\in [0,t],$ where $d_U$ is the interior
    distance on $U$;
   
   (ii) $TV_{s,t}(z)\subset\cC^u(\ve,\al)$ and
   $V_t(z)\in\cD_\ve^u(\Phi_\ve^tz,\al,\rho,\sqrt C)$ provided
   $\partial D\cap U^\ve_{D}(t,z,C\rho)=\emptyset$;
   
   (iii) For all $v\in V$ and $0\leq s\leq t$,
   \[
   |\pi_1\Phi^s_\ve v-\pi_1\Phi_\ve^sz|\leq Cc^{-1}\rho\ve\al^{-1}
   (1-\ve\al^{-1}-\ve\al^{-2})^{-1}e^{-\frac 12\ka(t-s)}.
   \]
    (iv) $c_{\rho,T}\leq m_D(V)J^u_\ve(t,z)\leq c_{\rho,T}^{-1}$ provided 
   $t\leq T/\ve$, where $m_D$ is the induced (not normalized) Riemannian 
   volume on $D.$
    \end{lemma}
    \begin{proof}
   (i)  Let $\gam$ be a smooth curve on 
   $V_t(z)$ connecting $a=\Phi_\ve^tz$ and $b=\Phi_\ve^tv.$ then
   $\tilde\gam=\Phi_\ve^{s-t}\gam$ is a smooth curve on $V_{s,t}(z)$ connecting
   $\Phi_\ve^sz$ and $\Phi_\ve^sv.$ Since $T\tilde\gam\subset TV_s(z)\subset
   \cC^u(\ve,\al)$ then by (\ref{1.3.3}), length$(\gam)\geq 
   e^{\ka(t-s)/2}$length$(\tilde\gam)$ if $t-s\geq t_1.$ Then for such $t$
   and $s,$
   \begin{equation*}
   d_{V_s}(\Phi^s_\ve v,\Phi_\ve^s z)\leq\,\mbox{length}(\tilde\gam)\leq
   e^{-\ka(t-s)/2}\mbox{length}(\gam).
   \end{equation*}
   Observe that (i) is nontrivial only for large $t-s$, so minimizing in
   $\gam$ in the above inequality  we derive the assertion (i).
   
   Next, we derive (ii). Its first part follows from (\ref{1.3.2}).
   By the definition of $V$, $d_{V_t(z)}(w,z)\leq C\rho$ for any
    $w\in\partial V_t(z)$. It remains to show that 
    $d_{V_t(z)}(w,z)\geq\sqrt C\rho$ for any $w\in\partial V_t(z).$ Indeed, 
    suppose $d_{V_t(z)}(w,z)<\sqrt C\rho.$  Set
    $$
    d_1=\sup_{-t_1\leq s\leq 0,v\in\cX\times \cW}\|D_v\Phi_\ve^s\|.
    $$
    Next, we conclude via perturbation arguments that
    $d_1>0$ provided $\ve$ is small enough.
    Let $w=\Phi_\ve^tv.$ It follows from Lemma \ref{lem1.3.1} that
    $d_{V_s(z)}(\Phi_\ve^sv,\Phi_\ve^sz)<d_1\sqrt C\rho$ for all
     $s\in[0,t].$ Hence, if $\sqrt C\geq d_1$ then $v\not\in\partial V,$
     and so $w\not\in\partial V_t(z).$
     
     In order to derive (iii) take an arbitrary smooth curve $\gam$ on 
     $V_{s,t}(z)$ connecting $\Phi_\ve^sv$ and $\Phi_\ve^sz$. Then 
    $\frac {d\gam(s)}{ds}=\dot\gam(s)\in\cC^u_{\gam(s)}(\ve,\al).$ It follows
    that if $\gam(s)=(\gam^\cX(s),\gam^\bfM(s))$ with $\gam^\cX(s)\in\cX$ and
    $\gam^\bfM(s)\in\bfM$ then $\dot\gam^\bfM(s)=\dot\gam^{0s}(s)+
    \dot\gam^u(s)$, $\|\dot\gam^\cX(s)\|\leq\ve\al^{-1}\|\dot\gam^u(s)\|,$
    $\|\dot\gam^{0s}(s)\|\leq\ve\al^{-2}\|\dot\gam^u(s)\|,$ and so
    \[
    |\|\dot\gam(s)|\|\geq\|\dot\gam^u(s)\|-\|\dot\gam^{0s}(s)\|-
    \|\dot\gam^\cX(s)\|\geq\|\dot\gam^u(s)\| (1-\ve\al^{-1}-\ve\al^{-2}).
    \]
    Hence
    \[
    \|\dot\gam^\cX(s)\|\leq\ve\al^{-1}(1-\ve\al^{-1}-\ve\al^{-2})^{-1}
    |\|\dot\gam(s)|\|,
    \]
    and so
    \begin{equation*}
    |\pi_1\Phi_\ve^sv-\pi_1\Phi_\ve^sz|\leq\ve\al^{-1}(1-\ve\al^{-1}-
    \ve\al^{-2})^{-1}\mbox{length}(\gam).
    \end{equation*}
    Minimizing the right hand side here over such $\gam$ we obtain (iii) 
    using (i).
     Finally, (iv) follows from (\ref{1.3.5}), (i), (ii), and the H\" older
     continuity of $\vf^u$ (as a function on $\cX\times W).$
      \end{proof}

    For each $y\in\La_x$ and $\vrho>0$ small enough set $W_x^s(y,\vrho)=
    \{\tilde y\in\cW:\, d_\bfM(F_x^ty, F_x^t\tilde y)\leq\vrho\,\,\,\forall
    t\geq 0\}$ and  $W_x^u(y,\vrho)=\{\tilde y\in\cW:\, d_\bfM(F_x^ty,
    F_x^t\tilde y)\leq\vrho\,\,\,\forall t\leq 0\}$ which are local stable
    and unstable manifolds for $F_x$ at $y.$ According to \cite{HPPS}
    and \cite{Ro} these families can be included into continuous families
    of $n^s$ and $n^u-$dimensional stable and unstable $C^1$ discs 
    $W_x^s(y,\vrho)$ and $W^u_x(y,\vrho),$ respectively, defined for all
    $y\in\cW$ and such that $W_x^s(y,\vrho)$ is tangent to $\Gam^s_x$,
    $W_x^u(y,\vrho)$ is tangent to $\Gam^u_x$,
    $F^t_xW^s_x(y,\vrho)\subset W^s_x(F^t_xy,\vrho)$, and $W^u_x(y,\vrho)
    \supset F_x^{-t}W^u_x(F^t_xy,\vrho).$ Actually, as we noted it above if
    $y\in\cW$ then $y$ belongs to a stable manifold $W^s_x(\tilde y)$ of
    some $\tilde y\in\La_x$ and we choose $W^s_x(y,\vrho)$ to be the subset
    of $W^s_x(\tilde y).$  
   
   \begin{lemma}\label{lem1.3.3} For any $0<\rho_1<\rho_0$ small enough and a 
   continuous function $g$ on $\cX\times \cW$ uniformly in 
   $D\in\cD^u_\ve(z,\al,\rho,C)$, $x'\in\cX,\, z\in\cX\times\cW$
    and $\rho\in[\rho_1,\rho_0]$,
   \begin{equation}\label{1.3.6}
   \lim_{t\to\infty}\frac 1t\log\int_D\exp\left(\int_0^tg(x',F^r_{\pi_1z}
   \pi_2v)dr\right)dm_D(v)=P_{\pi_1z}(g(x',\cdot)+\vf_{\pi_1z}^u)
   \end{equation}
   where $m_D$ is the induced Riemannian volume on $D.$
   \end{lemma}
   \begin{proof} Set $\tilde W=\pi_2D,$ $x=\pi_1z$ and $y\in\pi_2z.$ By
   standard transversality considerations we can define a one-to-one map
   $\tilde\pi:\tilde W\rightarrow\tilde\pi\tilde W\subset W^u_x(y,\tilde C\rho)$
   by $\tilde\pi(\tilde w)=w\in F^\tau_xW_x^s(\tilde w,\tilde C\rho)$
   provided $\al,\rho,|\tau|>0$ are small and $\tilde C>0$ is sufficiently
   large. By the absolute continuity of the stable foliation arguments
   (see, for instance, \cite{Ma}, Section 3.3) we conclude that $\tilde\pi$
   and its inverse have bounded Jacobians. It follows that it suffices to
   establish (\ref{1.3.6}) for $D=\{ x\}\times W$ where $W=W^u_x(y,\gam)$
   uniformly in $\gam\in[\gam_0,\gam_0^{-1}],\,\gam_0>0.$
   
   Set $W_\tau=F_x^\tau W,$ $W_{r,q}=\cup_{r\leq \tau\leq q}F_x^\tau W$ and 
   \[
   I^V_{x,x'}(t)=\int_V\exp\big(\int_0^tg(x',F_x^\tau v)d\tau\big)dm_V(v).
   \]
   Then
   \begin{eqnarray*}
   &I_{x,x'}^{W_\tau}(t)=\int_{W_\tau}\exp\big(\int_0^tg(x',F_x^\te v)d\te
   \big)dm_{W_\tau}(v)\\
   &=\int_W\exp\big(\int_\tau^{t+\tau}g(x',F^\te_xv)d\te\big)J^u_x(\tau,v)
   dm_W(v),
   \end{eqnarray*}
   and so 
   \begin{equation*}
   re^{-2(\tau+r)\| g\|}I^W_{x,x'}(t)\leq I_{x,x'}^{W_{\tau,\tau+r}}(t)\leq
   re^{(\tau+r)(2\| g\|+\|\vf^u\|)}I^W_{x,x'}(t)
   \end{equation*}
   where $\|\cdot\|$ is the supremum norm on $\cX\times\cW.$
   Since $\cW\subset\cup_{v\in\La_x}W^s_x(v)$ by \cite{BR} then given
   $\eta>0$ there exist $\gam(\eta),\tau(\eta)>0,$ $v\in\La_x,$ and
   $U\subset\cup_{-r\leq\te\leq r}F^\te_xW^u_x(v)$ such that for any
   $\gam\leq\gam(\eta)$ we can define a one-to-one map $\phi:\, U\to
   W_{\tau(\eta),\tau(\eta)+r}$ by $\phi(w)\in W^s_x(w,\eta).$ By standard
   absolute continuity of the stable foliation considerations (see
    \cite{Ma}, Section 3.3) it follows that $\phi$ and its
   inverse have bounded Jacobians which together with the above arguments
   yield that it suffices to prove Lemma \ref{lem1.3.3} only when $y\in
   \La_x,$ and so (see \cite{BR}), $W=W^u_x(y,\gam)\subset\La_x.$ We
   observe also that without loss of generality we can assume $\gam$ 
   to be sufficiently small since we can always cover $W^u_x(y,\gam)$
   by $W^u_x(y_i,\tilde\gam),\, i=1,...,k$ with $y_1=y,$ $k=k(\gam,\tilde
   \gam)$ and small $\tilde\gam\leq\gam,$ so proving Lemma \ref{lem1.3.3}
   for all such $W^u_x(y_i,\tilde\gam)$ will imply it for $W^u_x(y,\gam)$ 
   itself.
   
   So now assume that $W=W^u_x(y,\gam)\subset\La_x$ and we claim that
   for any $\eta>0$ there exists $\tau(\eta,\gam)>0$ such that $W_{0,\tau}$
   forms an $\eta-$net in $\La_x$ for any $\tau\geq\tau(\eta,\gam).$
   Indeed, by topological transitivity there exists $v\in\La_x$ whose
   orbit is dense in $\La_x,$ and so by standard ergodicity considerations
   with respect to any ergodic invariant measure with full support on
   $\La_x$ (take, for instance, the SRB measure) we conclude that 
    for any $\tau$ already
   $\{ F_x^rv,\, r\geq\tau\}$ is dense in $\La_x$. Then by transversality
   of $W^s_x$ and $\cup_\te F_x^\te W^u_x$
   there exists $r>0$ such that $F^r_xv\in W_x^s(w,\gam)$ for some
   $w\in W,$ and so the forward orbit of $w$ is dense in $\La_x$, whence
   our claim holds true. By compactness and structural stability
   considerations it follows that we can choose the same $\tau(\eta,\gam)$
   for all $y\in\cW$ and $x\in\cX.$
   
   For any set $V\subset\La_x$ put $U_V(t,y,\zeta)=\{ v\in V:\,
   d(F^r_xy,F^r_xv)\leq\zeta\,\,\,\,\forall r\in[0,t]\}.$ Recall, that a
   finite set $E\subset\La_x$ is called $(\zeta,t)-$separated for the
   flow $F^t_x$ if $y,\tilde y\in E,\, y\ne\tilde y$ implies $\tilde y
   \not\in U_{\La_x}(t,y,\zeta).$ A set $E\subset\La_x$ is called  
   $(\zeta,t)-$spanning if for any $y\in\La_x$ there is $\tilde y\in E$
    such that $y\in U_{\La_x}(t,\tilde y,\zeta)$. Let $W_{0,\tau}$ be an
    $\eta-$net in $\La_x$ and $E$ be a maximal $(\zeta,t)-$separated subset
    of $W_{0,\tau}$. Then $U_{W_{0,\tau}}(t,y,\zeta/2),\, y\in E$ are
    disjoint sets. By transversality of $W_{0,\tau}$ and $W^s_x$ there 
    exists $C_1>0$
    such that for any $y\in\La_x$ we can find $v(y)\in W_{0,\tau}$ such
    that $y\in W^s_x(v(y),C_1\eta),$ and so for some $w(y)\in E,$
    $y\in U_{\La_x}(t,w(y),C_2(\zeta+\eta))$ with some $C_2>0$ large
    enough but independent of $x,y,\zeta,\eta.$ Hence, $E$ is 
    $(C_2(\zeta+\eta),t)-$spanning, and so $W_{0,\tau}\subset\cup_{y\in E}
    U_{W_{0,\tau}}(t,y,C_2(\zeta+\eta))$. Assume, first, that $g=g(x',v)$ 
    in (\ref{1.3.6}) is H\" older continuous in $v.$ Then by standard
    volume lemma arguments (see \cite{BR}) we obtain for $V= 
    U_{W_{0,\tau}}(t,y,\gam),\,\gam>0$ and $y\in W_{0,\tau}$ that
    \begin{eqnarray*}
    &\big\vert\log\int_V\exp\big(\int^t_0g(x',F_x^rv)dr\big)dm_V(v)\\
    &-\int_0^t\big(g(x',F^r_xy)+\vf_x^u(F^r_xy)\big)dr\big\vert\leq C(\gam)
    \nonumber\end{eqnarray*}
    where $C(\gam)>0$ does not depend on $x,x',y,t.$ Now (\ref{1.3.6})
    follows from the above integral estimates
    and the uniform in $(\gam,t)-$separated and $(\gam,t)-$spanning 
    sets approximation of the topological pressure (see, for instance,
    \cite{Bo1} and \cite{Fra}). For a general continuous $g$ approximate
    it uniformly by H\" older continuous functions and (\ref{1.3.6}) will
    follow in this case again. The limit (\ref{1.3.6}) is uniform in $x'$
    and in $z$ since $P_{F_{\pi_1z}}(g(x',\cdot)+\vf_{\pi_1z}^u)$ uniformly
    continuous in $x'$ (easy to see) and it is uniformly continuous in $z$
    (see \cite{Co}) and, furthermore, it follows from Lemma 5.1 from \cite{Co}
    that the family $\frac 1t\log\int_D\exp\left(\int_0^tg(x',F^r_{\pi_1z}
   \pi_2v)dr\right)dm_D(v),\, t\geq 1$ is equicontinuous in $z$.
   \end{proof}

   \begin{proposition}\label{prop1.3.4} 
   For any $\rho,C,b>0$ with $C$ large and $C\rho$ small enough
   there exists a positive function $\zeta_{b,\rho,T}(\Del,s,\ve)$ such that
   \begin{equation}\label{1.3.7}
    \limsup_{\Del\to 0}\limsup_{\ve\to 0}\limsup_{s\to\infty}
    \zeta_{b,\rho,T}(\Del,s,\ve)=0
    \end{equation}
    and for any $x,x'\in\cX,\, y\in \cW,\, t\geq t_1,$ $\tau\leq\frac T\ve
    -t$, $\be\in\bbR^d,\,|\be|\leq b$, $D\in\cD_\ve^u((x,y),\al,\rho,C),$
   $z\in D$ and $V=U^\ve_D(t,z,C\rho)$ satisfying $V\cap\partial D=\emptyset$ 
   we have
     \begin{eqnarray}\label{1.3.8}
    &\,\,\,\,\bigg\vert\frac 1\tau\log\int_V\exp\langle\be,\int_t^{t+\tau}B(x',
    Y^\ve_v(s))ds\rangle dm_D(v)+\frac 1\tau\log J^u_\ve(t,z)\\
    &-P_{\pi_1z_t}(\langle\be, B(x',\cdot)\rangle+\vf^u_{\pi_1z_t})\bigg
    \vert\leq\zeta_{b,\rho,T}(\ve\tau,\min (\tau,(\log\frac 1\ve)^\la),\ve)
    \nonumber
    \end{eqnarray}
    where $\langle\cdot,\cdot\rangle$ is the inner product, $z_t=\Phi_\ve^tz,$
     $\la\in(0,1)$, and $m_D$ is the induced Riemannian volume on $D.$
    \end{proposition}
     \begin{proof} 
    By (\ref{1.1.1}) and (\ref{1.2.15}) for any $w,\tilde w\in\bbR^d\times\cW,$
    \begin{eqnarray*}
    &d(\Phi_\ve^sw,\Phi^s\tilde w)\leq d(w,\tilde w)+\int_0^s\| b(\Phi_\ve^uw)
    -b(\Phi^u\tilde w)\|du\\
    &\leq d(w,\tilde w)+K\int_0^sd(\Phi_\ve^uw,\Phi^u\tilde w)du,
    \end{eqnarray*}
    where, $d=dist$ and, recall, $\Phi=\Phi_0.$ Then, by Gronwall's inequality
    \begin{equation*}
    d(\Phi_\ve^sw,\Phi^s\tilde w)\leq e^{Ks}d(w,\tilde w).
    \end{equation*}
    Hence, 
    \begin{equation*}
    d_{\bfM}(\pi_2\Phi_\ve^sw,F^s_{\pi_1z_r}\pi_2w)\leq d(\Phi^s_\ve w,\Phi^s
    (\pi_1z_r,\pi_2w))\leq e^{Ks}|\pi_1w-\pi_1z_r|.
    \end{equation*}
    
    Recall, that by Lemma \ref{lem1.3.2}(iii) for any $w\in V_r(z),\, r\leq t$,
    \begin{equation*}
    |\pi_1z_r-\pi_1w|\leq Cc^{-1}\rho\ve\al^{-1}(1-\ve\al^{-1}-
    \ve\al^{-2})^{-1}.
    \end{equation*}
   
     Set 
    \begin{equation*}
    I^\ve_{x'}(v,r,q)=\exp\big\langle\be,\int_r^qB(x',\pi_2(\Phi_\ve^sv))ds
    \big\rangle.
    \end{equation*}
    Then
    \begin{equation}\label{1.3.9}
    \int_{V_{0,r}(z)}I^\ve_{x'}(v,r,r+s)dm_D(v)=\int_{V_r(z)}I^\ve_{x'}(w,0,s)
    J_\ve^{\Xi_w}(-r,w)dm_{V_r(z)}(w)
    \end{equation}
   where $V_{0,r}(z)=U^\ve_D(r,z,C\rho)$, $\Xi_w$ is the tangent space to 
   $V_r(z)$ at $w$ and $m_{V_r(z)}$
  is the induced Riemannian volume on $V_r(z).$ By (\ref{1.3.4}), (\ref{1.3.5}),
  Lemma \ref{lem1.3.2}(i), and the H\" older continuity of the function 
  $\vf^u,$
    \begin{equation}\label{1.3.10}
    C_2^{-1}(1+C_1\ve)^{-r}\leq J_\ve^{\Xi_w}(-r,w)J^u_\ve(r,z)\leq
     C_2(1-C_1\ve)^{-r}
     \end{equation}
     for some $C_2>0$ independent of $\ve,\, r,\, z\in D$ and $w\in V_r(z).$
     Since $V_r(z)\in\cD_\ve^u(z_r,\al,\rho,\sqrt C)$ by Lemma 
    \ref{lem1.3.2}(ii), it follows from (\ref{1.2.15}) and the above estimates
     that
     \begin{eqnarray}\label{1.3.11}
     &\left(\nu_{b}(\ve,s)\right)^{-1}
     \leq\int_{V_r(z)}I^\ve_{x'}(w,0,s)dm_{V_r(z)}(w)\\
     &\times\left(\int_{V_r(z)}\exp\big\langle\be,\int^s_0 B(x',
     F^\sig_{\pi_1z_r}\pi_2w)d\sig\big\rangle dm_{V_r(z)}(w)\right)^{-1}
     \leq \nu_{b}(\ve,s)
     \nonumber\end{eqnarray}
     where
     \begin{equation*}
     \nu_{b}(\ve,s)=2+C_3\exp(C_3be^{Ks}\ve ),
     \end{equation*}
      $C_3>0$ is a constant independent of $\ve,\rho,b,r,z,x'$ and 
     $\be\in\bbR^d$ with $|\be|\leq b$.
     
     Next, choose $\la\in(0,1)$ and $\te(\ve)\in[(\log\frac 1\ve)^\la,
     2(\log\frac 1\ve)^\la]$ so that $n=\tau/\te(\ve)$ is an integer.
     If $n\leq 1$ then (\ref{1.3.8}) follows from (\ref{1.3.9})--(\ref{1.3.11})
     and Lemma \ref{lem1.3.3}. Now, let $n>1,\, k<n$ and $v\in V_{0,t}(z)$. 
      Then by (\ref{1.2.15}) and Lemma \ref{lem1.3.2}(i) for any 
     $w\in V_{0,t+k\te(\ve)}(v),$
     \begin{eqnarray*}
     &C_4^{-1}\leq I^\ve_{x'}(w,t,t+(k+1)\te(\ve))\\
     &\times\bigg( I^\ve_{x'}(v,t,t+k\te(\ve))I^\ve_{x'}(w,t+k\te(\ve),t+
     (k+1)\te(\ve))\bigg)^{-1}\leq C_4 \nonumber
     \end{eqnarray*}
     and
     \begin{equation}\label{1.3.12}
     C_4^{-1}\leq I^\ve_{x'}(v,t,t+k\te(\ve))\left(I^\ve_{x'}(w,t,t+k\te(\ve))
     \right)^{-1}\leq C_4
     \end{equation}
     where $C_4=C_4(b)= e^{4bKC\rho c^{-1}\ka^{-1}}$. Integrating the 
     inequalities above we obtain
     \begin{eqnarray}\label{1.3.13}
     &\,\,\,\,\,\,\,\,\,\,C_4^{-1}\leq \int_{V_{0,t+k\te(\ve)}(v)}I^\ve_{x'}
     (w,t,t+(k+1)\te(\ve))dm_D(w)\bigg( I^\ve_{x'}(v,t,t+k\te(\ve))\\
     &\times\int_{V_{0,t+k\te(\ve)}(v)}I^\ve_{x'}(w,t+k\te(\ve),t+(k+1)\te(\ve))
     dm_D(w)\bigg)^{-1}\leq C_4.\nonumber
     \end{eqnarray}
      From the estimates (\ref{1.3.9})--(\ref{1.3.11}) together with Lemma 
      \ref{lem1.3.3} we conclude that for some $C_5>0$ independent of $t,k,
      \ve,\rho,v,$ and $x'$,
     \begin{eqnarray}\label{1.3.14}
     &C_5^{-1}e^{-\te(\ve)\eta_{b,\rho}(\ve)}(\nu_{b}\big(\ve,\te(\ve))
     \big)^{-1}\leq\int_{V_{0,t+k\te(\ve)}(v)}
     I^\ve_{x'}(w,t+k\te(\ve),\\
     &t+(k+1)\te(\ve))dm_D(w)J^u_\ve(t+k\te(\ve),v)
     \exp\big(-\te(\ve)P_{F_{\pi_1v_{t+k\te(\ve)}}}(\langle\be, B(x',
     \cdot)\rangle\nonumber\\
     &+\vf^u_{\pi_1v_{t+k\te(\ve)}})\big)\leq e^{\te(\ve)\eta_{b,\rho}(\ve)}
     \nu_{b}(\ve,\te(\ve))\nonumber
     \end{eqnarray}
     where $v_s=\Phi_\ve^sv$ and $\eta_{b,\rho}(\ve)>0$, 
     $\eta_{b,\rho}(\ve)\to 0$ as $\ve\to 0.$ 
     Observe that by (\ref{1.1.1}), (\ref{1.2.15}) and Lemma 
     \ref{lem1.3.2}(iii),
     \[
     |\pi_1v_{t+k\te(\ve)}-\pi_1z_t|\leq K\ve\tau+Cc^{-1}\rho\ve\al^{-1}
     (1-\ve\al^{-1}-\ve\al^{-2})^{-1},
     \]
     and so setting $P=P_{F_{\pi_1z_t}}(\langle\be, B(x',\cdot)\rangle
     +\vf^u_{\pi_1z_t})$
     we obtain by (\ref{1.2.15}) and \cite{Co} (see also \cite{KKPW} and 
     \cite{Ru}) that
     \begin{equation}\label{1.3.15}
     \big\vert P-P_{F_{\pi_1v_{t+k\te(\ve)}}}(\langle\be, B(x',\cdot)
     \rangle +
     \vf^u_{\pi_1v_{t+k\te(\ve)}})\big\vert\leq C_6\ve\tau
     \end{equation}
     for some $C_6=C_6(b)>0$ independent of $v,k\leq n,\ve,t,z,x',$ and 
     $\be\in\bbR^d$ with $|\be|\leq b$ provided, say, $\tau\geq 1$ which we
      can assume without loss of generality.
     
     A finite set $E\subset D$ will be called
     $(s,\gam,\ve,D)-$separated if $v_i,v_j
     \in E,$ $v_i\ne v_j$ implies that $v_i\not\in U^\ve_D(s,v_j,\gam).$
     Let $E_k$ be a maximal $(t+k\te(\ve),C\rho,\ve,D)-$separated set
     in $D$ and define
     \[
     E_k^U=\{ v\in E_k:\, U^\ve_D(t+k\te(\ve),v,C\rho)\cap U\ne\emptyset\}.
     \]
     Then for $k\geq 1$,
     \begin{equation}\label{1.3.16}
     U^\ve_D(t,z,\gam+a_kC\rho)\supset\cup_{v\in E_k^{U_D^\ve(t,z,\gam)}}
     U^\ve_D(t+k\te(\ve),v,C\rho)\supset U^\ve_D(t,z,\gam)
     \end{equation}
     where $a_k=c^{-1}e^{-\frac 12\ka k\te(\ve)}$ and the left hand side 
     of (\ref{1.3.16}) follows from Lemma \ref{lem1.3.2}(i) assuming that $\ve$
     is small enough. Observe also that $U_D^\ve(t+k\te(\ve),v,C\rho/2)$ are
     disjoint for different $v\in E_k.$ For $k=1,2,...,n-1$
     set $V(k,\rho)=U^\ve_D(t,z,C\rho(1+\sum^{n-1}_{j=k}a_j))$ and
     $V(-k,\rho)=U^\ve_D(t,z,C\rho(1-\sum^{n-1}_{j=k}a_j))$ with 
     $V(n,\rho)=V=V_{0,t}(z).$ Then by (\ref{1.3.12})--(\ref{1.3.16}) and Lemma
     \ref{lem1.3.2}(iii),
     \begin{eqnarray*}
&\int_{V(k+1,\rho)}I^\ve_{x'}(v,t,t+(k+1)\te(\ve))dm_D(w)\\
&\leq\sum_{v\in E_{k}^{V(k+1,\rho)}}\int_{V_{0,t+k\te(\ve)}(v)}
I^\ve_{x'}(w,t,t+(k+1)\te(\ve))dm_D(w)\nonumber\\
&\leq C_4\sum_{v\in E_{k}^{V(k+1,\rho)}}I^\ve_{x'}(v,t,t+k\te(\ve))
\int_{V_{0,t+k\te(\ve)}(v)}I^\ve_{x'}(w,t+k\te(\ve),t\nonumber\\
&+(k+1)\te(\ve))dm_D(w)\leq C_5C_4e^{\te(\ve)(\eta_{b,\rho}(\ve)+C_6\ve\tau+P)}
 \nu_b(\ve,\te(\ve))\sum_{v\in E_{k}^{V(k+1,\rho)}}\nonumber\\
 &\big(J^u_\ve(t+k\te(\ve),v)\big)^{-1}
 I^\ve_{x'}(v,t,t+k\te(\ve))\nonumber\\
 &\leq C_5C_4^2c_{\rho/2,T}^{-1}
 e^{\te(\ve)(\eta_{b,\rho}(\ve)+C_6\ve\tau+P)}\nu_b(\ve,\te(\ve))\nonumber\\
 &\times\sum_{v\in E_{k}^{V(k+1,\rho)}}\int_{U_D^\ve(t+k\te(\ve),v,C\rho/2)}
     I^\ve_{x'}(w,t,t+k\te(\ve))dm_D(w)\nonumber\\
     &\leq C_5C_4^2c_{\rho/2,T}^{-1} e^{\te(\ve)(\eta_{b,\rho}(\ve)+C_6\ve\tau
     +P)}\nu_b(\ve,\te(\ve))\int_{V(k,\rho)}I^\ve_{x'}(w,t,t+k\te(\ve))dm_D(w).
  \nonumber
 \end{eqnarray*}
 Similarly, we obtain
 \begin{eqnarray*}
  &\int_{V(-(k+1),\rho)}I^\ve_{x'}(w,t,t+(k+1)\te(\ve))dm_D(w)\\
  &\geq C_5^{-1}C_4^{-2}c_{\rho,T}e^{-\te(\ve)(\eta_{b,\rho/2}(\ve)+C_6\ve\tau
  +P)}(\nu_b(\ve,\te(\ve))^{-1}\nonumber\\
  &\times\int_{V(-k,\rho)}I^\ve_{x'}(w,t,t+k\te(\ve))dm_D(w).
  \nonumber\end{eqnarray*}
  Emloying these estimates recursively for $k=n-1, 
  n-2,...,1$ and estimating $\int_{V(\pm 1,\rho)}I^\ve_{x'}(w,t,t+\te(\ve))
  dm_D(w)$ by (\ref{1.3.14}) with $k=0$ and with $V(\pm 1,\rho)$ in place
  of $V_{0,t}(v)$  we derive (\ref{1.3.8}) with
  $\zeta_{b,\rho,T}(\Del,s,\ve)=s^{-1}\log\big(C_5C^2_4(c^{-1}_{\rho,T}+
  c^{-1}_{\rho/2,T})\big)+\eta_{b,\rho}(\ve)+\eta_{b,\rho/2}(\ve)+
  2(\te(\ve))^{-1}\log\nu_b(\ve,\te(\ve))+C_6\Del$
  provided $n\geq 1.$ 
  \end{proof}

    Next, under Assumption \ref{ass1.2.6} we derive a volume lemma 
    \index{volume lemma} type 
  assertion (see \cite{BR}) which will be needed in Sections \ref{sec1.7} and
  \ref{sec1.8} and which will hold true on any time intervals and not just on 
  time intervals of order $1/\ve$ as in Lemma \ref{lem1.3.2}(iv). In order to 
  do so we will consider a subset of embedded $C^2$ discs from $\cD_\ve^u(z,\al,
  \rho,C)$ taking special care of their $C^2$ bounds.
    
    Namely, let Exp$_y:\, T_y\bfM\to\bfM$ be the exponential map which is a 
    diffeomorphism of $V^y_\del=\{\xi\in T_y\bfM:\,\|\xi\|<\del\}$ onto the 
    open $\del-$neighborhood $U_\del(y)$ of $y$ provided $\del>0$ is small
    enough. Given $x,\tilde x\in\bbR^d,\, y\in\bfM$ and $\xi\in T_y\bfM$ set
    \begin{equation*}
    \chi_{x,y}(\tilde x,\xi)=(x+\tilde x,\,\mbox{Exp}_y\xi)
    \end{equation*}
    which is a diffeomorphism of $\bbR^d\times V_\del^y$ onto $\bbR^d\times
    U_\del(y)$. Let $D\in\cD_\ve^u(z,\al,\rho,C)$, $z=(x,y)$, $y\in\cW$ be an
    embedded $C^2$ disc. Assuming that $C^2\rho<\del$ we can define 
    $\hat D=\chi^{-1}_{x,y}(D)$ which is a $C^2$ hypersurface in 
    $\{\tilde x:\,|\tilde x\|<\del\}\times V^y_\del$. If $\del$ is sufficiently 
    small then the tangent subbubndle $T\hat D$ over $\hat D$ still stays
    close to $\Gam^u_x(y)$, and so we can represent $\hat D$ as a parametric
    set $(\eta,\vf(\eta),x(\eta))$ where $\eta\in\Gam^u_x(y)$, $\vf(\eta)\in
    \Gam_x^{0s}(y)$ and $x(\eta)\in\bbR^d$. We will write that
    $D\in\hat\cD_\ve^u(z,\al,\rho,C,L)$ if the parametric representation
    of the corresponding $\hat D$ as above satisfies
    \begin{equation*}
    \max_{i,j,k,l}\max\big(\big\vert\frac {\partial^2\vf_i(\eta)}
    {\partial\eta_k\partial\eta_l}\big\vert,\,\frac {\partial^2x_j(\eta)}
    {\partial\eta_k\partial\eta_l}\big\vert\big)\leq L .
    \end{equation*}
    
    \begin{lemma}\label{lem1.3.5} There exists $t_1\geq t_0$ such that for any
    $t\geq t_1$ we can choose $\del>0$ small enough and $L>0$ large enough
    so that if $D\in\hat\cD^u_\ve(z,\al,\rho,C,L)$ and $C^2\rho<\del$ then
    \begin{equation*}
    \{ v\in\Phi_\ve^tD:\, d_{\Phi_\ve^tD}(\Phi_\ve^tz,v)\leq C^2\rho\}
    \in\hat\cD_\ve^u(\Phi_\ve^tz,\al,\rho,C,L).
    \end{equation*}
    \end{lemma}
    \begin{proof}  Since the differential $D_0$Exp$_y$ of the exponential map
    at zero is the identity map it follows from the definition of 
    $\cD_\ve^u(z,\al,\rho,C)$ that 
    \begin{equation*}
    \max_{i,j,k}\max\big(\big\vert\frac {\partial\vf_i(\eta)}
    {\partial\eta_k}\big\vert,\,\frac {\partial x_j(\eta)}
    {\partial\eta_k}\big\vert\big)\leq c(\ve,\del)
    \end{equation*}
    where $c(\ve,\del)\to 0$ (uniformly in all $D$ as above) as 
    $\ve,\del\to 0$. Let $z=(x,y)$ and set
    \begin{equation*}
    f^t_{x,y,\ve}=\chi^{-1}_{\Phi_\ve^t(x,y)}\circ\Phi_\ve^t\circ\chi_{x,y}
    \end{equation*}
    which for each fixed $t>0$ and a sufficiently small $\del>0$ (depending
    on $t$) defines a diffeomorphism of $\bbR^d\times V_\del^y$ onto its
    image. By (\ref{1.3.2}) the tangent subbundle over $\Phi_\ve^tD$ is 
    contained in $\cC^u(\ve,\al)$, and so for small $\del>0$ the tangent
    subbundle $T(f^t_{x,y,\ve}\hat D)$ over $f^t_{x,y,\ve}\hat D$ stays close 
    to $\Gam^u_{\tilde x}(\tilde y)$ where $\tilde x=\pi_1\Phi_\ve^t(z)$ and
     $\tilde y=\pi_2\Phi_\ve^t(z)$. Hence, we can represent $f^t_{x,y,\ve}
     (\hat D)$ in a parametric form $\big(\tilde\eta,\tilde\vf(\tilde\eta),
     \tilde x(\tilde\eta)\big)$ where $\tilde\eta\in\Gam^u_{\tilde x}
     (\tilde y),\,\tilde\vf(\tilde\eta)\in\Gam_{\tilde x}^{0s}(\tilde y)$
     and $\tilde x(\tilde\eta)\in\bbR^d$. Fix some $t>t_0$ so that
     (\ref{1.3.3}) holds true. Write $f^t_{x,y,\ve}(\eta,\vf,x)=(\tilde\eta,
     \tilde\vf,\tilde x)$, so that, in particular,
     \begin{equation*}
     f^t_{x,y,\ve}(\eta,\vf(\eta),x(\eta))=\big(\tilde\eta,
     \tilde\vf(\tilde\eta),\tilde x(\tilde\eta)\big).
     \end{equation*}
     Then
     \begin{equation}\label{1.3.17}
     \tilde\eta=A\eta+a_{\ve,\del}(\eta,\vf,x),\,\tilde\vf=B\vf+b_{\ve,\del}
     (\eta,\vf,x),\,\tilde x=x+c_{\ve,\del}(\eta,\vf,x)
     \end{equation}
     where $\eta\in\Gam^u_x(y)$, $\vf\in\Gam_x^{0s}(y)$, $x\in\bbR^d$,
     $A$ is an $n^u\times n^u$--matrix, $B$ is an 
     $n^{0s}\times n^{0s}$--matrix with $n^{0s}=n_\bfM-n^u$ and 
     \begin{equation}\label{1.3.18}
     \| a_{\ve,\del}(\eta,\vf,x)\|_{C^1}+\| b_{\ve,\del}(\eta,\vf,x)\|_{C^1}
     +\| c_{\ve,\del}(\eta,\vf,x)\|_{C^1}\leq c(\ve,\del)
     \end{equation}
     for all $(\eta,\vf)\in V^y_\del$ and $|x|<\del$ where $c(\ve,\del)\to 0$
     as $\ve,\del\to 0$. By (\ref{1.2.15}), (\ref{1.3.1}) and Assumption
     \ref{ass1.2.6} it follows that there exists a constant $R>0$ such that
     for any $\xi\in\Gam_x^{0s}(y)$, $x\in\bbR^d$, $y\in\bfM$,
     \begin{equation}\label{1.3.19}
     \| DF_x^t\xi\|\leq R\|\xi\|.
     \end{equation}
     By (\ref{1.3.1}) we can choose $t>0$ large enough and then $\ve$ and $\del$
     small enough so that for all $\eta\in\Gam_x^u(y)$,
     \begin{equation}\label{1.3.20}
     \| A\eta\|\geq (1+R)\|\eta\|.
     \end{equation}
     Now $t$ is fixed and we can choose $\ve$ and $\del$ so small that
     (\ref{1.3.19}) implies that,
     \begin{equation}\label{1.3.21}
     \| B\|\leq 1+R.
     \end{equation}
     In order to shorten notations for every vector function $f\big(\zeta)=
     (f_1(\zeta),...,f_l(\zeta)\big)$, $\zeta=(\zeta_1,...,\zeta_k)$ we
     denote by $\frac {\partial f}{\partial\zeta}$ the Jacobi matrix
     $(\partial f_i(\zeta)/\partial\zeta_j)$ and by $\frac {\partial^2f}
     {\partial\zeta^2}$ we denote the collection $(\partial^2 f_i(\zeta)/
     \partial\zeta_j\partial\zeta_k)$. We set also 
     \begin{equation*}
     \|\frac {\partial f}{\partial\zeta}\|=\max_{i,j}|\frac 
     {\partial f_i(\zeta)}{\partial\zeta_j}|\,\,\mbox{and}\,\,
     \|\frac {\partial^2f}{\partial\zeta^2}\|=\max_{i,j,k}|\frac 
     {\partial^2f_i}{\partial\zeta_j\partial\zeta_k}|.
     \end{equation*}
    Observe that by Assumptions \ref{ass1.2.1}, \ref{ass1.2.2} and 
    \ref{ass1.2.6}
    for any $t>0$ there exists $\hat R=\hat R_t>0$ such that
    \begin{equation}\label{1.3.22}
    \sup_{\ve\leq 1}\sup_{|u|\leq t}\|\Phi_\ve^u\|_{C^2}\leq\hat R
    \end{equation}
    and
    \begin{equation}\label{1.3.23}
    \max\big(\| a_{\ve,\del}(\eta,\vf,x)\|_{C^2},\,\| b_{\ve,\del}
    (\eta,\vf,x)\|_{C^2},\,\|c_{\ve,\del}(\eta,\vf,x)\|_{C^2}\big)
    \leq 2\hat R+1.
    \end{equation}
    It follows by (\ref{1.3.17})--(\ref{1.3.23}) (with natural product 
    notations) that
    \begin{equation*}
    \big\|\frac {\partial^2\tilde\vf}{\partial\tilde\eta^2}\big\|=
    \big\|\frac {\partial^2\tilde\vf}{\partial\eta^2}\big(\frac {\partial\eta}
    {\partial\tilde\eta}\big)^2+\frac {\partial\tilde\vf}{\partial\eta}
    \frac {\partial^2\eta}{\partial\tilde\eta^2}\big\|\leq
    \big\|\frac {\partial^2\tilde\vf}{\partial\eta^2}\big\|(1+R)^{-2}+
    c(\ve,\del)\hat R
    \end{equation*}
    and
    \begin{eqnarray*}
    &\big\|\frac {\partial^2\tilde\vf}{\partial\eta^2}\big\|=\big\|
    B\frac {\partial^2\vf}{\partial\eta^2}+\frac {\partial^2b_{\ve,\del}}
    {\partial\eta^2} +2\frac {\partial^2b_{\ve,\del}}{\partial\eta\partial\vf}
    \frac {\partial\vf}{\partial\eta}+2\frac {\partial^2b_{\ve,\del}}
    {\partial\eta\partial x}\frac {\partial x}{\partial\eta}
    \frac {\partial^2b_{\ve,\del}}{\partial\vf^2}\big(\frac {\partial\vf}
    {\partial\eta}\big)^2\\
    &+2\frac {\partial^2b_{\ve,\del}}
    {\partial\vf\partial x}
    \frac {\partial\vf}{\partial\eta}\frac {\partial x}{\partial\eta}
    +\frac {\partial^2b_{\ve,\del}}{\partial x^2}\big(\frac {\partial x}
    {\partial\eta}\big)^2\big\|
    \leq (1+R)L+(2\hat R+1)(1+c(\ve,\del))^2.\nonumber
    \end{eqnarray*}
    Similarly,
    \begin{equation*}
    \big\|\frac {\partial^2\tilde x}{\partial\eta^2}\big\|
     \leq (1+R)L+(2\hat R+1)(1+R)^{-2}(1+c(\ve,\del))^2+c(\ve,\del)\hat R.
     \end{equation*}
    Choosing $L\geq R^{-1}(2\hat R+R+2)$ we obtain that if 
    \[
    \max\big\{\big\|\frac {\partial^2\vf}{\partial\eta^2}\big\|,\,
    \big\|\frac {\partial^2x}{\partial\eta^2}\big\|\big\}\leq L
    \]
    then
    \[
    \max\big\{\big\|\frac {\partial^2\tilde\vf}{\partial\tilde\eta^2}\big\|,\,
    \big\|\frac {\partial^2\tilde x}{\partial\tilde\eta^2}\big\|\big\}
    \leq L
    \]
    and the assertion of Lemma \ref{lem1.3.5} follows.
     \end{proof}
     The main purpose of the previous result is to derive the following 
     volume lemma type assertion which plays an essential role in Section 
     \ref{sec1.6}.
     
     \begin{lemma}\label{lem1.3.6} For any $\be\in(0,C^2\rho)$ there exists
     $c_\be>0$ such that if $D\in\hat D^u_\ve(z,\al,\rho,C,L)$ and $L$
     is large enough then for any $t>0$ and $v,w\in D$ satisfying
     $w\in U^\ve_D(t,v,\be)\subset D$,
     \begin{equation}\label{1.3.24}
     c_\be\leq m_D\big(U^\ve_D(t,v,\be)\big)J_\ve^{T_wD}(t,w)\leq c_\be^{-1}.
     \end{equation}
     \end{lemma}
     \begin{proof} Set $V_{s,t}=\Phi_\ve^sU_D^\ve(t,v,\be)$ and $V_t=V_{t,t}$.
     Similarly to Lemma \ref{lem1.3.2}(ii), $V_t\in\cD_\ve^u(\Phi_\ve^tv,\al,
     \be C^{-1},\sqrt C)$, and so by uniformity considerations there exists
     $\tilde c_\be>0$ independent of $v,t$ and $D$ as above such that
     \begin{equation}\label{1.3.25}
     \tilde c_\be\leq m_{V_t}(V_t)=\int_{U^\ve_D(t,v,\be)}J^{T_wD}_\ve(t,w)
     dm_D(w)\leq\tilde c_\be^{-1}.
     \end{equation}
     Choose $l\in\bbN$ so that $t_2=t/l\in[t_1,2t_1)$ and set $w_k=
     \Phi_\ve^{kt_2}w,\, v_k=\Phi_\ve^{kt_2}v$. Then for any $w\in 
     U^\ve_D(t,v,\be)$,
     \begin{equation}\label{1.3.26}
     J_\ve^{T_wD}(t,w)=\prod_{k=0}^{l-1}J_\ve^{T_{w_k}V_{kt_2,t}}(t_2,w_k)
     \end{equation}
     and by Lemma \ref{lem1.3.2}(i),
     \begin{equation}\label{1.3.27}
     d_{V_{kt_2,t}}(w_k,v_k)\leq c^{-1}\be e^{-\frac 12\ka(l-k)t_2}.
     \end{equation}
     By (\ref{1.3.5}), (\ref{1.3.22}), (\ref{1.3.26}), and (\ref{1.3.27}) 
     together
     with Lemma \ref{lem1.3.5} we conclude that there exists a constant 
     $\tilde C>0$ such that
     \begin{equation}\label{1.3.28}
     \big\vert\ln J_\ve^{T_{w_k}V_{kt_2,t}}(t_2,w_k)-
     \ln J_\ve^{T_{v_k}V_{kt_2,t}}(t_2,v_k)\big\vert\leq\tilde 
     Ce^{-\frac 12\ka(l-k)t_2}.
     \end{equation}
     Now (\ref{1.3.24}) follows from (\ref{1.3.25}), (\ref{1.3.26}), and 
     (\ref{1.3.28}) with 
     \[
     c_\be=\tilde c_\be\exp\big(-2\tilde C(1-e^{-\frac 12\ka t_2})^{-1}\big).
     \]
     \end{proof}

      \section{Large deviations: preliminaries}\label{sec1.4}\setcounter{equation}{0}

We will need the following version of general large deviations bounds when
usual assumptions hold true with errors. An upper bound similar to (\ref{1.4.3})
below appeared previously in \cite{Ki7}.
For simplicity we will formulate the result for $\bbR^d-$valued random 
vectors though the same arguments work for random variables with values in a 
Banach space. The proof is a strightforward modification of the standard one 
(cf. \cite{Ki1}) but still we exhibit it here for readers' convenience.

\begin{lemma}\label{lem1.4.1} Let $H=H(\be)$, $\eta=\eta(\be)$ be
uniformly bounded on compact sets functions on $\bbR^d$ and 
$\{\Xi_\tau,\,\tau\geq 1\}$ be a family of $\bbR^d-$valued random
vectors on a probability space $(\Om,\cF,P)$ such that $|\Xi_\tau|\leq C
<\infty$ with probability one for some constant $C$ and all $\tau\geq 1$.
For any $a>0$ and $\al,\be_0\in\bbR^d$ set
\begin{equation}\label{1.4.1}
L_a^{\be_0}(\al)=\sup_{\be\in\bbR^d,|\be+\be_0|\leq a}(\langle\be,\al\rangle-
H(\be)),\,\,L_a(\al)=L_a^0(\al),\,\, L(\al)=L_\infty(\al).
\end{equation}
(i) For any $\la,a>0$ there exists $\tau_0=\tau(\la, a,C)$ such that
whenever for some $\tau\ge\tau_0$, $\be_0\in\bbR^d$ and each $\be\in\bbR^d$
with $|\be+\be_0|\leq a$,
\begin{equation}\label{1.4.2}
H_\tau(\be)=\tau^{-1}\log Ee^{\tau \langle\be,\Xi_\tau\rangle}
\leq H(\be)+\eta(\be)
\end{equation}
then for any compact set $\cK\subset\bbR^d$,
\begin{equation}\label{1.4.3}
P\{\Xi_\tau\in \cK\}\leq\exp\left(-\tau (L_a^{\be_0}(\cK)-\eta_a^{\be_0}-
\la|\be_0|-\la)\right)
\end{equation}
where
\begin{equation}\label{1.4.4}
\eta_a^{\be_0}=\sup\{\eta(\be):\,|\be+\be_0|\leq a\}\,\,\mbox{and}\,\,
L_a^{\be_0}(\cK)=\inf_{\al\in \cK}L_a^{\be_0}(\al).
\end{equation}
(ii) Suppose that $\al_0\in\bbR^d$, $0<a\leq\infty$ and there exists 
$\be_0\in\bbR^d$ such that $|\be_0|\leq a$ and
\begin{equation}\label{1.4.5}
H(\be_0)=\langle\be_0,\al_0\rangle -L_a(\al_0).
\end{equation}
If (\ref{1.4.2}) holds true then for any $\del>0$,
\begin{equation}\label{1.4.6}
P\{ |\Xi_\tau-\al_0|\leq\del\}\leq\exp\left( -\tau(L_a(\al_0)-\eta(\be_0)-
\del|\be_0|)\right).
\end{equation}
(iii) Assume that $\al_0,\be_0\in\bbR^d$ satisfy (\ref{1.4.5}). For any
$\la,a>0$ there exists $\tau_0=\tau(\la,a,C)$ such that whenever for some
$\tau\geq\tau_0$ and each $\be\in\bbR^d$ with $|\be|\leq a$ the inequality
(\ref{1.4.2}) holds true together with
\begin{equation}\label{1.4.7}
\tau^{-1}\log Ee^{\tau \langle\be,\Xi_\tau\rangle}\geq H(\be)-\eta(\be)
\end{equation}
then for any $\gam,\del>0,\,\gam\leq\del$,
\begin{eqnarray}\label{1.4.8}
&P\{ |\Xi_\tau-\al_0|<\del\}\geq\exp\left( -\tau(L(\al_0)+\eta(\be_0) +
\gam|\be_0|)\right)\\
&\times\left( 1-\exp\big(-\tau(\tilde L_a^{\be_0}(\cK_{\gam,C}(\al_0))-\eta_a
-\eta(\be_0)-\la|\be_0|-\la)\big)\right)\nonumber
\end{eqnarray}
where 
\[
\tilde L_a^{\be_0}(\al)=L_a(\al)-\langle\be_0,\al\rangle +H(\be_0),
\]
$\tilde L_a^{\be_0}(\cK)=\inf_{\al\in \cK}\tilde L_a^{\be_0}(\al),$ $\eta_a=
\eta_a^0$, $\cK_{\gam,C}(\al_0)=\overline{U_C(0)}\setminus U_\gam(\al_0),$
 $U_\gam(\al)=\{\tilde\al:\,|\tilde\al-\al|<\gam\}$ and $\bar U$ denotes 
 the closure of $U$.
\end{lemma}
\begin{proof}(i) By (\ref{1.4.1}) for any $\al\in \cK_C=\cK\cap\overline{U_C(0)}$
and $\la>0$ there exists $\be_{\la}(\al)\in\bbR^d$ such that
\begin{equation}\label{1.4.9}
|\be_\la(\al)+\be_0|\leq a\,\,\mbox{and}\,\,
\langle\be_{\la}(\al),\al\rangle -H(\be_{\la}(\al))>L_a^{\be_0}(\al)-\la/2.
\end{equation}
Set $\gam_{a,\la}(\al)=\frac \la 2\min(1,a^{-1})$ and
 cover the compact set $\cK_C$ by open balls $U_{\gam_{a,\la}}(\al),\,
 \al\in \cK_C.$ Let $U_{\gam_{a,\la}(\al_1)},...,U_{\gam_{a,\la}(\al_n)}$
 be a finite subcover with a minimal number $n$ of elements. Observe that
 $n$ does not exceed the maximal number of points in $\overline{U_C(0)}$ with
 pairwise distances at least $\frac 12\gam_{a,\la}$ and the latter number
 depends only on $C,a$ and $\la$. By (\ref{1.4.2}) and (\ref{1.4.9}) for each 
 $i=1,...,n$,
 \begin{eqnarray*}
 &e^{\tau\eta_a^{\be_0}}\geq E\bbI_{\Xi_\tau\in U_{\gam_{a,\la}(\al_i)}
 (\al_i)}e^{\tau(\langle\be_\la(\al_i),\Xi_\tau\rangle -H(\be_\la(\al_i)))}\\
 &\geq e^{-\tau\la/2}E\bbI_{\Xi_\tau\in U_{\gam_{a,\la}(\al_i)}
 (\al_i)}e^{\tau\big(\langle\be_\la(\al_i),\al_i\rangle-|\be_0|\la/2-
 H(\be_\la(\al_i))\big)}\\
 &e^{\tau\big(L_a^{\be_0}(\al_i)-|\be_0|\la/2-\la/2\big)}
 P\{\Xi_\tau\in U_{\gam_{a,\la}(\al_i)}(\al_i)\}.
 \end{eqnarray*}
 Since $L_a^{\be_0}(\al_i)\geq L_a^{\be_0}(\cK)$ and $|\Xi_\tau|\leq C$ then
 summing these inequalities in $i=1,...,n$ we obtain
 \begin{equation}\label{1.4.10}
  P\{\Xi_\tau\in \cK\}=P\{\Xi_\tau\in \cK_C\}\leq e^{-\tau\big(L_a^{\be_0}(\cK)
  -\eta_a^{\be_0}-|\be_0|\la/2-\tau^{-1}\log n-\la/2\big)}.
  \end{equation}
  Since $n$ is bounded by a number depending only on $\la,a,$ and $C$ we
  can choose $\tau_0=\tau_0(\la,a,C)$ so that $\tau_0^{-1}\log n\leq\la/2$
  which together with (\ref{1.4.10}) yield (\ref{1.4.3}).
  
  (ii) By (\ref{1.4.2}) and (\ref{1.4.5}),
  \begin{equation*}
 e^{\tau\eta(\be_0)}\geq E\bbI_{\Xi_\tau\in U_{\del}(\al_0)}
 e^{\tau\big(\langle\be_0,\Xi_\tau\rangle -H(\be_0)\big)}\geq
 e^{\tau\big(L_a(\al_0)-\del |\be_0|\big)}
 P\{\Xi_\tau\in U_\del(\al_0)\}
 \end{equation*}
 and (\ref{1.4.6}) follows.
 
 (iii) By (\ref{1.4.5}) and (\ref{1.4.7}) for any $\gam\leq\del$,
 \begin{eqnarray}\label{1.4.11}
  &P\{ |\Xi_\tau-\al_0|<\del\}\geq P\{ |\Xi_\tau-\al_0|<\gam\}\\
  &=E^{\be_0}_\tau\bbI_{|\Xi_\tau-\al_0 |<\gam}e^{-\tau\big(\langle\be_0,
  \Xi_\tau\rangle -H_\tau(\be_0)\big)}\nonumber\\
  &\geq e^{-\tau\big(L(\al_0)+|\be_0|\gam +\eta(\be_0\big)}
  P^{\be_0}_\tau\{ |\Xi_\tau-\al_0|<\gam\}\nonumber
   \end{eqnarray}
   where $E^{\be_0}_\tau$ is the expectation with respect to the
   probability measure $P^{\be_0}$ on $(\Om,\cF)$ such that
   \begin{equation*}
   \frac {dP^{\be_0}_\tau}{dP}=e^{\tau(\langle\be_0,\Xi_\tau\rangle
   -H_\tau(\be_0))}.
   \end{equation*}
   Now by (\ref{1.4.2}) and (\ref{1.4.5}) for any $\be\in\bbR^d$ with
   $|\be+\be_0|\leq a$ we obtain that
   \begin{equation}\label{1.4.12}
   \tau^{-1}\log E_\tau^{\be_0}e^{\tau\langle\be ,\Xi_\tau
   \rangle }=H_\tau(\be+\be_0)-H_\tau(\be_0))
   \leq\tilde H^{\be_0}(\be)+\tilde\eta^{\be_0}(\be)
   \end{equation}
   where $\tilde H^{\be_0}(\be)=H(\be+\be_0)-H(\be_0)$ and $\tilde\eta^{\be_0}
   (\be)=\eta(\be+\be_0)+\eta(\be_0)$. Observe that
   \begin{equation}\label{1.4.13}
   \sup_{\be\in\bbR^d,|\be+\be_0|\leq a}\big(\langle\be,\al\rangle-
   \tilde H^{\be_0}(\be)\big)=L_a(\al)-\langle\be_0,\al\rangle +H(\be_0)=
   \tilde L_a^{\be_0}(\al).
   \end{equation}
   Thus, applying (i) on the probability space $(\Om,\cF,P_\tau^{\be_0})$
   we derive that
   \begin{equation}\label{1.4.14}
   P^{\be_0}_\tau\{ |\Xi_\tau-\al_0 |\geq\gam\}\leq 
   \exp\big(-\tau\big(\tilde L_a^{\be_0}(\cK_{\gam,C}(\al_0))-\eta_a-
   \eta(\be_0)-\la|\be_0|-\la\big)\big)
   \end{equation}
   provided $\tau\geq\tau_0$ for a sufficiently large $\tau_0=\tau_0(\la,a,C)$.
   This together with (\ref{1.4.11}) yield (\ref{1.4.8}). 
   \end{proof}
   
   \begin{lemma}\label{lem1.4.2} Let $S_n,\, n=1,2,...$ be a nondecreasing sequence of lower
   semicontinuous functions on a metric space $M$ and let $S=\lim_{n\to\infty}
   S_n.$ Assume that $S$ is also lower semicontinuous and for any compact set
   $\cK\subset M$ denote
   \[
   S_n(\cK)=\inf_{\gam\in \cK}S_n(\gam)\,\,\mbox{and}\,\, 
   S(\cK)=\inf_{\gam\in \cK}S(\gam).
   \]
   Then
   \begin{equation}\label{1.4.15}
   \lim_{n\to\infty}S_n(\cK)=S(\cK).
   \end{equation}
   \end{lemma}
   \begin{proof} By the lower semicontinuity of $S_n$ and $S$ and by compactness
    of $\cK$ it follows that there exist $\hat\gam_n,\hat\gam\in \cK$ such that
    $S_n(\hat\gam_n)=S_n(\cK)$ and $S(\hat\gam)=S(\cK)$. Passing if needed to a
    subsequence assume that $\hat\gam_n\to\tilde\gam\in \cK$ as $n\to\infty$.
    Since
    \begin{equation}\label{1.4.16}
    S_n(\cK)=S_n(\hat\gam_n)\leq S_n(\hat\gam)
    \end{equation}
    then
    \begin{equation}\label{1.4.17}
    \limsup_{n\to\infty}S_n(\cK)\leq S(\hat\gam)=S(\cK).
    \end{equation}
    Assume now that $S(\cK)<\infty$. Since
    \[
    S(\cK)=S(\hat\gam)\leq S(\tilde\gam)
    \]
    then for any $\ve>0$ there exists $n(\ve)$ such that
    \begin{equation}\label{1.4.18}
    S(\hat\gam)\leq S_{n(\ve)}(\tilde\gam)+\ve .
    \end{equation}
    By the lower semicontinuity of $S_{n(\ve)}(\gam)$ it follows that for
    $m\geq n(\ve)$ large enough
    \begin{equation}\label{1.4.19}
    S(\hat\gam)\leq S_{n(\ve)}(\hat\gam_m)+2\ve\leq S_m(\hat\gam_m)+2\ve
    \end{equation}
    where we use also that $S_m,m=1,2,...$ is a nondecreasing sequence.
    Since (\ref{1.4.19}) holds true for any $m\geq n(\ve)$ large enough and
    for each $\ve>0$ we can pass there to the limit so that, first, 
    $m\to\infty$ and then $\ve\to 0$ yielding that
    \[
    S(\cK)\leq\liminf_{m\to\infty}S_n(\cK)
    \]
    which together with (\ref{1.4.17}) give (\ref{1.4.15}) under the condition
    $S(\cK)<\infty$. If $S(\cK)=\infty$ then $S(\tilde\gam)=\infty$ and for any
    $A>0$ there exists $n(A)$ such that $S_n(\tilde\gam)>A$ for any $n\geq
    n(A)$. By the lower semicontinuity of $S_n$ we conclude that
     $S_n(\hat\gam_m)>A$ for $m\geq n$ large enough which implies that
     $S_m(\hat\gam_m)>A$ for all sufficiently large $m$. Hence
     \begin{equation}\label{1.4.20}
     \liminf_{m\to\infty}S_m(\cK)=\liminf_{m\to\infty}S_m(\hat\gam_m)>A
     \end{equation}
     and since $A$ is arbitrary the left hand side of (\ref{1.4.20}) equals
     infinity, i.e. again (\ref{1.4.15}) holds trues with both parts of it
     being equal $\infty$.   
     \end{proof}
   
   In the next section we will employ the following general result which 
   will enable us to subdivide time into small intervals freezing the slow
   variable on each of them so that the estimate (\ref{1.3.8}) of 
   Proposition \ref{prop1.3.4} becomes sufficiently precise and, on the other
   hand, we will not change much the corresponding functionals $S_{0T}$
   appearing in required large deviations estimates. This result is
   certainly not new, it is cited in \cite{Ve3} as a folklore fact and
   a version of it can be found in \cite{Kr}, p.67 but for readers convenience
   we give its proof here.
   
   \begin{lemma}\label{lem1.4.3}
   Let $f=f(t)$ be a measurable function on $\bbR^1$ equal zero outside
   of $[0,T]$ and such that $\int_0^T|f(t)|dt<\infty$. For each positive
   integer $m$ and $c\in[0,T]$ define $f_m(t,c)=f([(t+c)\Del^{-1}]\Del-c)$
   where $\Del=T/m$ and $[\cdot]$ denotes the integral part. Then there
   exists a sequence $m_i\to\infty$ such that for Lebesgue almost all
   $c\in[0,T]$,
   \begin{equation}\label{1.4.21}
   \lim_{i\to\infty}\int_0^T|f(t)-f_{m_i}(t,c)|dt=0.
   \end{equation}
   \end{lemma}
   \begin{proof}
   For each $\del>0$ there exists a $C^1$ function $g$ on $\bbR^1$ equal
   zero outside of $[0,T]$ and such that
   \begin{equation}\label{1.4.22}
   \int_0^T|g(t)-f(t)|dt\leq\del/T.
   \end{equation}
   Define $g_m(t,c)$ as above with $g$ in place of $f$. Then
    \begin{eqnarray}\label{1.4.23}
  &\,\,\,\,\,\,\,\int_0^Tdc\int_0^T|g_m(t,c)-f_m(t,c)|dt\leq
  \int_0^Tdc\sum_{i=0}^\infty|g(i\Del-c)-f(i\Del-c)|\Del\\
  &=\sum_{i=1}^m\Del\int_0^{i\Del}|g(u)-f(u)|du=
  \sum_{i=1}^m\Del\sum_{k=0}^{i-1}\int_{k\Del}^{(k+1)\Del}|g(u)-f(u)|du
  \nonumber\\
  &=\Del\sum_{k=0}^{m-1}(m-k)\int_{k\Del}^{(k+1)\Del}|g(u)-f(u)|du\leq
  T\int_0^T|g(u)-f(u)|du\leq\del.\nonumber
  \end{eqnarray}
  We have also
  \begin{equation}\label{1.4.24}
  \int_0^T\int_0^T|g_m(t,c)-g(t)|dtdc\leq\Del\sup_{0\leq t\leq T}g'(t).
  \end{equation}
  Since
  \[
  |f(t)-f_m(t,c)|\leq|f(t)-g(t)|+|g(t)-g_m(t,c)|+|g_m(t,c)-f_m(t,c)|
  \]
  it follows from (\ref{1.4.22})--(\ref{1.4.24}) that
  \[
  \lim_{m\to\infty}\int_0^T\int_0^T|f(t)-f_m(t,c)|dtdc=0.
  \]
  This together with the Chebyshev inequality and the Borel--Cantelli lemma
  yield (\ref{1.4.21}) for some sequence $m_i\to\infty$ and Lebesgue almost
  all $c\in[0,T]$.
  \end{proof}

   \section{Large deviations: Proof of Theorem 1.2.3}\label{sec1.5}
\setcounter{equation}{0}

\begin{lemma}\label{lem1.5.1} Let $x_i,\tilde x_i\in\cX,\, i=0,1,...,N,$
$0=t_0<t_1<...<t_{N-1}<t_N=T,$ $\Del=\max_{0\leq i\leq N-1}(t_{i+1}-t_i),$
$\xi_i=(x_i-x_{i-1})(t_i-t_{i-1})^{-1}$, $n(t)=\max\{ j\geq 0:\, t\geq t_j\},$
$\psi(t)=\tilde x_{n(t)}$, $v\in\cX\times\bfM$,
\[
\Xi_j^\ve(v,x)=(t_j-t_{j-1})^{-1}\int_{t_{j-1}}^{t_j}B(x,Y^\ve_v(s/\ve))ds,
\]
and for $t\in[0,T]$,
\begin{equation}\label{1.5.1}
Z_{v,x}^{\ve,\psi}(t)=x+\int_0^tB(\psi(s),Y_v^\ve(s/\ve))ds.
\end{equation}
Then
\begin{eqnarray}\label{1.5.2}
&\big\vert\Xi_j^\ve(v,x_{j-1})-(t_j-t_{j-1})^{-1}(Z^\ve_v(t_j)-
Z^\ve_v(t_{j-1}))\big\vert\\
&\leq K\big\vert Z^\ve_v(t_{j-1})-x_{j-1}\big\vert
+\frac 12K^2(t_j-t_{j-1}),\nonumber
\end{eqnarray}
\begin{eqnarray}\label{1.5.3}
&\sup_{0\leq s\leq t}\big\vert Z_{v,x}^{\ve,\psi}(s)-\psi(s)\big\vert\leq
|x-x_0|+\max_{0\leq j\leq n(t)}|x_j-\tilde x_j|\\
&+K\Del+n(t)\Del\max_{1\leq j\leq n(t)}\big\vert\Xi_j^\ve(v,\tilde x_{j-1})
-\xi_j\big\vert\nonumber
\end{eqnarray}
and
\begin{equation}\label{1.5.4}
\sup_{0\leq s\leq t}\big\vert Z^\ve_v(s)-Z_{v,x}^{\ve,\psi}(s)\big\vert\leq
e^{Kt}\big(|\pi_1v-x|+Kt\sup_{0\leq s\leq t}\big\vert Z_{v,x}^{\ve,\psi}(s)-
\psi(s)\big\vert\big)
\end{equation}
where, recall, $Z^\ve_v(s)=X_v^\ve(s/\ve)$ and $\pi_1v=z\in\cX$ if
$v=(z,y)\in\cX\times\bfM$.
\end{lemma}
\begin{proof} By (\ref{1.2.12}),
\begin{eqnarray*}
&\big\vert (t_j-t_{j-1})\Xi_j^\ve(v,x_{j-1})-(Z^\ve_v(t_j)-Z^\ve_v(t_{j-1}))
\big\vert\leq\int_{t_{j-1}}^{t_j}\big\vert B(x_{j-1},Y^\ve_v(\frac s\ve))\\
&-B(Z_v^\ve(s),Y^\ve_v(\frac s\ve))\big\vert ds
\leq\int_{t_{j-1}}^{t_j}\big(\big\vert B(x_{j-1},Y^\ve_v(\frac s\ve))-
B(Z_v^\ve(t_{j-1}),Y^\ve_v(\frac s\ve))\big\vert\\
&+\big\vert B(Z_v^\ve(t_{j-1}),Y^\ve_v(\frac s\ve))-
B(Z_v^\ve(s),Y^\ve_v(\frac s\ve))\big\vert\big) ds
\leq K(t_j-t_{j-1})|Z^\ve_v(t_{j-1})\\
&-x_{j-1}|+K\int_{t_{j-1}}^{t_j}
|Z^\ve_v(s)-Z^\ve_v(t_{j-1})|ds\leq K(t_j-t_{j-1})|Z^\ve_v(t_{j-1})-x_{j-1}|\\
&+K^2\int_{t_{j-1}}^{t_j}(s-t_{j-1})ds\leq
 K(t_j-t_{j-1})|Z^\ve_v(t_{j-1})-x_{j-1}|+\frac 12K^2(t_j-t_{j-1})^2
 \end{eqnarray*}
 and (\ref{1.5.2}) follows.
 
 Observe, that
 \begin{eqnarray*}
 &Z_{v,x}^{\ve,\psi}(s)-\psi(s)=x-x_0+(x_{n(s)}-\tilde x_{n(s)})\\
 &+\sum_{j=1}^{n(s)}(t_j-t_{j-1})\big(\Xi_j^\ve(v,\tilde x_{j-1})-\xi_j\big)
 +\int_{t_{n(s)}}^sB\big(\tilde x_{n(s)},Y^\ve_v(\frac u\ve)\big)du
 \end{eqnarray*}
 and (\ref{1.5.3}) follows in view of (\ref{1.2.15}). Next, by (\ref{1.2.15}),
 \begin{eqnarray*}
 &\big\vert Z_{v}^{\ve}(s)-Z_{v,x}^{\ve,\psi}(s)\big\vert\leq |\pi_1v-x|+
 \int_0^s\big\vert B(Z_{v}^{\ve}(u),Y^\ve_v(\frac u\ve))\\
 &-B(Z_{v,x}^{\ve,\psi}(u),Y^\ve_v(\frac u\ve))\big\vert du
 +\int_0^s\big\vert B(Z_{v,x}^{\ve,\psi}(u),Y^\ve_v(\frac u\ve))-
 B(\psi(u),Y^\ve_v(\frac u\ve))\big\vert du\\
 &\leq |\pi_1v-x|+K\int_0^s\big\vert Z_{v,x}^{\ve,\psi}(u)-\psi(u)
 \big\vert du+K\int_0^s
 \big\vert Z_{v}^{\ve}(u)- Z_{v,x}^{\ve,\psi}(u)\big\vert du
 \end{eqnarray*}
 and (\ref{1.5.4}) follows from the Gronwall inequality \index{Gronwall
 inequality}.
 \end{proof}
 
 For any $x',x''\in\cX$ and $\be,\xi\in\bbR^d$ set
 \[
 L_b(x',x'',\xi)=\sup_{\be\in\bbR^d,|\be|\leq b}\big(\langle\be,\xi\rangle-
 H(x',x'',\be)\big),
 \]
 and $L_b(x,\xi)=L_b(x,x,\xi)$ with $H(x',x'',\be)$ given by (\ref{1.2.10}).
 \begin{proposition}\label{prop1.5.2}
 Let $x_j,t_j,\xi_j,N,\Del,T$ and $\Xi^\ve_j$ be the same as in Lemma 
 \ref{lem1.5.1} and assume that
 \begin{equation}\label{1.5.5}
 \hat\Del=\min_{0\leq i\leq N-1}(t_{i+1}-t_i)\geq\Del/3.
 \end{equation}
 Fix also $\rho>0$ so that Proposition \ref{prop1.3.4} holds true.
 
 (i) There exist $\del_0>0,\ve_0(\Del)>0$ and $C_T(b)>0$ independent of
 $x,x_j,\tilde x_j,\xi_j$ such that if $\del\leq\del_0$ and $\ve\leq
 \ve_0(\Del)$ then for any $b>0$,
 \begin{eqnarray}\label{1.5.6}
 &m\big\{ y\in\cW:\,\max_{1\leq j\leq N}\big\vert\Xi^\ve_j((x,y),\tilde
 x_{j-1})-\xi_j|<\del\big\}\\
 &\leq\exp\big\{-\frac 1\ve\big(\sum_{j=1}^N(t_j-t_{j-1})L_b(\tilde x_{j-1},
 \xi_j)-\eta_{b,T}(\ve,\Del)-C_T(b)(d+\del)\big)\big\}\nonumber
 \end{eqnarray}
 where $d=|x-x_0|+\max_{0\leq j\leq N}|x_j-\tilde x_j|$, $\eta_{b,T}(\ve,
 \Del)$ does not depend on $x,x_j,\tilde x_j,\xi_j$ and
 \begin{equation}\label{1.5.7}
 \lim_{\Del\to 0}\limsup_{\ve\to 0}\eta_{b,T}(\ve,\Del)=0.
 \end{equation}
 In particular, if for each $j=1,...,N$ there exists $\be_j\in\bbR^d$ such
 that
 \begin{equation}\label{1.5.8}
 L(\tilde x_j,\xi_j)=\langle\be_j,\xi_j\rangle-H(\tilde x_j,\be_j)
 \end{equation}
 and
 \begin{equation}\label{1.5.9}
 \max_{1\leq j\leq N}|\be_j|\leq b<\infty
 \end{equation}
 then (\ref{1.5.6}) holds true with $L(\tilde x_j,\xi_j)$ in place of 
 $L_b(\tilde x_j,\xi_j)$, $j=1,...,N$.
 
 (ii) For any $b,\la,\del,q>0$ there exist $\Del_0=\Del_0(b,\la,\del,q)>0$
 and $\ve_0=\ve_0(b,\la,\del,q,\Del)$, the latter depending also on $\Del>0$, 
 such that if $\xi_j$ and $\be_j$ satisfy (\ref{1.5.8}) and (\ref{1.5.9}),
 $\max_{1\leq j\leq N}|\xi_j|\leq q$, $\Del<\Del_0$ and $\ve<\ve_0$ then
 \begin{eqnarray}\label{1.5.10}
 &m\big\{ y\in\cW:\,\max_{1\leq j\leq N}\big\vert\Xi^\ve_j((x,y),\tilde
 x_{j-1})-\xi_j|<\del\big\}\\
 &\geq\exp\big\{-\frac 1\ve\big(\sum_{j=1}^N(t_j-t_{j-1})L(\tilde x_{j-1},
 \xi_j)+\eta_{b,T}(\ve,\Del)+C_T(b)d+\la\big)\big\}\nonumber
 \end{eqnarray}
 with some $C_T(b)>0$ depending only on $b$ and $T$.
 \end{proposition}
 \begin{proof} (i) Assuming that $\rho$ is small and $C\geq 2$ is large so that 
 $C^6\rho$ is still small, we consider for each $x\in\cX_T$ and $y\in\cW$
 closed discs $D_0\in\cD^u_\ve((x,y),\al,\rho,C^3)$
  and $D\in\cD^u_\ve((x,y),\al,\rho,C)$ with $D_0\supset D$. 
  For each small $r\geq 0$ set 
 \[
 D(r)=\{ v\in D_0:\,\inf_{\tilde v\in D}d_{D_0}(v,\tilde v)\leq r\}\,\,
 \mbox{and}
 \]
 \[
 D(-r)=\{ v\in D:\,\inf_{\tilde v\in D_0\setminus D}d_{D_0}(v,\tilde v)\geq r\}.
 \]
 Then $D(r)\cap\partial D_0=\emptyset$ provided $r=r(\rho)< C^2\rho(C-1)$. For
 any pair of compact sets
 $\tilde D\subset\hat D\subset\bbR^d\times\bfM$ and $\vrho>0$ a finite set
 $G\subset\tilde D$ will be called $(s,\vrho,\ve,\tilde D,\hat D)$-separated 
 \index{ $(\cdot,\cdot,\cdot,\cdot,\cdot)$-separated set} if 
 $v_i,v_j\in G,$ $v_i\ne v_j$ implies that 
 $v_i\not\in U^\ve_{\hat D}(s,v_j,\vrho).$ Choose a maximal
   $(t_{n-1}\ve^{-1},C\rho,\ve,\tilde D,D_0)$-separated set $G_{n-1}(\tilde D)$
  in $\tilde D\subset D_0$ (where maximal means that the set cannot be enlarged
  still remaining $(\cdot,\cdot,\cdot,\cdot,\cdot)$-separated). Then
   \[
   \cup_{v\in G_{n-1}(\tilde D)}U^\ve_{D_0}(t_{n-1}\ve^{-1},v,C\rho)\supset
   \tilde D
   \]
   and, by Lemma \ref{lem1.3.2}(i) for small $\ve$, $n>1$, and $v\in\tilde D$,
   \[
   U^\ve_{D_0}(t_{n-1}\ve^{-1},v,C\rho)\subset\tilde D(\ve).
   \]
   Set
   \[
   \Gam_{\tilde D}^j(r)=\big\{ v\in\tilde D:\,\big\vert\Xi_j^\ve(v,
   \tilde x_{j-1})-\xi_j\big \vert<r\big\}
   \]
   and
   \[
   G^\Gam_{n-1}(\tilde r,r)=\big\{ v\in G_{n-1}(D(\tilde r)):\, U^\ve_{D_0}
   (t_{n-1}\ve^{-1},v,C\rho)\cap\big(\bigcap_{j=1}^{n-1} 
   \Gam_{D(\tilde r)}^j(r)\big)\ne\emptyset\big\}
   \]
   assuming that $D(\tilde r)\subset D_0$.
   Then for $\tilde r\geq 0$, $\tilde r< r(\rho)=C^2\rho(C-1)$,
   \begin{eqnarray}\label{1.5.11}
 &m_{D_0}\big\{ v\in D(\tilde r):\,\max_{1\leq j\leq n}
 \big\vert\Xi^\ve_j(v,\tilde x_{j-1})-\xi_j|<r\big\}\\
 &=m_{D_0}\big(\bigcap_{j=1}^{n}
 \Gam_{D(\tilde r)}^j(r)\big)\nonumber\\
 &\leq\sum_{v\in G^\Gam_{n-1}(\tilde r,r)}m_{D_0}\big(
 U^\ve_{D_0}(t_{n-1}\ve^{-1},v,C\rho)\cap
 \Gam_{D(\tilde r)}^n(r)\big).\nonumber
 \end{eqnarray}
 By Lemma \ref{lem1.3.2}(i) if $n>1$ and $\ve$ is small enough then
 $d(v',v)\leq\ve$ for any $v'\in U^\ve_{D_0}(t_{n-1}\ve^{-1},v,C\rho)$
 and using, in addition,  Assumption \ref{ass1.2.2} and the inequality
 (\ref{1.5.3}) we obtain that for any $j\leq n-1$,
 \begin{eqnarray}\label{1.5.12}
 &\big\vert\Xi^\ve_j(v,x)-\Xi^\ve_j(v',x)\big\vert\\
 &\leq K(t_j-t_{j-1})c^{-1}
 C\rho\int_{t_{j-1}}^{t_j}e^{-\frac \ka{2\ve}(t_{n-1}-s)}ds\leq 4K\Del^{-1}
 c^{-1} C\rho\ka^{-1}\ve.\nonumber
 \end{eqnarray}
 Hence, if $v\in G_{n-1}^\Gam(\tilde r,r)$ then for $C_1=8Kc^{-1}C\rho\ka^{-1}$
 and $\tilde r<r(\rho)-\ve$,
  \begin{equation}\label{1.5.13}
  U^\ve_{D_0}(t_{n-1}\ve^{-1},v,C\rho)\subset
  \bigcap_{j=1}^{n-1}\Gam_{D(\tilde r+\ve)}^j(r+C_1\ve\Del^{-1}),
  \end{equation}
  provided $\ve$ is small enough, and so, by (\ref{1.5.3}) and (\ref{1.5.4}),
   \begin{eqnarray}\label{1.5.14}
   &|\pi_1v_{t_{n-1}\ve^{-1}}-\tilde x_{n-1}|\leq d_{n-1}=e^{Kt_{n-1}}
   \sup_{v\in D_0}|\pi_1v-x_0|\\
   &+(e^{Kt_{n-1}}Kt_{n-1}+1)\big(\max_{0\leq j\leq n-1}|x_j-\tilde x_j|+
   K\Del+(n-1)\Del(r+C_1\ve\Del^{-1})\big)\nonumber
   \end{eqnarray}
   where we set $v_s=\Phi^s_\ve v$.
   Since $H(x',x'',\be)$ is (Lipschitz) continuous in $\be$ there exists
   $\be^{(a)}_n(x',x'')\in\bbR^d$ such that
   \begin{equation}\label{1.5.15}
   |\be_n^{(a)}(x',x'')|\leq a\,\,\mbox{and}\,\, L_a(x',x'',\xi_n)=\langle
   \be_n^{(a)}(x',x''),\xi_n\rangle -H(x',x'',\be_n^{(a)}(x',x'')).
   \end{equation}
   Let $v\in G^\Gam_{n-1}(r)$ and $\be^{(a)}_n=\be^{(a)}_n(\pi_1v_{t_{n-1}
   \ve^{-1}},\tilde x_{n-1}).$ Since $H(x',x'',\be)$ is Lipschitz continuous
   (and even $C^1$) in $x'$ and $x''$ (see \cite{Co}) it follows from
   (\ref{1.5.14}) that
   \begin{equation}\label{1.5.16}
   \big\vert H(\pi_1v_{t_{n-1}\ve^{-1}},\tilde x_{n-1},\be_n^{(a)})-
   H(\tilde x_{n-1},\be_n^{(a)})\big\vert\leq C(a)d_{n-1}
   \end{equation}
   where $C(a)>0$ depends only on $a$. Since 
   $U^\ve_{D(\tilde r)}(t_{n-1}\ve^{-1},v,C\rho)\cap\partial D_0=\emptyset$
   provided $v\in G^\Gam_{n-1}(\tilde r,r),\, n>1,\,\tilde r<r(\rho)-\ve$
    we derive from Lemma \ref{lem1.3.2}(iv), Proposition \ref{prop1.3.4}, and 
   Lemma \ref{lem1.4.1}(i) that for such $v,n,\tilde r,\ve$ and any $a>0$, 
   \begin{eqnarray}\label{1.5.17}
   &m_{D_0}\big(U^\ve_{D_0}(t_{n-1}\ve^{-1},v,C\rho)
   \cap\Gam_{D(\tilde r)}^n(r)\big)\\
   &\leq m_{D_0}\big(U^\ve_{D_0}
   (t_{n-1}\ve^{-1},v,C\rho)\big)\nonumber \\
   &\times\exp\big(-\frac {(t_n-t_{n-1})}{\ve}
   (L_a(\tilde x_{n-1},\xi_n)-\tilde\eta_{a,T}(\ve,\Del)-C(a)d_{n-1}-ra)\big)
   \nonumber\end{eqnarray}
   where $\tilde\eta_{a,T}(\ve,\Del)\to 0$ as, first, $\ve\to 0$ and then
   $\Del\to 0$.
   
   Since $U^\ve_{D_0}(t_{n-1}\ve^{-1},v,\frac 12C\rho)$ are 
   disjoint for different $v\in G_{n-1}(D(\tilde r))$ we obtain from 
   (\ref{1.5.13}) and Lemma \ref{lem1.3.2}(iv) that
    \begin{eqnarray}\label{1.5.18}
    &\sum_{v\in G_{n-1}^\Gam(\tilde r,r)}m_{D_0}
    \big(U^\ve_{D_0}
    (t_{n-1}\ve^{-1},v,C\rho)\big)\\
    &\leq c^{-1}_{\frac 12\rho,T}c^{-1}_{\rho,T}
    \sum_{v\in G_{n-1}^\Gam(\tilde r,r)}m_{D_0}
    \big(U^\ve_{D_0}
    (t_{n-1}\ve^{-1},v,\frac 12 C\rho)\big)\nonumber\\
    &\leq c^{-1}_{\frac 12\rho,T}c^{-1}_{\rho,T}
    m_{D_0}\big(\bigcup_{v\in G_{n-1}^\Gam(\tilde r,r)}
    U^\ve_{D_0}(t_{n-1}\ve^{-1},v,\frac 12C\rho)\big)\nonumber \\
    &\leq c^{-1}_{\frac 12\rho,T}c^{-1}_{\rho,T}m_{D_0}\big(
    \bigcap_{j=1}^{n-1}\Gam^j_{D(\tilde r+\ve)}(r+C_1\ve\Del^{-1})\big).
    \nonumber\end{eqnarray}
    Employing (\ref{1.5.11}), (\ref{1.5.17}) and (\ref{1.5.18}) for $n=N,N-1,...,
    2$ with $r=\del+C_1\ve\Del^{-1}, \del+2C_1\ve\Del^{-1},...,
    \del+(N-1)C_1\ve\Del^{-1}$ and $\tilde r=\ve,2\ve,...,(N-1)\ve$, 
    respectively, and using only (\ref{1.5.17}) for $n=1$ we derive that
    \begin{eqnarray}\label{1.5.19}
  &m_{D_0}\big\{ v\in D:\,\max_{1\leq j\leq N}
 \big\vert\Xi^\ve_j(v,\tilde x_{j-1})-\xi_j|<\del\big\}\\  
 &\leq\exp\big\{-\frac 1\ve\big(\sum_{j=1}^N(t_j-t_{j-1})L_a(\tilde x_{j-1},
 \xi_j)-\eta_{a,\rho,T}(\ve,\Del)-C(a,T)(d+\del)\big)\big\}\nonumber
 \end{eqnarray}
 provided $\del+2C_1\ve T\Del^{-2}\leq\rho$ and $\ve T\Del^{-1}<r(\rho)$
 with $\eta_{a,\rho,T}(\ve,\Del)$
 satisfying (\ref{1.5.7}) and with the same $d$ as in (\ref{1.5.6}).
 
 Let $D_x(r,w)$ be a ball on $W^u_x(w,\vrho)$ centered at $w$ and having radius
  $Cr$, $\rho\leq r\leq 2\rho<\vrho$ in the interior metric on $W^u_x(w,\vrho)$ 
  (which, recall, is a semi-invariant extension of the family of local unstable
   manifolds on $\La_x$-- see Section 3 and \cite{Ro}). Then 
   $D_x(r,w)\in\cD^u_\ve((x,w),\al,\rho,C)$ if $C\geq 2$.
 Recall, that if $\rho$ is small enough then the extended local unstable 
 and stable discs $W^u_x(w,r(\rho))$ and $W^s_x(w,r(\rho))$ are defined for all 
 $w\in\cW$ and, in fact, by (\ref{1.2.4}), the compactness arguments and by
 \cite{Ro} such discs can be defined for all $w$ from a small neighborhood 
 $U$ of $\bar\cW$ which is still contained in the basin of attraction of
 each $\La_z$. For each $w\in\bar W$ set
 \[
 Q_x(w,\rho)=\bigcup\big\{ D_x(r(\rho),F_x^r\tilde w):\, |r|\leq C\rho,\,
 \tilde w\in W_x^s(w,C\rho)\big\}
 \]
 and assume that $\rho$ is small enough so that $Q_x(w,\rho)\subset U$. Then
 (\ref{1.5.19}) together with the Fubini theorem yield (\ref{1.5.6}) with the
 box $Q_x(w,\rho)$ in place of the whole $\cW$. Relying on the transversality
 of unstable and weakly stable submanifolds together with compactness arguments
 we conclude that there exist an integer $n_\rho$ depending only on $\rho$
 such that $\cW$ can be covered by $n_\rho$ boxes $Q_x(w_i,\rho),\, i=1,2,...,
 n_\rho$ which yields now (\ref{1.5.6}) in the required form.
 
 (ii) We start proving (\ref{1.5.10}) by using (\ref{1.5.12}) in order to 
 conclude similarly to (\ref{1.5.13}) that if $n>1$,
 $v\in G^\Gam_{n-1}(\tilde r-\ve,r-C_1\ve\Del^{-1})$, and $\tilde v\in
 U^\ve_{D_0}(t_{n-1}\ve^{-1},v,C\rho)$ then $\tilde v\in
 \cap_{j=1}^{n-1}\Gam^j_{D(\tilde r)}(r)$. Hence,
  \begin{eqnarray}\label{1.5.20}
 &m_{D_0}\big\{ v\in D(\tilde r):\,\max_{1\leq j\leq N}
 \big\vert\Xi^\ve_j(v,\tilde x_{j-1})-\xi_j|<r\big\}\\
 &\geq m_{D_0}\bigg(
  \Gam_{D(\tilde r)}^n(r)\cap\big(\bigcup_{v\in G^\Gam_{n-1}
  (\tilde r-\ve,r-C_1\ve\Del^{-1})} U^\ve_{D_0}
  (t_{n-1}\ve^{-1},v,C\rho)\big)\bigg)\nonumber \\
  &\geq\sum_{v\in G^\Gam_{n-1}(\tilde r-\ve,r-C_1\ve\Del^{-1})}m_{D_0}
  \big(\Gam_{D(\tilde r)}^n(r)\cap 
  U^\ve_{D_0}(t_{n-1}\ve^{-1},v,\frac 12C\rho)\big)\nonumber
 \end{eqnarray}
 where the last inequality holds true since $U^\ve_{D_0}
 (t_{n-1}\ve^{-1},v,\frac 12C\rho)$ are disjoint for different
  $v\in G_{n-1}^\Gam(\tilde r,r-C_1\ve\Del^{-1})$. Using (\ref{1.5.16}), 
  Lemma \ref{lem1.3.2}(iv),
  Proposition \ref{prop1.3.4}, and Lemma \ref{lem1.4.1}(iii) we obtain that
  for any $v\in G^\Gam_{n-1}(\tilde r-\ve,r-C_1\ve\Del^{-1})$, $\vsig\le\del$, 
  $\sig>0$ and $b\geq\max_{1\leq j\leq N}|\be_j|$,
   \begin{eqnarray}\label{1.5.21}
   &m_{D_0}\big(\Gam_{D(\tilde r)}^n(r)\cap 
  U^\ve_{D_0}(t_{n-1}\ve^{-1},v,\frac 12C\rho)\big)\\
  &\geq m_{D_0}\big(U^\ve_{D_0}(t_{n-1}\ve^{-1},v,
  \frac 12C\rho)\big)g_{n,b}(\ve,\Del,\vsig,\sig)
  \exp\big(-\frac {(t_n-t_{n-1})}{\ve}L(\tilde x_{n-1},\xi_n)\big)\nonumber
  \end{eqnarray}
  where
  \begin{eqnarray*}
  &g_{n,b}(\ve,\Del,\vsig,\sig)=\exp\bigg(-\frac {(t_n-t_{n-1})}{\ve}
  \big(\tilde\eta_{b,T}(\ve,\Del)+C_T(b)d_{n-1}+\vsig b\big)\bigg)\\
  &\times\bigg(1-\exp\big(-\frac {(t_n-t_{n-1})}{\ve}(d(b)-\tilde\eta_{b,T}
  (\ve,\Del)-\sig b-\sig)\big)\bigg),
  \end{eqnarray*}
  \[
  d(b)=\min_{1\leq j\leq N}\tilde L_b^{\be_j}(\tilde x_{j-1},\cK_{\vsig,C}
  (\xi_j)),\, \tilde L^\be_b(x,\cK)=\inf_{\al\in \cK}\tilde L_b^\be(x,\al),
  \]
  \[
  \cK_{\vsig,C}(\al)=\bar U_C(0)\setminus U_\vsig(\al),\,
  \tilde L_b^\be(x,\al)=L_b(x,\al)-\langle\be,\al\rangle +H(x,\be),\,
  C_T(b)>0,
  \]
  and $\tilde\eta_{b,T}(\ve,\Del)\to 0$ as, first, $\ve\to 0$ and then 
  $\Del\to 0$.
  
  By Lemma \ref{lem1.3.2}(iv) and the definitions of $\Gam^j$ and $G^\Gam$,
   \begin{eqnarray}\label{1.5.22}
    &\sum_{v\in G_{n-1}^\Gam(\tilde r-\ve,r-C_1\ve\Del^{-1})}m_{D_0}
    \big(U^\ve_{D_0}(t_{n-1}\ve^{-1},v,\frac 12C\rho)\big)\\
    &\geq c_{\frac 12\rho,T}c_{\rho,T}
    \sum_{v\in G_{n-1}^\Gam(\tilde r-\ve,r-C_1\ve\Del^{-1})}m_{D_0}
    \big(U^\ve_{D_0}(t_{n-1}\ve^{-1},v,C\rho)\big)\nonumber\\
    &\geq c_{\frac 12\rho,T}c_{\rho,T}
    m_{D_0}\big(\bigcap_{j=1}^{n-1}
    \Gam^j_{D(\tilde r-\ve)}(r-C_1\ve\Del^{-1})\big).
    \nonumber\end{eqnarray}
  Employing (\ref{1.5.20})--(\ref{1.5.22}) for $n=N,N-1,...,2$ with $r=\del$,
  $\del-C_1\ve\Del^{-1},$...,$\del-(N-2)C_1\ve\Del^{-1}$ and $\tilde r=0,
  -\ve,-2\ve,...,-(N-2)\ve$, respectively, and using only (\ref{1.5.21}) for
  $n=1$ we derive that
   \begin{eqnarray}\label{1.5.23}
    &m_{D_0}\big\{ v\in D:\,\max_{1\leq j\leq N}
 \big\vert\Xi^\ve_j(v,\tilde x_{j-1})-\xi_j|<\del\big\}\\  
 &\geq\exp\bigg(-\frac 1\ve\big(\sum_{j=1}^N(t_j-t_{j-1})L(\tilde x_{j-1},
 \xi_j)+C(\rho,\del)\ve\Del^{-1}\big)\bigg)\nonumber \\
 &\times\prod_{n=1}^N g_{n,b}(\ve,\Del,\vsig,\sig)\nonumber
 \end{eqnarray}
 for some $C(\rho,\del)>0$ provided, say, $NC_1\ve\Del^{-1}\leq 
 2TC_1\ve\Del^{-2}\leq\frac {\del}2$ and $T\ve\Del^{-1}<C\rho/2$. 
 Since $H(x,\be)$ is differentiable in $\be$
 (see \cite{Co}) then
 \begin{equation*}
 \tilde L(\tilde x_j,\al)=L(\tilde x_j,\al)-\langle\be_j,\xi_j\rangle +
 H(\tilde x_j,\be_j)>0
 \end{equation*}
 for any $\al\ne\xi_j$ (see Theorems 23.5 and 25.1 in \cite{Roc}), and so
 by the lower semicontinuity of $L(x,\al)$ in $\al$ (and, in fact, also
 in $x$), 
 \[
 \tilde L^{\be_j}(\tilde x_{j-1},\cK_{\vsig,C}(\xi_j))=\inf_{\al\in 
 \cK_{\vsig,C}(\xi_j)}\tilde L^{\be_j}(\tilde x_{j-1},\al)>0.
 \]
 This together with Lemma \ref{lem1.4.2} yield that $d(b)$ appearing in the
 definition of $g_{n,b}(\ve,\Del,\vsig,\sig)$ is positive provided $b$ is
 sufficiently large. In fact, it follows from the lower semicontinuity
 of $L(x,\al)$ that $d(b)$ is bounded away from zero by a positive
 constant independent of $\tilde x_j$ and $\xi_j$, $j=1,...,N$ if these
 points vary over fixed compact sets and (\ref{1.5.8}) together with
 (\ref{1.5.9}) hold true. Now, given $\la>0$ choose, first, sufficiently
 large $b$ as needed and then subsequently choosing small $\sig$ and $\vsig$,
 then small $\Del$, and, finally, small enough $\ve$ we end up with an
 estimate of the form
 \begin{equation}\label{1.5.24}
 g_{n,b}(\ve,\Del,\vsig,\sig)\geq\exp\big(-\frac {(t_n-t_{n-1})}{\ve}
 (\eta_{b,\rho,T}(\ve,T)+C_T(b)d+\la)\big)
 \end{equation}
 where $C_T(b)>0$ and $\eta_{b,\rho,T}(\ve,T)$ satisfies (\ref{1.5.7}).
 Finally, (\ref{1.5.10}) follows from (\ref{1.5.23}), (\ref{1.5.24}) and the
 Fubini theorem (similarly to (i)).
 \end{proof}
 
 Next, we pass directly to the proof of Theorem \ref{thm1.2.3} starting with 
 the lower bound. Some of the details below are borrowed from \cite{Ve3}
 but we believe that our exposition and the way of proof are more precise, 
 complete
 and easier to follow. Assume that $S_{0T}(\gam)<\infty$, and so that $\gam$ is
 absolutely continuous, since there is nothing to prove otherwise. Then by 
 (\ref{1.2.13}), $L(\gam_s,\dot{\gam}_s)<\infty$ for Lebesgue almost all 
 $s\in[0,T].$ By (\ref{1.2.15}) and (\ref{1.3.6}),
 \begin{equation}\label{1.5.25}
 H(x,\be)\leq K|\be|,
 \end{equation}
 and so if $L(\gam_s,\dot{\gam}_s)<\infty$ it follows from (\ref{1.2.12}) that
 $|\dot{\gam}|\leq K$. Suppose that $\cD(L_s)=\{\al:\, L(\gam_s,\al)<\infty\}
 \ne\emptyset$ and let ri$\cD(L_s)$ be the interior of $\cD(L_s)$ in its affine
 hull (see \cite{Roc}). Then either ri$\cD(L_s)\ne\emptyset$ or $\cD(L_s)$ (by
 its convexity) consists of one point and recall that $\dot{\gam}_s\in\cD(L_s)$
 for Lebesgue almost all $s\in[0,T]$.  By (\ref{1.2.10}) and (\ref{1.5.25}),
\begin{equation}\label{1.5.26}
0=H(\gam_s,0)=\inf_{\al\in\bbR^d}L(\gam_s,\al).
\end{equation}
This together with the nonnegativity and lower semi-continuity of 
$L(\gam_s,\cdot)$ yield that there exists $\hat\al_s$ such that 
$L(\gam_s,\hat\al_s)=0$ and by a version of the measurable selection (of the
implicit function) theorem (see \cite{CV}, Theorem III.38), $\hat\al_s$ can be
chosen to depend measurably in $s\in[0,T]$. Of course, if ri$\cD(L_s)=\emptyset$
then $\cD(L_s)$ contains only $\hat\al_s$ and in this case $\hat\al_s=
\dot{\gam}_s$ for Lebesgue almost all $s\in[0,T]$. Taking $\al_s=\hat\al_s$ 
and $\be_s=0$ we obtain
\begin{equation}\label{1.5.27}
L(\gam_s,\al_s)=\langle\be_s,\al_s\rangle-H(\gam_s,\be_s).
\end{equation}

Observe that $\ell(s,\al)=L(\gam_s,\al)$ is measurable as a function of $s$
and $\al$ since it is obtained via (\ref{1.2.12}) as a supremum in one argument 
of a family of continuous functions, and so this supremum can be taken there
over a countable dense set of $\be$'s. Hence, the set $A=\{ (s,\al):\,
s\in[0,T],\,\al\in\cD(L_s)\}=\ell^{-1}[0,\infty)$ is measurable, and so the set
$B=A\setminus\{(s,\dot{\gam_s}),\, s\in[0,T]\}$ is measurable, as well. Its
projection $V=\{ s\in[0,T]:\, (s,\al)\in B\,\mbox{for some}\,\al\in\bbR^d\}$
on the first component of the product space is also measurable and $V$ is the
set of $s\in[0,T]$ such that $\cD(L_s)$ contains more than one point. Employing
Theorem III.22 from \cite{CV} we select $\bar\al_s\in\bbR^d$ measurably in
$s\in V$ and such that $(s,\bar\al_s)\in B$. By convexity and lower 
semicontinuity of $L(\gam_s,\cdot)$ it follows from Corollary 7.5.1 in
\cite{Roc} that 
\begin{equation}\label{1.5.28}
L(\gam_s,\dot{\gam}_s)=\lim_{p\uparrow\infty}L(\gam_s,\al^{(p)}_s)\,\,
\,\mbox{where}\,\,\,\al^{(p)}_s=(1-p^{-1})\dot{\gam}_s+p^{-1}\bar\al_s.
\end{equation}
For each $\del>0$ set 
\[
n_\del(s)=\min\{ n\in\bbN:\, |L(\gam_s,\dot{\gam}_s)-
L(\gam_s,\al^{(n)}_s)|+|\dot{\gam}_s-\al^{(n)}_s|<\del\}.
\]
 Then, clearly, $n_\del(s)$ is a measurable function of $s$, and so 
 $\al_s=\al_s^{(\del)}=\al^{(n_\del(s))}_s$ and $L(\gam_s,\al_s)$ are
  measurable in $s$, as well. By Theorems 23.4 and 23.5 from \cite{Roc}
  for each $\al_s=\al_s^{(\del)}$ there exists $\be_s=\be_s^{(\del)}\in\bbR^d$
  such that (\ref{1.5.27}) holds true. Given $\del',\la>0$ take $\del=\min(\del',
  \la/3)$ and for $s\in[0,T]\setminus V$ set $\al_s=\hat\al_s$. Then
\begin{equation}\label{1.5.29}
\int_0^T\big\vert L(\gam_s,\dot{\gam}_s)-L(\gam_s,\al_s)\big\vert ds<\la/3\,\,\,
\mbox{and}\,\,\int_0^T|\dot{\gam}_s-\al_s|ds<\del'.
\end{equation}

For each $b>0$ set $\al_s^b=\al_s$ if the corresponding $\be_s$ in (\ref{1.5.27})
satisfies $|\be_s|\leq b$ and $\al_s^b=\hat\al_s$, otherwise. Note, that 
(\ref{1.5.27}) remains true with $\al_s^b$ in place of $\al_s$ with $\be_s=0$ if 
$\al_s^b=\hat\al_s$. As observed above $|\al|\leq K$ whenever $L(z,\al)<\infty$,
and so $|\hat\al_s|\leq K$ for Lebesgue almost all $s\in[0,T]$. We recall also
that $|\dot{\gam}_s-\al_s|<\del$ and $\dot{\gam}_s\leq K$ for Lebesgue
almost all $s\in[0,T]$. Since $S_{0T}(\gam)<\infty$, $|L(\gam_s,\dot{\gam}_s)-
L(\gam_s,\al_s)|<\del$, and $L(\gam_s,\al_s^b)\uparrow L(\gam_s,\al_s)$ as
$b\uparrow\infty$ for Lebesgue almost all $s\in[0,T]$, we conclude from 
(\ref{1.5.29}) and the above observations that for $b$ large enough
\begin{equation}\label{1.5.30}
\int_0^T\big\vert L(\gam_s,\al_s)-L(\gam_s,\al_s^b)\big\vert ds<\la/3\,
\mbox{and}\,\int_0^T|\al_s-\al_s^b|ds<\del'.
\end{equation}
Next, we apply Lemma \ref{lem1.4.3} to conclude that there exists a sequence
$m_j\to\infty$ such that for each $\Del_j=T/m_j$ and Lebesgue almost all $c\in[0,T)$,
\begin{equation}\label{1.5.31}
\int_0^T\big\vert L(\gam_s,\al_s^b)-L(\gam_{q_j(s,c)},\al_{q_j(s,c)}^b)\big
\vert ds<\la/3\,\mbox{and}\,\int_0^T|\al_s^b-\al_{q_j(s,c)}^b|ds<\del'.
\end{equation}
where $q_j(s,c)=[(s+c)\Del^{-1}_j]\Del_j-c$, $[\cdot]$ denotes the integral part
and we assume $L(\gam_s,\al_s^b)=0$ and $\al_s^b=0$ if $s<0$.

Choose $c=c_j\in[\frac 13\Del_j,\frac 23\del_j]$ and set $\hat\gam_s=x+\int_0^s
\al^b_{q_j(u,c)}du$, $\psi_s=\gam_{q_j(s,c)}$ where $\gam_u=\gam_0$ if $u<0$,
$x_0=\tilde x_0=x$, $x_N=\hat\gam_T,$ $\tilde x_N=\gam_T$ and $x_k=\hat\gam_
{k\Del_j-c}$, $\tilde x_k=\gam_{k\Del_j-c}$ for $k=1,...,N-1$ and $\xi_k=
\al^b_{(k-1)\Del_j-c}$ for $k=1,2,...,N$ where $N=\min\{ k:k\Del_j-c>T\}$.
Since $|\dot{\gam}_s|\leq K$ for Lebesgue almost all $s\in[0,T]$ then
$\bfr_{0T}(\gam,\psi)\leq K\Del_j$ and, in addition, $\bfr_{0T}(\gam,\hat\gam)
\leq 3\del'$ by (\ref{1.5.29})--(\ref{1.5.31}). This together with (\ref{1.5.3}) and
(\ref{1.5.4}) yield that for $v=(x,y)$,
\begin{eqnarray}\label{1.5.32}
&\,\,\,\,\,\,\,\bfr_{0T}(Z^\ve_v,\gam)\leq\bfr_{0T}(Z^\ve_v,\psi)+
K\Del_j\leq(KTe^{KT}+1)\bfr_{0T}(Z^{\ve,\psi}_v,\psi)+K\Del_j\\
&\leq(KTe^{KT}+1)\big(3\del'+K\Del_j+(T+1)\max_{1\leq k\leq N}\big\vert
\Xi^\ve_k(v,\tilde x_{k-1})-\xi_k\big\vert\big)+K\Del_j \nonumber
\end{eqnarray}
provided $\Del_j\leq 1$ where $Z_v^{\ve,\psi}$ and $\Xi^\ve_k(v,x)$ are the same
as in Lemma \ref{lem1.5.1}, the latter is defined with $t_k=k\Del_j-c$, $k=1,...,N-1$
and $t_N=T$. Choose $\del'$ so small and $m_j$ so large that
\[
(KTe^{KT}+1)\big(3\del'+K\Del_j+(T+1)\del'\big)+K\Del_j<\del
\]
then by (\ref{1.5.32}),
\begin{equation}\label{1.5.33}
\big\{ y\in\cW:\,\bfr_{0T}(Z^\ve_{x,y},\gam)<\del\big\}\supset\big\{ y\in\cW:\,
\max_{1\leq k\leq N}\big\vert\Xi^\ve_k(v,\tilde x_{k-1})-\xi_k\big\vert 
<\del'\big\}.
\end{equation}
By (\ref{1.5.29})--(\ref{1.5.31}),
\begin{equation}\label{1.5.34}
\sum_{k=1}^N(t_k-t_{k-1})L(\tilde x_{k-1},\xi_k)\leq S_{0T}(\gam)+\la
\end{equation}
and by the construction above the conditions of the assertion (ii) of Proposition
\ref{prop1.5.2} hold true, so choosing $m_j$ sufficiently large we derive 
(\ref{1.2.16})
(with $2\la$ in place of $\la$) from (\ref{1.5.10}), (\ref{1.5.33}) and
 (\ref{1.5.34}) provided $\ve$ is small enough.

Next, we pass to the proof of the upper bound (\ref{1.2.17}). Assume that 
(\ref{1.2.17}) is not true, i.e. there exist $a,\la,\del>0$ and $x\in\cX_T$ 
such that for some sequence $\ve_k\to 0$ as $k\to\infty$,
\begin{equation}\label{1.5.35}
m\big\{ y\in\cW:\,\bfr_{0T}\big( Z^{\ve_k}_{x,y},\Psi^a_{0T}(x)\big)\geq 3
\del\big\}>\exp\big(-\frac 1{\ve_k}(a-\la)\big).
\end{equation}
Since $\| B(x,y)\|\leq K$ by (\ref{1.2.15}) all paths of $Z^{\ve}_{x,y}(t),\, 
t\in[0,T]$ and of $Z^{\ve,\psi}_{v,x}(t),\, t\in[0,T]$ given by (\ref{1.5.1}) 
(the latter for any measurable $\psi$) belong to a compact set 
$\tilde \cK^x\subset C_{0T}$ which consists of curves starting at $x$ and 
satisfying the Lipschitz condition with the constant
$K$. Let $\tilde U^x_{\rho}$ denotes the open $\rho$-neighborhood of the 
compact set
$\Psi^a_{0T}(x)$ and $\cK^x_\rho=\tilde \cK^x\setminus\tilde U^x_{\rho}$ . 
For any small
$\del'>0$ choose a $\del'$-net $\gam_1,...,\gam_n$ in $\cK^x_{2\del}$ where 
$n=n(\del')$. Since
\begin{equation*}
\big\{ y\in\cW:\,\bfr_{0T}\big( Z^{\ve_k}_{x,y},\Psi^a_{0T}(x)\big)\geq 3\del
\big\}\subset\bigcup_{n\geq j\geq 1}\big\{ y\in\cW:\,\bfr_{0T}
\big( Z^{\ve_k}_{x,y},\gam_j)\leq\del'\big\}
\end{equation*}
then there exists $j$ and a subsequence of $\{\ve_k\}$, for which we use the 
same notation, such that
\begin{equation}\label{1.5.36}
m\big\{ y\in\cW:\,\bfr_{0T}( Z^{\ve_k}_{x,y},\gam_j)\leq\del'\big\}>
n^{-1}\exp\big(-\frac 1{\ve_k}(a-\la)\big).
\end{equation}
Denote such $\gam_j$ by $\gam^{\del'}$, choose a sequence $\del_l\to 0$ and set
$\gam^{(l)}=\gam^{\del_l}$. Since $\cK^x_{2\del}$ is compact there exists a 
subsequence $\gam^{(l_j)}$ converging in $C_{0T}$ to $\hat\gam\in \cK^x_{2\del}$ 
which together with (\ref{1.5.36}) yield
\begin{equation}\label{1.5.37}
\limsup_{\ve\to 0}\ve\ln m\big\{y\in\cW:\,\bfr_{0T}(Z^{\ve}_{x,y},\hat\gam)\leq
\del'\big\}>-a+\la
\end{equation}
for all $\del'>0$.

We claim that (\ref{1.5.37}) contradicts (\ref{1.5.2}) and the assertion (i) of
Proposition \ref{prop1.5.2}. Indeed, set
\[
S^\psi_{b,0T}(\gam)=\int_0^TL_b(\psi(s),\dot{\gam}(s))ds\,\,\mbox{and}\,\,
S_{b,0T}(\gam)=S^\gam_{b,0T}(\gam).
\]
By the monotone convergence theorem
\begin{equation}\label{1.5.38}
S^\psi_{b,0T}(\gam)\uparrow S^\psi_{0T}(\gam)\,\,\mbox{and}\,\,
S_{b,0T}(\gam)\uparrow S_{0T}(\gam)\,\,\mbox{as}\,\, b\uparrow\infty.
\end{equation}
Similarly to our remark (before Assumption \ref{ass1.2.2}) in Section 
\ref{sec1.2}
it follows from the results of Section 9.1 of \cite{IT} that the functionals
$S^\psi_{b,0T}(\gam), S^\psi_{0T}(\gam)$ and $S_{b,0T}(\gam),S_{0T}(\gam)$ are
lower semicontinuous in $\psi$ and $\gam$ (see also Section 7.5 in \cite{FW}).
This together with (\ref{1.5.38}) enable us to apply Lemma \ref{lem1.4.2} in 
order to conclude that
\begin{equation}\label{1.5.39}
\lim_{b\to\infty}S_{b,0T}(\cK_\del^x)=S_{0T}(\cK_\del^x)=\inf_{\gam\in 
\cK^x_\del}S_{0T}(\gam)>a
\end{equation}
where $S_{b,0T}(\cK_\del^x)=\inf_{\gam\in \cK^x_\del}S_{b,0T}(\gam)$. The last
inequality in (\ref{1.5.39}) follows from the lower semicontinuity of $S_{0T}$.
Thus we can and do choose $b>0$ such that
\begin{equation}\label{1.5.40}
S_{b,0T}(\cK_\del^x)>a-\la/8.
\end{equation}

By the lower semicontinuity of $S_{b,0T}^\psi(\gam)$ in $\psi$ there exists a
function $\del_\la(\gam)>0$ on $\cK^x_\del$ such that for each $\gam\in 
\cK^x_\del$,
\begin{equation}\label{1.5.41}
S_{b,0T}^\psi(\gam)>a-\la/4\,\,\mbox{provided}\,\,\bfr_{0T}(\gam,\psi)<\del_\la
(\gam).
\end{equation}
Next, we restrict the set of functions $\psi$ to make it compact. Namely, we 
allow from now on only functions $\psi$ for which there exists $\gam\in 
\cK^x_\del$ such that either $\psi\equiv\gam$ or $\psi(t)=\gam(kT/m)$ for
 $t\in[kT/m,(k+1)T/m)$, $k=0,1,...,m-1$ and $\psi(T)=\gam(T)$ where $m$ is
 a positive integer. It is easy to see that the set of such functions $\psi$ 
 is compact with respect to the uniform convergence topology in $C_{0T}$ and 
 it follows that $\del_\la(\gam)$ in (\ref{1.5.41}) constructed with such $\psi$
  in mind is lower semicontinuous in $\gam$. Hence
\begin{equation}\label{1.5.42}
\del_\la=\inf_{\gam\in \cK^x_\del}\del_\la(\gam)>0.
\end{equation}

Now take $\hat\gam$ satisfying (\ref{1.5.37}) and for any integer $m\geq 1$ set 
$\Del=\Del_m=T/m$, $x_k=x_k^{(m)}=\hat\gam(k\Del)$, $k=0,1,...,m$ and 
$\xi_k=\xi_k^{(m)}=\Del^{-1}\big(\hat\gam(k\Del)-\hat\gam((k-1)\Del)\big),$ 
$k=1,...,m$. Define a piecewise linear 
$\chi_m$ and a piecewise constant $\psi_m$ by
\begin{equation}\label{1.5.43}
\chi_m(t)=x_k+\xi_k\Del\,\,\mbox{and}\,\,\psi_k(t)=x_k\,\,\mbox{for}\,\, 
t\in[k\Del,(k+1)\Del)
\end{equation}
and $k=0,1,...,m-1$ with $\chi_m(T)=\psi_m(T)=\hat\gam(T)$. Since $\hat\gam$ is
 Lipschitz continuous with the constant $K$ then
\begin{equation}\label{1.5.44}
\bfr_{0T}(\chi_m,\psi_m)\leq K\Del\,\,\mbox{and}\,\,\bfr_{0T}(\hat\gam,\psi_m)
\leq K\Del.
\end{equation}
If $m$ is large enough and $\ve>0$ is sufficiently small then
\begin{equation}\label{1.5.45}
\Del<K^{-1}\min(\del/2,\del_\la)\,\,\mbox{and}\,\,\eta_{b,T}(\ve,\Del)<\la/8
\end{equation}
where $\eta_{b,T}(\ve,\Del)$ is the same as in (\ref{1.5.6}). Since 
$\hat\gam\in\cK^x_{2\del}$
it follows from (\ref{1.5.44}) and (\ref{1.5.45}) that $\chi_m\in \cK^x_{\del}$
 and by (\ref{1.5.41})
and the first inequality in (\ref{1.5.45}) we obtain that
\begin{equation}\label{1.5.46}
S^{\psi_m}_{b,0T}(\chi_m)=\Del\sum_{k=0}^{m-1}L_b(x_k,\xi_k)>a-\frac \la{4}.
\end{equation}
Hence, by (\ref{1.5.6}) and the second inequality in (\ref{1.5.45}) for all 
$\ve$ small enough,
\begin{equation}\label{1.5.47}
 m\big\{ y\in\cW:\,\max_{1\leq k\leq m}\big\vert\Xi^\ve_k((x,y),
 x_{k-1})-\xi_k|<\rho\big\}\leq e^{-\frac 1\ve(a-\la/2)}
 \end{equation}
provided $C_T(b)\rho<\la/8$ (taking into account that $x_0=x$). By (\ref{1.5.2}) 
and the definition of vectors $\xi_k$ for any $v\in\cW$,
\begin{eqnarray}\label{1.5.48}
&\big\vert\Xi^\ve_k(v,x_{k-1})-\xi_k\big\vert\leq\big\vert\Xi^\ve_k(v,x_{k-1})
-\Del^{-1}
\big(Z_v^\ve(k\Del)-Z_v^\ve((k-1)\Del)\big)\big\vert\\
&+2\Del^{-1}\bfr_{0T}(Z^\ve_v,\hat\gam)\leq(K+2\Del^{-1})\bfr_{0T}(Z^\ve_v,
\hat\gam)+\frac 12K^2\Del.\nonumber
\end{eqnarray}
Therefore,
\begin{eqnarray}\label{1.5.49}
&\big\{ y\in\cW:\,\bfr_{0T}(Z^\ve_{x,y},\hat\gam)\leq\del'\big\}\subset\big
\{ y\in\cW:\,\\
&\max_{1\leq k\leq m}\big\vert\Xi^\ve_k((x,y), x_{k-1})-\xi_k|\leq (K+2\Del^{-1})\del'+
\frac 12K^2\Del\big\}.\nonumber
\end{eqnarray}
Choosing, first, $m$ large enough so that $\Del$ satisfies (\ref{1.5.45}) with
 all sufficiently
small $\ve$ and also that $8C_T(b)K^2\Del<\la$, and then choosing $\del'$ so 
small that $16C_T(b)(K+2\Del^{-1})\del'<\del$, we conclude that (\ref{1.5.47}) 
together with (\ref{1.5.49})
contradicts (\ref{1.5.37}), and so the upper bound (\ref{1.2.17}) holds true. 
Since $S_{0T}(\gam)=0$ if and only if $\gam=\gam^u$ satisfying (\ref{1.2.14})
the estimate (\ref{1.2.18}) follows from (\ref{1.2.17}) and the lower 
semicontinuity of the functional $S_{0T}$, completing the proof of Theorem 
\ref{thm1.2.3}.\qed

\section{Further properties of $S$-functionals}\label{sec1.6}
\setcounter{equation}{0}

In this section we study essential properties of the functionals $S_{0T}$
which will be needed in the proofs of Theorems \ref{thm1.2.5} and 
\ref{thm1.2.7}
in the next sections. We will start with the following general fact which do 
not require specific conditions of Theorems \ref{thm1.2.5} and \ref{thm1.2.7}.
\begin{lemma}\label{lem1.6.1} There exists $r>0$ such that if $x\in\bar\cX$
then any $\mu_{xx}$ from the space $\cM_x$ of $F^t_x$-invariant probability
measures on $\La_x$ can be included into a weakly continuous in $z$ family
$\mu_{xz}\in\cM_z,\,|z-x|<r$ (considered in the space of probability measures
on $\overline {\cW}$) for which $\bar B_{\mu_{xz}}(z)=\int B(x,y)d\mu_{xz}(y)$ 
is $C^1$ in $z$ and the entropy $h_{\mu_{xz}}(F^1_z)$ is continuous in 
$z$ as $|z-x|<r$. Furthermore, there exists $C>0$ such that
\begin{equation}\label{1.6.1}
|\bar B_{\mu_{xz_1}}(z_1)-\bar B_{\mu_{xz_2}}(z_2)|<C|z_1-z_2|\,\,
\mbox{whenever}\,\,x,z_1,z_2\in\bar\cX,\, |z_i-x|<r,\, i=1,2
\end{equation}
and for any $\al>0$ there exists $\be>0$ such that if $x,z\in\bar\cX,\, 
|z-x|<\be$ then
\begin{equation}\label{1.6.2}
|h_{\mu_{xz}}(F^1_z)-h_{\mu_{xx}}(F^1_x)|<\al\,\,\mbox{and}\,\,
|I_z(\mu_{xz})-I_x(\mu_{xx})|<\al.
\end{equation}
\end{lemma}
\begin{proof} The following argument (whose ingredients appear already in
\cite{LMM}, \cite{KKPW}, and \cite{Co}) was 
indicated to me by A.Katok. If $r$ is small enough the structural stability
theorem for Axiom A flows obtained in \cite{Ro} can be applied in order to 
compare $F^1_x$ and $F^1_z$ but here we will need its more recent form 
derived in \cite{LMM}, \cite{KKPW}, and \cite{Co} which yields a 
homeomorphism $u_{xz}:\La_x\to\La_z$ and a continuous function $c_{xz}$ on 
$\La_z$ both with $C^1$ dependence on $z$ and such that the conjugate flow
$\tilde F^t_z=u_{xz}F^t_xu^{-1}_{xz}$ satisfies
\[
\frac {d\tilde F^t_zy}{dt}=c_{xz}(\tilde F^t_zy)b(z,\tilde F^t_zy)
\]
where $u_{xx}$ is the identity map on $\La_x$ and $c_{xx}\equiv 1$. 
By the standard direct verification we see that $\mu=u_{xz}\mu_{xx}$  is an 
$\tilde F_z$-invariant probability measure. It is known (see, for instance,
 \cite{To}, Theorem 4.2) that then the probability measure $\mu_{xz}$ on 
 $\La_z$ defined by its Radon--Nikodim derivative
\[
\frac {d\mu_{xz}}{d\mu}(y)=c_{xz}(y)\big(\int_{\La_z}c_{xz}d\mu\big)^{-1}
\]
is $F^t_z$-invariant. In our case this can be seen easily since for any
$C^1$ function $q$ on $\La_z$,
\begin{eqnarray*}
&\frac d{dt}\int_{\La_z}q\circ F^t_zd\mu_{xz}\big\vert_{t=0}=
\big(\int_{\La_z}c_{xz}d\mu\big)^{-1}\int_{\La_z}c_{xz}(b(z,\cdot),
\nabla q)d\mu\\
&=\big(\int_{\La_z}c_{xz}d\mu\big)^{-1}
\frac d{dt}\int_{\La_z}q\circ\tilde F^t_zd\mu\big\vert_{t=0}=0
\end{eqnarray*}
where the last equality holds true by $\tilde F^t_z$-invariance of $\mu$.

Now
\begin{eqnarray*}
&\bar B_{\mu_{xz}}(z)=\int_{\La_z}B(z,y)d\mu_{xz}(y)=
\big(\int_{\La_z}c_{xz}d\mu\big)^{-1}\int_{\La_z}B(z,y)c_{xz}(y)d\mu(y)\\
&\big(\int_{\La_x}c_{xz}(u_{xz}y)d\mu_{xx}(y)\big)^{-1}
\int_{\La_x}B(z,u_{xz}y)c_{xz}(u_{xz}y)d\mu_{xx}(y).
\end{eqnarray*}
This together with (\ref{1.2.15}) yield the differentiability of 
$\bar B_{\mu_{xz}}(z)$ in $z$ taking into account that $c_{xz}$ and $u_{xz}$ 
are $C^1$ in $z$ (see \cite{Co}) and since the proof of this fact relies on a
version of the implicit function theorem (see  \cite{KKPW}) which provides 
derivatives in $z$ uniformly in $x\in\cX$ whenever $|z-x|<r$ and $r$ is small 
enough we derive also (\ref{1.6.1}). 
Next, clearly, 
$h_{\mu_{xx}}(F^1_x)=h_\mu(\tilde F^1_z)$. If we knew that $\mu_{xx}$ were
ergodic then, of course, $\mu$ would be ergodic, as well, and it would follow 
from Theorem 10.1 in \cite{To} that 
\[
h_{\mu_{xz}}(F^1_z)=h_{\mu_{xx}}(F^1_x)\int_{\La_x}c_{xz}(u_{xz}y)d\mu_{xx}(y)
\]
which would yield the differentiability of $h_{\mu_{xz}}(F^1_z)$ in $z$. 
In the general case we obtain from \cite{To} that
\[
h_{\mu_{xz}}(F^1_z)\inf_{y\in\La_z}c_{xz}(y)\leq h_\mu(\tilde F^1_z)\leq
h_{\mu_{xz}}(F^1_z)\sup_{y\in\La_z}c_{xz}(y),
\]
and so
\[
\big\vert h_{\mu_{xz}}(F^1_z)-h_{\mu_{xx}}(F^1_x)\big\vert\leq h_{\mu_{xz}}
(F^1_z)\max\big(|\sup_{y\in\La_z}c_{xz}(y)-1|,\,|1-\inf_{y\in\La_z}c_{xz}(y)|
\big).
\]
Since by Ruelle's inequality (see, for instance \cite{KH}),
\[
h_{\mu_{xz}}(F^1_z)\leq\sup_{y\in\La_z}|\vf^u_z(y)|
\]
we derive both the continuity of $h_{\mu_{xz}}(F^1_z)$ in $z$ and the
first part of (\ref{1.6.2}). The second part of (\ref{1.6.2}) follows from its 
first part in view of (\ref{1.2.8}) taking into account that the function 
$\vf^u_x(y)$ defined by (\ref{1.2.5}) is H\" older continuous in $y$ and 
uniformly Lipschitz continuous (even $C^1$) in $x$ (see \cite{Co}) and 
that $B(x,y)$ is Lipschitz continuous in both variables (see (\ref{1.2.15})).
\end{proof}

The following result gives, in particular, sufficient conditions for a set 
to be an $S$-compact.
\begin{lemma}\label{lem1.6.2} (i) There exists $C>0$ and for each $x\in\cX$ 
where the vector field $B$ is complete there exists $r=r(x)>0$ such that if 
$|z_1-x|<r$ and $|z_2-x|<r$ then we can construct $\gam\in C_{0t}$ with 
$t\leq C|z_1-z_2|$ satisfying 
\[
\gam_0=z_1,\,\,\gam_{t}=z_2\,\,\mbox{and}\,\, S_{0t}(\gam)\leq C|z_1-z_2|.
\]
It follows that $R(\tilde z,z)$ and $R(z,\tilde z)$ are locally Lipschitz
continuous in $z$
belonging to the open $r$-neighborhood of $x$ when $\tilde z$ is fixed.

(ii) Let $\cO\subset\cX$ be a compact $\Pi^t$-invariant 
set which either contains a dense in $\cO$ orbit of $\Pi^t$ or $R(x,z)=0$ for
any pair $x,z\in\cO$. Suppose that $B$ is complete at each point of $\cO$. 
Then $\cO$ is an $S$-compact.

(iii) Assume that for any $\eta>0$ there exists $T(\eta)>0$ such that for
each $x\in\cO$ its orbit $\{\Pi^tx,\, t\in[0,T(\eta)]\}$ of length $T(\eta)$
forms an $\eta$-net in $\cO$ or, equivalently, that $\Pi^t$ is a minimal flow 
on $\cO$. Suppose that $B$ is complete at a point of $\cO$. Then $\cO$ is an
 $S$-compact.
\end{lemma}
\begin{proof} (i) Fix some $x\in\cX$. In view of the ergodic decomposition (see,
for instance, \cite{KH}) any $\mu\in\cM_x$ can be represented as an integral
over the space of ergodic measures from $\cM_x$. Using the specification
(see \cite{Bo1} and \cite{Fra}) any ergodic $\mu\in\cM_x$ can be approximated
(in the weak sense) by $F^t_x$-invariant measures sitting on its periodic
orbits, i.e. by measures of the form $\mu_\vf=\frac 1{t_\vf}\int_0^{t_\vf}
\del_{F^s_xy}ds$ where $\vf=\{ F^s_xy,\, 0\leq s\leq t_\vf\}$, $F_x^{t_\vf}y
=y$ is a periodic orbit of $F^t_x$ with a period $t_\vf$. This is done in a
standard way by choosing a generic point of an ergodic measure $\mu$, i.e.
a point $w$ which satisfies $\lim_{t\to\infty}t^{-1}\int_0^tg(F^s_xw)ds=
\int gd\mu$ for any continuous function $g$ on $\La_x$, and then approximating
the orbit of $w$ by periodic orbits of $F^t_x$ using the specification
theorem (see Theorem 3.8 in \cite{Bo1}). It is well known (see \cite{Bo1})
that there are countably many periodic orbits of $F^t_x$ which together
with the above discussion yield that the closed convex hull $\Gam_x^{(0)}$
of the set $\{\bar B_{\mu_\vf}(x):\,\vf$ is a periodic orbit of $F^t_x\}
\subset\bbR^d$ coincides with $\Gam_x=\{\bar B_\mu(x):\,\mu\in\cM_x\}$.

Now assume that $B$ is complete at $x$. Then $\{\al\Gam_x,\,\al\in[0,1]\}=
\{\al\Gam_x^{(0)},\,\al\in[0,1]\}$ contains an open neighborhood of 0 in
$\bbR^d$. But then we can find a simplex $\Del_x$ with vertices in 
$\Gam_x^{(0)}$ such that $\{\al\Del_x,\,\al\in[0,1]\}$ contains an open
neighborhood of 0 in $\bbR^d$ and for some periodic orbits $\vf_1,...,
\vf_k$ of $F^t_x,\, k\geq d+1$,
\[
\Del_x=\{\sum^k_{i=1}\la_i\bar B_{\mu^{(i)}}(x):\,\sum^k_{i=1}\la_i=1,\,
\la_i\geq 0\,\forall i\}
\]
where we denote $\mu^{(i)}=\mu_{\vf_i}$. By compactness of $\Del_x$ it
follows also that
\[
\mbox{dist}(\Del_x,0)=d_x>0.
\]
Now, set $\mu^{(i)}_{xx}=\mu^{(i)}$, $i=1,...,k$ and include each 
$\mu^{(i)}_{xx}$ into the weakly continuous in $z$ families $\mu_{xz}^{(i)}$
constructed in Lemma \ref{lem1.6.1} for $z$ in some neighborhood of $x$. If
$|z-x|\leq r(x)$ and $r(x)$ is small enough each simplex 
\[
\Del_z=\{\sum^k_{i=1}\la_i\bar B_{\mu^{(i)}_{xz}}(z):\,\sum^k_{i=1}\la_i=1,\,
\la_i\geq 0\,\forall i\}
\]
intersects and not at 0 with any ray emanating from $0\in\bbR^d$ or, in other
words, $\{\al\Del_z,\,\al\in[0,1]\}$ contains an open neighborhood of 0 in 
$\bbR^d$ and, moreover,
\[
\mbox{dist}(0,\Del_z)\geq\frac 12d_x.
\]
Since all $\sum_{i=1}^k\la_i\mu^{(i)}_{xz}$ are $F^t_z$-invariant probability 
measures provided $\sum^k_{i=1}\la_i=1,\,\la_i\geq 0$ we conclude that for
any $z$ in the $r(x)$-neighborhood of $x$ and any vector $\xi$ there exists
an $F^t_z$-invariant probability measure $\mu_\xi$ such that 
$\bar B_{\mu_{\xi}}(z)$ has the same direction as $\xi$ and
\[
K\geq |\bar B_{\mu_\xi}(z)|\geq\frac 12d_x
\]
where $K$ is the same as in (\ref{1.2.15}). It follows that any two points
$z_1$ and $z_2$ from the open $r(x)$-neighborhood of $x$ can be connected
by a curve $\gam$ lying on the interval connecting $z_1$ and $z_2$ with
$K\geq |\dot {\gam}_s^{(1)}|\geq\frac 12d_x$, i.e. $\gam_0=z_1,\,\gam_t=z_2$
with some $t\in[K^{-1}|z_1-z_2|,2d^{-1}_x|z_1-z_2|]$ and by (\ref{1.2.9}),
\[
S_{0t}(\gam)\leq 2d_x^{-1}|z_1-z_2|\sup_{z\in\bar\cX,y\in\La_z}|\vf_z^u(y)|.
\]
In view of the triangle inequality for $R$ what we have proved yields the
continuity of $R(\tilde z,z)$ and $R(z,\tilde z)$ in $z$ belonging to the
open $r(x)$-neighborhood of $x$ when $\tilde z$ is fixed. Covering $\bar\cX$
by $r(x)$-neighborhoods of points $x\in\bar\cX$ and choosing a finite 
subcover we obtain (i) with the same constant for all $\bar\cX$.

Next, we derive the sufficient conditions of (ii) for the $S$-compactness. 
First, observe that both assumptions there imply that for any $\eta>0$ there 
exist $t_\eta>0$ and $\gam^\eta\in C_{0t_\eta}$ such that $\gam^\eta$ form an 
$\eta/4C$-net in $\cO$ and $S_{0t_\eta}(\gam^\eta)<\eta/4$ where $C$ is the
same as in (i). Indeed, if there exists 
a dense orbit of $\Pi^t$ in $\cO$ then a sufficiently long piece of this orbit 
will work as such $\gam$ with its $S$-functional equal 0. If $R(x,z)=0$ for any
$x,z\in\cO$ then we can choose an $\eta/4C$-net $x_1,...,x_n$ in $\cO$ and then
construct curves $\gam^{(i)}$ such that $\gam_0^{(i)}=x_i,\,\gam_{t_i}^{(i)}
=x_{i+1},\, i=1,...,n-1$ with $S_{0t_i}(\gam^{(i)})<\eta/4n$. Taking 
$\gam_t^\eta=\gam^{(i)}_{t-\sum_{1\leq j\leq i-1}t_j}$ for $t\in
[\sum_{1\leq j\leq i-1}t_j,\, \sum_{1\leq j\leq i}t_j]$ we obtain the required
curve. Now, for each $x\in\cO$ let $U_x$ be the open $r(x)$-neighborhood
of $x$ in $\bbR^d$ where the construction of the part (i) can
be implemented. Since $\cO$ is compact we can choose from the cover
$\{ U_x,\, x\in\cO\}$ of $\cO$ a finite subcover $\cU=\{ U_{x_1},...
,U_{x_\ell}\}$ of $\cO$. For any positive $\eta$ such that $\eta/4C$
is less than the
Lebesgue number (see \cite{Wa}) of $\cU$ we construct $\gam^\eta$ as above
and then for any $z\in\cO$ there is $i$ and $\tilde z\in\gam^\eta$ such that
$z,\tilde z\in U_{x_i},\, |z-\tilde z|\leq\eta/4C$, and so by the assertion (i)
 we can connect $z$ and $\tilde z$ by a curve $\gam^z\in 
C_{0t_z}$ with $t_z\leq\eta/4$ and $S_{0t_z}(\gam^z)\leq\eta/4$. It follows
 that any two points $x,z\in\cO$ can be connected by a curve $\gam\in C_{0t}$
 with $t\in[0,\, t_\eta+\eta/2]$ and $S_{0t}(\gam)\leq 3\eta/4$. Now set
 $\cO_\rho=\{ z:\,\mbox{dist}(z,\cO)\leq\rho\}$ and suppose that $\cO_{\rho_0}
 \subset\cup_{1\leq i\leq\ell}U_{x_i}$. Let $\eta/4C<\rho_0$ be smaller than
 the Lebesgue number of the cover $\{ U_{x_i},...,U_{x_\ell}\}$ of 
 $\cO_{\rho_0}$ and set $U_\eta=\{ z:\,$dist$(z,\cO)<\eta/4C\}$. Then for any
 $z\in U_\eta$ there exists $x\in\cO$ with $|z-x|<\eta/4C$, and so $x,z\in 
 U_{x_i}$ for some $i$. Hence, by (i) there exists
 a curve $\tilde\gam\in C_{0\tilde t}$ connecting $x$ with $z$ and such that 
 $\tilde t\leq\eta/4$ and $S_{0\tilde t}(\tilde\gam)\leq\eta/4$. By above
 we can connect any $\tilde x\in\cO$ with $x$ by a curve $\gam\in C_{0t}$
 with $t\in[0,\, t_\eta+\eta/2]$ and $S_{0t}(\gam)\leq 3\eta/4$ and then
 using $\tilde\gam$ we arrive at a combined curve connecting $\tilde x$ with
 $z$ and satisfying the conditions required to ensure that $\cO$ is an 
 $S$-compact by taking $T_\eta=t_{\eta}+\eta$.
 
 (iii) Now assume that for any $\eta>0$ and each $x\in\cO$ its piece of the
 $\Pi^t$-orbit of length $T(\eta)$ forms an $\eta$-net in $\cO$ and suppose
 that $B$ is complete at $x_0\in\cO$. Set $L=\sup_{-1\leq t\leq 1}\sup_{z\in
 \bar V}\| D_z\Pi^t\|$ where $D_z\Pi^t$ is the differential of $\Pi^t$ at $z$.
 Let $x\in\cO$ and $z\in\cX$ with dist$(z,\cO)<\eta L^{-T(\eta/3C)}/3C$
 where $\eta<Cr(x_0)$. Then for some $\tilde z\in\cO$, $|z-\tilde z|
 <\eta L^{-T(\eta/3C)}/3C$ and $|\Pi^{-s}\tilde z-x_0|<\eta/3C<\frac 
 13r(x_0)$,  $|\Pi^{-s}\tilde z-\Pi^{-s}z|<\eta/3C<\frac 13r(x_0)$ for some
 $s\in[0,T(\eta/3C)]$, and so $|\Pi^{-s}z-x_0|<2\eta/3C<r(x_0)$. In addition,
 for any $\eta<3Cr(x_0)$ there exists $t(\eta)>0$ so that $|\Pi^{t(\eta)}x-
 x_0|\leq\eta/3C$ with $t(\eta)\in[0,T(\eta/3C)]$. Now, by the assertion (i)
  we can connect $\Pi^{t(\eta)}x$ with $x_0$ by a curve $\gam^{(1)}\in
 C_{0t_1}$ with $t_1\leq\eta/3$ and $S_{0t_1}(\gam^{(1)})\leq\eta/3$, then 
 connect $x_0$ with $\Pi^{-s}z$ by a curve $\gam^{(2)}\in C_{0t_2}$ with
 $t_2\leq 2\eta/3$ and $S_{0t_2}(\gam^{(2)})\leq 2\eta/3$. Finally, we can
 connect $x$ with $z$ by the curve $\gam\in C_{0,t(\eta)+t_1+t_2+s}$ with
 $S_{0,t(\eta)+t_1+t_2+s}(\gam)\leq\eta$ and such that $\gam_t=\Pi^tx$
 for $t\in[0,t(\eta)]$, $\gam_t=\gam^{(1)}_{t-t(\eta)}$ for $t\in[t(\eta),
 t(\eta)+t_1]$, $\gam_t=\gam^{(2)}_{t-t(\eta)-t_1}$ for $t\in[t(\eta)+t_1,
 t(\eta)+t_1+t_2]$, and $\gam_t=\Pi^{t-t(\eta)-t_1-t_2-s}z$ for 
 $t\in[t(\eta)+t_1+t_2, t(\eta)+t_1+t_2+s]$ yielding that $\cO$ is an 
 $S$-compact.
 \end{proof}

The following assertion which relies on Lemma \ref{lem1.6.1} will be also useful
in our analysis.
\begin{lemma}\label{lem1.6.3} For any $\eta>0$ and $T>0$ there exists $\zeta>0$
such that if $\gam\in C_{0T},\,\gam\subset\cX$, $S_{0T}(\gam)<\infty$, 
$\gam_0=x_0$, and $|z_0-x_0|<\zeta$ then we can find $\tilde\gam\in C_{0T}$, 
$\tilde\gam\subset\cX$ with $\tilde\gam_0=z_0$ satisfying
\begin{equation}\label{1.6.3}
\bfr_{0T}(\gam,\tilde\gam)<\eta\,\,\mbox{and}\,\,|S_{0T}(\tilde\gam)-
S_{0T}(\gam)|<\eta.
\end{equation}
\end{lemma}
\begin{proof} By (\ref{1.2.13}) and the lower semicontinuity of the functionals
$I_z(\nu)$ there exist measures $\nu_t\in\cM_{\gam_t},\,
t\in[0,T]$ such that $\dot\gam_t=\bar B_{\nu_t}(\gam_t)$ for Lebesgue almost
all $t\in[0,T]$ and $I_{\gam_t}(\nu_t)=L(\gam_t,\dot\gam_t)$ for Lebesgue 
almost all $t\in[0,T]$. Recall also that $\dot\gam_t$ is measurable in $t$. 
Introduce the (measurable) map
$q:\,[0,T]\times\cP(\bar\cW)\to\bbR\cup\{\infty\}\times\bbR^d$ defined
by $q(t,\nu)=\big( I_{\gam_t}(\nu),\bar B_\nu(\gam_t)\big)$. Recall that
$\dot\gam_t$ is measurable in $t$, and so another map $r:\,[0,T]\to
\bbR\cup\{\infty\}\times\bbR^d$ defined by $r(t)=\big(L(\gam_t,\dot\gam_t\big),
\dot\gam_t\big)$ is also measurable in $t\in[0,T]$. Then $q(t,\nu_t)=r(t)$
 and it follows from the measurable selection in the implicit function theorem 
(see \cite{CV}, Theorem III.38) that measures $\nu_t$ satisfying this 
condition can be chosen to depend measurably on $t\in[0,T]$. 

Now, given $\eta>0$ we pick up a small $\zeta>0$ which will be specified 
later on and employ Lemma \ref{lem1.4.3} in the same way as in 
(\ref{1.5.31}) together with (\ref{1.2.9}), (\ref{1.2.11}), and (\ref{1.2.13})
in order to conclude that for all $n\in\bbN$ large enough 
there exists $t^{(n)}_1\in[0,T/n)$ such that if $t^{(n)}_{j+1}=t_1^{(n)}+
jn^{-1}T,\, j=1,2,...,n-1,\, t^{(n)}_{n+1}=T$ then
\begin{eqnarray}\label{1.6.4}
&\int_0^{t_1^{(n)}}\big\vert\bar B_{\nu_s}(\gam_s)\big\vert ds+
\sum_{j=1}^n\int_{t_j^{(n)}}^{t_{j+1}^{(n)}}\big\vert\bar B_{\nu_s}(\gam_s)
-\bar B_{\nu_{t_j^{(n)}}}(\gam_{t_j^{(n)}})\big\vert ds\\
&+S_{0t_1^{(n)}}(\gam)+\sum_{j=1}^n\int_{t_j^{(n)}}^{t_{j+1}^{(n)}}\big\vert
I_{\gam_s}(\nu_s)-I_{\gam_{t_j^{(n)}}}(\nu_{t_j^{(n)}})\big\vert ds<\zeta.
\nonumber\end{eqnarray}
Set $\dot\psi_s^{(n)}=0$ for $s\in[0,t^{(n)}_1)$ and $\dot\psi_s^{(n)}=
\bar B_{\nu_{t_j^{(n)}}}(\gam_{t_j^{(n)}})$ for $s\in[t_j^{(n)},
t^{(n)}_{j+1}),\, j=1,...,n$. Then $\psi_t^{(n)}=\gam_0+
\int_0^t\dot\psi_s^{(n)}ds,\, t\in[0,T]$ defines a polygonal line such that
\begin{equation*}
\bfr_{0T}(\gam,\psi^{(n)})<\zeta.
\end{equation*}

Next, set $\tilde\gam_t=z_0$ for all $t\in[0,t_1^{(n)}]$ and continue the
construction of $\tilde\gam$ in the following recursive way. Suppose that
$\tilde\gam_t$ is already defined for all $t\in[0,t_j^{(n)}]$ and some 
$j\geq 1$. Denote $x_j=\gam_{t_j^{(n)}}$, $y_j=\psi^{(n)}_{t_j^{(n)}}$,
$z_j=\tilde\gam_{t_j^{(n)}}$ and
suppose that $|z_j-x_j|<r-KTn^{-1}$ where $K$ is the same as in (\ref{1.2.15})
and $r$ comes from Lemma \ref{1.6.1}. For $t\in[t_j^{(n)},t_{j+1}^{(n)}]$
define $\tilde\gam_t$ as the integral curve starting at $z_j$ of the vector
field $\tilde B(z)=\bar B_{\mu_{x_jz}}(z),\, |z-x_j|<r$ with $\mu_{x_jz}\in
\cM_z$ obtained in Lemma \ref{lem1.6.1} for $\mu_{x_jx_j}=\nu_{t_j^{(n)}}$,
i.e. $\tilde\gam_t$ is the solution of the equation
\[
\tilde\gam_t=z_j+\int_{t_j^{(n)}}^t\tilde B(\tilde\gam_s)ds.
\]
This definition is legitimate since in view of (\ref{1.2.15}) and our assumption
 on $z_j$ the curve $\tilde\gam_t,\, t\in[t_j^{(n)},t_{j+1}^{(n)}]$ does not
 exit the $r$-neighborhood of $x_j$. By (\ref{1.6.1}) and the above
 for all $t\in[t_j^{(n)},t_{j+1}^{(n)}]$,
 \[
 \big\vert\tilde B(\tilde\gam_t)-\bar B_{\nu_{t_j^{(n)}}}(x_j)\big\vert <
 C(|z_j-x_j|+KTn^{-1})<C(|z_j-y_j|+\zeta+KTn^{-1}),
 \]
 and so
 \begin{equation*}
 |z_{j+1}-y_{j+1}|\leq\sup_{t\in[t_j^{(n)},t_{j+1}^{(n)}]}|
 \tilde\gam_t-\psi_t^{(n)}|\leq|z_j-y_j|(1+CTn^{-1})+\zeta T/n+CKT^2n^{-2}.
 \end{equation*}
 Assuming that $|z_0-x_0|<\zeta$ with $\zeta$ small enough and since
 $x_0=y_0$ we obtain successively from here that for all $j=1,2,...,n$,
 \[
 |z_j-y_j|\leq (1+CTn^{-1})^n(2\zeta+KTn^{-1})\leq e^{CT}(2\zeta+
 KTn^{-1}),
 \]
 which enables us to continue our construction recursively for $j=1,2,...,n$
 if $\zeta$ and $n^{-1}$ are small enough yielding also that
 \[
 \bfr_{0T}(\tilde\gam,\psi^{(n)})\leq e^{CT}(2\zeta+KTn^{-1}).
 \]
 Hence, the first part of (\ref{1.6.3}) follows provided $\zeta$ and $n^{-1}$
 are sufficiently small. 
 
 Next, observe that 
 \[
 \sup_{t\in[t_j^{(n)},t_{j+1}^{(n)}]}|\tilde\gam_t-x_j|\leq |z_j-y_j|
 |y_j-x_j|+KTn^{-1}\leq KTn^{-1}(1+e^{CT})+\zeta(1+2e^{CT})
 \]
 and the right hand side here can be made as small as we wish choosing $\zeta$
 small and $n$ large. Hence, by Lemma \ref{lem1.6.1} we can make 
 \[
 \max_{0\leq j\leq n}\sup_{t\in[t_j^{(n)},t_{j+1}^{(n)}]}\big\vert
 I_{\tilde\gam_t}(\mu_{x_j\tilde\gam_t})-I_{x_j}(\nu_{t_j^{(n)}})\big\vert
 <\eta/2T
 \]
 which together with (\ref{1.6.4}) yield the second part of (\ref{1.6.3}).
 \end{proof}
 
 The following result will enable us to control the time which the slow 
 motion can spend away from the $\om$-limit set of the averaged motion.
 \begin{lemma}\label{lem1.6.4} Let $G\subset\cX$ be a compact set not 
    containing entirely any forward semi-orbit of the flow $\Pi^t$.
    Then there exist positive
     constants $a=a_G$ and $T=T_G$ such that for any $x\in G$ and $t\geq 0$,
     \begin{equation*}
     \inf\big\{ S_{0t}(\gam):\,\gam\in C_{0t}\,\,\mbox{and}\,\,\gam_s\in G\,\,
     \mbox{for all}\,\, s\in[0,t]\big\}\geq a[t/T]
     \end{equation*}
     where $[c]$ denotes the integral part of $c$.
     \end{lemma}
     \begin{proof} For each $x\in G$ set $\sig_x=\inf\{ t\geq 0:\, \Pi^tx
     \not\in G\}$. By the assumption of the lemma $\sig_x<\infty$ for
     each $x\in G$ and it follows from continuous dependence of solutions
     of (\ref{1.1.6}) on initial conditions that $\sig_x$ is upper semicontinuous.
     Hence, $\tilde T+\sup_{x\in G}\sig_x<\infty$. Set $T=\tilde T+1$ and
     $\Gam=\{\gam\in C_{0T}:\,\gam_s\in G$ for all $s\in[0,T]\}$. Since
     no $\gam\in\Gam$ can be a solution of the equation (\ref{1.2.14}) then
     $S_{0T}(\gam)>0$ for any $\gam\in\Gam$. The set $\Gam$ is closed with
      respect to the uniform convergence and since the functional $S_{0T}$
      is lower semicontinuous we obtain that
      \[
      \inf_{\gam\in\Gam}S_{0T}(\gam)=a>0.
      \]
      This together with (\ref{1.2.13}) yield the assertion of Lemma 
      \ref{lem1.6.4}.
      \end{proof}

Untill now we have not used specific assumptions of Theorem \ref{thm1.2.5} 
but some of them will be needed for the following auxiliary result.
\begin{lemma}\label{lem1.6.5} 
Let $V$ be a connected open set with a piecewise smooth boundary and assume
that (\ref{1.2.20}) holds true. Then the function $R_\partial(x)$ is upper
semicontinuous at any $x_0\in V$ for which $R_\partial(x_0)<\infty$. 
Let $\cO\subset V$ be an $S$-compact. 

(i) Then for each $z\in\bar V$ the function $R(x,z)$ takes on the same
value $R^\cO(z)$ for all $x\in\cO$, and so $R_\partial(x)$ takes on the same
value $R_\partial$ for all $x\in\cO$ and the set $\partial_{\min}(x)=\{ z\in
\partial V:\, R(x,z)=R_\partial\}$ coincides with the same (may be empty)
set $\partial_{\min}$ for all $x\in\cO$. Furthermore, for each $\del>0$ 
there exists $T(\del)>0$ such that for any $x\in\cO$ we can construct 
$\gam^x\in C_{0t_x}$ with $t_x\in(0,T(\del)]$ satisfying
\begin{equation}\label{1.6.5}
\gam^x_0=x,\,\,\gam^x_{t_x}\in\partial V\,\,\mbox{and}\,\, S_{0t_x}(\gam^z)
\leq R_\partial+\del.
 \end{equation}
 
 (ii) Suppose that $R_\partial<\infty$ and dist$(\Pi^tx,\cO)
 \leq d(t)$ for some $x\in V$ and $d(t)\to 0$ as $t\to\infty$. Then
 $R_\partial(x)\leq R_\partial$ and for any $\del>0$ there exist $T_{\del,d}
 >0$ (depending only on $\del$ and the function $d$ but not on $x$) and
 $\hat\gam^x\in C_{0s_x}$ with $s_x\in(0,T_{\del,d}]$ satisfying
 \begin{equation}\label{1.6.6}
\hat\gam^x_0=x,\,\,\hat\gam^x_{s_x}\in\partial V\,\,\mbox{and}\,\, S_{0s_x}
(\hat\gam^x)\leq R_\partial+\del.
 \end{equation}
 In particular, if $R_\partial <\infty$ then $R_\partial(x)<\infty$ and
 if $\cO$ is an $S$-attractor of the flow $\Pi^t$ then $R_\partial(x)<\infty$
 for all $x\in V$.
 
 (iii) Suppose that for any open set $U\supset\cO$ the compact set
 $\bar V\setminus U$ does not contain entirely any forward semi-orbit of
 the flow $\Pi^t$. Then the function $R^\cO(z)$ is lower semicontinuous
 in $z\in\bar V$, $R^\cO(z)\to 0$ as dist$(z,\cO)\to 0$, and $\partial_{\min}$
 is a nonempty compact set.
\end{lemma}
\begin{proof} Let $R_\partial(x_0)<\infty$
for some $x_0\in V$. Then for any $\al>0$ there exist $T>0$ and $\gam\in 
C_{0T}$ such that $\gam_0=x_0,\,\gam_T\in\partial V$ and $S_{0T}(\gam)\leq
 R_\partial(x_0)+\al$. By Lemma \ref{lem1.6.3} for any $\eta>0$ we can choose 
 $\zeta>0$ so that if $z_0\in V,\, |z_0-x_0|<\zeta$ then there exists 
 $\tilde\gam\in C_{0T}$ such that $\tilde\gam_0=z_0,\,\bfr_{0T}(\tilde\gam,
 \gam)<\eta$ and $|S_{0T}(\tilde\gam)-S_{0T}(\gam)|<\eta$.
 Let $\mu_{\gam_Tz}\in\cM_z,\, |z-\gam_T|<r$ be measures obtained in 
 Lemma \ref{lem1.6.1} for $\mu_{\gam_T\gam_T}=\nu$ with $\nu$ satisfying
 the second part of (\ref{1.2.20}). Since the boundary $\partial V$ is
 piecewise smooth it follows from the continuous dependence 
 of solutions of ordinary differential equations on initial conditions that
 for all small $\eta>0$ there exists $t(\eta)\to 0$ as $\eta\to 0$ such 
 that if $\psi_t,\, t\in[0,t(\eta)]$ is an integral curve of the vector
  field $\bar B_{\mu_{\gam_Tz}}(z),\, |z-\gam_T|<r$ with $\psi_0=\tilde z,
  \, |\tilde z-\gam_T|<\eta$ then $\psi_{t(\eta)}\not\in V$. Since
  $|\tilde\gam_T-\gam_T|<\eta$ we can define $\tilde\gam_t=\psi_{t-T}$
  for $t\in[T,T+t(\eta)]$. Now, $\tilde\gam_{T+t(\eta)}\not\in V$ and by
  (\ref{1.2.9}),
  \[
  \big\vert S_{0,T+t(\eta)}(\tilde\gam)-S_{0T}(\tilde\gam)\big\vert
  \leq t(\eta)\sup_{x\in\bar\cX,y\in\La_x}|\vf^u_x(y)|.
  \]
  Thus we can choose $\eta$ so small that $R_\partial(z_0)\leq R_\partial(x_0)
  +2\al$ and the upper semicontinuity of $R_\partial(x)$ at $x_0$ follows.
  
  From now on till the end of the proof of this lemma we assume that $\cO$
  is an $S$-compact and prove, first, the assertion (i). It follows from
  the definition of an $S$-compact that $R(x_1,x_2)=0$ for any pair 
  $x_1,x_2\in\cO$, and so $R(x_1,z)=R(x_2,z)$ for any such $x_1,x_2$ and
  each $z\in\bar V$. It follows that $R_\partial(x)$ takes on the same 
  value $R_\partial$ for all $x\in\cO$ and all sets $\partial_{\min}(x),
  x\in\cO$ coincide with some, may be empty, set $\partial_{\min}$. Fix
  $x_0\in\cO$. Then for each $\del>0$ there exists $t^{(0)}_\del>0$
  and $\gam^{(0)}\in C_{0t_\del^{(0)}}$ such that
  \[
  \gam_0^{(0)}=x_0,\,\,\gam^{(0)}_{t^{(0)}_\del}\in\partial V\,\,\mbox{and}
  \,\, S_{0t_\del^{(0)}}(\gam^{(0)})\leq R_\partial+\del/2.
  \]
  By the definition of an $S$-compact there exists $T_{\del/2}>0$ such
  that for any $z\in\cO$ we can construct $\gam^{(z,\del)}\in 
  C_{0t_\del^{(z)}}$ with $t_\del^{(z)}\in[0,T_{\del/2}]$ satisfying
  \[
  \gam_0^{(z,\del)}=z,\,\,\gam^{(z,\del)}_{t^{(z)}_\del}=x_0\,\,\mbox{and}
  \,\, S_{0t_\del^{(z)}}(\gam^{(z,\del)})\leq\del/2.
  \]
  Defining $\gam^z$ by $\gam^z_t=\gam_t^{(z,\del)}$ for $t\in[0,t_\del^{(z)}]$
  and $\gam^z_t=\gam^{(0)}_{t-t_\del^{(z)}}$ for $t\in[t_\del^{(z)},
  t_\del^{(z)}+t_\del^{(0)}]$ we obtain a curve satisfying (\ref{1.6.5})
  with $T(\del)=t_\del^{(0)}+T_{\del/2}$.
  
  Next, we prove (ii) assuming that $R_\partial <\infty$ and that
   dist$(\Pi^tx,\cO)\leq d(t)$ for some
   $x\in V$ with $d(t)\to 0$ as $t\to\infty$. By (i), for
   any $\eta>0$ there exists $T_\eta>0$ such that for any $z\in\cO$ we
   can construct $\gam^z\in C_{0t_z}$ with $t_z\in(0,T_\eta]$ and $\gam^z
   \in C_{0t_z}$ satisfying (\ref{1.6.5}) with $\del=\eta$. For such $\eta$
   and $T_\eta$ choose $\zeta$ by Lemma \ref{lem1.6.3} so that if $|\tilde x-z|
   <\zeta$ and $z\in\cO$ then in the same way as at the beginning of the
   proof of this lemma we can construct
   $\tilde\gam\in C_{0,t_z+t(\eta)}$ with $t(\eta)\to 0$ as $\eta\to 0$ such
    that
    \[
    \tilde\gam_0=\tilde x,\,\,\tilde\gam_{t_z+t(\eta)}\in\partial V\,\,
  \mbox{and}\,\, |S_{0,t_z+t(\eta)}(\tilde\gam)-S_{0t_z}(\gam^z)|\leq
  \eta+Ct(\eta).
  \]
  Pick up $\tilde t=\tilde t(d,\zeta)$ so that $d(\tilde t)<\zeta$. Then 
  $|\tilde x-z|<\zeta$ for $\tilde x=\Pi^{\tilde t}x$ and some $z\in\cO$. 
  Now construct as above $\tilde\gam$ for such $z$ and define 
  $\hat\gam^x\in C_{0s_x}$ with $s_x=\tilde t+t_z+t(\eta)$ setting 
  \[
  \hat\gam_t^x=\Pi^tx\,\,\mbox{for}\,\, t\in[0,\tilde t]\,\,\mbox{and}\,\,
  \hat\gam^x_t=\tilde\gam(t-\tilde t)\,\,\mbox{for}\,\, t\in[\tilde t,
  \tilde t+t_z+t(\eta)].
  \]
  Then $S_{0s_x}(\hat\gam^x)\leq R_\partial+2\eta+Ct(\eta)$ and $s_x\in
  (0,T_\eta+\tilde t+t(\eta)]$. Choosing $\eta$ so small that $2\eta +Ct(\eta)
  \leq\del$ and then taking $T_{\del,d}=T_\eta+\tilde t+t(\eta)$ we conclude
  that $\hat\gam^x$ satisfies (\ref{1.6.6}). Since $\eta$ is arbitrary we obtain
  that $R_\partial(x)\leq R_\partial$. If $\cO$ is an $S$-attractor whose 
  basin contains $\bar V$ then we can choose $d(t)\to 0$ as $t\to 0$ which
  in view of the continuous dependence of $\Pi^tx$ on $x$ will be the same
  for all $x\in\bar V$ (though for this lemma $d(t)$ as above depending 
  on $x$ would suffice, as well), so our conditions are satisfied now for all 
  $x\in V$. Hence, in this case $R_\partial(x)$ is finite in 
  the whole $V$, completing the proof of (ii).
  
   Finally, we prove (iii). Recall, that by the definition of an $S$-compact 
  $\cO$ it follows that $R(x,z)=0$ whenever $x,z\in\cO$. For all $\eta>0$
  let $U_\eta\supset\cO$ be open sets appearing in the definition of an
  $S$-compact. If $x\in\cO$ and $z\in U_\eta$ then $R(x,z)\leq\eta$.
   Hence, if dist$(z_n,\cO)\to 0$
  as $n\to\infty$ then $R(x,z_n)\to 0$. Now, let $z_0\in\bar V\setminus\cO$
  and $z_n\to z_0$ as $n\to\infty$. For each $\del\geq 0$ set $\cO_\del=
  \{ z\in V:\,$dist$(z,\cO)\leq\del\}$ and let $\del(\eta)=\frac 12
  \inf\{ |x-z|:\, x\in\cO,\, z\in\bar V\setminus U_\eta\}$. Without loss of
  generality we will assume that $z_i\not\in U_{\eta_0}$ for some $\eta_0
  >0$ and all $i=0,1,2,\dots$. Fix $x\in\cO$. By the definition of the 
  function $R$ for any $\zeta>0$ we can choose $t_{n,\zeta}>0$ and 
  $\gam^{(n,\zeta)}\in C_{0t_{n,\zeta}}$, $n=0,1,2,...$ such that
  \begin{equation}\label{1.6.7}
  \gam_0^{(n,\zeta)}=x,\,\gam^{(n,\zeta)}_{t_{n,\zeta}}=z_n\,\,\mbox{and}\,\,
  S_{0t_{n,\zeta}}(\gam^{(n,\zeta)})\leq R(x,z_n)+\zeta.
  \end{equation}
  For each $\eta\leq\eta_0$ set
  \[
  s_n=s_{n,\eta,\zeta}=\sup\{ t\geq 0:\,\gam_t^{(n,\zeta)}\in
  \cO_{\del(\eta)}\}.
  \]
  Consider $\tilde\gam^{(n)}\in C_{0,t_{n,\zeta}-s_n}$ defined by 
  $\tilde\gam^{(n)}_t=\gam^{(n,\zeta)}_{t+s_n}$ for $t\in[0,t_{n,\zeta}-s_n]$
  which stays in $\bar V\setminus$int$\cO_{\del(\eta)}$ (where int$G$ means
  the interior of a set $G$), and so by Lemma \ref{lem1.6.4} we conclude that
  \[
  t_{n,\zeta}-s_n\leq a^{-1}_\eta(R(x,z_n)+1)
  \]
  provided, say, $\zeta\leq 1/2$ where $a_\eta>0$ depends only on $\eta$.
  In order to verify the lower semicontinuity of $R(x,z)$ at $z=z_0$ we
  have only to consider the case
  \[
  \liminf_{n\to\infty}R(x,z_n)=A<\infty,
  \]
  and so we can assume that $R(x,z_n)\leq 2A$ for all $n=0,1,2,\dots$.
  Passing to a subsequence and denoting its members by the same letters we 
  can assume also that
  \[
  \lim_{n\to 0}R(x,z_n)=A.
  \]
  The curves $\tilde\gam^{(n)}$ are Lipschitz continuous with a constant
  $K$ from (\ref{1.2.15}), and so this sequence is relatively compact. Hence,
  we can choose a uniformly converging subsequence and denoting, again, its
  members by the same letters we obtain now that 
  \[
  \tilde\gam^{(n)}\to\tilde\gam^{(0)}\,\,\mbox{as}\,\, n\to\infty
  \]
  where $\tilde\gam^{(0)}\in C_{0t_0}$ with $t_0\in(0,a^{-1}_\eta(A+1)],\,$
  dist$(\tilde\gam_0^{(0)},\cO)=\del(\eta)$ and $\tilde\gam^{(0)}_{t_0}=
  z_0$. Each curve $\tilde\gam^{(n)},\, n=0,1,2,...$ can be extended to a
  curve in $C_{0,T}$ with $T=a_\eta^{-1}(2A+1)$ and the same $S$-functional
  by adding to one of its ends a piece of the orbit of the flow $\Pi^t$.
  Hence, we can rely on the lower semicontinuity of the functional $S_{0T}$
  in order to derive from (\ref{1.6.7}) that
  \[
  S_{0t_0}(\tilde\gam^{(0)})\leq A+\zeta.
  \]
  By the definition of an $S$-compact there exists $\hat\gam\in C_{0r}$ 
  with $r\in[0,T_{2\eta}]$ such that $\hat\gam_0=x,\,\hat\gam_r=\tilde\gam_0$
  and $S_{0r}(\hat\gam)\leq 2\eta$. It follows that
  \[
  R(x,z_0)\leq A+2\eta+\zeta
  \]
  and since $\eta$ and $\zeta$ can be chosen arbitrarily small we conclude 
  that $R(x,z_0)\leq A$ obtaining the lower semicontinuity of $R(x,z)$ at
  $z=z_0$. Finally, the lower semicontinuity of $R(x,z)$ in $z\in\partial V$
  for a fixed $x\in\cO$ implies that $\partial_{\min}(x)$ is nonempty and
  compact and since $\partial_{\min}(x)$ is the same for all $x\in\cO$ by
  (i), the proof of Lemma \ref{lem1.6.5} is complete.
  \end{proof}

\section{"Very long" time behavior: exits from a domain}\label{sec1.7}
\setcounter{equation}{0}

In this section we derive Theorems \ref{thm2.2.5} relying on 
certain "Markov property type" \index{Markov property}
arguments which are substantial modifications
of the corresponding arguments from Sections 4 and 5 of \cite{Ki2}. In this
and the following section in order to simplify notations we will write
$\cD^u_\ve(z,\al,\rho,C)$ for $\hat\cD^u_\ve(z,\al,\rho,C,L)$ (both introduced
in Section \ref{sec1.3}) with some large $L$ so that appropriate discs on
(extended) unstable leaves $W^u_x$ and all their $\Phi_\ve^s$-iterates belong
to this set. We start with the following result which will not only yield 
Theorem \ref{thm1.2.5} but also will play an important role in the proof of 
Theorem \ref{thm1.2.7} in the next section. 
\begin{proposition}\label{prop1.7.1} Let $V$ be a connected open set with
a piecewise smooth boundary $\partial V$ such that $\bar V=V\cup\partial V
\subset\cX$. Assume that for each $z\in\partial V$ there exist $\iota=
\iota(z)>0$
and an $F^t$-invariant probability measure $\nu$ on $\La_z$ so that
\begin{equation}\label{1.7.1}
z+s\bar B_\nu(z)\in\bbR^d\setminus\bar V\,\,\mbox{for all}\,\, s\in(0,\iota],
\end{equation}
i.e. $\bar B_\nu(z)\ne 0$ and it points out into the exterior of $\bar V$.

(i) Suppose that for some $A_1,T>0$ and any $z\in\bar V$ there exists 
$\vf^z\in C_{0T}$ such that for some $t=t(z)\in(0,T]$,
\begin{equation}\label{1.7.2}
\vf_0^z=z,\,\vf^z_t\not\in V\,\,\mbox{and}\,\, S_{0t}(\vf^z)\leq A_1.
\end{equation}
Then for each $x\in V$,
\begin{equation}\label{1.7.3}
\limsup_{\ve\to 0}\ve\log\int_{\cW}\tau^\ve_{x,y}(V)dm(y)\leq A_1
\end{equation}
and for any $\al>0$ there exists $\la(\al)=\la(x,\al)>0$ such that for
all small $\ve>0$,
\begin{equation}\label{1.7.4}
m\big\{ y\in\cW:\,\tau^\ve_{x,y}(V)\geq e^{(A_1+\al)/\ve}\big\}\leq 
e^{-\la(\al)/\ve}.
\end{equation}

(ii) Assume that there exists an open set $G$ such that $V$ contains its 
closure $\bar G$ and the intersection of $\bar V\setminus G$ with the 
$\om$-limit set of the flow $\Pi^t$ is empty. Let $\Gam$
be a compact subset of $\partial V$ such that
\begin{equation}\label{1.7.5}
\inf_{x\in G,z\in\Gam}R(x,z)\geq A_2
\end{equation}
for some $A_2>0$. Then for some $T>0$ and any $\be>0$ there exists 
$\la(\be)>0$ such that for each $x\in V$ and any small $\ve>0$,
\begin{eqnarray}\label{1.7.6}
&m\big\{ y\in\cW:\, Z^\ve_{x,y}(\tau^\ve_{x,y}(V))\in\Gam,\,
\tau^\ve_{x,y}(V)\leq e^{(A_2-\be)/\ve}\big\}\\
&\leq m\big\{ y\in\cW:\, Z^\ve_{x,y}(\tau^\ve_{x,y}(V))\in\Gam,\,
\tau^\ve_{x,y}(V)<T\big\}+e^{-\la(\be)/\ve}.\nonumber
\end{eqnarray}
Suppose that for some $x\in V$,
\begin{equation}\label{1.7.7}
a(x)=\inf_{t\geq 0}\mbox{dist}(\Pi^tx,\partial V)>0.
\end{equation}
Then $R_\partial(x)>0$ and for each $T>0$ there exists $\hat\la(T)=
\hat\la(T,x)>0$ such that for all small $\ve>0$,
\begin{equation}\label{1.7.8}
m\{ y\in\cW:\,\tau^\ve_{x,y}(V)<T\}\leq e^{-\hat\la(T)/\ve}
\end{equation}
and if the set $\Gam$ from (\ref{1.7.5}) coincides with the whole
$\partial V$ then
\begin{equation}\label{1.7.9}
\liminf_{\ve\to 0}\ve\log\int_\cW\tau^\ve_{x,y}(V)dm(y)\geq A_2.
\end{equation}
The corresponding to (\ref{1.7.3}), (\ref{1.7.4}), (\ref{1.7.8}) and 
(\ref{1.7.9})
assertions hold true also when $\cW$ and $m$ in these estimates are replaced 
by a disc $D\in\cD^u_\ve(z,\al,\rho,C)$ with $\pi_1z=x$ and by $m_D$, 
respectively ((see (\ref{1.7.21}) and\ref{1.7.22}), (\ref{1.7.34}), 
(\ref{1.7.36}) and (\ref{1.7.37}) below).
\end{proposition}
\begin{proof}
Observe that applying to (\ref{1.5.19}) and (\ref{1.5.23}) the 
arguments which were used in order to derive Theorem \ref{thm1.2.3} from 
Proposition \ref{prop1.5.2} and the latter from Proposition \ref{prop1.3.4} and 
Lemma \ref{lem1.4.1} we obtain that (\ref{1.2.16}) and (\ref{1.2.17}) can be 
written 
for any disc $D\in\cD^u_\ve(z,\al,\rho,C),\, z=(x,y)$ in place of the whole 
$\cW$, namely, for any $\gam\in C_{0T}$ with $\gam_0=x$, $\del,\la,a>0$, and 
$\ve$ small enough 
 \begin{equation}\label{1.7.10}
    m_D\left\{ v\in D:\,\bfr_{0T}(Z^\ve_v,\gam)<\del\right\}\geq
    \exp\left\{-\frac 1\ve(S_{0T}(\gam)+\la)\right\}
    \end{equation}
    and
    \begin{equation}\label{1.7.11}
    m_D\left\{ v\in D:\,\bfr_{0T}(Z^\ve_v,\Psi^a_{0T}(x))\geq\del\right\}
    \leq\exp\left\{-\frac 1\ve(a-\la)\right\}
    \end{equation}
    which holds true in the same sense as (\ref{1.2.16})--(\ref{1.2.17}) and
    (\ref{1.7.10})--(\ref{1.7.11}) are uniform in $D$ as above.
    
    In order to prove (i) we observe, first, that the assumption (\ref{1.7.1})
    above together with Lemma \ref{lem1.6.2}(i) and the compactness of
     $\partial V$
    considerations enable us to extend any $\vf^z,\, z\in V$ slightly so that 
    it will exit some fixed neighborhood of $V$ with only slight increase in
    its $S$-functional. Hence, from the beginning we assume that for each
    $\be>0$ there exists $\del=\del(\be)>0$ such that for any $z\in V$ we can
    find $T>0$, $\vf^z\in C_{0T}$ and $t=t(z)\in(0,T]$ satisfying
    \[
    \vf_0^z=z,\,\vf^z_t\not\in V_\del\,\,\mbox{and}\,\, S_{0t}(\vf^z)\leq
    A_1+\be
    \] 
    where $V_\del=\{ x:\,$dist$(x,V)\leq\del\}$. It follows that for any
    $x\in V,\, n\geq 1$, and $D\in\cD^u_\ve\big((x,w),\al,\rho,C\big)$,
    \begin{eqnarray}\label{1.7.12}
    &\quad\,\,\,\big\{ v\in D:\,\tau^\ve_v(V)>nT\big\}=\big\{ v\in D:\,
    Z^\ve_v(t)\in V,\,\forall\,t\in[0,nT]\big\}\\
    &=\big\{ v\in D:\,\tau^\ve_{\Phi_\ve^{kT/\ve}v}(V)>T,\,\,
    \forall\, k=0,1,...,n-1\big\}\subset G^\ve_{n,\del}\nonumber\\
    &\stackrel{\mbox{def}}{=}\big\{ v\in D:\,\Phi_\ve^{kT/\ve}v
    \not\in\bigcup_{z\in V}A_{\del}^{\ve,T}(\vf^z),
    \forall\, k=0,1,...,n-1\big\}\nonumber
    \end{eqnarray}
    where for any subset $H\subset C_{0T}$ and $c>0$, 
    \[
    A_{c}^{\ve,T}(H)=\{ v\in V\times\cW:\,\bfr_{0T}(Z^\ve_v,H)<c\}.
    \]
    
     For $k=1,2,...$ define
    \[
    Q^\ve_{k,\del}=\big\{ v\in D:\, U^\ve_D(kT/\ve,v,\del/4)\cap 
    G^\ve_{k,\del}\ne\emptyset\big\}
    \]
    and 
    \[
    R^\ve_{k,\del}=\bigcup_{v\in Q^\ve_{k,\del}}U^\ve_D(kT/\ve,v,\del/4)
    \]
    which are, clearly, compact sets satisfying
    \begin{equation*}
    G^\ve_{k,\del}\subset Q^\ve_{k,\del}\subset R^\ve_{k,\del}
    \subset G^\ve_{k,\del/2}.
    \end{equation*}
    Let $E_k$ be a maximal $\big(kT/\ve,\del/2,\ve,Q^\ve_{k,\del},D\big)$-
    separated set in $Q^\ve_{k,\del}$. Then 
    \begin{equation}\label{1.7.13}
    \bigcup_{v\in E_k}U_D(kT/\ve,v,\del/4)\subset R^\ve_{k,\del}\subset
    \bigcup_{v\in E_k}U_D(kT/\ve,v,3\del/4)
    \end{equation}
    and the left hand side of (\ref{1.7.13}) is a disjoint union. This together
      with Lemma \ref{lem1.3.6} give
    \begin{eqnarray}\label{1.7.14}
    &m_D(R^\ve_{k,\del})\leq\sum_{v\in E_k}m_D\big(U_D(kT/\ve,v,3\del/4)
    \big)\\
    &\leq c^{-1}_{3\del/4}c^{-1}_{\del/4}
    \sum_{v\in E_k}m_D\big(U_D(kT/\ve,v,\del/4)\big).\nonumber
    \end{eqnarray}
    
    By Lemma \ref{lem1.3.2}(ii),
    \begin{equation}\label{1.7.15}
    D_k(v)=\Phi_\ve^{kT/\ve}U_D(kT/\ve,v,\del/4)\in\cD^u_\ve
    (\Phi_\ve^{kT/\ve}v,\al,\frac {\del}{4\sqrt C},\sqrt C).
    \end{equation}
    Clearly, for any $v\in E_k$,
    \begin{eqnarray}\label{1.7.16}
    &\Gam_k(v)=\big\{ w\in U_D(kT/\ve,v,\del/4):\\
    &\Phi_\ve^{kT/\ve}w\in
    \bigcup_{z\in V}A_{\del/2}^{\ve,T}(\vf^z)\big\}\subset 
    R^\ve_{k,\del}\setminus R^\ve_{k+1,\del}.\nonumber
    \end{eqnarray}
    In view of (\ref{1.7.15}) we can apply (\ref{1.7.10}) which together with 
    the
    choice of curves $\vf^x$ yield that for any $\la>0$ and $\ve$ small
    enough,
    \begin{eqnarray}\label{1.7.17}
    &m_{D_k(v)}\big(\Phi_\ve^{kT/\ve}\Gam_k(v)\big)
    \geq m_{D_k(v)}\big\{w\in D_k(v):\,\\
    &\bfr_{0,T}(Z^\ve_w,\,\vf^{z_k^\ve(v)})<\del/2\big\} 
    \geq\exp\big(-\frac 1\ve(A_1+\be+\la)\big)\nonumber
    \end{eqnarray}
    where $z_k^\ve(v)=\pi_1(\Phi_\ve^{kT/\ve}v)$.
    By Lemma \ref{lem1.3.6} it follows that
    \begin{equation}\label{1.7.18}
    m_D(\Gam_k(v))\geq c(\del)m_D\big( U_D(kT/\ve,v,\del/4))
    \exp\big(-\frac 1\ve(A_1+\be+\la)\big)
    \end{equation}
    for some $c(\del)>0$. Since $U_D(kT/\ve,v,\del/4)$
    are disjoint for different $v\in E_k$ we derive from (\ref{1.7.14}), 
    (\ref{1.7.16}) and (\ref{1.7.18}) that
    \begin{eqnarray}\label{1.7.19}
    &\,\,\,\,\,m_D(R^\ve_{k,\del})-m_D(R^\ve_{k+1,\del})\geq m_D\big(
    \bigcup_{v\in E_k} \Gam_k(v)\big)=\sum_{v\in E_k}m_D(\Gam_k(v))\\
    &\geq c(\del)\exp\big(-\frac 1\ve(A_1+\be+\la)\big)\sum_{v\in E_k}
    U_D(kT/\ve,v,\del/4)\nonumber \\
    &\geq\tilde c_\del m_D(R^\ve_{k,\del})\exp\big(-\frac 1\ve(A_1
    +\be+\la)\big)\nonumber
    \end{eqnarray}
    where $\tilde c_\del=c(\del)c_{3\del/4}c_{\del/4}$. Applying (\ref{1.7.19}) 
    for $k=1,2,...,n-1$ we obtain that
    \begin{equation}\label{1.7.20}
    m_D(G^\ve_{n,\del})\leq m_D(R^\ve_{n,\del})\leq\bigg( 1-\tilde c_\del 
    \exp\big(-\frac 1\ve(A_1+\be+\la)\big)\bigg)^nm_D(D).
    \end{equation}
    This together with (\ref{1.7.12}) yield that for any $\be>0$ there exists
    $c(\be)>0$ such that for all small $\ve>0$,
    \begin{equation}\label{1.7.21}
    m_D\big\{ v\in D:\,\tau_v^\ve(V)>e^{(A_1+\be)/\ve}\big\}
    <e^{-c(\be)/\ve}.
    \end{equation}
    
     Observe that by (\ref{1.7.12}) and (\ref{1.7.20}),
    \begin{eqnarray}\label{1.7.22}
    &\int_D\tau_v^\ve(V)dm_D(v)\leq\sum_{n=0}^\infty (n+1)T
    \big( m_D\big\{ v\in D:\,\\
    &\tau_v^\ve(V)>nT\big\}- m_D\big\{ v\in D:\,
    \tau_v^\ve(V)>(n+1)T\big\}\big) \nonumber \\
    &=T\sum_{n=0}^\infty m_D\big\{ v\in D:\,
    \tau_v^\ve(V)>nT\big\}\nonumber \\
    &\leq Tm_D(D)
    \tilde c_\del^{-1}\exp\big(\frac 1\ve(A_1+\be+\la)\big).
    \nonumber\end{eqnarray}
    In the same way as at the end of the proof of Proposition 
    \ref{prop1.5.2}(i) we fix now an initial point $x=Z_v^\ve(0)\in\cX$ and 
    choose discs $D$ to be small
    balls on the (extended) local unstable manifolds $W^u_x(w,\vrho),\,
    w\in\cW$ which by means of the Fubini theorem and compactness arguments
    enable us to extend (\ref{1.7.21}) and (\ref{1.7.22}) to the case when 
    $m_D$
    is replaced by $m$ and $D$ by $\cW$ yielding (\ref{1.7.3}) and (\ref{1.7.4})
    since $\be$ and $\la$ in (\ref{1.7.22}) can be chosen arbitrarily small as
    $\ve\to 0$. 
    
     Next, we derive the assertion (ii). Let $t>0$ and $n$ be 
    the integral part of $t/T$ where $T>0$ will be chosen later. Let, 
    again, $D\in\cD^u_\ve\big((x,w),\al,\rho,C\big)$ and $x\in V$. Then
    \begin{eqnarray}\label{1.7.23}
    & m_D\{ v\in D:\, Z^\ve_v(\tau^\ve_v)\in\Gam,\,\tau_v^\ve(V)<t\}\\
     &\leq m_D\{ v\in D:\, Z^\ve_v(\tau^\ve_v(V))\in\Gam,\,
     \tau^\ve_v(V)<(n+1)T\}\nonumber\\
     &=\sum_{k=0}^nm_D\{ v\in D:\, Z^\ve_v(\tau^\ve_v(V))
     \in\Gam,\, kT\leq\tau^\ve_v(V)<(k+1)T\}.\nonumber 
    \end{eqnarray}
    Let $K$ be the intersection of the $\om$-limit set of the flow $\Pi^t$
    with $\bar V$. Then $K$ is a compact set and by our assumption 
    $K\subset G$. Hence, 
    \[
    \del=\frac 13\inf\{ |x-z|:\, x\in K,\,\, z\in\bar V\setminus G\}>0
    \]
    and if we set $U_\eta=\{ z\in V:\,\mbox{dist}(z,K)<\eta\}$ then
    $U_{3\del}\subset G$. Now suppose that 
    $kT\leq\tau^\ve_{x,w}(V)<(k+1)T$ for some $k\geq 1$ and
    $Z^\ve_{x,w}(\tau^\ve_{x,w}(V))\in\Gam$ with $x\in V$ and $w\in\cW$.
    Then either there is $t_1\in[(k-1)T,kT]$ such that $Z^\ve_{x,w}(t)\in
    \bar V\setminus U_{2\del}$ for all $t\in[t_1,t_1+T]$ or there exist
    $t_2,t_3>0$ such that $(k-1)T\leq t_2<t_3<(k+1)T$ and $Z^\ve_{x,w}(t_2)
    \in U_{2\del}$ while $Z^\ve_{x,w}(t_3)\in\Gam$. Set $\cT_z=
    \{\gam\in C_{0,2T}:\,\gam_0=z$ and either there is $t_1\in[0,T]$ so that
    $\gam_t\in\bar V\setminus U_{2\del}$ for all $t\in[t_1,t_1+T]$ or 
    $\gam_{t_2}\in U_{2\del}$ and $\gam_{t_3}\in\Gam$ for some
    $0\leq t_2<t_3<2T\}$. Then for any $k\geq 1$,
    \begin{eqnarray}\label{1.7.24}
    &\{ v\in D:\, Z^\ve_v(\tau^\ve_v(V))\in\Gam,\,
    kT\leq\tau^\ve_v(V)<(k+1)T\}\\
    &\subset\big\{ v\in D:\, Z^\ve_v(\tau^\ve_v(V))\in\Gam,\,
    \Phi_\ve^{(k-1)T/\ve}v\in A_0^{\ve,2T}(\cT_{z_{k-1}^\ve(v)})\big\}\nonumber
    \end{eqnarray}
    where $z_k^\ve(v)=\pi_1(\Phi_\ve^{kT/\ve}v)$ and $A^{\ve,T}_0(H)=\{
    w\in V\times\cW:\, Z^\ve_w\in H\}$.
    
      Let $D\subset D_0\in\cD^u_\ve\big((x,y),\al,\rho,C^3\big),\, x\in V$ 
    (later both discs will be small balls on $W^u_x(y,C^6\rho)$) assuming 
    that $\rho$ is small and $C\geq 2$ is large so that $C^6\rho$ is still 
   small. Choose a maximal $\big((k-1)T/\ve,C\rho,\ve,D_0,D_0\big)$-separated 
   set $\tilde E_{k-1}$ in $D_0$ and let
     \[
     E_{k-1}=\big\{ v\in\tilde E_{k-1}:\, U^\ve_{D_0}((k-1)T/\ve,v,C\rho)
     \cap D\ne\emptyset\big\}.
     \]
     Then for $\ve$ small enough,
     \begin{equation}\label{1.7.25}
     D_0\supset \bigcup_{v\in E_{k-1}}U^\ve_{D_0}
     \big((k-1)T/\ve,v,C\rho\big)\supset D
     \end{equation}
     and for any $v,w\in E_{k-1},\,v\ne w$,
     \begin{equation}\label{1.7.26}
      U^\ve_{D_0}\big((k-1)T/\ve,v,C\rho/2\big)\cap
     U^\ve_{D_0}\big((k-1)T/\ve,w,C\rho/2\big)=\emptyset.
      \end{equation} 
     If $v\in E_{k-1},\, w\in U^\ve_{D_0}\big((k-1)T/\ve,v,C\rho\big)$ and
    $\Phi_\ve^{(k-1)T/\ve}w\in A_0^{\ve,2T}(\cT_{z^\ve_{k-1}(w)})$ then by
    Lemma \ref{lem1.3.2}(iii), $|z^\ve_{k-1}(w)-z^\ve_{k-1}(v)|$ is of order
    $\ve$, and so for each $\eta>0$ if $\ve$ is small enough then 
    $\Phi_\ve^{(k-1)T/\ve}w\in
    A_\eta^{\ve,2T}(\cT_{z^\ve_{k-1}(v)})$. For each $q>0$
    set $\cT_z^q=\{\gam\in C_{0,2T}:\gam_0=z\,\,\mbox{and}\,\,\bfr_{0,2T}
    (\gam,\cT_z)\leq q\}$ and suppose that for some $\eta>0$ there is
    $d_\eta\geq 0$ so that
    \begin{equation}\label{1.7.27}
    \inf_{z\in V}\inf_{\gam\in\cT_z^{2\eta}}S_{0,2T}(\gam)> d_\eta .
    \end{equation}
    Then $\cT_z^{2\eta}\cap\Psi^{d_\eta}_{0,2T}(z)=\emptyset$, 
    where $\Psi^a_{0,t}(z)$ is the same as in Theorem \ref{thm1.2.3}, and so
    \begin{equation}\label{1.7.28}
    \cT^\eta_z\subset\big\{\gam\in C_{0,2T}:\,\gam_0=z\,\,\mbox{and}\,\,
    \bfr_{0,2T}(\gam,\Psi^{d_\eta}_{0,2T}(z))\geq\eta\big\}.
    \end{equation}
    Hence,
    \begin{equation}\label{1.7.29}
    A_\eta^{\ve,2T}(\cT_z)\subset\big\{ (z,w)\in V\times\cW:\,\bfr_{0,2T}
    (Z^\ve_{(z,w)},\Psi^{d_\eta}_{0,2T}(z))\geq\eta\big\}.
    \end{equation}
    By Lemma \ref{lem1.3.2}(ii),
    \begin{equation*}
    D_{k-1}(v)=\Phi_\ve^{(k-1)T/\ve}U^\ve_{D_0}\big((k-1)T/\ve,v,C\rho\big)
    \in\cD^u_\ve\big(\Phi_\ve^{(k-1)T/\ve}v,\al,\rho,\sqrt C\big),
    \end{equation*}
    and so applying (\ref{1.7.11}) to $D_{k-1}(v)$ we obtain from (\ref{1.7.27})--
    (\ref{1.7.29}) that for any $\be>0$ and sufficiently small $\ve$ uniformly
    in discs $D_{k-1}(v)$ as above,
    \begin{equation*}
    m_{D_{k-1}(v)}\big(A_\del^{\ve,2T}(\cT_{z^\ve_{k-1}(v)})\big)\leq
    \exp(-(d_\eta-\be)/\ve).
    \end{equation*}
    This together with Lemma \ref{lem1.3.6} yield that for each $v\in E_{k-1}$,
    \begin{eqnarray*}
    &m_{D_0}\big\{\tilde v\in U^\ve_{D_0}\big((k-1)T/\ve,v,C\rho\big):\,
    \Phi_\ve^{(k-1)T/\ve}\tilde v\in A_\eta^{\ve,2T}(\cT_{z_{k-1}^\ve(v)})
    \big\}\\
    &\leq\tilde Ce^{-(d_\eta-\be)/\ve}m_{D_0}
    \big(U^\ve_{D_0}((k-1)T/\ve,v,C\rho)\big)\nonumber
    \end{eqnarray*}
    for some $\tilde C>0$ depending only on $C\rho$. Combining this with
     (\ref{1.7.24})--(\ref{1.7.26}) and Lemma \ref{lem1.3.6}
     we obtain that for any $k\geq 1$,
     \begin{equation}\label{1.7.30}
     m_D\{ v\in D:\, Z^\ve_v(\tau^\ve_v)\in\Gam,\,
     kT\leq\tau^\ve_v<(k+1)T\}\leq\hat Ce^{-(d_\eta-\be)/\ve}
     \end{equation}
     for some $\hat C>0$ depending only on $C\rho$.
     
     Next, we will specify $d_\eta$ in (\ref{1.7.27}) choosing 
     $\eta\leq\frac 12\del$. For each $z\in V$ we can write
    \begin{equation}\label{1.7.31}
    \cT_z^{2\eta}\subset\tilde\cT_z^\eta\cup\hat\cT_z^\eta
    \end{equation}
    where $\tilde\cT_z^\eta=\{\gam\in C_{0,2T}:\,\gam_0=z,\,
    \gam_{t_2}\in U_{3\del}$ and $\gam_{t_3}\in\Gam_{2\eta}$ for some 
    $0\leq t_2<t_3<2T\}$ with $\Gam_{r}=\{z:$ dist$(z,\Gam)\leq r\}$ and 
    $\hat\cT_z=\{\gam\in C_{0,2T}:\,\gam_0=z$ and there is $t_1\in[0,T]$ so
    that $\gam_t\in V_{2\eta}\setminus U_\del$ for all $t\in[t_1,t_1+T]\}$.
    By (\ref{1.7.5}) and the lower semicontinuity of the functional $S_{0,2T}$
    it follows that for any $\zeta>0$ we can choose $\eta>0$ small enough
    so that 
    \begin{equation}\label{1.7.32}
    \inf_{z\in V}\inf_{\gam\in\tilde\cT_z^\eta}S_{0,2T}(\gam)> A_2-\zeta.
    \end{equation}
    
    Since $\bar V\setminus U_\del$ is disjoint with the $\om$-limit set of the
    flow $\Pi^t$ and the latter is closed then if $\eta$ is sufficiently small
    $V_{2\eta}\setminus U_\del$ is also disjoint with this $\om$-limit set
    and, in particular, it does not contain any forward semi-orbit
    of $\Pi^t$. Hence we can apply Lemma \ref{lem1.6.4} which in view of
    (\ref{1.2.13}) implies that there exists $a>0$ such that for all small
     $\eta>0$,
    \begin{equation}\label{1.7.33}
    \inf_{z\in V}\inf_{\gam\in\hat\cT_z}S_{0,2T}(\gam)> aT
    \end{equation}
    which is not less than $A_2$ if we take $T=A_2/a$. Now, (\ref{1.7.32}) and 
    (\ref{1.7.33}) produce (\ref{1.7.27}) with $d=A_2-\zeta$, and so 
    (\ref{1.7.30}) follows with such $d_\eta$.
    This together with (\ref{1.7.23}) yield that for any $\be>0$ we can choose
    sufficiently small $\zeta,\la >0$ and then $\eta>0$ so that for all $\ve$
    small enough
    \begin{eqnarray}\label{1.7.34}
    &m_D\big\{ v\in D:\, Z^\ve_v(\tau^\ve_v)\in\Gam,\,
    \tau^\ve_v(V)\leq e^{(A_2-\be)/\ve}\big\}\\
    &\leq m_D\{ v\in D:\, Z^\ve_v(\tau^\ve_v)\in\Gam,\,\tau_v^\ve(V)<T\}
    +e^{-\la/2\ve}.\nonumber
    \end{eqnarray}
    
    Now assume that (\ref{1.7.7}) holds true for some $x\in V$. 
    Recall, that $S_{0T}(\gam)=0$ implies that $\gam$ is a piece of an
    orbit of the flow $\Pi^t$. Since no $\gam\in C_{0T}$ satisfying
    \begin{equation}\label{1.7.35}
    \gam_0=x\,\,\mbox{and}\,\,\inf_{t\in[0,T]}\,\mbox{dist}
    (\gam_t,\partial V)\leq a(x)/2
    \end{equation}
    can be such piece of an orbit we conclude by the lower semicontinuity of
    $S_{0T}$ that $S_{0T}(\gam)> c(x)$ whenever (\ref{1.7.35}) holds true 
    for some $c(x)>0$ independent of $\gam$ (but depending on $x$).
     Hence, by (\ref{1.7.11}), 
    \begin{eqnarray}\label{1.7.36}
    &m_D\{ v\in D:\,\tau^\ve_v<T\}\leq m_D\big\{ v\in D:\,\bfr_{0T}
    \big(Z^\ve_v, \Psi^{c(x)}_{0T}(x)\big)\\
    &\geq a(x)/2\big\}\leq\exp(-c(x)/2\ve)\nonumber
    \end{eqnarray}
    provided $\ve$ is small enough and (\ref{1.7.8}) follows. 
    Observe also that any $\gam\in C_{0t}$
    with $\gam_0=x\in V$ and $\gam_t\in\partial V$ should contain a piece
    which either belongs to some $\tilde\cT_z^\eta$ or $\hat\cT_z^\eta$, 
    as above, or to satify (\ref{1.7.35}). By (\ref{1.7.32}), (\ref{1.7.33}), 
    and 
    the above remarks it follows that $S_{0t}(\gam)\geq q(x)$ for such $\gam$
    where $q(x)>0$ depends only on $x$, and so $R_\partial(x)\geq q(x)$.
    
    Finally, similarly to Proposition \ref{prop1.5.2} we fix $x=Z^\ve_v(0)\in V$,
    choose discs $D$ and $D_0$ to be small balls on the (extended) local
    unstable manifolds $W^u_x(w,q),\, w\in\cW,q>0$ and using the Fubini theorem 
    we extend (\ref{1.7.34}) and (\ref{1.7.36}) to the
    case when $D$ and $m_D$ are replaced by $\cW$ and $m$, respectively,
    yielding (\ref{1.7.6}). If $\Gam=\partial V$ then by (\ref{1.7.6}) and 
    (\ref{1.7.8}),
    \begin{eqnarray}\label{1.7.37}
    &\int_{\cW}\tau^\ve_{x,y}(V)dm(y)\geq e^{(A_2-\be)/\ve}\, m\big\{ y\in\cW:
    \, \tau^\ve_{x,y}(V)\geq e^{(A_2-\be)/\ve}\big\}\\
    &\geq e^{(A_2-\be)/\ve}(1-e^{-\la(\be)/\ve}-e^{-\hat\la(T)/\ve})\nonumber
    \end{eqnarray}
    and, since $\be>0$ is arbitrary, (\ref{1.7.9}) follows completing the
    proof of Proposition \ref{prop1.7.1}.
    \end{proof}

   Now we will derive Theorem \ref{thm1.2.5} from Proposition \ref{prop1.7.1}.
   Assume, first, that $R_\partial <\infty$. Then by Lemma \ref{1.6.4}, 
   $R_\partial(x)$ is finite and continuous in the whole $V$. Moreover,
   since $\cO$ is an $S$-attractor the conditions of Lemma \ref{lem1.6.5}
   are satisfied with some $d(t)\to 0$ as $t\to\infty$ the same for all
   points of $V$ which yields the conditions of Proposition \ref{prop1.7.1}(i)
   with $A_1=R_\partial +\del$ for any $\del>0$. Hence, (\ref{1.7.3}) and 
   (\ref{1.7.4})
   hold true with $A_1=R_\partial$. Since $\cO$ is an $S$-attractor of the
   flow $\Pi^t$ and its basin contains $\bar V$ then the intersection of
   $\bar V\setminus\cO$ with the $\om$-limit set of $\Pi^t$ is empty.
   By the definition of an $S$-attractor for any $\eta>0$ there exists
   an open set $U_\eta\supset\cO$ such that $R(x,z)\leq\eta$ whenever
   $x\in\cO$ and $z\in U_\eta$. Hence, by the triangle inequality for
   the function $R$ and Lemma \ref{lem1.6.5} for any set 
   $\Gam\subset\partial V$,
   \begin{equation}\label{1.7.38}
   \inf_{z\in U_\zeta,\tilde z\in\Gam}R(z,\tilde z)\geq\inf_{\tilde z\in\Gam}
   R^\cO(\tilde z)-\eta.
   \end{equation}
   If $\Gam=\partial V$ then by Lemma \ref{lem1.6.5} the right hand side of 
   (\ref{1.7.38}) equals
   $A_2=R_\partial-\eta$. Assuming that $R_\partial<\infty$ we can apply
   Proposition \ref{prop1.7.1}(ii) with such $A_2$ yielding (\ref{1.7.6}), 
   (\ref{1.7.8}) and since $\eta>0$ is arbitrary (\ref{1.2.21}) and 
   (\ref{1.2.22})
   follow in this case. If $R_\partial=\infty$ then (\ref{1.2.22}) is trivial 
   and by (\ref{1.7.38}), $R(z,\tilde z)=\infty$ for any $z\in U_\zeta$ and
   $\tilde z\in\partial V$, and so we can apply Proposition \ref{prop1.7.1}(ii)
   with any $A_2$ which sais that the left hand side in (\ref{1.7.9}) equals
   $\infty$, and so (\ref{1.2.21}) holds true in this case, as well.
   
    Next, we establish (\ref{1.2.23}). For small $\del,\be>0$ and large $T>0$
  which will be specified later on set $\Gam_1=\{ v\in V_\del\times\cW:\,
  Z^\ve_v(T)\in V\setminus U_{\del/2}(\cO)\}$, $\Gam_2=\{ v\in U_{\del/2}(\cO)
  \times\cW:\,\tau^\ve_v(U_\del(\cO))\leq e^{\be/\ve}\}$ and $t_\ve=
  T+e^{\be/\ve}$. Then
  \begin{eqnarray}\label{1.7.39}
  &\Te^\ve_v((n+1)t_\ve\wedge\tau^\ve_v(V))-
  \Te^\ve_v(nt_\ve\wedge\tau^\ve_v(V))\\
  &\leq T+t_\ve\big(\bbI_{\Gam_1}(\Phi_\ve^{t_\ve n/\ve}v)+
  \bbI_{V_\del\times\cW\setminus\Gam_1}(\Phi_\ve^{t_\ve n/\ve}v)\bbI_{\Gam_2}
  (\Phi_\ve^{t_\ve n/\ve}v)\big).
  \nonumber\end{eqnarray}
  If $\del$ is sufficiently small then $V_\del$ is still contained in the
  basin of $\cO$ with respect to the flow $\Pi^t$, and so we can choose
  $T$ (depending only on $\del$) so that
  \[
  \Pi^TV_\del\subset U_{\del/4}(\cO).
  \]
  Then for some $a>0$,
  \[
  \inf\big\{ S_{0T}(\gam):\,\gam\in C_{0T},\,\gam_0\in V_\del,\,\gam_T\not\in
  U_{\del/3}(\cO)\big\}>a,
  \]
  and so if $\gam_0\in V_\del$ and $\gam_T\not\in U_{\del/2}(\cO)$ then 
  dist$(\gam,\Psi^a_{0T}(z))\geq\del/6$ for any $z\in V_\del$. Relying on
  (\ref{1.7.11}) we obtain that for any $\tilde D\in\cD^u_\ve(z,\al,\rho,
  \sqrt C)$ with $\tilde D\subset V_\rho\times\cW$,
  \begin{equation}\label{1.7.40}
  m_{\tilde D}(\Gam_1\cap\tilde D)\leq e^{-a/2\ve}
  \end{equation}
  provided $\ve$ is small enough. Next, the same arguments which yield 
  (\ref{1.7.34}) and (\ref{1.7.36}) enable us to conclude that if $\be>0$ is
  small enough then for any $\tilde D\in\cD^u_\ve(z,\al,\rho,\sqrt C)$
  with $\tilde D\subset U_{\del/2}(\cO)\times\cW$,
  \begin{equation}\label{1.7.41}
   m_{\tilde D}(\Gam_2\cap\tilde D)\leq e^{-\be/\ve}.
  \end{equation}
   Now let $D\in\cD^u_\ve(z,\al,\rho,C)$, $z=(x,y)$, $D\subset V_\del\times
   \cW$ and $D\subset D_0\in\cD^u_\ve(z,\al,\rho,C^3)$ with $\rho$ small and
   $C$ large so that $C^6\rho$ is still small. Let $E_n^{(1)}$ and $E_n^{(2)}$
   be maximal $(nt_\ve\ve^{-1},C\rho,\ve,D,D_0)-$ and 
   $((nt_\ve+T)\ve^{-1},C\rho,\ve,D,D_0)-$separated sets, respectively. Then
   \[
   \cup_{v\in E^{(1)}_n}U^\ve_D(nt_\ve\ve^{-1},v,C\rho)\supset D
   \]
   and since the last union is contained in a small neighborhood of $D$ and
   $U^\ve_D(nt_\ve\ve^{-1},v,C\rho/2)$ are disjoint for different 
   $v\in E^{(1)}_n$ we obtain using Lemma \ref{lem1.3.6} that
   \[
   \sum_{v\in E_n^{(1)}}m_D\big(U^\ve_D(nt_\ve\ve^{-1},v,C\rho)\big)\leq 
   2c^{-1}_{C\rho}c^{-1}_{C\rho/2}m_D(D).
   \]
   Similarly,
   \[
   \cup_{v\in E^{(2)}_n}U^\ve_D((nt_\ve+T)\ve^{-1},v,C\rho)\supset D
   \]
   and 
   \[
   \sum_{v\in E_n^{(2)}}m_D\big(U^\ve_D((nt_\ve+T)\ve^{-1},v,C\rho)\big)\leq 
   2c^{-1}_{C\rho}c^{-1}_{C\rho/2}m_D(D).
   \]
  Since by Lemma \ref{lem1.3.2}(ii),
  \[
  \Phi^t_\ve U^\ve_D(t,v,C\rho)\in\cD(\Phi^t_\ve v,\al,\rho,\sqrt C)
  \]
  we can apply (\ref{1.7.39})--(\ref{1.7.41}) together with Lemma \ref{lem1.3.6}
  (similarly to the proof of (\ref{1.7.30}))
  in order to conclude that for sufficiently small $\be$ and any much smaller
  $\ve$,
  \[
  \int_D\big(\Te^\ve_v((n+1)t_\ve\wedge\tau^\ve_v(V))-
  \Te^\ve_v(nt_\ve\wedge\tau^\ve_v(V))\big)dm_D(v)\leq t_\ve e^{-\be/\ve}(T+1).
  \]
  Choosing discs $D$ to be small balls on the (extended) local unstable 
  manifolds $W^u_x(w,\vrho),\, w\in\cW$ together with the Fubini theorem 
  we extend this estimate to
  \begin{equation}\label{1.7.42}
  \int_\cW\big(\Te^\ve_v((n+1)t_\ve\wedge\tau^\ve_v(V))-
  \Te^\ve_v(nt_\ve\wedge\tau^\ve_v(V))\big)dm(v)\leq\tilde Ct_\ve e^{\be/\ve}
  \end{equation}
  for some $\tilde C>0$ depending on $\del$ but independent of $n$ and $\ve$. 
  Finally, (\ref{1.2.22}) and (\ref{1.7.42}) together with the Chebyshev 
  inequality
  yield that for $n(\ve)=[e^{(R_\partial+\be/4)/\ve}t_\ve^{-1}]$, each $x\in V$,
  a small $\be>0$ and any much smaller $\ve>0$,
  \begin{eqnarray}\label{1.7.43}
  &m\big\{ w\in\cW:\,\Te^\ve_{x,w}(\tau^\ve_{x,w}(V))\geq e^{-\be/4\ve}
  \tau^\ve_{x,w}(V)\big\}\\
  &\leq m\big\{ w\in\cW:\,\Te^\ve_{x,w}((n(\ve)+1)t_\ve)
  \geq e^{-\be/4\ve}e^{(R_\partial-\be/4)/\ve}\big\}\nonumber\\
  &+m\big\{ w\in\cW:\,
  \tau_{x,w}^\ve(V)<e^{(R_\partial-\be/4)/\ve}\,\,\mbox{or}\,\,
  \tau_{x,w}^\ve(V)>e^{(R_\partial+\be/4)/\ve}\big\}\nonumber\\
  &\leq\tilde Ce^{-\be/4\ve}\big(1+e^{-(R_\partial+\be/4)/\ve}
  (T+e^{\be/\ve})\big)+e^{-\la(\be/4)/\ve}.\nonumber
  \end{eqnarray}
   Since $R_\partial>0$ and we can choose $\be$ to be arbitrarily small, 
 (\ref{1.7.43}) yields (\ref{1.2.23}).
 
 In order to complete the proof of Theorem \ref{thm1.2.5} we have to derive
  (\ref{1.2.24}). If $\partial_{\min}=\partial V$ then 
    there is nothing to prove, so we assume that $\partial_{\min}$ is a
    proper subset of $\partial V$ and in this case, clearly,
    $R_\partial<\infty$. Since $\Gam=\{ z\in\partial V:\,\mbox{dist}
    (z,\partial_{\min})\geq\del\}$ is compact and disjoint with 
    $\partial_{\min}$ which is also compact then by the lower semicontinuity
    of $R^\cO(z)$ established in Lemma \ref{lem1.6.5}(iii) it follows that
    $R^\cO(z)\geq R_\partial +\be$ for some $\be>0$ and all $z\in\Gam$.
    Then by (\ref{1.7.38}), $R(z,\tilde z)\geq R_\partial+\be/2$ for any
    $z\in U_{\be/2}$ and $\tilde z\in\Gam$. Hence, applying Proposition 
    \ref{prop1.7.1} we obtain that
    \[
  m\big\{ y\in\cW:\,\tau^\ve_{x,y}(V)\geq e^{(R_\partial+\frac 13\be)/\ve}
  \big\}\leq e^{-\la/\ve}
  \]
  and
  \[
  m\big\{ y\in\cW:\, Z^\ve_{x,y}(\tau^\ve_{x,y}(V))\in\Gam,\,
  \tau^\ve_{x,y}(V)\leq
  e^{(R_\partial+\frac 13\be)/\ve}\big\}< e^{-\la/\ve}
  \]
  for some $\la>0$ and all $\ve$ small enough yielding (\ref{1.2.24})
  and completing the proof of Theorem \ref{thm1.2.5}. 
   \qed

 \section{Adiabatic transitions between basins of attractors}\label{sec1.8}
\setcounter{equation}{0}
    In this section we will prove Theorem \ref{thm1.2.7} relying, again, on 
    Proposition \ref{prop1.7.1} together with "Markov property type"
    arguments and at the end of the proof we will apply even some rough "strong
 Markov property type" \index{strong Markov property}
 arguments in order to deal with subsequent transitions
 between basins of attractors. In view of (\ref{1.2.27}) and Lemma 
 \ref{lem1.6.2}{i}
 any curve $\gam\in C_{0t}$ starting at $\gam_0=x\in V_{j_1}$ and ending at
    $\gam_t=z\in\cap_{1\leq i\leq k}\partial V_{j_i},\, k\leq\ell$ can be
    extended into each $V_{j_i},\, i=1,...,k$ with arbitrarily small increase
    in its $S$-functional. Hence, 
    \begin{equation}\label{1.8.1}
    R^{(i)}_\partial=\min_{j\ne i}R_{ij}
    \end{equation}
    where $R^{(i)}_\partial=\inf\{ R(x,z):\, x\in\cO_i,\, z\in\partial V_i\}$.
    Let $Q$ be an open ball of radius at least $r_0$ centered at the origin
    of $\bbR^d$. By Assumption \ref{ass1.2.6} the slow motion $Z^\ve_{x,y}$
    cannot exit $Q$ provided $x\in Q$ and $y\in\cW$. Furthermore, it is clear
    that $Q$ contains the $\om$-limit set of the averaged flow $\Pi^t$.
    Assumption \ref{ass1.2.6} enables us to deal only with restricted basins
    $V^Q_i=V_i\cap Q$ and though the boundaries $\partial V^Q$ of $V_i^Q$
    may include now parts of the boundary $\partial Q$ of $Q$ it makes no
    difference since $Z^\ve$ cannot reach $\partial Q$ if it starts in $Q$.
    Set $V^{(i)}=Q\setminus\cup_{j\ne i}U_\del(\cO_j)$ where $\del>0$ is
    small enough. We claim that in view of (\ref{1.2.27}) each $V_i$ satisfies 
  conditions of Proposition \ref{prop1.7.1}(i) for any $\be>0$ with $A_1=
  R_\partial^{(i)}+\be$ and some $T=T_\be$ depending on $\be$. Indeed, set 
 \[
 \partial(\eta)=\{ v\in Q:\,\mbox{dist}(v,\cup_{1\leq j\leq\ell}\partial V_j)
 \leq\eta\},\,\eta>0.
 \]
 In view of (\ref{1.2.9}) and (\ref{1.2.27}) there exists $L>0$ such that if
 $\eta$ is small enough and $z\in\partial(\eta)$ we can construct a curve
 $\vf^z\in C_{0,L\eta}$ with $S_{0,L\eta}(\vf^z)\leq\tilde L\eta$, 
 $\vf^z_0=z,\,\vf_t^z\in V_j\setminus\partial(\eta)$ for some $t\in[0,L\eta]$
  and $j=1,...,\ell$ where $\tilde L=L\sup_{x,y}|\vf^u_x(y)|$.
   Since $V_j$ is the basin of $\cO_j$ there exists 
 $T=T_{\eta,\del}$ such that $\Pi^T\vf^z_t\in U_\del(\cO_j)$ and extending
 $\vf^z$ by the piece of the orbit of $\Pi^t$ we obtain a curve $\tilde\vf^z
 \in C_{0,L\eta+T}$ starting at $z$, entering $U_\del(\cO_j)$ and satisfying
 $S_{0,L\eta+T}(\tilde\vf^z)\leq\tilde L\eta$. Hence, for $z\in\partial(\eta)$
 the condition (\ref{1.7.2}) holds true with $V=V^{(i)}$ and $A_1=\tilde L\eta$.
 Since the $\om$-limit set of the flow $\Pi^t$ is contained in $Q\cap
 \big(\cup_{1\leq j\leq\ell}(\partial V_j\cup\cO_j)\big)$ it follows from
 Assumption \ref{ass1.2.6} and compactness considerations that there exists
 $\tilde T=\tilde T_{\eta,\delta}$ such that for any $z\in Q\setminus V_i$
 we can find $t_z\in[0,\tilde T]$ with $\Pi^{t_z}z\in\partial(\eta)\cup
 \big(\cup_{j\ne i}U_\del(\cO_j)\big)$. If $\Pi^{t_z}z\in
 \cup_{j\ne i}U_\del(\cO_j)$ then we take $\vf^z_t=\Pi^tz,\, t\in[0,\tilde T]$
 to satisfy (\ref{1.7.2}) for $V=V^{(i)}$ and $A_1=0$. If $\Pi^{t_z}z\in
 \partial(\eta)$ then we extend the curve $\vf^z_t=\Pi^tz,\, t\in[0,t_z]$
 as in the above argument which yields a curve $\tilde\vf^z$ starting at
 $z$, ending in some $U_\del(\cO_j),\, j\ne i$ and having its $S$-functional
 not exceeding $\tilde L\eta$. Finally, in the same way as in the proof of 
 Theorem \ref{thm1.2.5} for any $\be>0$ there exists $\hat T=\hat T_{\eta,\del,
 \be}$ such that whenever $z\in V_i(\eta)=V_i\cap Q\setminus\partial(\eta)$
 we can construct $\vf^z\in C_{0\hat T}$ such that (\ref{1.7.2}) holds true 
 with $V=V_i(\eta)$ and $A_1=R^{(i)}_\partial +\be/2$ and, moreover,
 dist$(\vf_t^z,V_j)\leq\eta$ for some $t\leq\hat T$ and $j\ne i$ with $R_{ij}
 =R^{(i)}_\partial$. Then in the same way
 as above we can extend $\vf^z$ to some $\tilde\vf^z\in C_{\hat T+\tilde T}$
 so that $\tilde\vf_t^z\in U_\del(V_j)$ for some $j$ as above,
 $t\leq\hat T+\tilde T$ and $S_{0,\hat T+\tilde T}(\tilde\vf^z)
 \leq R^{(i)}_\partial+\be/2+\tilde L\eta$ which gives (\ref{1.7.2}) for all 
 $z\in V=V^{(i)}$ with $A_1=R^{(i)}_\partial +\be$ provided $\eta$ is small 
 enough. Hence, Proposition \ref{prop1.7.1}(i) yields the estimates 
 (\ref{1.7.3}) 
 and (\ref{1.7.4}) for $\tau^\ve_{x,y}(i)$ in place of $\tau^\ve_{x,y}(V)$ with 
   $A_1=R_\partial^{(i)}$. In order to obtain the corresponding bounds in the
   other direction observe that in view of (\ref{1.2.27}),
   \begin{equation}\label{1.8.2}
   R^{(i)}_\partial(\del)=\inf\{ R(x,z):\, x\in\cO_i,\, z\not\in V_i(\eta)\}
   \to R^{(i)}_\partial\,\,\mbox{as}\,\,\del\to 0.
   \end{equation}
    Since $\overline {V_i(\eta)}$ is contained in the basin of $\cO_i$ we
    can apply to $V_i(\eta)$ the same estimates as in Theorem \ref{thm1.2.5}
    which together with (\ref{1.8.2}) and the fact that the exit time of
    $Z^\ve$ from $V_i(\eta)$ is smaller than its exit time from $V_i$ provide
    the remaining bounds yielding (\ref{1.2.28}) and (\ref{1.2.29}).
    
    Next, we derive (\ref{1.2.30}) similarly to (\ref{1.2.23}) but taking into
    account that $\cup_{1\leq j\leq\ell}\partial V_j$ may contain parts of
    the $\om$-limit set of the flow $\Pi^t$ which allows the slow motion
    $Z^\ve$ to stay long time near these boundaries. Still, set
    \[
    \te^\ve_v=\inf\{ t\geq 0:\, Z^\ve_v(t)\in\cup_{1\leq j\leq\ell}U_{\del/3}
    (\cO_j)\}.
    \]
    Using the same arguments as above we conclude that for any $\eta>0$
    there exists $T=T_{\eta,\del}$ such that whenever $z\in Q$ we can 
    construct $\vf^z\in C_{0T}$ with $\vf^z_0=z,\,\vf_T^z\in\cup_{1\leq j
    \leq\ell}U_\del(\cO_j)$ and $S_{0T}(\vf^z)\leq\eta$. This together with
    (\ref{1.7.21}) and Assumption \ref{ass1.2.6} yield that for any disc
    $D\in\cD^u_\ve(z,\al,\rho,C),\, D\subset Q\times\cW$,
    \[
    m_D\{ v\in D:\,\te^\ve_v\geq e^{2\eta/\ve}\}\leq e^{-\la(\eta)/\ve}
    \]
    for some $\la(\eta)=\la(x,\eta)>0$ and all small $\ve$. Set 
    \[
    \Gam_1=\big\{ v\in Q\times\cW:\, Z^\ve_v(e^{2\eta/\ve})\in Q\setminus
    \cup_{1\leq j\leq\ell}U_{\del/2}(\cO_j)\big\},
    \]
    \[
    \Gam_2=\big\{ v\in\cup_{1\leq j\leq\ell}U_{\del/2}(\cO_j):\,\tau^\ve_v
    \big(\cup_{1\leq j\leq\ell}U_\del(\cO_j)\big)\leq e^{\be/\ve}\big\}
    \]
    and $t_\ve=e^{2\eta/\ve}+e^{\be/\ve}$ where $\eta$ is much smaller than 
    $\be$. Then proceeding similarly to the proof of (\ref{1.2.23}) as in
    (\ref{1.7.40})--(\ref{1.7.43}) above we arrive at (\ref{1.2.30}).
     
    Next, we obtain (\ref{1.2.31}) relying on additional assumptions specified
    in the statement of Theorem \ref{thm1.2.7}. Let $V^Q_i$ be the same as 
    above and $\partial^{(i)}_0(x)=\{ z\in\partial V_i^Q:\, 
  R(x,z)=R_\partial^{(i)}\}$. Since $\cO_i$ is an $S$-attractor it follows from
  Lemma \ref{lem1.6.5}(i) that $R(x,z)$ and $\partial_0^{(i)}(x)$ coincide 
  with 
  the same function $R^{\cO_i}(z)$ and the same (in general, may be empty) set 
   $\partial_0^{(i)}$, respectively, for all $x\in\cO_i$. By Lemma 
   \ref{lem1.6.2}(i), our assumption that $B$ is complete on $\partial V_{i}$
   implies that $R^{\cO_{i}}(z)$ is continuous in a neighborhood of
   $\partial V_{i}$, and so $\partial_0^{(i)}$ is a nonempty compact set.
   Since we assume that $\io(i)\ne i$ is the unique index $j$ for which 
   $R_{ij}=R_{i\io(i)}=R_\partial^{(i)}$ then by (\ref{1.2.27}),
    \begin{equation*}
    \min_{j\ne i,\io(i)}\inf_{z\in\partial_0^{(i)}}\mbox{dist}
    (z,\partial V_j)>0.
    \end{equation*}
    Observe that if $\tilde\cO\subset\partial V_{i}$ is an $S$-compact then
    either $\tilde\cO\subset\partial_0^{(i)}$ or $\tilde\cO\cap
    \partial_0^{(i)}=\emptyset$. Denote by $L_\Pi$ the $\om$-limit set of
    the averaged flow $\Pi^t$. Since $L_\Pi\cap\partial V_{i}$ consists
    of a finite number of $S$-compacts it follows that
    \[
    \inf\{ |z-\tilde z|:\, z\in L_\Pi\cap\partial_0^{(i)},\,\tilde z\in
    L_\Pi\setminus\partial_0^{(i)}\}>0.
    \]
    By the continuity of $R^{\cO_{i}}(z)$ in $z\in\partial 
    V_{i}$ there exists $a>0$ such that 
    \[
    \inf\big\{  R^{\cO_{i}}(z):\, z\in\big(\cup_{j\ne i,\io(i)}(\partial 
    V_{i}\cap\partial V_j)\big)\cup\big((L_\Pi\setminus\partial_0^{(i)})
    \cap\partial V_{i}\big)\big\}\geq R_\partial^{(i)}+9a.
    \]
   These considerations enable us to construct a connected open set $G$
   with a piecewise smooth boundary $\partial G$ such that
   \[
   \bar G\subset V_{i}\cup(V_{\io(i)}\setminus\cO_{\io(i)})
   \cup\big((\partial V_{i}\cap
   \partial V_{\io(i)})\setminus (L_\Pi\setminus\partial_0^{(i)})\big)
   \]
   and for $\Gam=\partial G\setminus U_\del(\cO_{\io(i)})$ and some 
   $a(\del)>0$,
   \begin{equation}\label{1.8.3}
   \inf_{z\in\Gam}R^{\cO_{i}}(z)\geq R_\partial^{(i)}+8a
   \end{equation}
   provided $a\leq a(\del)$. The idea of this construction is that if
   $Z^\ve_{x,y}(\tau^\ve_{x,y}(i))\not\in V_{\io(i)}$ then the slow motion 
   should exit $G$ through the part $\Gam$ of its boundary. Somewhat
   similarly to the proof of Proposition \ref{prop1.7.1}(ii) we will show
   that for "most" initial conditions $y$ this can only occur after the
   time $\exp\big((R_\partial^{(i)}+2a)/\ve\big)$ and, on the
   other hand, we conclude from (\ref{1.2.29}) that for "most" initial 
   conditions $y$ the exit time $\tau^\ve_{x,y}(i)$ does not exceed
   $\exp\big((R_\partial^{(i)}+a)/\ve\big)$.
   
   Let $U_0$ be a sufficiently small open neighborhood of $\partial_0^{(i)}$ 
   so that, in particular,
   \[
   \sup_{z\in U_0}R^{\cO_{i}}(z)\leq R_\partial^{(i)}+a
   \]
   and set
   \[
   \tau^\ve_{x,y}(G)=\inf\{ t\geq 0:\, Z^\ve_{x,y}(\tau^\ve_{x,y}(G))
   \not\in G\}.
   \]
   For each disc $D\in\cD^u_\ve((x,w),\al,\rho,C)$ we can write
   \begin{eqnarray}\label{1.8.4}
   &\,\,\,\,\,\,\,\,\big\{ v\in D:\,\tau^\ve_v(G)\leq e^{(R_\partial^{(i)}
   +a)/\ve},\,
   Z^\ve_v(\tau^\ve_v(G))\in\Gam\big\}\subset\bigcup_{0\leq n\leq n(\ve)+1}
   \big(A^{(1)}_D(n)\cup\\
   &\bigcup_{(n-1)t_\ve\leq k\leq(n+1)t_\ve}\big( A^{(2)}_D(k)+A^{(3)}_D(k)
   +\bigcup_{k-2t_\ve\leq m\leq k-2T} A^{(4)}_D(m)\cap A^{(5)}_D(k)\big)\big)
   \nonumber\end{eqnarray}
   where $t_\ve=e^{\be/\ve}$ for some small $\be>0$, $n(\ve)=
   \big[e^{(R_\partial^{(i)}+a-\be)/\ve}\big]$, 
   $A^{(1)}_D(n)=\{ v\in D:\, Z^\ve_v(t)\in G\setminus\big( U_\eta(\cO_{i})
   \cup U_\del(\cO_{\io(i)})\big)\,\,\mbox{for all}\,\, t\in[(n-1)t_\ve,nt_\ve]
   \big\}$ for a sufficiently small $\eta>0$, $A^{(2)}_D(k)=\big\{ v\in D:\,
   \exists t_1,t_2\,\mbox{with}\, k\leq t_1<t_2<k+3T,\, Z^\ve_v(t_1)
   \in U_\eta(\cO_{i}),\, Z^\ve_v(t_2)\in\Gam\big\}$,\,
   $A^{(3)}_D(k)=\{ v\in D:\, Z^\ve_v(t)\in G\setminus\big(U_0\cup
    U_\eta(\cO_{i})\cup U_\del(\cO_{\io(i)})\big)\,\,\mbox{for all}\,\, 
    t\in[k,k+T]\big\}$,\, $A^{(4)}_D(m)=\big\{ v\in D:\,
   \exists t_1,t_2\,\mbox{with}\, m\leq t_1<t_2<m+T,\, Z^\ve_v(t_1)
   \in U_\eta(\cO_{i}),\, Z^\ve_v(t_2)\in U_0\big\}$,\, and
   $A^{(5)}_D(k)=\big\{ v\in D:\,\exists t_3,t_4\,\mbox{with}\, 
   k\leq t_3<t_4<k+T,\, Z^\ve_v(t_3)\in U_0,\, Z^\ve_v(t_4)\in\Gam\big\}$.
   Observe that $G\setminus\big(U_\eta(\cO_{i})\cup U_\del(\cO_{\io(i)})\big)$
   satisfies conditions of Proposition \ref{prop1.7.1}(i) with arbitrarily 
   small $A_1$, so similarly to (\ref{1.7.21}) (and taking into account
   Lemma \ref{lem1.3.6}) we can estimate
   \begin{equation}\label{1.8.5}
   m_D(A^{(1)}_D(n))\leq\exp(-\frac 12e^{\be/\ve}).
   \end{equation}
   Similarly to the proof of Proposition \ref{prop1.7.1}(ii) we obtain also 
   that
   \begin{equation}\label{1.8.6}
   \max\big(m_D(A^{(2)}_D(k)),m_D(A^{(3)}_D(k))\big)\leq 
   e^{-(R_\partial^{(i)}+3a)/\ve}
   \end{equation}
   where we, first, choose $\eta$ small and then $T$ large enough. 
   
   Next, we estimate $m_D\big(A^{(4)}_D(m)\cap A^{(5)}_D(k)\big)$ for 
   $m\leq k-2T$ by the following Markov property type argument. Let
   $D\subset D_0\in\cD^u_\ve((x,y),\al,\rho,C^3)$ and choose a maximal
   $(k/\ve,C\rho,\ve,D,D_0)$-separated set $E$ in $D$. Let
   $\tilde E=\big\{ v\in E:\, U^\ve_{D_0}(k/\ve,v,C\rho)\cap A^{(4)}_D(m)
   \ne\emptyset\}$, $\tilde A=\cup_{v\in\tilde E}U^\ve_{D_0}(k/\ve,v,C\rho)$ 
   and $\cT^\ve=\big\{\gam\in C_{0T}:\,\exists t_1,t_2\,\mbox{with}\, 
   0\leq t_1<t_2\leq T,\,\gam_{t_1}\in U_{2\eta}(\cO_{i})\,\mbox{and dist}
   (\gam_{t_2},U_0)\leq\eta\big\}$. Assume that dist$(z,\bbR^d\setminus U_0)
   \leq\eta$ for any $z\in\partial V_i$. By Lemma \ref{lem1.6.2}(i), $R(x,z)$ 
   is continuous in $z$ when $z$ belong to a sufficiently small neighborhood
   of $\partial V_{i}$ which together with the definition of $S$-compacts
   yields that $S_{0T}(\gam)\geq R^{(i)}_\partial-3a/2$
   for any $\gam\in\cT^\ve$, provided $\eta$ is small enough. Observe that
   $Z^\ve_{\Phi_\ve^{m/\ve}v}\in\cT^\ve$ for any $v\in\tilde A$, provided
   $\ve$ is sufficiently small. These together with the arguments similar
   to the proof of Proposition \ref{prop1.7.1}(ii) yield the estimate
   \begin{equation}\label{1.8.7}
   m_D(\tilde A)\leq e^{-(R^{(i)}_\partial-2a)/\ve}
   \end{equation}
   for all $\ve$ small enough. Since $U^\ve_{D_0}(k/\ve,v,\frac 12C\rho)$
   are disjoint for different $v\in E$ we obtain by Lemma \ref{lem1.3.6},
   \begin{equation}\label{1.8.8}
   m_D(\tilde A)\geq\sum_{v\in\tilde E}U^\ve_{D_0}(k/\ve,v,\frac 12C\rho)\geq
   \tilde c\sum_{v\in\tilde E}U^\ve_{D_0}(k/\ve,v,C\rho)
   \geq\tilde cm_D(\tilde A).
   \end{equation}
    where $\tilde c=c_{C\rho}c_{C\rho/2}$. In a similar way we obtain that
    for each disc $\tilde D=\Phi_\ve^{k/\ve}U^\ve_{D_0}(k/\ve,v,C\rho)$,
    \begin{equation}\label{1.8.9}
    m_{\tilde D}\big(\tilde D\cap A^{(5)}_D(k)\big)\leq e^{-6a/\ve}
    \end{equation}
    provided $\ve$ is small enough. By (\ref{1.8.7})--(\ref{1.8.9}) together
    with Lemma \ref{lem1.3.6},
    \begin{equation}\label{1.8.10}
    m_D\big(A^{(4)}_D(m)\cap A^{(5)}_D(k)\big)\leq 
    e^{-(R^{(i)}_\partial+3a)/\ve}
    \end{equation}
    provided $m\leq k-2T$ and $\ve$ is small enough. Summing in $m,\,k$ and 
    $n$ we obtain from (\ref{1.8.4})--(\ref{1.8.6}) and (\ref{1.8.10}) that
    for a small $\be$ and all sufficiently small $\ve$,
     \begin{equation}\label{1.8.11}
     m_D\big\{ v\in D:\,\tau^\ve_v(G)\leq e^{(R^{(i)}_\partial+a)/\ve},\,
     Z^\ve_v(\tau^\ve_v(G))\in\Gam\big\}\leq e^{-a/\ve}.
     \end{equation}
     Taking discs $D$ to be small balls on the (extended) local unstable
      manifolds $W^u_x(w,q)$ and using the Fubini theorem as before we 
      obtain (\ref{1.8.11}) for $m$ and $\cW$ in place of $m_D$ and $D$,
       respectively. On the other hand, employing Proposition \ref{prop1.7.1}(i)
       we derive that
       \[
        m\big\{ v\in\cW:\,\tau^\ve_v(G)>e^{(R^{(i)}_\partial+a)/\ve}\big\}
        \leq e^{-\la /\ve}
     \]
     for some $\la>0$ and all $\ve$ small enough which together with 
     (\ref{1.8.11}) considered for $m$ and $\cW$ in place of $m_D$ and $D$ 
     yield (\ref{1.2.31}). 
     
     In order to complete the proof of Theorem \ref{thm1.2.7} it remains to
     derive (\ref{1.2.32}) and (\ref{1.2.33}). Both statements hold true for 
     $n=1$ in view of (\ref{1.2.29}) and (\ref{1.2.31}) but, in fact, we will
      use them as the induction base with $m_D$ and $D$ in place of $m$ 
      and $\cW$ where $D\in\cD^u_\ve(z,\al,\rho,C)$ and $D\subset(Q\cap V_i)
      \times\cW$ which holds true in view of (\ref{1.8.4})--(\ref{1.8.11})
      together with the corresponding form of Proposition \ref{prop1.7.1}.
      For such $D$ and $n\in\bbN$ set
      \[
      H_D(n,\al)=\big\{ v\in D:\, \Sigma^\ve_i(k,-\al)\leq\tau_v(i,k)\leq
      \Sigma^\ve_i(k,\al)\,\,\,\forall k\leq n\big\}
      \]
      and
      \[
      G_D(n)=\big\{ v\in D:\, Z^\ve_v(\tau_v(i,k))\in V_{\io_k(i)}
      \,\,\,\forall k\leq n\big\}.
      \]
      As the induction hypotesis we assume that for any $\al>0$ there exist
      $\la(\al)>0$ and $\la>0$ such that for all small $\ve$,
      \begin{equation}\label{1.8.12}
      m_D\big(H_D(n,\al)\big)\geq m_D(D)-ne^{-\la(\al)/\ve}\,\,\mbox{and}\,\,
      m_D\big(G_D(n)\big)\geq m_D(D)-ne^{-\la/\ve}.
      \end{equation}
      Set $a=\del/4K$ where $K$ is the same as in (\ref{1.2.15}) so that if
      \[
      \Gam_D(l)=\Gam_D(l,n,\al)=\big\{ v\in D:\, (l-1)a\leq\tau^\ve_v(i,n)
      <la\big\}\cap G_D(n)\cap H_D(n,\al)
      \]
      then 
      \begin{equation}\label{1.8.13}
      3\del/4\leq\mbox{dist}\big(Z^\ve_v(t),\cO_{\io_n(i)}\big)\leq 5\del/4
      \,\,\mbox{for all}\,\, v\in\Gam_D(l)\,\,\mbox{and}\,\, t\in[(l-1)a,la].
      \end{equation}
      Choose also $N=N_a$ so that for any $t\geq (N-1)a$,
      \begin{equation}\label{1.8.14}
      \Pi^tU_{2\del}(\cO_j)\subset U_{\del/4}(\cO_j)\,\,\mbox{for each}\,\,
      j=1,...,\ell.
      \end{equation}
      
      Let $E_l$ be a maximal $\big(la/\ve,C\rho,\ve,\overline {\Gam_D(l)},
      D_0\big)$-separated set where $D\subset D_0\in\cD_\ve^u(z,\al,
      \rho,C^3)$ as before. Set
      \[
      \Gam^U_D(l)=\cup_{v\in E_l}U^\ve_{D_0}(la/\ve,v,C\rho),
      \]
      then for $\ve$ small enough,
      \[
      D_0\supset\Gam^U_D(l)\supset\Gam_D(l).
      \]
      We claim that there exists $\be>0$ such that if $C\rho\leq\del/4$ then
      for all small $\ve$,
      \begin{equation}\label{1.8.15}
      m_D\big(\Gam^U_D(l)\cap\bigcup_{j=l+N}^\infty\Gam^U_D(j)\big)\leq
      e^{-\be/\ve}m_D(\Gam_D^U(l)).
      \end{equation}
      Indeed, let $\cT=\big\{\gam\in C_{0,Na}:\,\gam_0\in U_{2\del}
      (\cO_{\io_n(i)}),\,\gam_{t_1}\in U_{3\del/4}(\cO_{\io_n(i)})\,\,\mbox
      {for some}\,\,t_1\in[0,Na]\,\,\mbox{and}\,\,\gam_t\not\in U_{\del/2}
      (\cO_{\io_n(i)})\,\,\mbox{for all}\,\, t\in[0,Na]\big\}$.
      Then by (\ref{1.8.14}) and the lower semicontinuiti of the functional
      $S_{0,Na}$ we obtain that
      \begin{equation}\label{1.8.16}
      \inf\big\{ S_{0,Na}(\gam):\,\gam\in\cT\big\}=\eta>0.
      \end{equation}
      Since by Lemma \ref{lem1.3.2}(ii) and (iii) for any $w\in\tilde D(t,v)=
      \Phi_\ve^tU^\ve_{D_0}(t,v,C\rho)$ the distance $|\pi_1w-
      \pi_1\Phi_\ve^tv|$ has the order of $\ve$ we conclude from (\ref{1.8.13})
      and (\ref{1.8.16}) that for each $v\in\overline {\Gam_D(l)}$,
      \begin{equation*}
      \bfr_{0,Na}\big(Z^\ve_{\tilde v},\Psi^\eta_{0,Na}(\tilde z)\big)\geq
      \del/8\,\,\mbox{for any}\,\,\tilde v\in\tilde D(la/\ve,v)\cap
      \Phi_\ve^{la/\ve}\cup_{j=l+N}^\infty\Gam_D^U(j).
      \end{equation*}
      Hence, by (\ref{1.7.11}) and Lemma \ref{lem1.3.6} it follows that for any 
      $v\in\overline{\Gam_D(l)}$ and all $\ve$ small enough,
      \begin{equation*}
      m_{D_0}\big(U^\ve_{D_0}(la/\ve,v,C\rho)\cap\bigcup_{j=l+N}^\infty
      \Gam_D^U(j)\big)\leq e^{-\eta/2\ve}m_{D_0}\big(U^\ve_{D_0}
      (la/\ve,v,C\rho)\big)
      \end{equation*}
      and since $U^\ve_{D_0}(la/\ve,v,C\rho/2)$ are disjoint for different
      $v\in E_l$ we apply Lemma \ref{1.3.6} once more and obtain (\ref{1.8.15}).
      
      Set $Q_D(n)=H_D(n,\al)\cap G_D(n)\setminus H_D(n+1,\al)\cap G_D(n+1)$.
      Applying (\ref{1.8.12}) with $n=1$ to each $\tilde D=\tilde D(la/\ve,v),\,
      v\in E_l$ and using Lemma \ref{lem1.3.6} we derive also that
      \begin{equation}\label{1.8.17}
      m_D\big(\Gam_D^U(l)\cap Q_D(n)\big)\leq e^{-\be/\ve}m_{D_0}
      (\Gam^U_D(l))
      \end{equation}
      for some $\be>0$ and all small $\ve$. By (\ref{1.8.17}) we can write
      \begin{eqnarray}\label{1.8.18}
      &m_D(Q_D(n))=\sum_{\Sigma^\ve_i(n,-\al)\leq l\leq \Sigma^\ve_i(n,\al)}
      m_D\big(
      \Gam_D(l)\cap Q_D(n)\big)\leq\\
      &\sum_{\Sigma^\ve_i(n,-\al)\leq l\leq \Sigma^\ve_i(n,\al)}
      m_D\big(\Gam_D^U(l)
      \cap Q_D(n)\big)\leq e^{-\be/\ve}\sum_{1\leq l\leq
       \Sigma^\ve_i(n,\al)}m_{D_0}\big(\Gam_D^U(l)\big).\nonumber
       \end{eqnarray}
       Observe that for any finite measure $\mu$, measurable sets
       $A_1,A_2,...$ and integers $k,N>0$,
       \begin{equation}\label{1.8.19}
       \sum_{l=1}^{kN}\mu(A_i)=\sum_{j=1}^N\big(\mu(\cup_{i=0}^{k-1}A_{j+ik})
       +\sum_{i=0}^{k-2}\mu(A_{j+iN}\cap\bigcup_{r=i+1}^{k-1}A_{j+rN})\big)
       \end{equation}
       which follows applying $\mu(A\cup\tilde A)=\mu(A)+\mu(\tilde A)-
       \mu(A\cap\tilde A)$ to $A=A_{j+iN}$ and $\tilde A=\cup_{r=i+1}^{k-1}
       A_{j+rN}$ for $i=0,1,...,k-2$. Applying (\ref{1.8.19}) for $\mu=
       m_{D_0}$ and $A_l=\Gam^U_D(l)$ it follows from (\ref{1.8.15}) that
       \begin{equation*}
       \sum_{l=1}^{kN}m_{D_0}\big(\Gam_D^U(l)\big)\leq Nm_{D_0}(D_0)+
       e^{-\be/\ve}\sum_{l=1}^{kN}m_{D_0}\big(\Gam_D^U(l)\big),
       \end{equation*}
       i.e. for any $k\in\bbN$,
       \begin{equation*}
       \sum_{l=1}^{kN}m_{D_0}\big(\Gam_D^U(l)\big)\leq (1-e^{-\be/\ve})^{-1}
       Nm_{D_0}(D_0).
       \end{equation*}
       This together with (\ref{1.8.12}) and (\ref{1.8.18}) complete the
       induction step and proves (\ref{1.2.32}) and (\ref{1.2.33}) for $D$ and
       $m_D$ in place of $\cW$ and $m$. Finally, as before we complete the 
       proof of Theorem \ref{thm1.2.7} by choosing discs $D$ to be small balls 
       on the (extended) local unstable manifolds $W^u_x(w,\vrho),\, w\in
       \cW$ which together with the Fubini theorem enables us to extend the
       estimates to $\cW$ and $m$ as required in (\ref{1.2.32}) and 
       (\ref{1.2.33}).
      \qed

    \section{Averaging in difference equations}\label{sec1.9}\setcounter{equation}{0}

For readers convenience we start this section with the setup and necessary 
 technical results from \cite{Ki5} refering there for the corresponding proofs.
 These results are similar to Section \ref{sec1.3} and we refer the reader also
 to \cite{Ki7} where more details of proofs can be found than in \cite{Ki5} 
 and though \cite{Ki7} deals only with the continuous time case the
 corresponding discrete time proofs can be obtained, essentially, by 
 simplification. We will discuss below mainly the Axiom A case since the
 corresponding proofs for expanding transformations can be obtained,
 essentially, by simplification of the same arguments, roughly speaking,
 by ignoring the stable direction.

As in Section \ref{sec1.3} we will use the representations 
$\xi=\xi^\cX+\xi^\cW$
of vectors $\xi\in T(\bbR^d\times\bfM)=\bbR^d\oplus T\bfM$, the norms
$|\|\xi|\|$ and the distances $d_\bfM$ and $d(\cdot,\cdot)$ on $\bfM$ and
 on $\bbR^d\times\bfM$, respectively. It is known (see \cite{HPPS}) that the 
 hyperbolic splitting $T_{\La_x} \bfM=\Gam_x^s\oplus\Gam_x^u$ over $\La_x$ 
 can be continuously extended to the splitting $T_\cW\bfM=
 \Gam_x^s\oplus\Gam_x^u$ over $\cW$ which is forward invariant with respect to
  $DF_x$ and satisfies exponential estimates (\ref{1.3.1}) with a uniform in 
 $x\in\cX$ exponent $\ka>0.$  Moreover, by \cite{Ru} (see also \cite{Co}) we 
can choose these extensions so that $\Gam_x^s(w)$ and $\Gam_x^u(w)$ will be 
H\" older continuous in $w$ and $C^1$ in $x$ in the corresponding Grassmann 
bundle. Actually, since $\cW$ is contained in the basin of each attractor 
$\La_x$, any point $w\in\cW$ belongs to the stable manifold $W_x^s(v)$ of some
 point $v\in\La_x$ (see \cite{BR}), and so we 
choose naturally $\Gam_x^s(w)$ to be the tangent space to $W^s_x(v)$ at $w.$
Now each vector $\xi\in T_{x,w}(\cX\times \cW)=T_x\cX\oplus T_w\cW$ can be 
represented uniquely in the form $\xi=\xi^\cX +\xi^s+\xi^u$ with $\xi^\cX\in
 T_x\cX$, $\xi^s\in\Gam^s_x(w)$, and $\xi^u\in\Gam^u_x(w)$.
For each small $\ve,\al>0$ set $\cC^u(\ve,\al)=\{\xi\in T(\cX\times \cW):\,
\|\xi^{s}\|\leq\ve\al^{-2}\|\xi^u\|\,\mbox{and}\, \|\xi^\cX\|\leq\ve\al^{-1}
\|\xi^u\|\}$ and $\cC^u_{x,w}(\ve,\al)=\cC^u(\ve,\al)\cap T_{x,w}(\cX\times \cW)$
which are cones around $\Gam^u$ and $\Gam^u_x(w),$ respectively. Similarly,
we define $\cC^s(\ve,\al)=\{\xi\in T(\cX\times \cW):\,
\|\xi^{u}\|\leq\ve\al^{-2}\|\xi^s\|\,\mbox{and}\, \|\xi^\cX\|\leq\ve\al^{-1}
\|\xi^s\|\}$ and $\cC^u_{x,w}(\ve,\al)=\cC^u(\ve,\al)\cap T_{x,w}(\cX\times \cW)$
which are cones around $\Gam^s$ and $\Gam^s_x(w),$ respectively. The
corresponding version of Lemma \ref{lem1.3.1} is proved in \cite{Ki5} and
the discrete time versions of Lemmas \ref{lem1.3.2} and \ref{lem1.3.3} follow
in the same way as in \cite{Ki7}. Let, again, $\cD_\ve^u(z,\al,\rho,C)$ be the
set of all $C^1$ embedded $n_u-$dimensional closed discs $D\subset\cX\times \cW$
 such that $z\in D$, $TD\subset\cC^u(\ve,\al)$ and if $v\in\partial D$ then
  $C\rho\leq d_D(v,z)\leq C^2\rho$. For
  $D\in\cD^u_\ve(z,\al,\rho,C)$ and $z=(x,y)\in D\subset\cX\times \cW$ set
   $U^\ve_{D}(n,z,\vrho)=\{\tilde z\in D:\,\max_{0\leq k\leq n}
    d_D(\Phi_\ve^kz,\Phi_\ve^k\tilde z)\leq\vrho\}$ and
   let $\pi_1:\,\cX\times\cW\to\cX$ and $\pi_2:\,\cX\times\cW\to\cW$ be
   natural projections on the first and second factors, respectively. The
   same proof as in \cite{Ki7} yields the following discrete time version of
   Proposition \ref{prop1.3.4}.
   
   \begin{proposition}\label{prop1.9.1} 
   For any $\rho,C,b>0$ with $C$ large and $C\rho$ small enough
   there exists a positive function $\zeta_{b,\rho,T}(\Del,s,\ve)$ satisfying
   (\ref{1.3.7}) such that for any $x,x'\in\cX,\, y\in \cW,\, n\geq n_1,$ 
   $k\leq\frac T\ve -n$, $\be\in\bbR^d,\,|\be|\leq b$, $D\in\cD_\ve^u((x,y),
   \al,\rho,C),$ $z\in D$ and $V=U^\ve_D(t,z,C\rho)\subset D$ we have
     \begin{eqnarray}\label{1.9.1}
    &\,\,\,\,\bigg\vert\frac 1k\log\int_V\exp\langle\be,\sum_{j=n}^{n+k-1}
    B(x',Y^\ve_v(j))\rangle dm_D(v)+\frac 1k\log J^u_\ve(n,z)\\
    &-P_{F_{\pi_1z_n}}(\langle\be, B(x',\cdot)\rangle+\vf^u_{\pi_1z_n})\bigg
    \vert\leq\zeta_{b,\rho,T}(\ve k,\min (k,(\log\frac 1\ve)^\la),\ve)
    \nonumber
    \end{eqnarray}
    where $z_n=\Phi_\ve^nz$ and $\la\in(0,1)$.
    \end{proposition}
 
  Next, observe that the results of Section \ref{sec1.4} above are so general
  that they work both for the continuous and the discrete time case. Now, we
  will discuss the discrete time version of Lemma \ref{1.5.1}.
  
  \begin{lemma}\label{lem1.9.1} Let $x_i,\tilde x_i\in\cX,\, i=0,1,...,N,$
$0=t_0<t_1<...<t_{N-1}<t_N=T,$ $\Del=\max_{0\leq i\leq N-1}(t_{i+1}-t_i),$
$\xi_i=(x_i-x_{i-1})(t_i-t_{i-1})^{-1}$, $n(t)=\max\{ j\geq 0:\, t\geq t_j\},$
$\psi(t)=\tilde x_{n(t)}$ and
\begin{equation}\label{1.9.2}
\Xi_j^\ve(v,x)=([\ve^{-1}t_j]-[\ve^{-1}t_{j-1}])^{-1}
\sum_{[\ve^{-1}t_{j-1}]\leq k\leq [\ve^{-1}t_j]}B(x,Y^\ve_v(k)).
\end{equation}
Set for $t=\ve n\in [0,T],\, n\in\bbN$,
\begin{equation}\label{1.9.3}
Z_{v,x}^{\ve,\psi}(t)=x+\ve\sum_{0\leq k<n}B(\psi(\ve k),Y_v^\ve(k))
\end{equation}
and for $t\in[n\ve,(n+1)\ve),$
\begin{equation}\label{1.9.4}
Z_{v,x}^{\ve,\psi}(t)=(n+1-t/\ve)Z_{v,x}^{\ve,\psi}(\ve n)+
(t/\ve-n)Z_{v,x}^{\ve,\psi}(\ve(n+1)).
\end{equation}
Then
\begin{eqnarray}\label{1.9.5}
&\,\,\,\,\,\big\vert\Xi_j^\ve(v,x_{j-1})-(t_j-t_{j-1})^{-1}(Z^\ve_v(t_j)-
Z^\ve_v(t_{j-1}))\big\vert\leq K\big\vert Z^\ve_v(t_{j-1})-x_{j-1}\big\vert\\
&+\frac 12K^2(t_j-t_{j-1})+K\ve(4+KT+KT^2+T|\pi_1v|+|x_{j-1}|),\nonumber
\end{eqnarray}
\begin{eqnarray}\label{1.9.6}
&\sup_{0\leq s\leq t}\big\vert Z_{v,x}^{\ve,\psi}(s)-\psi(s)\big\vert\leq
|x-x_0|+\max_{0\leq j\leq n(t)}|x_j-\tilde x_j|\\
&+2\ve(K+\max_{1\leq i\leq N}|\xi_i|)+K\Del+n(t)\Del
\max_{1\leq j\leq n(t)}\big\vert\Xi_j^\ve(v,\tilde x_{j-1})
-\xi_j\big\vert\nonumber
\end{eqnarray}
and
\begin{equation}\label{1.9.7}
\sup_{0\leq s\leq t}\big\vert Z^\ve_v(s)-Z_{v,x}^{\ve,\psi}(s)\big\vert\leq
e^{Kt}\big(4K\ve+|\pi_1v-x|+Kt\sup_{0\leq s\leq t}\big\vert Z_{v,x}^{\ve,\psi}
(s)-\psi(s)\big\vert\big)
\end{equation}
where, recall, $Z^\ve_v(s)=X_v^\ve(s/\ve)$ and $\pi_1v=z\in\cX$ if
$v=(z,y)\in\cX\times\bfM$.
\end{lemma}
\begin{proof} The proof of (\ref{1.9.5}) and (\ref{1.9.6}) is strightforward
using the definitions (\ref{1.9.2})--(\ref{1.9.4}) in the same way as the proof
of (\ref{1.5.2}) and (\ref{1.5.3}) only the integrals in the latter case should be
replaced by the corresponding sums in the former one. The estimate (\ref{1.9.7})
follows in the same way as (\ref{1.5.4}) only the use of the standard Gronwall 
inequality in the latter proof should be replaced by the discrete time version
of the Gronwall inequality as in Lemma 4.20 of \cite{El}.
\end{proof} 

Now the proof of the discrete time version of Proposition \ref{prop1.5.2}
and of the remaining part of the proof of large deviations bounds (\ref{1.2.16})
and (\ref{1.2.17}) for the discrete time case proceeds almost verbatim as the
corresponding continuous time proofs in Section \ref{sec1.5}. Observe that in
the discrete time case the functionals $S_{0T}$ are given again by (\ref{1.2.13})
with $I_x(\nu)$ defined by (\ref{1.2.8}) where $F^1_x=F_x$ and $\vf^u(x)$ is
given by (\ref{1.2.34}). The property of $I$-functionals described in Lemma
\ref{lem1.6.1} follows directly in the discrete time case via conjugation since
we do not have to deal with the time change here. Other auxiliary results of
Section \ref{sec1.6} are derived in the discrete time case exactly in the same
 way as there. The proof of the
discrete time versions of Theorems \ref{thm1.2.5} and \ref{thm1.2.7} under the
corresponding assumptions goes through exactly in the
same way as its continuous time counterpart in Section \ref{sec1.6} yielding
the assertion of Theorem \ref{thm1.2.10}.  \qed

Next, we exhibit computations demonstrating a discrete time version of
Theorem \ref{thm1.2.7} for simple examples. The maps $F_x$ in both examples
 have the form $F_xy=3y+x\,(\mbox{mod}\, 1)$ where $x\in\bbR^1$ and 
 $y\in [0,1]$ but by identifying the end points of the unit interval we 
 view $F_x$ as expanding maps of the circle $\bbT^1$. The function
 $B$ from (\ref{1.1.10}) is given in the first example by
\[
 B(x,y)=x(x^2-4)(1-x^2)+50\sin 2\pi y.
\]
 Hence, we are dealing here with the
maps $\Phi_\ve:\,\bbR^1\times\bbT^1\to\bbR^1\times\bbT^1$ 
defined by
\[
\Phi_\ve(x,y)=\big( x+\ve(x(x^2-4)(1-x^2)+50\sin 2\pi y),\, 3y+x\,
(\mbox{mod}\, 1)\big).
\] 
All maps $F_x$ preserve the normalized Lebesgue measure $\mbox{Leb}$ on
$\bbT^1$ and it is the SRB measure $\SRBx$ for each $F_x$ in this simple
case. The averaged equation (\ref{1.1.11}) for $\brZ(t)=\brX^\ve(t/\ve)$
 has here the form
\[
\frac {d\brZ(t)}{dt}=\brB(\brZ(t)),\, 
\]
where $\brB(x)=x(x^2-4)(1-x^2)$. The one dimensional vector field $\brB(x)$ has 
three attracting fixed points $\cO_1=2,\cO_2=0,\cO_3=-2$ and
two repelling fixed points 1 and $-1$. In order to apply the discrete time
version of Theorem \ref{thm1.2.7} (i.e. Theorem \ref{thm1.2.10}) to this
 example
we have to verify that $B$ is complete at the fixed points $-2,-1,0,1,2$ of 
the averaged system. Since at these points $F_x$ coincides with the map
$y\to 3y$(mod 1) we can take the periodic orbits $1/8,3/8$ and $5/8,7/8$
of the latter and notice that the average of $\sin 2\pi y$ along the former
is $1/\sqrt 2$ and along the latter $-1/\sqrt 2$ which yields completness
of $B$ at zeros of $\bar B$.
 
 According to the corresponding part of Theorem \ref{thm1.2.10} which is a 
 discrete time version of Theorem \ref{thm1.2.7} the transitions between 
 $\cO_1,\cO_2,$ and $\cO_3$ 
 are determined by  $R_{ij},\, i,j=1,2,3$ which are obtained via the 
 functionals $S_{0t}(\gam)$ given by (\ref{1.2.13}). Even here these 
 functionals
 are not easy to compute though their main ingredients the functionals 
 $I_x(\nu)$ from (\ref{1.2.8}) are given now by the simple formula
 \[
   I_x(\nu)=\left\{\begin{array}{ll}
  \ln 3-h_\nu(F^1_x) &\mbox{if $\nu$ is $F_x$-invariant}\\
  \infty &\mbox{otherwise}
  \end{array}\right.
  \]
  and the set of $F_x$-invariant measures can be reasonably described since
  all $F_x$'s are conjugate to the simple map $y\to 3y\,(\mbox{mod}\, 1)$.
  We plot below the histogram of a single orbit of the slow motion
  $X^\ve_{x,y}(n),\, n=0,1,2,...,10^9$ with $\ve=10^{-3}$ 
  and the initial values $x=0,\, y=0.001$.
  The histogram shows that most of the points of the orbit stay near the
  attractors $\cO_1,\,\cO_2$ and $\cO_3$
  and $X^\ve_{x,y}(n)$ hops between basins of attraction of these points. 
  The form of the histogram indicates (according to Theorem \ref{thm1.2.7}) 
  the equality $R_{21}=R_{23}$ and in this case Theorem \ref{thm1.2.7} (or its
  discrete time version) cannot specify whether the slow motion exits from
  the basin of $\cO_2$ to the basin of $\cO_1$ or to the basin of $\cO_3$.
  Observe that Theorem \ref{thm1.2.7} is an asymptotical as 
   $\ve\to 0$ result and it takes an exponential in $1/\ve$ time for a typical
  orbit to exit from the basin of one attractor and to hop to the basin of 
  another one. Hence, the computations should be done for small $\ve$ and 
  exponentially long in $1/\ve$ orbits which is time consuming, so we put a
  big coefficient in front of $\sin$ which makes this exponent smaller. 
  Of course, it is hard to be absolutely sure that $\ve$ in our computations is 
  small enough and the number of iterates is large enough to demonstrate 
  faithfully the real situation in this case but we found that our histograms 
  are rather robust, for instance, their shapes have the same form for 
  $\ve=10^{-3}$ when the number of iterates ranges from $10^8$ to, at least, 
  $10^{11}$ and various initial conditions were checked, as well.

 \begin{figure} 
\centerline{\epsfig{file=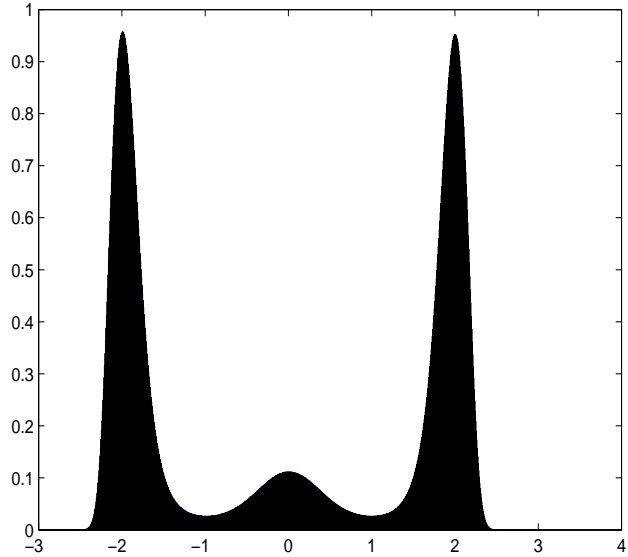, width=10cm, height=8.5cm}}
\caption{Symmetrical basins case}
\label{fig1.9.1:fig}
\end{figure}

Our second example differs from the first one only in $B$ which is given
now by
\[
 B(x,y)=x(x^2-4)(1-x)(1.5+x)+50\sin 2\pi y.
\]
Here the averaged system has the same attracting fixed points $\cO_1=2,
\cO_2=0,\cO_3=-2$ but one of two repelling fixed points moves from $-1$ to
$-3/2$. This makes the basin of attraction of $-2$ smaller while the left
interval of the basin of attraction of $0$ becomes larger. The latter leads
to the inequality $R_{23}>R_{21}$ which according to the discrete time version
of Theorem \ref{thm1.2.7} makes it more difficult for the slow motion to exit 
to 
the left from the basin of $\cO_2$ than to the right. As in the first example
in order to apply the latter result we have to check that $B$ is complete at
all zeros of $\bar B$ but since we did this already for all integer points it
remains to verify completness only for $x=-3/2$ which follows since
 $\sin 2\pi y$ equals 1 and $-1$ at two fixed points $1/4$ and $3/4$ of
 $F_{-3/2}$, respectively. In the histogram here
we plot $X^\ve_{x,y}(n),\, n=0,1,2,...,10^9$ with $\ve=10^{-3}$ 
and the initial values $x=-2,\, y=0.001$. In compliance with the discrete time
version of Theorem \ref{thm1.2.7} the histogram demonstrates that the slow
motion leaves the basin of $\cO_3$ and after arriving at the basin of $\cO_2$
it exits mostly to the basin of $\cO_1$, and so the slow motion hops mostly
between basins of $\cO_1$ and $\cO_2$ staying most of the time in small
neighborhoods of these points.

\begin{figure}
\centerline{\epsfig{file=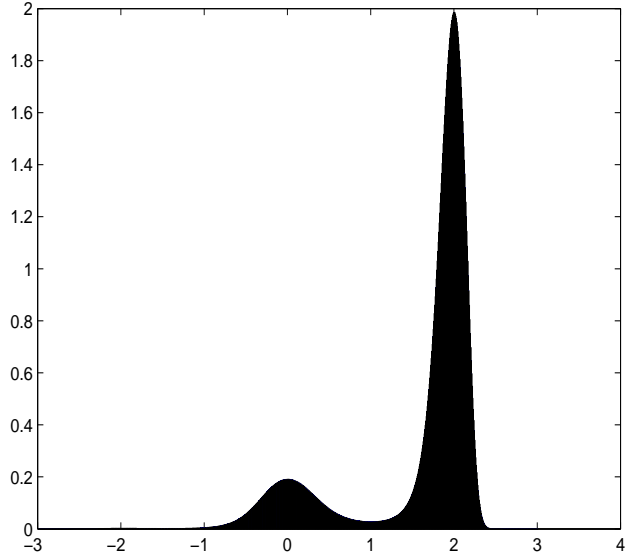, width=10cm, height=8.5cm}}
\caption{Asymmetrical basins case}
\label{fig1.9.2:fig}
\end{figure}

 \section{Extensions: stochastic resonance}\label{sec1.10}\setcounter{equation}{0}

The scheme for the stochastic resonance type phenomenon described below is a
slight modification of the model suggested by M.Freidlin (cf. \cite{Fre+}) and 
it can be demonstrated in the setup of three scale systems
 \begin{eqnarray}\label{1.10.1}   
&   \frac {dW^{\ve,\del}(t)}{dt}=\del\ve A(W^{\ve,\del}(t),X^{\ve,\del}(t),
Y^{\ve,\del}(t))\nonumber\\
&\frac {dX^{\ve,\del}(t)}{dt}=\ve B(W^{\ve,\del}(t),X^{\ve,\del}(t),
Y^{\ve,\del}(t))\\
&\frac {dY^{\ve,\del}(t)}{dt}=b(W^{\ve,\del}(t),X^{\ve,\del}(t),
Y^{\ve,\del}(t)),\nonumber
\end{eqnarray}
$W^{\ve,\del}=W^{\ve,\del}_{w,x,y}$, $X^{\ve,\del}=X^{\ve,\del}_{w,x,y}$, 
$Y^{\ve,\del}=Y^{\ve,\del}_{w,x,y}$ with initial conditions $W^{\ve,\del}(0)=w$,
$X^{\ve,\del}(0)=x$ and $Y^{\ve,\del}(0)=y$. 
We assume that $W^{\ve,\del}\in\bbR^l$, $X^{\ve,\del}\in\bbR^d$ while
 $Y^{\ve,\del}$ evolves on a compact $n_\bfM$-dimensional $C^2$ Riemannian 
 manifold 
$\bfM$ and the coefficients $A$, $B$, $b$ are bounded smooth vector fields on 
$\bbR^l$, $\bbR^d$ and $\bfM$, respectively, depending on other variables as 
parameters. The solution of (\ref{1.10.1}) determines the flow of 
diffeomorphisms 
$\Phi_{\ve,\del}^t$ on $\bbR^l\times\bbR^d\times\bfM$  acting by 
$\Phi_{\ve,\del}^t(w,x,y)=(W^{\ve,\del}_{w,x,y}(t), X^{\ve,\del}_{w,x,y}(t),
Y^{\ve,\del}_{w,x,y}(t))$. Taking $\ve=\del=0$ we arrive at the (unperturbed) 
flow $\Phi^t=\Phi^t_{0,0}$ acting by $\Phi^t(w,x,y)=(w,x,F^t_{w,x}y)$ 
where $F^t_{w,x}$ is another family of flows given by $F^t_{w,x}y=Y_{w,x,y}(t)$
with $Y=Y_{w,x,y}=Y^{0,0}_{w,x,y}$ which are solutions of
 \begin{equation}\label{1.10.2}
 \frac {dY(t)}{dt}=b(w,x,Y(t)),\,\, Y(0)=y.
 \end{equation}
 It is natural to view the flow $\Phi^t$ as describing an idealized
  physical system where parameters $w=(w_1,...,w_l)$, $x=(x_1,...,x_d)$ are
  assumed to be constants of motion while the perturbed flow 
  $\Phi^t_{\ve,\del}$ is regarded as describing a real system where evolution
  of these parameters is also taken into consideration but unlike the averaging
  setup (\ref{1.1.1}) we have now two sets of parameters moving with very 
  different speeds.

Set $\tilde W^{\ve,\del}(t)=W^{\ve,\del}(\frac t{\del\ve}),$
$\tilde X^{\ve,\del}(t)=X^{\ve,\del}(\frac t{\del\ve}),$
$\tilde Y^{\ve,\del}(t)=Y^{\ve,\del}(\frac t{\del\ve}),$ and pass from 
(\ref{1.10.1}) to the equations in the new time
\begin{eqnarray}\label{1.10.3}   
&   \frac {d\tilde W^{\ve,\del}(t)}{dt}= A(\tilde W^{\ve,\del}(t),
\tilde X^{\ve,\del}(t),\tilde Y^{\ve,\del}(t))\nonumber\\
&\frac {d\tilde X^{\ve,\del}(t)}{dt}=\del^{-1}B(\tilde W^{\ve,\del}(t),
\tilde X^{\ve,\del}(t),\tilde Y^{\ve,\del}(t))\\
&\frac {d\tilde Y^{\ve,\del}(t)}{dt}=(\del\ve)^{-1}b(\tilde W^{\ve,\del}(t),
\tilde X^{\ve,\del}(t),\tilde Y^{\ve,\del}(t)).\nonumber
\end{eqnarray}
 
Assume that the equation (\ref{1.10.1}) satisfy the assumptions similar to 
Assumptions \ref{ass1.2.1}, \ref{ass1.2.2}, \ref{ass1.2.6} together with other
corresponding conditions appearing in the setup of Theorem \ref{thm1.2.7}
(with $\bbR^l\times\bbR^d$ in place of $\bbR^d$), in particular, that 
$F^t_{w,x}y=Y^{0,0}_{w,x,y}(t),\, w\in\bbR^l,x\in\bbR^d$ form a compact set
of flows in the $C^2$ topology with $C^2$ dependence on $w,x$ and for all 
$w,x$ they are Axiom A flows in a neighborhood $\cW$ which contains a basic 
hyperbolic attractor $\La_{w,x}$ for $F^t_{w,x}$ and $\cW$ itself is contained 
in the basin of each $\La_{w,x}$. Set 
\begin{equation}\label{1.10.4}
\bar B_w(x)=\bar B(w,x)=\int B(w,x,y)d\SRBwx(y)
\end{equation}
where $\SRBwx$ is the SRB measure for $F^t_{w,x}$
and let $\bar X^{(w)}$ be the solution of the averaged equation
\begin{equation}\label{1.10.5} 
\frac {d\brX^{(w)}(t)}{dt}=\brB_w(\brX^{(w)}(t)).
\end{equation}
First, we apply averaging and large deviations estimates in averaging from
the previous section to two last equations in (\ref{1.10.3}) freezing the 
slowest
variable $w$ (i.e. taking for a moment $\del=0$). Namely, set $\hat X^\ve(t)=
X^{\ve,0}_{w,x,y}(t/\ve)$ and $\hat Y^\ve(t)=Y^{\ve,0}_{w,x,y}(t/\ve)$ so that
\begin{eqnarray}\label{1.10.6}
  &\frac {d\hat X^\ve(t)}{dt}= B(w,\hat X^\ve(t),\hat Y^\ve(t))\\
  &\frac {d\hat Y^\ve(t)}{dt}=\ve^{-1}b(w,\hat X^\ve(t),\hat Y^\ve(t)).\nonumber
 \end{eqnarray}
 Suppose for simplicity that $l=d=1$ (i.e. both $W^{\ve,\del}$ and
 $X^{\ve,\del}$ are one dimensional) and that the solution $\brX^{(w)}(t)$
 of (\ref{1.10.5}) has the limit set consisting of two attracting points
 $\cO_1$ and $\cO_2$, which for simplicity we assume to be independent
 of $w$, and a repelling fixed point $\cO_0^w$ depending on $w$ and separating
 their basins. As an example of $\brB$ we may have in mind
  $\brB_w(x)=(x-w)(1-x^2), -1<w<1$. Let $S^w_{0T}(\gam)$ be the 
 large deviations rate functional for the system (\ref{1.10.6}) defined in 
 (\ref{1.2.13}) and set for $i,j=1,2$,
 \begin{equation}\label{1.10.7}
 R_{ij}(w)=\inf\{ S^w_{0T}(\gam):\,\gam\in C_{0T},\,\gam_0=\cO_i,\,\gam_T
 =\cO_j,\, T\geq 0\}
 \end{equation}
 (cf. with $R_{ij}$ in Theorem \ref{thm1.2.7}).
 Set
 \begin{equation}\label{1.10.8}
 \bar A_i(w)=\int A(w,\cO_i,y)d\SRBwoi(y)
 \end{equation}
 and assume that for all $w$,
  \begin{equation}\label{1.10.9}
  \bar A_1(w)<0\quad\mbox{and}\quad\bar A_2(w)>0
  \end{equation}
  which means in view of the averaging principle (see Theorem \ref{thm1.2.3}
  and the following it discussion) that $W^{\ve,\del}_{w,x,y}(t)$ decreases 
  (increases) while $X^{\ve,\del}_{w,x,y}(t)$ stays close to $\cO_1$ 
  (to $\cO_2$) for "most" $y$'s with respect to the Riemannian volume on 
  $\bfM$ restricted to $\cW$.
  
  The following statement suggests a "nearly" periodic behavior of the
  slowest motion.
   \begin{conjecture}\label{conj1.10.1} Suppose that there exist strictly 
  increasing
  and decreasing functions $w_{-}(r)$ and $w_{+}(r)$, respectively, so that
  \[
  R_{12}(w_{-}(r))=R_{21}(w_{+}(r))=r
  \]
  and $w_{-}(\la)=w_{+}(\la)=w^*$ for some $\la>0$ while $w_{-}(r)<w^*<w_{+}(r)$
  for $r<\la$. Assume that $\del\to 0$ and $\ve\to 0$ in such a way that 
  \begin{equation}\label{1.10.10}
  \lim_{\ve,\del\to 0}\ve\ln(\del\ve)=-\rho>-\la.
  \end{equation}
  Then for any $w,x$ there exists $t_0>0$ so that the slowest motion
  $\tilde W^{\ve,\del}_{w,x,y}(t+t_0),\, t\geq 0$ converges weakly 
  (as $\ve,\del\to 0$ so that (\ref{1.10.10}) holds true) as a random 
  process on the probability space $(\cW,m_\cW)$ (where 
  $m_\cW$ is the normalized Riemannian volume on $\cW$) to a periodic
  function $\psi(t)$, $\psi(t+T)=\psi(t)$ with 
  \[
  T=T(\rho)=\int_{w_{-}(\rho)}^{w_{+}(\rho)}\frac {dw}{|\bar A_1(w)|}+
  \int_{w_{-}(\rho)}^{w_{+}(\rho)}\frac {dw}{|\bar A_2(w)|}.
  \]
  \end{conjecture}
  
   The argument supporting this conjecture goes as follows.
    Set $\check W^{\ve,\del}(t)=W^{\ve,\del}(t/\ve)$,
   $\check X^{\ve,\del}(t)=X^{\ve,\del}(t/\ve)$ and $\check Y^{\ve,\del}(t)
   =Y^{\ve,\del}(t/\ve)$ which satisfy
   \begin{eqnarray}\label{1.10.11}   
&   \frac {d\check W^{\ve,\del}(t)}{dt}=\del A(\check W^{\ve,\del}(t),
\check X^{\ve,\del}(t),\check Y^{\ve,\del}(t))\nonumber\\
&\frac {d\check X^{\ve,\del}(t)}{dt}=B(\check W^{\ve,\del}(t),\check 
X^{\ve,\del}(t),\check Y^{\ve,\del}(t))\\
&\frac {d\check Y^{\ve,\del}(t)}{dt}=\ve^{-1}b(\check W^{\ve,\del}(t),
\check X^{\ve,\del}(t),\check Y^{\ve,\del}(t)).\nonumber
\end{eqnarray}
Since $\check W^{\ve,\del}$ moves much slower than $\check X^{\ve,\del}$
we can freeze the former and in place of (\ref{1.10.11}) we can study 
(\ref{1.10.6}).
Applying the arguments of Theorem \ref{thm1.2.7} to the pair $\hat X,\hat Y$ 
from (\ref{1.10.6}) we conclude by (\ref{1.2.30}) that the intermediate
   motion $\tilde X^{\ve,\del}$ most of the time stays very close to either
   $\cO_1$ or $\cO_2$ before it exits from the corresponding basin, and so
   in view of an appropriate averaging principle (which follows, for
   instance, from Theorem \ref{1.2.3}) on bounded time intervals the slowest 
   motion $\tilde W^{\ve,\del}$ mostly stays close to the corresponding 
   averaged motion determined by the vector fields $\bar A_1$ and $\bar A_2$
   given by (\ref{1.10.4}). When $\tilde X^{\ve,\del}$ is close to $\cO_1$ 
    the slowest motion $\tilde W^{\ve,\del}$ decreases until $w=w_{-}(\rho)$ 
    where $R_{12}(w)=\rho$. In view of (\ref{1.2.29}) and the scaling 
    (\ref{1.10.10})
    between $\ve$ and $\del$, a moment later $R_{12}(w)$ becomes less
    than $\rho$ and $\tilde X^{\ve,\rho}$ jumps immediately close to 
    $\cO_2$. There $\bar A_2(w)>0$, and so $\tilde W^{\ve,\del}$ starts to 
    grow until it reaches $w=w_{+}(\rho)$ where $R_{21}(w)=\rho$. A moment
    later $R_{21}(w)$ becomes  smaller than $\rho$ and in view of 
    (\ref{1.2.29}) the intermediate motion
    $\tilde X^{\ve,\del}$ jumps immediately close to
    $\cO_1$. This leads to a nearly periodic behavior of $\tilde W^{\ve,\del}$.
    In order to make these arguments precise we have to deal here with
    an additional difficulty in comparison with the two scale setup considered
    in previous sections since now the large deviations $S$-functionals from 
    Theorem \ref{thm1.2.3} and the $R$-functions describing adiabatic 
    fluctuations
    and transitions of Theorems \ref{thm1.2.5} and \ref{thm1.2.7} depend on 
    another
    very slowly changing parameter. Still, the technique of 
    Sections \ref{sec1.7} and \ref{sec1.8} above applied on time intervals 
    where changes in the $w$-variable can be neglected should work here
     but the details of this approach have not been worked out yet. 
     
      On the other hand, when the fast motion $Y^{\ve,\del}$ does not depend
     on the slow motions, i.e. when the coefficient $b$ in (\ref{1.10.1}) 
     depend only on the coordinate $y$ (but not on $w$ and $x$), then the above 
     arguments can be made precise without much effort. Indeed, we can obtain
     estimates for transition times $\tau^\ve(1)$ and $\tau^\ve(2)$ of 
     $X^{\ve,\del}(t/\ve)$ between neighborhoods of $\cO_1$ and $\cO_2$ as
     in Theorem \ref{thm1.2.7} applying the latter to $\hat X^\ve$ and
     $\hat Y^\ve$ from (\ref{1.10.6}) with freezed $w$-variable. This is
     possible since the method of Proposition \ref{prop1.7.1} requires us
     to make large deviations estimates, essentially, only for probabilities
      $m\{ v\in D: kT\leq\tau_v(i)<(k+1)T\}$, i.e. on bounded time intervals,
      and then combine them with the Markov property type arguments. During
      such times the slowest motion $W^{\ve,\del}$ can move only a distance of
      order $\del T$. Thus freezing $w$ and using the Gronwall inequality
      for the equation of $X^\ve$ in order to estimate the resulting error
      we see that the latter is small enough for our purposes. Observe, that
      it would be much more difficult to justify freezing $w$ in the
      coefficient $b$ of $Y^\ve$, if we allow the latter to depend on $w$,
      since a strightforward application of the Gronwall inequality there
      would yield an error estimate of an exponential in $1/\ve$ order which 
      is comparable with $1/\del$. Still, it may be possible to take care 
      about the general case using methods of Sections \ref{sec1.4} and 
      \ref{sec1.5} since we produce large deviations estimates there by gluing 
      large deviations estimates on smaller time intervals where the 
      $x$-variable (and so, of course, $w$-variable) can be freezed. Next, set
      \[
      W_{w,y}^{\ve,\del,i}(t)=w+\del\ve\int_0^tA(W_w^{\del,\ve,i}(s),\cO_i,
      Y(s))ds
      \]
      where now $Y$ does not depend on $\ve$ and $\del$. Then by (\ref{1.10.1})
      together with the Gronwall inequality we obtain that
      \[
      |W_{w,x,y}^{\ve,\del}(t)-W_{w,y}^{\ve,\del,i}(t)|\leq L\del\ve
      e^{\del\ve Lt}\int_0^t|X^{\ve,\del}_{w,x,y}(s)-\cO_i|ds
      \]
      where $L$ is the Lipschitz constant of $A$. If $x$ belongs to the 
      basin $\cO_i$ then according to Theorem \ref{thm1.2.7} $\hat X^\ve$,
      and so also $X^{\ve,\del}$, stays most of the time near $\cO_i$ up
      to its exit from the basin of the latter which yields according to
      the above inequality that $W^{\ve,\del}$ stays close to $W^{\ve,\del,i}$
      during this time. But now we can employ the averaging principle for
      the pair $W^{\ve,\del,i}(t),Y(t)$ which sais that $W^{\ve,\del,i}(t)$
      stays close on the time intervals of order $1/\del\ve$ to the averaged
      motion $\bar W^{\ve,\del,i}_w(t)$ defined by
      \[
      \bar W^{\ve,\del,i}_w(t)=w+\int_0^t\bar A_i(\bar W^{\ve,\del,i}_w(s))ds
      \]
      and in view of (\ref{1.10.9}), $\bar W^{\ve,\del,1}_w(t)$ decreases while
      $\bar W^{\ve,\del,2}_w(t)$ increases which leads to the behavior
      described in Conjecture \ref{conj1.10.1}.
 
    A similar conjecture can be made under the corresponding conditions for
    the discrete time case determined by a three scale difference system of
    equations of the form
    \begin{eqnarray}\label{1.10.12}
  &W^{\ve,\del}(n+1)-W^{\ve,\del}(n)=\ve\del A(W^{\ve,\del}(n),X^{\ve,\del}(n),
  Y^{\ve,\del}(n)),\,\,\,W^{\ve,\del}(0)=w,\nonumber\\
 &\,\,\,\,\,\,\,\,\,\,\,\,\, X^{\ve,\del}(n+1)-X^{\ve,\del}(n)=
 \ve B(W^{\ve,\del}(n),
 X^{\ve,\del}(n),Y^{\ve,\del}(n)),\,\,\, X^{\ve,\del}(0)=x,\\
 &Y^{\ve,\del}(n+1)=F_{W^{\ve,\del}(n),X^{\ve,\del}(n)}
  Y^{\ve,\del}(n),\,\,\, Y^{\ve,\del}(0)=y\nonumber
 \end{eqnarray}
where $A$ and $B$ are smooth vector functions and $F_{w,x}:\, \bfM\to\bfM$
 is a smooth map (a diffeomorphism or an endomorphism). We obtain an example 
 where discrete time versions of conditions of Conjecture \ref{conj1.10.1} hold
 true setting, for instance, $A(w,x,y)=x\cos 2\pi w+\sin 2\pi y$, $B(w,x,y)
 =(x-w)(1-x^2)+\sin 2\pi y$ and $F_{w,x}y=3y+x+w$ (mod 1).

 \section{Young measures approach to averaging}\label{sec1.11}
\setcounter{equation}{0}

This section deals with the averaging principle and a bit with the 
corresponding large deviations in the sense of convergence of Young measures 
\index{Young measures} and I thank K.Gelfert for asking me about Young 
measures applications in averaging and for indicating to me the paper \cite{AG}.

Let $\mu$ belongs to the space $\cP(\bbR^d\times\bfM)$ of probability measures
on $\bbR^d\times\bfM$ and consider the Young measure (which is a map from a
measure space to a space of measures, see \cite{AG}) $\zeta^\ve$ from 
$([0,T]\times\bbR^d\times\bfM,\ell_T\times\mu)$ to $\cP(\bbR^d\times\bfM)$ 
defined by
\[
\zeta^\ve(t,x,y)=\del_{X^\ve_{x,y}(t/\ve),Y^\ve_{x,y}(t/\ve)}
\]
where $\ell_T$ is the Lebesgue measure on $[0,T]$, $\del_w$ is the unit mass
at $w$, and $X^\ve,\, Y^\ve$ are solutions of $(\ref{1.1.1})$ on the product
$\cP(\bbR^d\times\bfM)$. We assume that for all $x,z\in\bbR^d$ and $y,v\in
\bfM$ the coefficients $b$ and $B$ satisfy 
\begin{eqnarray}\label{1.11.1}
&|B(x,y)|+|b(x,y)|\leq K\,\,\mbox{and}\\
& |B(x,y)-B(z,v)|+|b(x,y)-b(z,v)|\leq K\big(|x-z|+d_\bfM (y,v)\big)
\nonumber\end{eqnarray}
for some $L>0$ independent of $x,y,z,v$. Of course, we could require the
Lipschitz continuity and the boundedness conditions (\ref{1.11.1}) only in
some open domain as in Section \ref{sec1.2} but we can always extend these 
vector fields to the whole $\bbR^d$ keeping these properties intact.

Suppose that $\mu\in\cP(\bbR^d\times\bfM)$ has a disintegration
\begin{equation}\label{1.11.2}
d\mu(x,y)=d\mu_x(y)d\la(x),\,\,\,\la\in\cP(\bbR^d)
\end{equation}
such that for each Lipschitz continuous function $g$ on $\bfM$ and any
$x,z\in\bbR^d$,
\begin{equation}\label{1.11.3}
|\int gd\mu_x-\int gd\mu_z|\leq K_{L(g)}|x-z|
\end{equation}
for some $K_L>0$ depending only on $L$ where $L(g)$ is both a Lipschitz
constant of $g$ and it also bounds $|g|$. Set
\begin{equation}\label{1.11.4}
\bar B(x)=\int B(x,y)d\mu_x(y)
\end{equation}
then by (\ref{1.11.1}) and (\ref{1.11.3}), $\bar B$ is bounded and Lipschitz 
continuous, and so there exists a unique solution $\bar X^\ve(t)=
\bar X^\ve_x(t)$ of (\ref{1.1.3}).  For any bounded continuous function $g$
on $\bbR^d\times\bfM$ define
\[
\cE^g_\ve(t,\del)=\big\{ (x,y)\in\bbR^d\times\bfM:\,\big\vert\frac 1t\int_0^t
g(x,Y^\ve_{x,y}(u))du-\bar g(x)\big\vert >\del\big\}
\]
where $\bar g(x)=\int g(x,y)d\mu_x(y)$.

By the definition (see \cite{AG}), the Young measures $\zeta^\ve$ converge
as $\ve\to 0$ to the Young measure $\zeta^0$ defined by
\[
\zeta^0(t,x,y)=\del_{\bar Z_x(t)}\times\mu_{\bar Z_x(t)}\in
\cP(\bbR^d\times\bfM),
\] 
$\bar Z_x(t)=\bar X^\ve_x(t/\ve)$,
if for any bounded continuous function $f$ on $[0,T]\times\bbR^d\times\bfM$,
\[
\int_0^Tf(s,\Phi_\ve^{s/\ve}(x,y))ds\to\int_0^T\bar f(s,\bar Z_x(s))ds\,\,
\mbox{as}\,\,\ve\to 0.
\]
The following result provides a verifiable (in some interesting cases) criterion
for even stronger convergence.

\begin{theorem}\label{thm1.11.1} Let $\mu\in\cP(\bbR^d\times\bfM)$ has the
 disintegration (\ref{1.11.2}) satisfying (\ref{1.11.3}). Then
 \begin{equation}\label{1.11.5}
 \lim_{\ve\to 0}\int_{\bbR^d}\int_\bfM\sup_{0\leq t\leq T}\big\vert\int_0^t
 \big(f(s,\Phi^{s/\ve}(x,y))-\bar f(s,\bar Z_x(s))\big)ds\big\vert d\mu_x(y)
 d\la(x)=0
 \end{equation}
for any bounded continuous function $f=f(t,x,y)$ on $[0,T]\times\bbR^d\times
 \bfM$ where $\bar f(t,x)=\int f(t,x,y)d\mu_x(y)$ if and only if for each
  $N\in\bbN$ and any finite collection $g_1,...,g_N$ of bounded Lipschitz
  continuous functions on $\bbR^d\times\bfM$ there exists an integer valued 
function $n=n(\ve)\to\infty$ as $\ve\to 0$ such that for any $\del>0$
and $l=1,...,N$,
    \begin{equation}\label{1.11.6}
    \lim_{\ve\to 0}\max_{0\leq j< n(\ve)}\mu\{ 
    \Phi_\ve^{-jt(\ve)}\cE^{g_l}_\ve(t(\ve),\del)\}=0,
    \end{equation}
    where $t(\ve)=\frac T{\ve n(\ve)}$ and, recall, $\Phi_\ve^t(x,y)=
    (X^\ve_{x,y}(t),Y^\ve_{x,y}(t))$. 
 \end{theorem}
 \begin{proof} First, we prove that (\ref{1.11.5}) implies (\ref{1.11.6}). Let
 $g_1,...,g_N$ be bounded Lipschitz continuous functions on $\bbR^d\times
 \bfM$ and set
 \begin{equation}\label{1.11.7}
 \rho_{x,y}^{\ve,l}(t)=\ve\int_0^t\big(g_l(X^\ve_{x,y}(s),Y^\ve_{x,y}(s))
 -\bar g_l(\bar X^\ve_x(s))\big)ds.
 \end{equation}
 If
 \[
 \rho_{x,y}^{\ve,l}=\sup_{0\leq t\leq T/\ve}|\rho_{x,y}^{\ve,l}(t)|
 \]
 then by (\ref{1.11.5}) for each $l=1,...,N$,
 \begin{equation}\label{1.11.8}
 \rho^\ve_l=\int_{\bbR^d}\int_\bfM\rho_{x,y}^{\ve,l}d\mu(x,y)\to 0\,\,\,
 \mbox{as}\,\,\ve\to 0.
 \end{equation}
 Choose an integer valued function $n(\ve)\to\infty$ as $\ve\to 0$ so that
  \begin{equation}\label{1.11.9}
  n(\ve)\max_{1\leq l\leq N}\rho^\ve_l\to 0\,\,\,\mbox{as}\,\,\ve\to 0
 \end{equation}
 and let $t(\ve)=T/\ve n(\ve)$. Set $x_k^\ve=X^\ve_{x,y}(kt(\ve))$, $y_k^\ve=
 Y^\ve_{x,y}(kt(\ve))$ and $\bar x_k^\ve=\bar X^\ve_{x}(kt(\ve))$, $k=0,1,...$.
 Then by (\ref{1.11.7}),
 \begin{equation}\label{1.11.10}
 \rho_{x,y}^{\ve,l}((j+1)t(\ve))-\rho_{x,y}^{\ve,l}(jt(\ve))=
 \ve\int_0^{t(\ve)}\big(g_l(X^\ve_{x^\ve_j,y^\ve_j}(u),Y^\ve_{x^\ve_j,y^\ve_j}
 (u))-\bar g_l(\bar X^\ve_{\bar x_j^\ve}(u))\big)du.
 \end{equation}
 By (\ref{1.11.1}),
 \begin{eqnarray}\label{1.11.11}
 &\ve\big\vert\int_0^{t(\ve)}\big(g_l(X^\ve_{x^\ve_j,y^\ve_j}(u),
 Y^\ve_{x^\ve_j,y^\ve_j}(u))-g_l(x^\ve_j,Y^\ve_{x^\ve_j,y^\ve_j}(u))\big)du
 \big\vert\\
 &\leq\ve L_l\int_0^{t(\ve)}|X^\ve_{x^\ve_j,y^\ve_j}(u)-x_j^\ve|du\leq
 L_lL(\ve n(\ve))^2\nonumber
 \end{eqnarray}
 where $L_l$ is the Lipschitz constant of $g_l$. Similarly, by (\ref{1.11.1})
 and (\ref{1.11.3}),
 \begin{equation}\label{1.11.12}
 \ve\big\vert\int_0^{t(\ve)}\big(\bar g_l(\bar X^\ve_{\bar x_j^\ve}(u))-
 \bar g_l(\bar x_j^\ve)\big)du\big\vert\leq (L_l+K_{L_l})K(\ve t(\ve))^2
 \end{equation}
 and
 \begin{equation}\label{1.11.13}
 |\bar g_l(\bar x^\ve_j)-\bar g_l(x^\ve_j)|\leq (L_l+K_{L_l})|\bar x_j^\ve-
 x^\ve_j|\leq (L_l+K_{L_l})\rho^\ve_{x,y}.
 \end{equation}
 It follows from (\ref{1.11.10})--(\ref{1.11.13}) that
 \begin{eqnarray}\label{1.11.14}
 &\big\vert\frac 1{t(\ve)}\int_0^{t(\ve)}g_l(x^\ve_j,Y^\ve_{x^\ve_j,y^\ve_j}(u))
 du-\bar g_l(x^\ve_j)\big\vert\\
 &\leq TK(2L_l+K_{L_l})/n(\ve)+(L_l+K_{L_l}+2T^{-1}
 n(\ve))\rho^\ve_{x,y}.\nonumber
 \end{eqnarray}
 Given $\del >0$ choose $\ve_\del>0$ such that for all $\ve\leq\ve_\del$ and
 $l=1,...,N$,
 \[
 TK(2L_l+K_{L_l})/n(\ve)\leq\del/2.
 \]
 Then by (\ref{1.11.14}),
 \begin{equation*}
 \Phi_\ve^{-jt(\ve)}\cE_\ve^{g_l}(t(\ve),\del)\subset A_\ve(\del)
 =\big\{ (x,y)\in\bbR^d\times\bfM:\, (L_l+K_{L_l}+2T^{-1}n(\ve))
 \rho^{\ve,l}_{x,y}>\del/2\big\}.
 \end{equation*}
 By Chebyshev's inequality
 \begin{equation}\label{1.11.15}
 \mu(A_\ve(\del))\leq\frac 2{\del}(L_l+K_{L_l}+2T^{-1}n(\ve))\rho_l^\ve.
 \end{equation}
 By (\ref{1.11.9}) the right hand side of (\ref{1.11.15}) tends to 0 as $\ve\to 0$
 yielding (\ref{1.11.6}).
 
 Next, we derive (\ref{1.11.5}) from (\ref{1.11.6}).
 Since $f$ in (\ref{1.11.5}) is a bounded function and $\la$ is
 a probability measure it is easy to see that it suffices to prove 
 (\ref{1.11.5}) when the integration in $x$ there is restricted to
 compact subsets of $\bbR^d$. But if we integrate in (\ref{1.11.5}) in $x$
 running over a compact set $G\subset\bbR^d$ then by (\ref{1.1.1}) and 
 (\ref{1.11.1}),
 \begin{equation}\label{1.11.16}
 \sup_{0\leq s\leq T}\mbox{dist}\big(X^\ve_{x,y}(s/\ve),G\big)\leq KT,
 \end{equation}
 i.e. both $Z^\ve_{x,y}(s)=X^\ve_{x,y}(s/\ve)$ and $\bar Z_x(s)$ belong to
 the $KT-$neighborhood $G_{KT}$ of the set $G$ when $x\in G$ and $s\in[0,T]$.
 On $[0,T]\times G_{KT}\times\bfM$ we can approximate $f$ uniformly by
 Lipschitz continuous functions. Thus, in place of (\ref{1.11.5}) it suffices
 to show that for any compact set $G\subset\bbR^d$ and a bounded Lipschitz
 continuous function $f$ on $[0,T]\times G_{KT}\times\bfM$ with a Lipschitz
 constant $L=L(f)$ in all variables,
  \begin{equation}\label{1.11.17}
 \lim_{\ve\to 0}\ve\int_G\int_\bfM\sup_{0\leq t\leq T/\ve}\big\vert\int_0^t
 \big(f(\ve s,X^\ve_{x,y}(s),Y^\ve_{x,y}(s))-\bar f(\ve s,\bar X^\ve_x(s))\big)
 ds\big\vert d\mu_x(y)d\la(x)=0.
 \end{equation}
 
  By (\ref{1.11.2}), (\ref{1.11.3}) and (\ref{1.11.4}),
 \begin{eqnarray}\label{1.11.18}
 &\ve\big\vert\int_0^t\big(f(\ve s,X^\ve_{x,y}(s),Y^\ve_{x,y}(s))-
 \bar f(\ve s,\bar X^\ve_x(s))\big)ds\big\vert\\
 &\leq\ve\big\vert\int_0^t\big(f(\ve s,X^\ve_{x,y}(s),Y^\ve_{x,y}(s))-
 \bar f(\ve s,X^\ve_x(s))\big)ds\big\vert\nonumber\\
 &+(L+K_L)T\sup_{0\leq s\leq T/\ve}|X^\ve_x(s)-\bar X^\ve_x(s)|.\nonumber
 \end{eqnarray}
 Since (\ref{1.11.6}) holds true also for $g=B$, it follows from Theorem 2.1
 of \cite{Ki7} that
 \begin{equation}\label{1.11.19}
 \lim_{\ve\to 0}\int_G\int_\bfM\sup_{0\leq s\leq T/\ve}|X^\ve_x(s)-
 \bar X^\ve_x(s)|d\mu(x,y)=0,
 \end{equation}
 and so we have only to deal with the first absolute value in the right hand
 side of (\ref{1.11.18}). As before set $x_k^\ve=X^\ve_{x,y}(kt(\ve))$, $y_k^\ve=
 Y^\ve_{x,y}(kt(\ve))$, $\bar x_k^\ve=\bar X^\ve_{x}(kt(\ve))$ and fix a large
 $N\in\bbN$. Let $l=l(j)=[\ve jt(\ve)N/T]=[jN/n(\ve)]$ then by (\ref{1.11.1}),
 (\ref{1.11.2}) and (\ref{1.11.3}),
 \begin{eqnarray}\label{1.11.20}
 &\quad\,\,\,\,\,\,\,\,\,\,\,\,\,\,\, \ve\big\vert\int_0^{t(\ve)}
 \big(f(\ve jt(\ve)+\ve u, X^\ve_{x^\ve_j,y^\ve_j}
 (u),Y^\ve_{x^\ve_j,y^\ve_j}(u))\\
 &-f(lT/N,x^\ve_j,Y^\ve_{x^\ve_j,y^\ve_j}(u))
 \big)ds\big\vert\leq LT^2/Nn(\ve)\nonumber\\
 &+L\ve\int_0^{t(\ve)}|X^\ve_{x^\ve_j,y^\ve_j}(u)-x^\ve_j|du
 \leq LT^2/Nn(\ve)+LT^2(1+K)(n(\ve))^{-2}\nonumber
 \end{eqnarray}
 and
 \begin{eqnarray}\label{1.11.21}
 &\ve\big\vert\int_0^{t(\ve)}\big(\bar f(\ve jt(\ve)+\ve u,
 X^\ve_{x^\ve_j,y^\ve_j}(u))-\bar f(lT/N,x_j^\ve))du\big\vert\\
 &\leq LT^2/Nn(\ve)+T^2(L+KL+KK_L)(n(\ve))^{-2}.\nonumber
 \end{eqnarray}
 Now using (\ref{1.11.20}), (\ref{1.11.21}) and assuming that $|f|\leq\hat L_f$
 for some constant $\hat L_f>0$ we obtain
 \begin{eqnarray}\label{1.11.22}
 &\ve\sup_{0\leq t\leq T/\ve}\big\vert\int_0^t\big(f(\ve s,X^\ve_{x,y}(s),
 Y^\ve_{x,y}(s))-\bar f(\ve s, X^\ve_{x,y}(s))\big)ds\big\vert\\
 &\leq 2\hat L_f\ve t(\ve)+\ve\sum_{j=0}^{n(\ve)-1}
 \big\vert\int_{jt(\ve)}^{(j+1)t(\ve)}\big(f(\ve s,X^\ve_{x,y}(s),
 Y^\ve_{x,y}(s))-\bar f(\ve s,\bar X^\ve_x(s))\big)ds\big\vert\nonumber\\
 &\leq 2\hat L_f\ve t(\ve)+\ve\sum_{j=0}^{n(\ve)-1}
 \big\vert\int_0^{t(\ve)}\big(f(\ve jt(\ve)+\ve s,X^\ve_{x^\ve_j,y^\ve_j}(s),
 Y^\ve_{x^\ve_j,y^\ve_j}(s))\nonumber\\
 &-\bar f(\ve jt(\ve)+\ve s, X^\ve_{x^\ve_j,y^\ve_j}(s))\big)ds\big\vert
 \nonumber\\
 &\leq 2LT^2/N+2\big(\hat L_fT+T^2(L+KL+KK_L)\big)/n(\ve)\nonumber\\
 &+\ve t(\ve)\sum_{l=0}^{N-1}\sum_{ln(\ve)/N\leq j<(l+1)n(\ve)/N,j\leq 
 n(\ve)}\big\vert\frac 1{t(\ve)}\int_0^{t(\ve)}f(lT/N,x^\ve_j,
 Y^\ve_{x^\ve_j,y^\ve_j}(s))ds\nonumber\\
 &-\bar f(lT/N,x_j^\ve)\big\vert\leq 2LT^2/N+2\big(\hat L_fT+T^2(L+KL+KK_L)\big)
 /n(\ve)\nonumber\\
 &+\ve t(\ve)n(\ve)\del+2\hat L_f\ve t(\ve)\sum_{l=0}^{N-1}
 \sum_{ln(\ve)/N\leq j<(l+1)n(\ve)/N,j
 \leq n(\ve)}\bbI_{\cE^{f_l}_\ve(t(\ve),\del)}(x_j^\ve,y_j^\ve)\nonumber
 \end{eqnarray}
 where $f_l(z,v)=f(lT/N,z,v).$ Integrating against $\mu$ both parts of 
 (\ref{1.11.22}) over $G\times\bfM$ we obtain
 \begin{eqnarray}\label{1.11.23}
 &\,\,\,\,\,\,\,\,\,\,\ve\int_G\int_\bfM\sup_{0\leq t\leq T/\ve}\big\vert
 \int_0^t\big(f(\ve s,X^\ve_{x,y}(s),Y^\ve_{x,y}(s))\\
 &-\bar f(\ve s, X^\ve_{x,y}(s))\big)ds\big\vert d\mu(x,y)
 \leq 2\big(\hat L_fT+T^2(L+KL+KK_L)\big)/n(\ve)\nonumber\\ 
 &+2LT^2/N+T\del+2\hat L_f\max_{0\leq l\leq N-1}\eta_l(\ve,\del)\nonumber
 \end{eqnarray}
 where 
 \[
 \eta_l(\ve,\del)=\max_{0\leq j\leq n(\ve)-1}\mu\big\{(G\times\bfM)\cap
 \Phi_\ve^{-jt(\ve)}\cE_\ve^{f_l}(t(\ve),\del)\big\}.
 \]
 By the assumption there exists an integer valued function $n(\ve)\to\infty$
 as $\ve\to 0$ such that (\ref{1.11.6}) holds true for all $g=f_0,f_1,...,
 f_{N-1}$
 and then $\max_{0\leq l\leq N-1}\eta_l(\ve,\del)\to 0$ as $\ve\to 0$. Hence,
 letting first $\ve\to 0$, then $\del\to 0$ and, finally, $N\to 0$ we obtain 
 (\ref{1.11.17}) in view of (\ref{1.11.18}) and (\ref{1.11.19}), completing the proof
  of Theorem \ref{thm1.11.1}.
 \end{proof}
    
Observe that (\ref{1.11.5}) holding true for all bounded continuous functions
is, in principle, stronger than the averaging principle in the form 
(\ref{1.11.19}) since the latter is equivalent to (\ref{1.11.5}) with $f=B$.
In fact, if we require (\ref{1.11.6}) only for one function $g=B$ then in
the same way as in the proof of Theorem \ref{thm1.11.1} above we conclude
that (\ref{1.11.6}) is equivalent to (\ref{1.11.19}) if
we consider the latter over all compacts $G\subset\bbR^d$ (which was proved
earlier in Theorem 2.1 of \cite{Ki7}). Still, the main
interesting classes of systems, we are aware of, for which (\ref{1.11.5}) holds
true are the same for which (\ref{1.11.19}) is satisfied though it is easy to
construct examples of (somewhat degenerate) right hand sides $b$ and $B$ in
(\ref{1.1.1}) for which (\ref{1.11.19}) holds true but (\ref{1.11.5}) fails 
(since
in the latter we require convergence for all functions $f$ and in the 
former only for $f=B$).

Set
\[
\cE^g_0(t,\del)=\big\{ (x,y)\in\bbR^d\times\bfM:\,\big\vert\frac 1t\int_0^t
g(x,F^u_xy)du-\bar g(x)\big\vert >\del\big\}
\] 
where, recall, $F^u_xy=Y_{x,y}(u)$ and $Y(u)$ satisfies $(\ref{1.1.2})$. In the
same way as Corollary 3.1 in \cite{Ki7} we obtain
\begin{corollary}\label{cor1.11.2}
Suppose that there exists an integer valued 
function $n=n(\ve)\to\infty$ as $\ve\to 0$ such that 
$t(\ve)=T(\ve n(\ve))^{-1}=o(\log(1/\ve))$ and for any $\del>0$ and
each bounded Lipschitz continuous function $g=g(x,y)$ on $\bbR^d\times\bfM$,
    \begin{equation}\label{1.11.24}
    \lim_{\ve\to 0}\max_{0\leq j< n(\ve)}\mu\{ 
    \Phi_\ve^{-jt(\ve)}\cE^g_0(t(\ve),\del)\}=0.
    \end{equation}
    Then (\ref{1.11.6}) is also satisfied, and so (\ref{1.11.5}) holds true.
    \end{corollary}
    
    In the same way as in \cite{Ki7} we obtain that (\ref{1.11.24}) holds true
    in the Anosov theorem setup when $\mu_x$ is an $F^t_x$-invariant measure
    which is ergodic for $\la$-almost all $x$, where $\la$ is the normalized 
    Lebesgue measure on a large compact in $\bbR^d$, and $\mu_x(U)=
    \int_Uq(x,y)dm(y)$ with $q(x,y)>0$ differentiable in $x$ and $y$.
     Furthermore, in the same way as in Theorem 2.4 of \cite{Ki7} or similarly
     to Theorem 2.4 of \cite{Ki9} we conclude that (\ref{1.11.6}) and 
     (\ref{1.11.24}) hold true under Assumptions \ref{ass1.2.1} and 
     \ref{ass1.2.2}.
     Moreover, employing the method of \cite{Ki9} the result can be extended
     to some partially hyperbolic systems.
     
     Observe that under Assumptions \ref{ass1.2.1} and \ref{ass1.2.2} we can
     obtain also large deviations bounds in the form (\ref{1.2.16}) and
     (\ref{1.2.17}) for
     \[
     \tilde Z^\ve_{x,y}(t)=\int_0^tf(s,X^\ve_{x,y}(s/\ve),Y^\ve_{x,y}(s/\ve))
     ds
     \]
     with the functional
     \begin{eqnarray*}
     &\tilde S_{0T}(\tilde\gam)=\inf\big\{ S_{0T}(\gam):\, S_{0T}(\gam)=
     \int_0^TI_{\gam_t}(\nu_t)dt,\\
     &\dot {\gam}_t=\bar B_{\nu_t}(\gam_t),\,\tilde\gam_t=\int_0^t\bar 
     f_{\nu_s}(s,\gam_s)ds\,\,\forall t\in[0,T]\big\},\,\,\,\bar f_\nu(s,x)=
     \int f(s,x,y)d\nu(y),
     \end{eqnarray*}
     where $f$ is a bounded Lipschitz continuous vector function. This can
     be done deriving first an estimate similar to Proposition \ref{prop1.3.4}
     for $f(\ve s,x',Y^\ve_v(s))$ in place of $B(x',Y^\ve_v(s))$ there, which
     should follow in the same way as the proof of Proposition 4.4 of 
     \cite{Ki7}, and proceeding similarly to Sections \ref{sec1.4} and 
     \ref{sec1.5} above. Of course, analogous results can be obtained in the 
     discrete time setup of difference equations (\ref{1.1.10}).

     \part[Markov Fast Motions]{Markov Fast Motions}\label{part2}
    \setcounter{equation}{0}

\section{Introduction}\label{sec2.1}\setcounter{equation}{0}

Many real systems can be viewed as a combination of slow and fast motions
which leads to complicated double scale equations. Already in the 19th 
century in applications to celestial mechanics it was well understood
(though without rigorous justification) that a good approximation of the
slow motion can be obtained by averaging its parameters in fast variables.
Later, averaging methods were applied in signal processing and, rather
recently, to model climate--weather \index{climate--weather system}
interactions (see \cite{Ha}, \cite{CMP},
\cite{Ha1} and \cite{Ki4}). The classical setup of averaging justified 
rigorously in \cite{BM} presumes that the fast motion does not depend on the 
slow one and most of the work on averaging treats this case only. On the other
hand, in real systems both slow and fast motions depend on each other which 
leads to the more difficult fully coupled case which we study here.
This setup emerges, in particular, in perturbations of Hamiltonian systems
which leads to fast motions on manifolds of constant energy and slow
motions across them.

It is natural to view double scale models as describing physical systems
 considered as perturbations of an idealized one which depends on parameters
 $x=(x_1,...,x_d)\in\bbR^d$ assumed to be constants (integrals) of motion. In 
 Part \ref{part2} we suppose that the evolution of this idealized sistem is described
 by certain family of Markov processes \index{Markov process}
 $Y_x(t)=Y_{x,y}(t)=Y_{x,y}(t,\om),\,
 Y_{x,y}(0)=y$ on a separable metric space $\bfM$. In the perturbed system 
 parameters start 
 changing slowly in time and we assume that the corresponding slow motion
 \index{slow motion} 
 $X^\ve(t)=X^\ve_{x,y}(t)=X^\ve_{x,y}(t,\om)$ is described by an ordinary 
 differential equations in $\bbR^d$ having the form
 \begin{equation}\label{2.1.1}
\frac {dX^\ve(t)}{dt}=\ve B(X^\ve(t),Y^\ve(t)),\,\, X^\ve(0)=x,\, Y^\ve(0)=y
\end{equation}
where $B:\bbR^d\times \bfM\to\bbR^d$ is Lipschitz continuous and 
the fast motion \index{fast motion}
$Y^\ve(t)=Y^\ve_{x,y}(t)$ evolves on $\bfM$, it depends, 
in general, on the
slow one and tends to $Y_{x,y}(t)$ as $\ve\to 0$. Usually, $Y^\ve(t)$ is
determined by certain equations, in general, coupled with (\ref{2.1.1}) which
means that their coefficients depend on the slow motion $X^\ve(t)$.

Assume that a nonrandom limit
\begin{equation}\label{2.1.2}
\bar B(x)=\lim_{T\to\infty}T^{-1}\int_0^TB(x,Y_{x,y}(t))dt
\end{equation}
exists in some sense, it "essentially" does not depend on $y$ and it depends
Lipschitz continuously on $x$. Then there exists a unique solution 
$\bar X^\ve=\bar X^\ve_x$ of the averaged equation \index{averaged equation}
\begin{equation}\label{2.1.3}
\frac {d\brX^\ve(t)}{dt}=\ve\brB(\brX^\ve(t)),\quad \brX^\ve(0)=x.
\end{equation}
The averaging principle \index{averaging principle} suggests that often
\begin{equation}\label{2.1.4}
\lim_{\ve\to 0}\sup_{0\leq t\leq T/\ve}|X^\ve_{x,y}(t)-\brX_x^\ve(t)|=0
\end{equation}
in some sense. If unperturbed motions $Y^\ve_{x,y}=Y_y$ do not depend on the
slow variables $x$ and $Y^\ve_{x,y}=Y_y$ then the averaged principle holds 
true under quite general circumstances but when the fast motion depends on
the slow one (coupled case) the situation becomes more complicated and
approximation of $X_{x,y}^\ve$ by $X^\ve_x$ in the weak or the average
sense was justified under some conditions in\cite{Kh2} and \cite{Ve1}.
An extension of the averaging principle in the sense of convergence of
Young measures \index{Young measures} is discussed in Section \ref{sec2.9}
below. 

In this work we are interested in large deviations bounds for probabilities 
that the time changed slow motion $Z^\ve(t)=X^\ve(t/\ve)$ belongs to various
 sets of curves which leads, in particular, to exponential bounds of the form
 \begin{equation}\label{2.1.5}
P\big\{ \sup_{0\leq t\leq T}|Z^\ve_{x,y}(t)-\brZ^\ve_x(t)|>
\del\}\leq e^{-\ka/\ve},\,\,\ka,\del>0
\end{equation}
where $\bar Z_x(t)=\bar X^\ve_x(t/\ve)$ satisfies
\begin{equation}\label{2.1.6}
\frac {d\brZ_x^\ve(t)}{dt}=\ve\brB(\brZ_x^\ve(t)),\quad \brZ_x^\ve(0)=x.
\end{equation}
When the fast motion do not depend on the slow one such results were obtained
in \cite{Fre} and \cite{FW} but the coupled case was dealt with much later
in \cite{Ve2} though (as we indicated this to the author) the proof there 
contained a vicious circle and
substantial gaps which, essentially, were fixed recently in \cite{Ve3}.
Still, \cite{Ve3} is rather difficult to follow and we find it useful to 
provide a precise and consistent exposition of this important result which 
also deals with a more general case including fast motions being random
evolutions whose extreme partial cases are diffusions and finite Markov
chains with continuous time. Moreover, we go beyond bounded time large 
deviations and describe the adiabatic behaviour of the slow motion $Z^\ve$ on
exponentially large in $1/\ve$ time intervals such as its exits from a domain
of attraction and transitions between attractors of the averaged system 
(\ref{2.1.6}). We observe that essentially the same proof yields the same
results for a bit more general case when both $B$ in \ref{2.1.1} and the
coefficients of the random evolutions in the next section depend also 
Lipschitz continuously on $\ve$.

We consider also the discrete time case where (\ref{2.1.1}) is replaced by a
difference equation of the form
\begin{equation}\label{2.1.7}
X^\ve(n+1)-X^\ve(n)=\ve B(X^\ve(n),Y^\ve(n)),\,\, X^\ve(0)=X^\ve_x=x
\end{equation}
where $B(x,y)$ is the same as in (\ref{2.1.1}) and the fast motion $Y^\ve(n)
=Y^\ve_{x,y}(n),\, n=0,1,...,, Y^\ve(0)=y$ is a perturbation of a family
$Y_{x,y}(n), n\geq 0$ of Markov chains parametrized by $x\in\bbR^d$. For
somewhat less general discrete time situation large diviations bounds were
obtained in \cite{GV} by a simpler approach but in our more general situation
we can rely only on methods similar to the continuous time case. Moreover,
unlike \cite{GV} we go farther and study also very long time "adiabatic"
behaviour of the slow motion similar to the continuous time case and
illustrate some of the results by computer simulations for simple models.

The strategy and many of arguments in Part \ref{part2} are rather similar to 
Part \ref{part1} where deterministic chaotic fast motions such as Anosov and
Axiom A systems were considered. Still, in view of the heavy dynamical
 systems background and machinery  Part \ref{part1} is hardly accessible
 for most of probabilists. By this reason we give full proofs here
refering to Part \ref{part1}  only for proofs of some general results on large 
deviations, rate functionals and some others which do not rely on the specific
 dynamical systems setup.

\section{Preliminaries and main results}\label{sec2.2}\setcounter{equation}{0}

We will assume that right hand side of (\ref{2.1.1}) is bounded and Lipschitz
continuous, i.e. for some $K>0$,
\begin{equation}\label{2.2.1}
\sup_{x,y}|B(x,y)|\leq K\,\,\mbox{and}\,\, |B(x,y)-B(z,v)|\leq K
\big(|x-z|+d_\bfM (y,v)\big)
\end{equation}
where $d_\bfM$ is the metric on $\bfM$. Our large deviations estimates will be 
derived
under the following general assumption on the fast motion which is satisfied,
as we explain it below, for random evolutions which are Markov processes with 
switching at random times between a finite number of diffusion processes.

\begin{assumption}\label{ass2.2.1} There exist a convex differentiable
in $\be$ and Lipschitz continuous in other variables function $H(x,x',\be)$ 
defined for all $\be\in\bbR^d$ and for $x,x'$ from the closure $\bar\cX$ of a 
relatively compact open connected set $\cX\subset\bbR^d$ and a positive
function $\zeta_{b,T}(\Del,s,\ve)$ satisfying
\begin{equation}\label{2.2.2}
\limsup_{\Del\to 0}\limsup_{\ve\to 0}\limsup_{s\to\infty}\zeta_{b,T}
(\Del,s,\ve)=0
\end{equation}
such that for all $t>0$, $x,x'\in\bar\cX,y\in\bfM$ and $|\be|\leq b$,
\begin{eqnarray}\label{2.2.3}
&\big\vert\frac 1t\log E\exp\langle\be,\int_0^tB\big(x',
Y^\ve_{x,y}(s)\big)ds\rangle\\
&-H(x,x',\be)\big\vert\leq\zeta_{b,T}\big(\ve t,\min(t,(\log 1/\ve)^\la),
\ve\big)\nonumber
\end{eqnarray}
where $\la\in(0,1)$ and $\langle\cdot,\cdot\rangle$ is the inner product.
\end{assumption}

Set
\begin{equation}\label{2.2.4}
 L(x,x',\al)=\sup_{\be\in\bbR^d}\big(<\al,\be>-H(x,x',\be)\big),
  \end{equation}
  $H(x,\be)=H(x,x,\be)$ and $L(x,\al)=L(x,x,\al)$. Since $H(x,h',0)=0$ then
  $L(x,x',\al)\geq 0$. In view of Assumption
  \ref{ass2.2.1} and standard convex analysis duality \index{convex analysis
  duality} results (see 
  \cite{AE} and \cite{Roc}) $L(x,x',\al)$ is (strictly) convex, lower 
  semicontinuous and we have also that
  \begin{equation}\label{2.2.5}
  H(x,x',\be)=\sup_{\al\in\bbR^d}\big(<\al,\be>-L(x,x',\al)\big).
  \end{equation}
  It follows also from Assumption \ref{ass2.2.1} that
  \begin{equation}\label{2.2.6}
  |H(x,x',\be)|\leq K|\be|.
  \end{equation}
  Since $H(x,h',0)=0$ by (\ref{2.2.6}) and $L(x,x',\al)$ is lower semicontinuous
  then it follows from (\ref{2.2.5}) that there exists a unique 
  $\al_{x,x'}\in\bbR^d$ such that 
  \begin{equation}\label{2.2.7}
  L(x,x',\al_{x,x'})=0.
  \end{equation}
  Set $\al_x=\al_{x,x}$. If $\al_x=\al(x)$ depends Lipschitz continuously
  in $x$ then we can define the averaged motion $\bar X^\ve=\bar X_x^\ve$ 
  in this general setup as the solution of the ordinary differential equation
  \begin{equation}\label{2.2.8}
  \frac {d\bar X^\ve(t)}{dt}=\al(\bar X^\ve(t))\,\, X^\ve(0)=x.
  \end{equation}
  Denote by $C_{0T}$ the space of continuous curves $\gam_t=\gam(t),\, t\in 
  [0,T]$ in $\cX$ which is the space of continuous maps of $[0,T]$ into $\cX.$
  For each absolutely continuous $\gam\in C_{0T}$ its velocity $\dot\gam_t$ 
  can be obtained as the almost everywhere limit of continuous functions
$n(\gam_{t+n^{-1}}-\gam_t)$ when $n\to\infty$. Hence $\dot\gam_t$ is 
measurable in $t$, and so we can set \index{$S$-functional}
  \begin{equation}\label{2.2.9}
  S_{0T}(\gam)=\int_0^TL(\gam_t,\dot{\gam}_t)dt
  \end{equation}
   Define the uniform metric on $C_{0T}$ by 
   \[\bfr_{0T}(\gam,\eta)=\sup_{0\leq t\leq T}|\gam_t-\eta_t|
   \]
   for any $\gam,\eta\in C_{0T}.$ Set $\Psi^a_{0T}(x)=\{\gam\in C_{0T}:\,
   \gam_0=x,\, S_{0T}(\gam)\leq a\}.$  
    Since $L(x,\al)$ is lower semicontinuous and convex in $\al$ and, in 
    addition, $L(x,\al)=\infty$ if $|\al|> K$ it follows that the conditions
    of Theorem 3 in 
   Ch.9 of \cite{IT} are satisfied as we can choose a fast growing minorant
   of $L(x,\al)$ required there to be zero in a sufficiently large ball and
   to be equal, say, $|\al|^2$ outside of it. As a result we conclude that
   $S_{0T}$ is lower semicontinuous functional on $C_{0T}$ with respect to
   the metric $\bfr_{0T}$, and so $\Psi^a_{0T}(x)$ is a closed set which plays 
   a crucial role in the large deviations \index{large deviations}
   arguments below. Set $\cX_t=
   \{ x\in\cX:\,\inf_{z\in\partial\cX}|x-z|\geq 2Kt\}$.
   \begin{theorem}\label{thm2.2.2} Suppose that (\ref{2.2.1}) and Assumption 
   \ref{ass2.2.1} hold true. Set $Z^\ve_{x,y}(t)=X^\ve_{x,y}(t/\ve)$ and let 
   $x\in\cX_T$. Then for any $a,\del,\la>0$ and every
    $\gam\in C_{0T},\,\gam_0=x$ there exists $\ve_0=\ve_0(x,\gam,a,\del,\la)
    >0$ such that for $\ve<\ve_0$ uniformly in $y\in\bfM$,
    \begin{equation}\label{2.2.10}
    P\left\{ \bfr_{0T}(Z^\ve_{x,y},\gam)<\del\right\}\geq
    \exp\left\{-\frac 1\ve(S_{0T}(\gam)+\la)\right\}
    \end{equation}
    and
    \begin{equation}\label{2.2.11}
    P\left\{ \bfr_{0T}(Z^\ve_{x,y},\Psi^a_{0T}(x))\geq\del\right\}
    \leq\exp\left\{-\frac 1\ve(a-\la)\right\}.
    \end{equation}
    \end{theorem}
    
    Next, let $V\subset\cX$ be a connected open set and put 
     $\tau^\ve_{x,y}(V)=\inf\{ t\geq 0:\, Z^\ve_{x,y}(t)
    \notin V\}$ where we take $\tau^\ve_{x,y}(V)=\infty$ if $X^\ve_{x,y}(t)\in
    V$ for all $t\geq 0.$ The following result follows directly from
    Theorem \ref{thm2.2.2}.
    \begin{corollary}\label{cor2.2.3} Under the conditions of Theorem 
    \ref{thm2.2.2} for any $T>0$ and $x\in V,$
    \begin{eqnarray*}
    &\lim_{\ve\to 0}\ve\log P\left\{ \tau^\ve_{x,y}(V)<T\right\}\\
    &=-\inf\left\{ S_{0t}(\gam):\,\gam\in C_{0T},\, t\in[0,T],\,\gam_0=x,\,
    \gam_t\not\in V\right\}\nonumber.
    \end{eqnarray*}
    \end{corollary}

    The main class of Markov processes satisfying our conditions which we
    have in mind consists of random evolutions \index{random evolutions}
    on $\bfM=M\times\{ 1,...,N\}$ where $M$ is a compact $n$-dimensional 
    $C^2$ Riemannian manifold and the unperturbed 
    parametric family of Markov processes $Y_{x,y}(t)$ is the pair 
    $Y_{x,v,k}(t)=(\hat Y_{x,v,k}(t),\nu_{x,v,k}(t))$ governed by the
    stochastic differential equations 
    \index{stochastic differential equation}
    \begin{equation}\label{2.2.12}
    d\hat Y_{x,v,k}(t)=\sig_{\nu_{x,v,k}(t)}\big(x,\hat Y_{x,v,k}(t)\big)dw_t
    +b_{\nu_{x,v,k}(t)}\big(x,\hat Y_{x,v,k}(t)\big)dt
    \end{equation}
    where $\hat Y_{x,v,k}(0)=v,\,\nu_{x,v,k}(0)=k$ and for all 
    $1\leq i,j\leq N,\, i\ne j$,
    \begin{equation}\label{2.2.13}
    P\big\{\nu_{x,v,k}(t+\Del)=j\big\vert\nu_{x,v,k}(t)=i,\hat Y_{x,v,k}(t)=w
    \big\}=q_{ij}(x,w)\Del+o(\Del)\,\,\mbox{as}\,
    \Del\downarrow 0. 
    \end{equation}
     We assume that $q_{kl}(x,w),\, k,l=1,...,N$ are bounded 
    positive $C^1$ functions, $\sig_k(x,v)\sig^*_k
    (x,v)=a_k(x,v)=\big( a^{ij}_k(x,v),\, i,j=1,...,n\big)$ is a $C^1$ field
    of positively definite symmetric matrices on $M$, $b_k(x,v)=
    \big(b^1_k(x,v),...,b_k^n(x,v)\big)$ is a $C^1$ vector field and all
    functions are defined and satisfy the above properties for $v\in M$ and
    $x$ belonging to an open neighborhood of $\bar\cX$. Here $w_t$ is the
    Brownian motion and the equation (\ref{2.2.12}) is written in local 
    coordinats. Observe that the existence and some properties of such Markov
    processes are discussed in \cite{Sk}. The generator \index{generator}
    $\cL^x$ of the Markov process $(\hat Y_x(t),\nu_x(t))$ is the
    operator acting on $C^2$ vector functions $f=(f_1,...,f_N)$ on $M$ by
    the formula
    \begin{equation}\label{2.2.14}
    (\cL^xf)_k(y)=\cL^x_kf_k(y)+\sum^N_{l=1}q_{kl}(x,y)\big(f_l(y)-f_k(y)\big)
    \end{equation}
    where $\cL^x_k$ is the elliptic second order differential operator
    \index{elliptic operator}
    \begin{equation}\label{2.2.15}
    \cL^x_k=\frac 12\langle a_k(x,\cdot)\nabla,\nabla\rangle +
    \langle b_k(x,\cdot),\nabla\rangle.
    \end{equation}
     Now, the perturbed fast motion $Y^\ve=(\hat Y^\ve,\nu^\ve)$ satisfies
    \begin{equation}\label{2.2.16}
    d\hat Y^\ve_{x,v,k}(t)=\sig_{\nu^\ve_{x,v,k}(t)}\big(X^\ve_{x,v,k}(t),
    \hat Y^\ve_{x,v,k}(t)\big)dw_t
    +b_{\nu^\ve_{x,v,k}(t)}\big(X^\ve_{x,v,k}(t),\hat Y^\ve_{x,v,k}(t)\big)dt,
    \end{equation}
    $X^\ve_{x,v,k}(0)=x,\, Y^\ve_{x,v,k}(0)=v,\,\nu^\ve_{x,v,k}(0)=k$ and
    \begin{eqnarray}\label{2.2.17}
    &P\big\{\nu^\ve_{x,v,k}(t+\Del)=j\big\vert\nu^\ve_{x,v,k}(t)=i,
    X^\ve_{x,v,k}(t)=z,\hat Y^\ve_{x,v,k}(t)=w\big\}\\
    &=q_{ij}(x,w)\Del+o(\Del)\,\,\mbox{as}\,\,\Del\downarrow 0\,\,\mbox{for 
    all}\,\,1\leq i,j\leq N,\, i\ne j\nonumber
    \end{eqnarray}
    where $X^\ve$ is given by (\ref{2.1.1}) with $B(x,y)=B(x,v,k)=B_k(x,v),\, 
    y=(v,k)$
    smoothly depending on $x$ and $v$, so that the triple $(X^\ve(t),Y^\ve(t),
    \nu^\ve(t))$ is a Markov processes. The following result which will
    be proved in Section \ref{sec2.3+} claims, in particular, that random 
    evolutions above satisfy Assumption \ref{ass2.2.1}
    
    \begin{proposition}\label{prop2.2.4} For the process $Y_x(t)=
    (\hat Y_x(t),\nu_x(t))$ defined by (\ref{2.2.12}) and (\ref{2.2.13}) the 
    limit
    \begin{equation}\label{2.2.18}
    H(x,x',\be)=\lim_{t\to\infty}\frac 1t\log E\exp\langle\be,\int_0^t
    B_{\nu_{x,v,k}(s)}(x',\hat Y_{x,v,k}(s))ds\rangle
    \end{equation}
    exists uniformly in $x,x'\in\bar\cX$, $y\in\bfM$ and $|\be|\leq b$, it is 
    strictly convex and differentiable in $\be$ and Lipschitz
    continuous in other variables, and it does not depend under our conditions 
    on $v$ and $k$. In this
    circumstances the function $L(x,x',\al)$ given by (\ref{2.2.4}) can be 
    represented in the explicit form 
    \begin{equation}\label{2.2.19}
    L(x,x',\al)=\inf\big\{ I_x(\mu):\,\sum_{k=1}^N\int_MB_k(x,v)d\mu_k(v)
    =\al\big\}
    \end{equation}
    where \index{$I$-functional}
    \begin{equation}\label{2.2.20}
    I_x(\mu)=-\inf_{u>0}\sum_{k=1}^N\int_M\frac {(\cL^xu)_k}{u_k}d\mu_k
    \end{equation}
    and the first infinum is taken over the set $\cP(\bfM)$ of probability
    measures on $\bfM$, i.e. over the vector measures $\mu=(\mu_1,...,\mu_N)$
    with $\sum_{k=1}^N\mu_k(M)=1$, and the second one is taken over positive 
    vector functions $u$ on $M$ belonging to the domain of the operator $\cL^x$.
    Clearly, $I_x(\mu)\geq 0$ and, furthermore, $I_x(\mu)=0$ if and only if
    $\mu$ is the invariant measure \index{invariant measure}
    $\mu^x=(\mu_1^x,...,\mu_N^x)$ of the Markov
     process $Y_x$ which is unique in our circumstances since the Doeblin 
     condition \index{Doeblin condition} (see \cite{Do})
     holds true here. The vector field $\bar B(x)=\int_\bfM
      B(x,y)d\mu^x(y)=\sum_{k=1}^N\int_MB_k(x,v)d\mu_k(v)$ is $C^1$ in $x$, 
      and so we can define the averaged motion \index{averaged motion}
       $\bar X^\ve=\bar X^\ve_x$ by
    \begin{equation}\label{2.2.21}
  \frac {d\bar X^\ve(t)}{dt}=\ve\bar B(\bar X^\ve(t)),\,\,\, X^\ve(0)=x.
  \end{equation}
  Hence, $S_{0T}(\gam)=0$ if and only if $\gam_t=\bar Z(t)=\bar X^\ve(t/\ve)$
  for all $t\in[0,T]$.
    The processes $Y^\ve$ given by (\ref{2.2.16}) and (\ref{2.2.17}) together
    with the function $H(x,x',\be)$ satisfy Assumption \ref{ass2.2.1}.
    \end{proposition}
   
     Clearly, if $N=1$ above then $Y^\ve$ becomes a
     diffusion process \index{diffusion process} and if all operators
     $\cL^x_k$ are just zero then we arrive to the case of continuous time
     Markov chains as fast motions \index{continuous time Markov chain} which
     also satify all our assumptions. We observe also that both
    Proposition \ref{prop2.2.4} and the results below can be extended to the
    case when $\cL^x_k$ are hypoelliptic operators satisfying natural 
    conditions so that we could rely, in particular, on results of Section 6.3
    from \cite{DS}.
    
    Suppose that the coefficients $\sig_k,b_k$ and $q_{ij}$ in (\ref{2.2.12})
    and (\ref{2.2.13}) do not depend on $x$. Then $Y^\ve_{x,y}(t)=Y_{x,y}(t)
    =Y_y(t)$ is an ergodic Markov process with the unique invariant measure
    $\mu$ and for any $y$ almost surely
    \begin{equation*}
    \lim_{T\to\infty}T^{-1}\int_0^TB(x,Y_y(t))dt=\bar B(x)=\int B(x,y)d\mu(y)
    \end{equation*}
    and by standard general results on the uncoupled averaging (see \cite{SV})
    it follows that for any $y$ almost surely
    \begin{equation}\label{2.2.22}
    \sup_{0\leq t\leq T}|X^\ve_{x,y}(t/\ve)-\bar X^\ve_x(t/\ve)|\to 0\,\,\,
    \mbox{as}\,\,\,\ve\to 0.
    \end{equation}
    In the fully coupled case (i.e. when $a_k,b_k,q_{ij}$ depend on $x$) 
    Theorem \ref{thm2.2.2} implies in the case of fast motions given by
    (\ref{2.2.16}) and (\ref{2.2.17}) that for each $\del>0$ there is 
    $\al(\del)>0$ such that for all small $\ve$,
    \begin{equation}\label{2.2.22+}
    P\{\sup_{0\leq t\leq T}|X^\ve_{x,y}(t/\ve)-\bar X^\ve_x(t/\ve)|>\del\}
    \leq e^{-\al(\del)/\ve}
    \end{equation}
    which means, in particular, that in this case we have in (\ref{2.2.22}) 
    convergence in probability. Examples from \cite{BK1} show that, in general,
    in the fully coupled setup we do not have convergence in (\ref{2.2.22})
    with probability one though in some
    cases such convergence can be derived from (\ref{2.2.22+}) if the
    derivatives of $X^\ve$ and $Y^\ve$ in $\ve$ grow subexponentially in 
    $1/\ve$ on time intervals of order $1/\ve$ (see Remark \ref{rem2.3.6}).
    
    In the following assertions we assume always that the fast 
    motions are obtained by means of (\ref{2.2.12}) and (\ref{2.2.13}) so that
    we could rely on (\ref{2.2.18})--(\ref{2.2.20}) though, in principle, it is
     possible to impose some general conditions on functions $L(x,\al)$ which
     would enable us to proceed with our arguments. 
    
    Precise large deviations bounds such as (\ref{2.2.10}) and (\ref{2.2.11})
     of Theorem \ref{thm2.2.2} are crucial in our study in Sections \ref{sec2.5}
     and \ref{sec2.6} of the "very long", i.e. exponential in 
    $1/\ve$, time "adiabatic" behaviour \index{adiabatic behavior} of the slow
     motion. Namely, we
    will describe such long time behavior of $Z^\ve$ in terms of the function 
     \[
     R(x,z)=\inf_{t\geq 0,\gam\in C_{0t}}\{ S_{0t}(\gam):\,\gam_0=x,\,
     \gam_t=z\}
     \]
     under various assumptions on the averaged motion $\brZ.$
     Observe that $R$ satisfies the triangle inequality $R(x_1,x_2)+
     R(x_2,x_3)\geq R(x_1,x_3)$ for any $x_1,x_2,x_3\in\cX$ and it determines
     a semi metric on $\cX$ which measures "the difficulty'" for the slow
     motion to move from point to point in terms of the functional $S$. 
     
     Introduce the averaged flow \index{averaged flow} $\Pi^t$ on $\cX_t$ by
     \begin{equation}\label{2.2.23}
     \frac {d\Pi^tx}{dt}=\bar B(\Pi^tx),\,\, x\in\cX_t
     \end{equation}
     where $\bar B(z)$ is the same as in (\ref{2.2.21}) and set $\bar B_\mu(z)=
     \int_\bfM B(z,y)d\mu(y)=\sum_{k=1}^N\int_MB_k(x,v)d\mu_k(v)$ for any 
     probability measure $\mu=(\mu_1,...,\mu_N)$ on $\bfM=M\times\{ 1,...,N\}$.
     Call a $\Pi^t$-invariant compact set $\cO\subset\cX$ an $S$-compact
     if for any $\eta>0$ there exist $T_\eta\geq 0$ and an open set 
     $U_\eta\supset\cO$ such that whenever $x\in\cO$ and $z\in U_\eta$ we can
     pick up $t\in[0,T_\eta]$ and $\gam\in C_{0t}$ satisfying
     \[
     \gam_0=x,\,\, \gam_t=z\,\,\mbox{and}\,\, S_{0t}(\gam)<\eta.
     \]
     It is clear from this definition that $R(x,z)=0$ for any pair points
     $x,z$ of an $S$-compact $\cO$ and by the above triangle inequality
     for $R$ we see that $R(x,z)$ takes on the same value when $z\in\cX$
     is fixed and $x$ runs over $\cO$.
     We say that the vector field $B$ on $\cX\times\bfM$ is complete 
     \index{complete} at 
     $x\in\cX$ if the convex set of vectors $\{\be\bar B_\mu(x):\,\be\in[0,1],
     \mu\in\cP(\bfM),\, I_x(\mu)<\infty\}$ contains an open neigborhood of 
     the origin in $\bbR^d$.
    It follows by Lemma \ref{lem1.6.2} in Part \ref{part1} that if 
    $\cO\subset\cX$ is a compact
     $\Pi^t$-invariant set such that $B$ is complete at each $x\in\cO$ and
      either $\cO$ contains a dense orbit of the flow $\Pi^t$ (i.e. $\Pi^t$ is
     topologically transitive \index{topologically transitive}
     on $\cO$) or $R(x,z)=0$ for any $x,z\in\cO$
      then $\cO$ is an $S$-compact. Moreover, to ensure that $\cO$ is an 
      $S$-compact it suffices to assume that $B$ is complete already at some 
      point of $\cO$ and the flow $\Pi^t$ on $\cO$ is minimal, \index{minimal}
      i.e. the $\Pi^t$-orbits of all points are dense in $\cO$ or, equivalently,
       for any $\eta>0$ there exists
       $T(\eta)>0$ such that the orbit $\{\Pi^tx,\, t\in[0,T(\eta)]\}$ of
       length $T(\eta)$ of each point $x\in\cO$ forms an $\eta$-net in $\cO$
       which is equivalent to minimality of the flow $\Pi^t$ on $\cO$
       (see \cite{Wa}).
       The latter condition obviously holds true when $\cO$ is a fixed point
       or a periodic orbit of $\Pi^t$ but also, more generally, when $\Pi^t$ 
       on $\cO$ is uniquely ergodic (see \cite{Wa}).
       
      A compact $\Pi^t$-invariant set $\cO\subset\cX$ is called an attractor
     (for the flow $\Pi^t$) if there is an open set $U\supset\cO$ and $t_U>0$
     such that
     \[
     \Pi^{t_U}\bar U\subset U\,\,\mbox{and}\,\,\lim_{t\to\infty}\,
     \mbox{dist}(\Pi^tz,\cO)=0\,\,\mbox{for all}\,\, z\in U.
     \]
     For an attractor $\cO$ the set $V=\{ z\in\cX:\,\lim_{t\to\infty}\,
     \mbox{dist}(\Pi^tz,\cO)=0\}$, which is clearly open, is called the
     basin (domain of attraction) of $\cO$. An attractor which is also
     an $S$-compact will be called an $S$-attractor \index{$S$-attractor}.
     
     In what follows we will speak about connected open sets $V$ with 
     piecewise smooth boundaries $\partial V$. The latter can be introduced in 
     various ways but it will be convenient here to adopt the definition from 
     \cite{Cow} saying that $\partial V$ is the closure of a finite union
     of disjoint, connected, codimension one, extendible $C^1$ (open or
     closed) submanifolds of $\bbR^d$ which are called faces of the boundary.
     The extendibility condition means that the closure of each face is a part 
     of a larger submanifold of the same dimension which coincides with the
     face itself if the latter is a compact submanifold. This enables us to 
     extend fields of normal vectors to the boundary of faces and to speak 
     about minimal angles between adjacent faces which we assume to be 
     uniformly bounded away from zero or, in other words, angles between 
     exterior normals to adjacent faces at a point of intersection of their 
     closures are uniformly bounded away from $\pi$ and $-\pi$. The following
     result which will be proved in Section \ref{sec2.5} describes exits of the
     slow motion from neighborhoods of attractors of the averaged motion.

     \begin{theorem}\label{thm2.2.5} Let $\cO\subset\cX$ be an $S$-attractor
      of the flow $\Pi^t$ whose basin contains the closure $\bar V$ of a 
     connected open set $V$ with a piecewise smooth boundary $\partial V$ 
     such that $\bar V\subset\cX$ and assume that for each $z\in\partial V$
     there exists $\vp=\vp(z)>0$ and a probability measure $\eta=\eta_z$ with
     $I_z(\eta_z)<\infty$ such that
     \begin{equation}\label{2.2.24}
     z+s\bar B(z)\in V\,\,\,\mbox{but}\,\,\, z+s\bar B_\eta(z)\in
     \bbR^d\setminus\bar V\,\,\mbox{for all}\,\, s\in(0,\vp],
     \end{equation}
     i.e. $\bar B(z)\ne 0,\,\bar B_\eta(z)\ne 0$ and the former vector
     points out into the interior while the latter into the exterior 
     of $V$. Set $R_{\partial}(z)=\inf\{ R(z,\tilde z):\,\tilde z\in
     \partial V\}$ and $\partial_{\min}(z)=\{\tilde z\in\partial V:\,
      R(z,\tilde z)=R_\partial(z)\}$. Then $R_\partial(z)$ takes on the 
      same value $R_\partial$ and $\partial_{\min}(z)$ coincides with
      the same compact nonempty set $\partial_{\min}$ for all $z\in\cO$
      while $R_\partial(x)\leq R_\partial$ for all $x\in V$. Furthermore,
      for any $x\in V$ uniformly in $y\in\bfM$,
     \begin{equation}\label{2.2.25}
     \lim_{\ve\to 0}\ve\log E\tau^\ve_{x,y}(V)=R_{\partial}>0
     \end{equation}
     and for each $\al>0$ there exists $\la(\al)=\la(x,\al)>0$ such that 
     uniformly in $y\in\bfM$ for all small $\ve>0$,
     \begin{equation}\label{2.2.26}
      P\big\{  e^{(R_\partial-\al)/\ve}>\tau^\ve_{x,y}(V)\,\,
     \mbox{or}\,\,\tau^\ve_{x,y}(V)>e^{(R_\partial+\al)/\ve}\big\}
     \leq e^{-\la(\al)/\ve}.
     \end{equation}
     Next, set
     \[
     \Te^\ve_v(t)=\Te_v^{\ve,\del}(t)=\int_0^t
     \bbI_{V\setminus U_\del(\cO)}(Z^\ve_v(s))ds
     \]
     where $U_\del(\cO)=\{ z\in\cX:$ dist$(z,\cO)<\del\}$ and $\bbI_\Gam(z)
     =1$ if $z\in\Gam$ and $=0$, otherwise. Then for any $x\in V$ and 
     $\del>0$ there exists $\la(\del)=\la(x,\del)>0$ such that 
     uniformly in $y\in\bfM$ for all small $\ve>0$,
     \begin{equation}\label{2.2.27}
     P\big\{ \Te^\ve_{x,y}(\tau^\ve_{x,y}(V))\geq e^{-\la(\del)/\ve}
     \tau^\ve_{x,y}(V)\big\}\leq  e^{-\la(\del)/\ve}.
     \end{equation}
     Finally, for every $x\in V$ and $\del>0$,
     \begin{equation}\label{2.2.28}
     \lim_{\ve\to 0}P\big\{ \mbox{dist}\big(Z^\ve_{x,y}
     (\tau^\ve_{x,y}(V)),\partial_{\min}\big)\geq\del\big\}=0
     \end{equation}
     provided $R_\partial <\infty$ and the latter holds true if and only if
      for some $T>0$ there exists 
     $\gam\in C_{0T},\,\gam_0\in\cO,\,\gam_T\in\partial V$ such that 
     $\dot {\gam}_t=\bar B_{\nu_t}(\gam_t)$ for Lebesgue almost all 
     $t\in[0,T]$ with $\nu_t\in\cM_{\gam_t}$ then $R_\partial <\infty$.
     \end{theorem}
     
     Theorem \ref{thm2.2.5} asserts, in particular, that typically the slow 
     motion $Z^\ve$ performs rare (adiabatic) fluctuations in the vicinity of 
     an $S$-attractor $\cO$ since it exists from any domain $U\supset\cO$ with 
     $\bar U\subset V$ for the time much smaller than $\tau^\ve(V)$ (as the 
     corresponding number $R_\partial=R_{\partial U}$ will be smaller) and by 
     (\ref{2.2.27}) it can spend in $V\setminus U_\del(\cO)$ only small 
     proportion of time which implies that $Z^\ve$ exits from $U$ and returns 
     to $U_\del(\cO)$ (exponentially in $1/\ve$) many times before it finally 
     exits $V$. We observe that in the much simpler uncoupled setup 
     corresponding results in the case of $\cO$ being an attracting point were 
     obtained for a continuous time Markov chain as a fast motion
      in \cite{Fre} but the proofs there 
     rely on the lower semicontinuity of the function $R$ which does not hold
     true in general, and so extra conditions like $S$-compactness of $\cO$
     or, more specifically, the completness of $B$ at $\cO$ should be assumed
     there, as well. It is important to observe that the intuition based on
     diffusion type small random perturbations of dynamical systems should be
     applied with caution to problems of large deviations in averaging since
     the $S$-functional of Theorem \ref{thm2.2.2} describing them is more
     complex and have rather different properties than the corresponding
     functional emerging in diffusion type random perturbations of dynamical
     systems (see \cite{FW}). The reason for this is the deterministic nature
     of the slow motion $Z^\ve$ which unlike a diffusion can move only with
     a bounded speed and, moreover, even in order to ensure its "diffusive
     like" local behaviour (i.e. to let it go in many directions) some extra
     nondegeneracy type conditions on the vector field $B$ are required.
     
     Our next result describes rare (adiabatic) transitions \index{adiabatic
     transitions} of the slow motion
     $Z^\ve$ between basins of attractors of the averaged flow $\Pi^t$ which 
     we consider now in the whole $\bbR^d$ and impose certain conditions on the
     structure of its $\om$-limit set.
     \begin{assumption}\label{ass2.2.6} Assumption \ref{ass2.2.1} holds true 
     for $\cX=\bbR^d$, the families $a^{ij}_k(x,\cdot),\, i,j=1,...,n\big)$ 
     and $b_k(x,\cdot)=\big(b^1_k(x,\cdot),...,b_k^n(x,\cdot)\big)$ of matrix
     and vector fields are compact sets in the $C^1$ topology, 
     \begin{equation}\label{2.2.29}
     \| B(x,y)\|_{C^2(\bbR^d\times\bfM)}\leq K
     \end{equation}
     for some $K>0$ independent of $x,y$ and there exists $r_0>0$ such that
     \begin{equation}\label{2.2.30}
     \big(x,B(x,y)\big)\leq -K^{-1}\,\,\mbox{for any}\,\, y\in\bfM
     \,\,\mbox{and}\,\, |x|\geq r_0.
     \end{equation}
     \end{assumption}
     
     The condition (\ref{2.2.30}) means that outside of some ball all vectors
     $B(x,y)$ have a bounded away from zero projection on the radial direction
     which points out to the origin. This condition can be weakened, for
     instance, it suffices that
     \[
     \lim_{d\to\infty}\inf\{ R(x,z):\,\mbox{dist}(x,z)\geq d\}=\infty
     \]
     but, anyway, we have to make some assumption which ensure that the slow
     motion stays in a compact
     region where really interesting dynamics takes place.
     
     Next, suppose that the $\om$-limit set of the averaged flow $\Pi^t$ is 
     compact and it consists of two parts, so that the first part is a finite 
     number of $S$-attractors $\cO_1,...,\cO_\ell$ whose basins
     $V_1,...,V_\ell$ have piecewise smooth boundaries $\partial V_1,...,
     \partial V_\ell$ and the remaining part of the $\om$-limit set is
     contained in $\cup_{1\leq j\leq\ell}\partial V_j$. We assume also
     that for any $z\in\cap_{1\leq i\leq k}\partial V_{j_i},\, k\leq\ell$
     there exist $\vp=\vp(z)>0$ and probability measures 
     $\eta_1,...,\eta_k$ such that $I(\eta_i)<\infty,\, i=1,...,k$ and
     \begin{equation}\label{2.2.31}
     z+s\bar B_{\eta_i}(z)\in V_{j_i}\,\,\mbox{for all}\,\, s\in(0,\vp]\,\,
     \mbox{and}\,\, i=1,...,k,
     \end{equation}
     i.e. $\bar B_{\eta_i}(z)\ne 0$ and it points out into the interior of
     $V_{j_i}$ which means that from any boundary point it is possible to
     go to any adjacent basin along a curve with an arbitrarily small
     $S$-functional. Let $\del>0$ be so small that the $\del$-neighborhood 
     $U_\del(\cO_i)=\{ z\in\cX:\,\mbox{dist}(z,\cO_i)<\del\}$
      of each $\cO_i$ is contained with its closure in the corresponding
      basin $V_i$. For any $x\in V_i$ set
     \[
     \tau^\ve_{x,y}(i)=\tau^{\ve,\del}_{x,y}(i)=\inf\big\{ t\geq 0:\,
      Z^\ve_{x,y}(t)\in\cup_{j\ne i}U_\del(\cO_j)\big\}.
     \]
     
      In Section \ref{sec2.6} we will derive the following result.
      
     \begin{theorem}\label{thm2.2.7} The function $R_{ij}(x)=\inf_{z\in V_j}
     R(x,z)$ takes on the same value $R_{ij}$ for all $x\in\cO_i,\, i\ne j$.
     Let $R_i=\min_{j\ne i,j\leq\ell}R_{ij}$. Then for any $x\in V_i$
     uniformly in $y\in\bfM$,
      \begin{equation}\label{2.2.32}
       \lim_{\ve\to 0}\ve\log E\tau^\ve_{x,y}(i)=R_i>0
     \end{equation}
      and for any $\al>0$ there exists $\la(\al)=\la(x,\al)>0$ such that for 
      all small $\ve>0$,
       \begin{equation}\label{2.2.33}
       P\big\{ e^{(R_i-\al)/\ve}>
       \tau^\ve_{x,y}(i)\,\,
       \mbox{or}\,\,\tau^\ve_{x,y}(i)>e^{(R_i+\al)/\ve}
       \big\}\leq e^{-\la(\al)/\ve}.
      \end{equation}
       Next, set
     \[
     \Te^{\ve,i}_v(t)=\Te_v^{\ve,i,\del}(t)=\int_0^t
     \bbI_{V_i\setminus U_\del(\cO_i)}(Z^\ve_v(s))ds.
     \]
      Then for any $x\in V_i$ and $\del>0$ there exists $\la(\del)=
      \la(x,\del)>0$ such that uniformly in $y\in\bfM$ for all small $\ve>0$,
     \begin{equation}\label{2.2.34}
     P\big\{ \Te^{\ve,i}_{x,y}(\tau^\ve_{x,y}(i))\geq 
     e^{-\la(\del)/\ve}\tau^\ve_{x,y}(i)\big\}\leq 
     e^{-\la(\del)/\ve}.
     \end{equation}
     Now, suppose that the vector field $B$ is complete on 
     $\partial V_{i}$ for some $i\leq\ell$ (which strengthens (\ref{2.2.31})
     there) and the restriction of the $\om$-limit set of $\Pi^t$ to 
     $\partial V_{i}$ consists of a finite number of $S$-compacts. Assume
     also that there is a unique index $\io(i)\leq\ell,\, \io(i)\ne i$ such 
     that $R_{i}=R_{i\io(i)}$. Then for any $x\in V_i$ there exists
      $\la=\la(x)>0$ such that uniformly in $y\in\bfM$ for all small 
     $\ve>0$,
      \begin{equation}\label{2.2.35}
      P\big\{ Z^\ve_{x,y}(\tau^\ve_{x,y}(i))\not\in V_{\io(i)}
      \big\}\leq e^{-\la/\ve}.
      \end{equation}
      Finally, suppose that the above conditions hold true for all 
      $i=1,...,\ell$. Define $\io_0(i)=i$, $\tau^\ve_v(i,1)=\tau^\ve_v(i)$
      and recursively, 
      \[
      \io_k(i)=\io(\io_{k-1}(i))\,\,\mbox{and}\,\,\tau^\ve_v(i,k)=
      \tau^\ve_v(i,k-1)+\tau^\ve_{v_\ve(k-1)}
      \big(j(v_\ve(k-1))\big),
      \]
      where $v_\ve(k)=\Phi_\ve^{\ve^{-1}\tau^\ve_v(i,k)}v$,
      $j((x,y))=j$ if $x\in V_j$, and set $\Sigma^\ve_i(k,a)=\sum_{l=1}^k
      \exp\big((R_{\io_{l-1}(i),\io_l(i)}+a)/\ve\big)$. Then for any $x\in V_i$
       and $\al>0$ there 
      exists $\la(\al)=\la(x,\al)>0$ such that uniformly in $y\in\bfM$
      for all $n\in\bbN$ and sufficiently small $\ve>0$,
       \begin{eqnarray}\label{2.2.36}
       &P\big\{ \Sigma^\ve_i(k,-\al)>
       \tau^\ve_{x,y}(i,k)\,\,\mbox{or}\,\,\\
       &\tau^\ve_{x,y}(i,k)>\Sigma^\ve_i(k,\al)\,\,\mbox{for some}\,\,
       k\leq n\big\}\leq ne^{-\la(\al)/\ve}\nonumber
      \end{eqnarray}
      and for some $\la=\la(x)>0$,
       \begin{equation}\label{2.2.37}
      P\big\{ Z^\ve_{x,y}(\tau^\ve_{x,y}(i,k))\not\in V_{\io_k(i)}
      \,\,\mbox{for some}\,\, k\leq n\big\}\leq ne^{-\la/\ve}.
      \end{equation}
      \end{theorem}
      
      Generically there exists only one index $\io(i)$ such that $R_i=
      R_{i\io(i)}$ and in this case Theorem \ref{thm2.2.7} asserts that
      $Z^\ve_{x,y},\, x\in V_i$ arrives (for "most" $y\in\cW$) at 
      $V_{\io(i)}$ after it leaves $V_i$. If $\cI(i)=\{ j:\, R_i=R_{ij}\}$
      contains more than one index then the method of the proof of Theorem
      \ref{thm2.2.7} enables us to conclude that in this case $Z^\ve_{x,y},\,
      x\in V_i$ arrives (for "most" $y\in\cW$) at $\cup_{j\in\cI(i)}V_j$
      after leaving $V_i$ but now we cannot specify the unique basin of
      attraction of one of $\cO_j$'s where $Z^\ve_{x,y}$ exits from $V_i$.
      If the succession function $\io$ is uniquely defined then it 
      determines an order of transitions of the slow motion $Z^\ve$
      between basins of attractors of $\bar Z$ and because of their finite
      number $Z^\ve$ passes them in certain cyclic order going around such
      cycle exponentially many in $1/\ve$ times while spending the total time
      in a basin $V_i$ which is approximately proportional to $e^{R_i/\ve}$.
      If there exist several cycles of indices $i_0,i_1,...,i_{k-1},i_k=i_0$ 
      where $i_j\leq\ell$ and $i_{j+1}=\io(i_j)$ then transitions between 
      different cycles may also be possible. In the uncoupled case with fast 
      motions being continuous time Markov chains a description of such 
      transitions via certain hierarchy of cycles appeared without a detailed 
      proof in \cite{Fre} and \cite{FW}. In our fully coupled setup the
       corresponding description does not seem to be different from the
       uncoupled situation since its justification relies only on the Markov
       property arguments and estimates of probabilities of transitions
       of $Z^\ve$ from $U_\del(\cO_i)$ to $U_\del(\cO_j)$.
     
     Set $\cI=\{ 1,...,\ell\}$. Following \cite{FW} we call a graph consisting
     of arrows $(k\to l)$ $(k\ne i,\, k,l\in\cI,\, k\ne l)$ an $i$-graph if
     every point $l\ne i$ is the origin of exactly one arrow and the graph
     has no circles. Let $G(i)$ be the set of all $i$-graphs. Next, choose
     $\del>0$ so small that $\overline {U_{2\del}(\cO_j)}\subset V_j,\, j=1,...,
     \ell$ and define stopping times $\sig_{x,y}^{\ve,\del}(0)=0$ and by
     induction for $k\geq 1$,
     \[
     \hat\sig_{x,y}^{\ve,\del}(k)=\inf\{ t\geq\sig_{x,y}^{\ve,\del}(k-1):\,
     Z^\ve_{x,y}(t)\not\in\cup_{1\leq j\leq\ell}U_{2\del}(\cO_j)\},
     \]
     \[
     \sig_{x,y}^{\ve,\del}(k)=\inf\{ t\geq\hat\sig_{x,y}^{\ve,\del}(k):\,
     Z^\ve_{x,y}(t)\in\cup_{1\leq j\leq\ell}U_{\del}(\cO_j)\}.
     \]
     Define the Markov chain 
     \[
     W^\ve_{x,y}(k)=\big(Z^\ve_{x,y}(\sig^{\ve,\del}_{x,y}(k)),\,
      Y^\ve_{x,y}(\ve^{-1}\sig^{\ve,\del}_{x,y}(k))\big)
     \]
     which evolves on the phase space $\cup_{1\leq j\leq\ell}\Gam_j$ where
     $\Gam_j=\partial U_\del(\cO_j)\times\bf M$.
     \begin{theorem}\label{thm2.2.8} Let $P(\cdot,\cdot)$ be the
     transition probability of the Markov chain $W^\ve$. Then for any $\be>0$
     there exist $\del_0,\ve_0>0$ such that if $\del<\del_0$ and $\ve<\ve_0$ 
     then
     \begin{equation}\label{2.2.38}
     e^{(-R_{ij}-\be)/\ve}\leq P\big((x,y),\Gam_j\big)\leq
     e^{(-R_{ij}+\be)/\ve}
     \end{equation}
     whenever $(x,y)\in\Gam_i$. Furthermore, if $\mu^\ve_W$ is an invariant
     measure of $W^\ve$ on $\cup_{1\leq j\leq\ell}\Gam_j$ then
     \begin{equation}\label{2.2.39}
     e^{-2\be(\ell -1)/\ve}\frac {Q_j}{Q_1+\cdots +Q_\ell}\leq
     \mu^\ve_W(\Gam_j)\leq e^{2\be(\ell-1)/\ve}\frac {Q_j}{Q_1+\cdots +Q_\ell}
     \end{equation}
     where
     \begin{equation}\label{2.2.40}
     Q_i=\sum_{g\in G(i)}\exp\big(-\ve^{-1}\sum_{(k\to l)\in g}R_{kl}\big).
     \end{equation}
     \end{theorem}
     
     Since total times spent by a Markov process in different sets are
     asymptotically proportional to masses given to these sets by corresponding
     invariant measures then Theorem \ref{thm2.2.8} (together with Theorem
     \ref{thm2.2.7}) yields actually that the slow motion $Z^\ve(t)$ spends in 
     a basin $V_j$ of the attractor $\cO_j$ a percentage of total time
     approximately proportional to $Q_j$ which will be illustrated by
     computational examples in Section \ref{sec2.7}.
      In fact, this description is effective 
     only if there is a unique $i_0$ and a graph $g\in G(i_0)$ such that
     $\sum_{(k\to l)\in g}R_{kl}$ is minimal possible among all such sums
     over all $j$-graphs. In this case the slow motion spends in $V_{i_0}$
     a proportion of time close to one.
      
     Next, we formulate our results for the discrete time case of difference 
     equations (\ref{2.1.7}).
 \begin{assumption}\label{ass2.2.9} There exist a convex differentiable in
 $\be$ and Lipschitz continuous in other variables function $H(x,x',\be)$ 
defined for all $\be\in\bbR^d$ and for $x,x'$ from the closure $\bar\cX$ of a 
relatively compact open connected set $\cX\subset\bbR^d$ and a positive
function $\zeta_{b,T}(\Del,s,\ve)$ satisfying (\ref{2.2.2})
such that for all $k>0$, $x,x'\in\bar\cX,y\in\bfM$ and $|\be|\leq b$,
\begin{eqnarray}\label{2.2.41}
&\big\vert\frac 1k\log E\exp\langle\be,\sum_{j=1}^kB\big(x',
Y^\ve_{x,y}(j)\big)\rangle\\
&-H(x,x',\be)\big\vert\leq\zeta_{b,T}\big(\ve k,\min(k,(\log 1/\ve)^\la),
\ve\big)\nonumber
\end{eqnarray}
where $\la\in(0,1)$ and $Y^\ve_{x,y}(n), n=0,1,2,...$ appears in (\ref{2.1.7}).
\end{assumption}

\begin{theorem}\label{thm2.2.10} Suppose that (\ref{2.2.1}) and Assumption 
\ref{ass2.2.9} are satisfied and that $X^\ve(n)=X^\ve_{x,y}(n),n=0,1,2,...$
are given by (\ref{2.1.7}).
 For $t\in [n,n+1]$ define $X^\ve(t)=(t-n)X^\ve(n+1)+(n+1-t)X^\ve(n)$ and set
 $Z^\ve(t)=Z^\ve_{x,y}(t)=X^\ve_{x,y}(t/\ve)$. Then Theorem \ref{thm2.2.2}
 and Corollary \ref{cor2.2.3} hold true with the corresponding functionals
 $S_{0T}$.
 \end{theorem}
 In a bit more restricted situation Theorem \ref{2.2.10} was proved by a
 simpler method in \cite{GV}.
 
 The main model of Markov chains \index{Markov chains}
 serving as fast motions $Y^\ve(n),n\geq 0$,
 we have in mind, is obtained in the following way. We start with a parametric
 family of Markov chains $Y_{x,y}(n),\, n\geq 0,\, Y_{x,y}(0)=y$ on a compact
 $C^2$ Riemannian manifold $M$ with transition probabilities $P^x(y,\Gam)=
 P^x_y\{ Y_{x,y}(1)\in\Gam\}$ having positive densities $p^x(y,z)=
 P^x(y,dz)/m(dz)$ with respect to the Riemannian volume $m$, so that $p^x(y,z)$ is $C^1$ in 
 $x$ and continuous in other variables. Next, we define $X^\ve(n)$ and 
 $Y^\ve(n)$ adding to (\ref{2.1.7}) another equation
 \begin{equation}\label{2.2.42}
 P\big\{ Y^\ve(n+1)\in\Gam\big\vert X^\ve(n)=x,\, Y^\ve(n)=y\big\}=P^x(y,\Gam).
 \end{equation}
 
 \begin{proposition}\label{prop2.2.11} Let $Y_{x,y}(n)$ be as above. Then the
 limit
 \begin{equation}\label{2.2.43}
 H(x,x',\be)=\lim_{k\to\infty}\frac 1k\log E\exp\langle\be,\sum_{j=1}^k
 B(x',Y_{x,y}(j))\rangle
 \end{equation}
 exists uniformly in $x,x'$ running over a compact set and in $y\in\bfM$
 and it satisfies conditions of Assumption \ref{ass2.2.9}. 
 In this circumstances the functionals $S_{0T}$ appearing in the large
 deviations estimates (\ref{2.2.10}) and (\ref{2.2.11}) again have the form
 (\ref{2.2.9}) with $L(x,\al)$ given by (\ref{2.2.19}) where now
 \begin{equation}\label{2.2.44}
 I_x(\mu)=-\inf_{u>0}\int_\bfM\log\frac {\int_\bfM p^x(y,v)u(v)dm(v)}{u(y)}
 d\mu(y).
 \end{equation}
 Clearly, $I_x(\mu)\geq 0$ and, furthermore, $I_x(\mu)=0$ if and only if
 $\mu$ is the invariant measure $\mu^x$ of the Markov chain $Y_x$ which is
 unique since the Doeblin condition (see \cite{Do}) holds true here. The 
 vector field $\bar B(x)=\int_\bfM B(x,y)d\mu^x(y)$ is $C^1$ in $x$, and so
 we can define uniquely the averaged motion $\bar X^\ve=\bar X^\ve_x$ by
 (\ref{2.2.21}) and, again, $S_{0T}(\gam)=0$ if and only if $\gam_t=
 \bar Z(t)=\bar X^\ve(t/\ve)$ for all $t\in[0,T]$. Furthermore,
 $Y^\ve(n),n\geq 0$ given by (\ref{2.2.42}) satisfies (\ref{2.2.41}).
 \end{proposition}
 The existence of the limit (\ref{2.2.43}) and its properties in our
 circumstances are well known 
 (see \cite{DV1}, \cite{DV2}, \cite{Ki1}, \cite{Ka}, \cite{HH}) and the fact 
 that (\ref{2.2.41}) holds true here will be explained at the beginning of 
 Section \ref{sec2.7}.
 
 \begin{theorem}\label{thm2.2.12} Let the fast motion $Y^\ve(n)=Y^\ve_{x,y}(n)$
 be constructed as above via (\ref{2.2.42}) then with the corresponding 
 definitions of $S$-compacts and under similar conditions the conclusions of
 Theorems \ref{thm2.2.5}, \ref{thm2.2.7} and \ref{thm2.2.8} remain true for the 
 corresponding slow motion $Z^\ve(t)$ defined in Theorem \ref{thm2.2.10}.
 \end{theorem}
 
 Observe, that we can easily produce a wide class of systems satisfying the 
 conditions of Theorems \ref{thm2.2.5}, \ref{thm2.2.7}, and \ref{thm2.2.8} or
 Theorem \ref{thm2.2.12} by setting $B(x,y)=\tilde B(x)+\hat B(x,y)$ so that
 $\int\hat B(x,y)d\mu^x(y)=0$ where $\mu^x$ is the unique invariant measure
 of $Y_x$ and the vector field $\tilde B$, which becomes now the averaged
 vector field $\bar B$, has an $\omega$-limit set satisfying conditions of 
 the above theorems. Simple examples of this construction will be exhibited
 in Section \ref{sec2.7} for which we also compute historgrams indicating
 proportions of time the slow motion spends near different attracting points
 of the averaged motion. We observe that the functional $S_{0T}$, which plays
 a crucial role in the above theorems, seems to be quite difficult to compute
 since this leads to difficult nonclassical variational problems.

 \section{Large deviations}\label{sec2.3}\setcounter{equation}{0}

We will need the following version of general large deviations bounds when
usual assumptions hold true with errors. The proof is a strightforward 
modification of the standard one (cf. \cite{Ki1}) and its details can be found 
in Part \ref{part1}, Lemma \ref{lem1.4.1}.

\begin{lemma}\label{lem2.3.1} Let $H=H(\be)$, $\eta=\eta(\be)$ be
uniformly bounded on compact sets functions on $\bbR^d$ and 
$\{\Xi_\tau,\,\tau\geq 1\}$ be a family of $\bbR^d-$valued random
vectors on a probability space $(\Om,\cF,P)$ such that $|\Xi_\tau|\leq C
<\infty$ with probability one for some constant $C$ and all $\tau\geq 1$.
For any $a>0$ and $\al,\be_0\in\bbR^d$ set
\begin{equation}\label{2.3.1}
L_a^{\be_0}(\al)=\sup_{\be\in\bbR^d,|\be+\be_0|\leq a}(\langle\be,\al\rangle-
H(\be)),\,\,L_a(\al)=L_a^0(\al),\,\, L(\al)=L_\infty(\al).
\end{equation}
(i) For any $\la,a>0$ there exists $\tau_0=\tau(\la, a,C)$ such that
whenever for some $\tau\ge\tau_0$, $\be_0\in\bbR^d$ and each $\be\in\bbR^d$
with $|\be+\be_0|\leq a$,
\begin{equation}\label{2.3.2}
H_\tau(\be)=\tau^{-1}\log Ee^{\tau \langle\be,\Xi_\tau\rangle}
\leq H(\be)+\eta(\be)
\end{equation}
then for any compact set $\cK\subset\bbR^d$,
\begin{equation}\label{2.3.3}
P\{\Xi_\tau\in \cK\}\leq\exp\left(-\tau (L_a^{\be_0}(\cK)-\eta_a^{\be_0}-
\la|\be_0|-\la)\right)
\end{equation}
where
\begin{equation}\label{2.3.4}
\eta_a^{\be_0}=\sup\{\eta(\be):\,|\be+\be_0|\leq a\}\,\,\mbox{and}\,\,
L_a^{\be_0}(\cK)=\inf_{\al\in \cK}L_a^{\be_0}(\al).
\end{equation}
(ii) Suppose that $\al_0\in\bbR^d$, $0<a\leq\infty$ and there exists 
$\be_0\in\bbR^d$ such that $|\be_0|\leq a$ and
\begin{equation}\label{2.3.5}
H(\be_0)=\langle\be_0,\al_0\rangle -L_a(\al_0).
\end{equation}
If (\ref{2.3.2}) holds true then for any $\del>0$,
\begin{equation}\label{2.3.6}
P\{ |\Xi_\tau-\al_0|\leq\del\}\leq\exp\left( -\tau(L_a(\al_0)-\eta(\be_0)-
\del|\be_0|)\right).
\end{equation}
(iii) Assume that $\al_0,\be_0\in\bbR^d$ satisfy (\ref{2.3.5}). For any
$\la,a>0$ there exists $\tau_0=\tau(\la,a,C)$ such that whenever for some
$\tau\geq\tau_0$ and each $\be\in\bbR^d$ with $|\be|\leq a$ the inequality
(\ref{2.3.2}) holds true together with
\begin{equation}\label{2.3.7}
\tau^{-1}\log Ee^{\tau \langle\be,\Xi_\tau\rangle}\geq H(\be)-\eta(\be)
\end{equation}
then for any $\gam,\del>0,\,\gam\leq\del$,
\begin{eqnarray}\label{2.3.8}
&P\{ |\Xi_\tau-\al_0|<\del\}\geq\exp\left( -\tau(L(\al_0)+\eta(\be_0) +
\gam|\be_0|)\right)\\
&\times\left( 1-\exp\big(-\tau(\tilde L_a^{\be_0}(\cK_{\gam,C}(\al_0))-\eta_a
-\eta(\be_0)-\la|\be_0|-\la)\big)\right)\nonumber
\end{eqnarray}
where 
\[
\tilde L_a^{\be_0}(\al)=L_a(\al)-\langle\be_0,\al\rangle +H(\be_0),
\]
$\tilde L_a^{\be_0}(\cK)=\inf_{\al\in \cK}\tilde L_a^{\be_0}(\al),$ $\eta_a=
\eta_a^0$, $\cK_{\gam,C}(\al_0)=\overline{U_C(0)}\setminus U_\gam(\al_0),$
 $U_\gam(\al)=\{\tilde\al:\,|\tilde\al-\al|<\gam\}$ and $\bar U$ denotes 
 the closure of $U$.
\end{lemma}
The proof of the following result is also standard and can be found in Part
\ref{part1}, Lemma \ref{lem1.4.2}.
\begin{lemma}\label{lem2.3.2} Let $S_n,\, n=1,2,...$ be a nondecreasing sequence
 of lower semicontinuous functions on a metric space $M$ and let 
 $S=\lim_{n\to\infty}S_n.$ Assume that $S$ is also lower semicontinuous and 
 for any compact set $\cK\subset M$ denote
   \[
   S_n(\cK)=\inf_{\gam\in \cK}S_n(\gam)\,\,\mbox{and}\,\, 
   S(\cK)=\inf_{\gam\in \cK}S(\gam).
   \]
   Then
   \begin{equation}\label{2.3.9}
   \lim_{n\to\infty}S_n(\cK)=S(\cK).
   \end{equation}
   \end{lemma}

   We will need also the following general result which 
   will enable us to subdivide time into small intervals freezing the slow
   variable on each of them so that the estimate (\ref{2.2.3}) of 
   Assumption \ref{ass2.2.1} becomes sufficiently precise and, on the other
   hand, we will not change much the corresponding functionals $S_{0T}$
   appearing in required large deviations estimates. This result is
   certainly not new, it is cited in \cite{Ve3} as a folklore fact and
   a version of it can be found in \cite{Kr}, p.67 while for a complete proof
   we refer the reader to Part \ref{part1}, Lemma \ref{lem1.4.3}.
   
   \begin{lemma}\label{lem2.3.3}
   Let $f=f(t)$ be a measurable function on $\bbR^1$ equal zero outside
   of $[0,T]$ and such that $\int_0^T|f(t)|dt<\infty$. For each positive
   integer $m$ and $c\in[0,T]$ define $f_m(t,c)=f([(t+c)\Del^{-1}]\Del-c)$
   where $\Del=T/m$ and $[\cdot]$ denotes the integral part. Then there
   exists a sequence $m_i\to\infty$ such that for Lebesgue almost all
   $c\in[0,T]$,
   \begin{equation}\label{2.3.10}
   \lim_{i\to\infty}\int_0^T|f(t)-f_{m_i}(t,c)|dt=0.
   \end{equation}
   \end{lemma}

Next we will need the following simple estimates whose proof uses the
Gronwall inequality and can be found in Part \ref{part1}, Lemma \ref{lem1.5.1}.

\begin{lemma}\label{lem2.3.4} Let $x_i,\tilde x_i\in\cX,\, i=0,1,...,N,$
$0=t_0<t_1<...<t_{N-1}<t_N=T,$ $\Del=\max_{0\leq i\leq N-1}(t_{i+1}-t_i),$
$\xi_i=(x_i-x_{i-1})(t_i-t_{i-1})^{-1}$, $n(t)=\max\{ j\geq 0:\, t\geq t_j\},$
$\psi(t)=\tilde x_{n(t)}$, $v\in\cX\times\bfM$,
\[
\Xi_j^\ve(v,x)=(t_j-t_{j-1})^{-1}\int_{t_{j-1}}^{t_j}B(x,Y^\ve_v(s/\ve))ds,
\]
and for $t\in[0,T]$,
\begin{equation}\label{2.3.11}
Z_{v,x}^{\ve,\psi}(t)=x+\int_0^tB(\psi(s),Y_v^\ve(s/\ve))ds.
\end{equation}
Then
\begin{eqnarray}\label{2.3.12}
&\big\vert\Xi_j^\ve(v,x_{j-1})-(t_j-t_{j-1})^{-1}(Z^\ve_v(t_j)-
Z^\ve_v(t_{j-1}))\big\vert\\
&\leq K\big\vert Z^\ve_v(t_{j-1})-x_{j-1}\big\vert
+\frac 12K^2(t_j-t_{j-1}),\nonumber
\end{eqnarray}
\begin{eqnarray}\label{2.3.13}
&\sup_{0\leq s\leq t}\big\vert Z_{v,x}^{\ve,\psi}(s)-\psi(s)\big\vert\leq
|x-x_0|+\max_{0\leq j\leq n(t)}|x_j-\tilde x_j|\\
&+K\Del+n(t)\Del\max_{1\leq j\leq n(t)}\big\vert\Xi_j^\ve(v,\tilde x_{j-1})
-\xi_j\big\vert\nonumber
\end{eqnarray}
and
\begin{equation}\label{2.3.14}
\sup_{0\leq s\leq t}\big\vert Z^\ve_v(s)-Z_{v,x}^{\ve,\psi}(s)\big\vert\leq
e^{Kt}\big(|\pi_1v-x|+Kt\sup_{0\leq s\leq t}\big\vert Z_{v,x}^{\ve,\psi}(s)-
\psi(s)\big\vert\big)
\end{equation}
where, recall, $Z^\ve_v(s)=X_v^\ve(s/\ve)$ and $\pi_1v=z\in\cX$ if
$v=(z,y)\in\cX\times\bfM$.
\end{lemma}

 For any $x',x''\in\cX$ and $\be,\xi\in\bbR^d$ set
 \[
 L_b(x',x'',\xi)=\sup_{\be\in\bbR^d,|\be|\leq b}\big(\langle\be,\xi\rangle-
 H(x',x'',\be)\big),
 \]
 and $L_b(x,\xi)=L_b(x,x,\xi)$ with $H(x',x'',\be)$ given by Assumption 
 \ref{ass2.2.1}.
 The following result is the crucial step in the proof of Theorem \ref{thm2.2.2}.
 \begin{proposition}\label{prop2.3.5}
 Let $x_j,t_j,\xi_j,N,\Del,T$ and $\Xi^\ve_j$ be the same as in Lemma 
 \ref{lem2.3.4} and assume that
 \begin{equation}\label{2.3.15}
 \hat\Del=\min_{0\leq i\leq N-1}(t_{i+1}-t_i)\geq\Del/3.
 \end{equation}
 
 (i) There exist $\del_0>0,\ve_0(\Del)>0$ and $C_T(b)>0$ independent of
 $x,y,x_j,\tilde x_j,\xi_j$ such that if $\del\leq\del_0$ and $\ve\leq
 \ve_0(\Del)$ then for any $b>0$,
 \begin{eqnarray}\label{2.3.16}
 &P\big(\max_{1\leq j\leq N}\big\vert\Xi^\ve_j((x,y),\tilde
 x_{j-1})-\xi_j|<\del\big)\\
 &\leq\exp\big\{-\frac 1\ve\big(\sum_{j=1}^N(t_j-t_{j-1})L_b(\tilde x_{j-1},
 \xi_j)-\eta_{b,T}(\ve,\Del)-C_T(b)(d+\del)\big)\big\}\nonumber
 \end{eqnarray}
 where $d=|x-x_0|+\max_{0\leq j\leq N}|x_j-\tilde x_j|$, $\eta_{b,T}(\ve,
 \Del)$ does not depend on $x,x_j,\tilde x_j,\xi_j$ and
 \begin{equation}\label{2.3.17}
 \lim_{\Del\to 0}\limsup_{\ve\to 0}\eta_{b,T}(\ve,\Del)=0.
 \end{equation}
 In particular, if for each $j=1,...,N$ there exists $\be_j\in\bbR^d$ such
 that
 \begin{equation}\label{2.3.18}
 L(\tilde x_j,\xi_j)=\langle\be_j,\xi_j\rangle-H(\tilde x_j,\be_j)
 \end{equation}
 and
 \begin{equation}\label{2.3.19}
 \max_{1\leq j\leq N}|\be_j|\leq b<\infty
 \end{equation}
 then (\ref{2.3.16}) holds true with $L(\tilde x_j,\xi_j)$ in place of 
 $L_b(\tilde x_j,\xi_j)$, $j=1,...,N$.
 
 (ii) For any $b,\la,\del,q>0$ there exist $\Del_0=\Del_0(b,\la,\del,q)>0$
 and $\ve_0=\ve_0(b,\la,\del,q,\Del)$, the latter depending also on $\Del>0$, 
 such that if $\xi_j$ and $\be_j$ satisfy (\ref{2.3.18}) and (\ref{2.3.19}),
 $\max_{1\leq j\leq N}|\xi_j|\leq q$, $\Del<\Del_0$ and $\ve<\ve_0$ then
 \begin{eqnarray}\label{2.3.20}
 &P\big(\max_{1\leq j\leq N}\big\vert\Xi^\ve_j((x,y),\tilde
 x_{j-1})-\xi_j|<\del\big)\\
 &\geq\exp\big\{-\frac 1\ve\big(\sum_{j=1}^N(t_j-t_{j-1})L(\tilde x_{j-1},
 \xi_j)-\eta_{b,T}(\ve,\Del)+C_T(b)d+\la\big)\big\}\nonumber
 \end{eqnarray}
 with some $C_T(b)>0$ depending only on $b$ and $T$.
 \end{proposition}
 \begin{proof} (i) Introduce the events
  \[
   \Gam^j(r)=\big\{ \big\vert\Xi_j(v,\tilde x_{j-1})-\xi_j\big \vert<r\big\},\,
   j=1,...,N
   \]
   so that we have
   \begin{equation}\label{2.3.21}
   P\big\{ \max_{1\leq j\leq n}\big\vert\Xi^\ve_j(v,\tilde x_{j-1})-\xi_j|
   <r\big\}=P\big(\cap_{j=1}^n\Gam^j(r)\big).
   \end{equation}
   Now for $v=(x,y)$ by the Markov property \index{Markov property}
   \begin{eqnarray}\label{2.3.22}
   &P\big(\cap_{j=1}^n\Gam^j(\del)\big)\\
   &=E\,\bbI_{\cap_{j=1}^{n-1}\Gam^j(\del)}P_{X^\ve_{x,y}(t_{n-1}\ve^{-1}),
   Y^\ve_{x,y}(t_{n-1}\ve^{-1})}\big\{\big\vert (t_n-t_{n-1})^{-1}\nonumber\\
   &\times\int_0^{t_n-t_{n-1}}B\big(\tilde x_{n-1},
   Y^\ve_{X^\ve_{x,y}(t_{n-1}\ve^{-1}),Y^\ve_{x,y}(t_{n-1}\ve^{-1})}(s/\ve)\big)
   ds-\xi_n\big\vert <\del\big\}.\nonumber
   \end{eqnarray}
   If $\om\in\cap_{j=1}^{n-1}\Gam^j(\del)$ then $X^\ve_{x,y}(t_{n-1}\ve^{-1},
   \om)=Z^\ve_{x,y}(t_{n-1},\om)$ in view of (\ref{2.3.13}) and (\ref{2.3.14}) 
   satisfies
   \begin{eqnarray}\label{2.3.23}
   &\big\vert X^\ve_{x,y}(t_{n-1}\ve^{-1},\om)-\tilde x_{n-1}\big\vert\leq
   d_{n-1}=(e^{Kt_{n-1}}Kt_{n-1}+1)\\
   &\times(|x-x_0|+\max_{0\leq j\leq n-1}|x_j-\tilde x_j|+K\Del+(n-1)\Del\del).
   \nonumber\end{eqnarray}
   Since $H(x',x'',\be)$ is Lipschitz continuous in $x'$ and $x''$ it follows
   from (\ref{2.3.22}) that
   \begin{equation}\label{2.3.24}
   \big\vert H\big(X^\ve_{x,y}(t_{n-1}\ve^{-1},\om),\tilde x_{n-1},
   \be_n^{(a)}\big)-H(\tilde x_{n-1},\be_n^{(a)})\big\vert\leq C(a)d_{n-1}
   \end{equation}
   provided $\om\in\cap_{j=1}^{n-1}\Gam^j(\del)$ where $C(a)>0$ depends
   only on $a$. In view of Assumption \ref{ass2.2.1} we can estimate from 
   above the probability in the right hand side of (\ref{2.3.22}) by means 
   of Lemma \ref{lem2.3.1}(i) which together with (\ref{2.3.24}) yield that
   \begin{eqnarray}\label{2.3.25}
   &P\big(\cap_{j=1}^n\Gam^j(\del)\big)\leq P\big(\cap_{j=1}^{n-1}
   \Gam^j(\del)\big)\\
    &\times\exp\big(-\frac {(t_n-t_{n-1})}{\ve}
   (L_a(\tilde x_{n-1},\xi_n)-\tilde\eta_{a,T}(\ve,\Del)-C(a)d_{n-1}-ra)\big)
   \nonumber\end{eqnarray}
   where $\tilde\eta_{a,T}(\ve,\Del)\to 0$ as, first, $\ve\to 0$ and then
   $\Del\to 0$. Applying (\ref{2.3.25}) for $n=N,N-1,...,2$ and estimating
   $P(\Gam^1(\del))$ by means of Lemma \ref{lem2.3.1}(i) we derive 
   (\ref{2.3.16}) in view of (\ref{2.3.21}).
   
   (ii) In order to obtain (\ref{2.3.20}) we rely on Assumption \ref{ass2.2.1}
   and Lemma \ref{lem2.3.1}(iii) estimating from below the probability 
   in the right hand side of (\ref{2.3.22}) which together with (\ref{2.3.23}) 
   yield
   \begin{eqnarray}\label{2.3.26}
   &P\big(\cap_{j=1}^n\Gam^j(\del)\big)\geq P\big(\cap_{j=1}^{n-1}
   \Gam^j(\del)\big)\\
    &\times\exp\big(-\frac {(t_n-t_{n-1})}{\ve}L_a(\tilde x_{n-1},\xi_n)\big)
    g_{n,b}(\ve,\del,\vsig,\sig)\nonumber
    \end{eqnarray}
     where
  \begin{eqnarray*}
  &g_{n,b}(\ve,\Del,\vsig,\sig)=\exp\bigg(-\frac {(t_n-t_{n-1})}{\ve}
  \big(\tilde\eta_{b,T}(\ve,\Del)+C_T(b)d_{n-1}+\vsig b\big)\bigg)\\
  &\times\bigg(1-\exp\big(-\frac {(t_n-t_{n-1})}{\ve}(d(b)-\tilde\eta_{b,T}
  (\ve,\Del)-\sig b-\sig)\big)\bigg),
  \end{eqnarray*}
  \[
  d(b)=\min_{1\leq j\leq N}\tilde L_b^{\be_j}(\tilde x_{j-1},\cK_{\vsig,C}
  (\xi_j)),\, \tilde L^\be_b(x,\cK)=\inf_{\al\in \cK}\tilde L_b^\be(x,\al),
  \]
  \[
  \cK_{\vsig,C}(\al)=\bar U_C(0)\setminus U_\vsig(\al),\,
  \tilde L_b^\be(x,\al)=L_b(x,\al)-\langle\be,\al\rangle +H(x,\be),\,
  C_T(b)>0,
  \]
  and $\tilde\eta_{b,T}(\ve,\Del)\to 0$ as, first, $\ve\to 0$ and then 
  $\Del\to 0$. Employing (\ref{2.3.26}) for $n=N,N-1,...,2$ and estimating
   $P(\Gam^1(\del))$ by means of Lemma \ref{lem2.3.1}(iii) we obtain from
   (\ref{2.3.21}) that
   \begin{eqnarray}\label{2.3.27}
   &P\big\{ \max_{1\leq j\leq n}\big\vert\Xi^\ve_j(v,\tilde x_{j-1})-\xi_j|
   <\del\big\}\\
   &\geq\exp\bigg(-\frac 1\ve\big(\sum_{j=1}^N(t_j-t_{j-1})L(\tilde x_{j-1},
 \xi_j)+C(\rho,\del)\ve\Del^{-1}\big)\bigg)\nonumber \\
 &\times\prod_{n=1}^N g_{n,b}(\ve,\Del,\vsig,\sig)\nonumber
 \end{eqnarray}
 for some $C(\rho,\del)>0$ provided, say, $NC_1\ve\Del^{-1}\leq 
 2TC_1\ve\Del^{-2}\leq\frac {\del}2$ and $T\ve\Del^{-1}<C\rho/2$. 
 Since $H(x,\be)$ is differentiable in $\be$ then
 \begin{equation*}
 \tilde L(\tilde x_j,\al)=L(\tilde x_j,\al)-\langle\be_j,\xi_j\rangle +
 H(\tilde x_j,\be_j)>0
 \end{equation*}
 for any $\al\ne\xi_j$ (see Theorems 23.5 and 25.1 in \cite{Roc}), and so
 by the lower semicontinuity of $L(x,\al)$ in $\al$ (and, in fact, also
 in $x$), 
 \[
 \tilde L^{\be_j}(\tilde x_{j-1},\cK_{\vsig,C}(\xi_j))=\inf_{\al\in 
 \cK_{\vsig,C}(\xi_j)}\tilde L^{\be_j}(\tilde x_{j-1},\al)>0.
 \]
 This together with Lemma \ref{lem2.3.2} yield that $d(b)$ appearing in the
 definition of $g_{n,b}(\ve,\Del,\vsig,\sig)$ is positive provided $b$ is
 sufficiently large. In fact, it follows from the lower semicontinuity
 of $L(x,\al)$ that $d(b)$ is bounded away from zero by a positive
 constant independent of $\tilde x_j$ and $\xi_j$, $j=1,...,N$ if these
 points vary over fixed compact sets and (\ref{2.3.18}) together with
 (\ref{2.3.19}) hold true. Now, given $\la>0$ choose, first, sufficiently
 large $b$ as needed and then subsequently choosing small $\sig$ and $\vsig$,
 then small $\Del$, and, finally, small enough $\ve$ we end up with an
 estimate of the form
 \begin{equation}\label{2.3.28}
 g_{n,b}(\ve,\Del,\vsig,\sig)\geq\exp\big(-\frac {(t_n-t_{n-1})}{\ve}
 (\eta_{b,\rho,T}(\ve,T)+C_T(b)d+\la)\big)
 \end{equation}
 where $C_T(b)>0$ and $\eta_{b,\rho,T}(\ve,T)$ satisfies (\ref{2.3.17}).
 Finally, (\ref{2.3.20}) follows from (\ref{2.3.27}) and (\ref{2.3.28}). 
 \end{proof}
 
 The remaining part of the proof of Theorem \ref{thm2.2.2} contains mostly
 some convex analysis arguments and it repeats almost verbatim the corresponding 
 part of the proof of Theorem \ref{thm1.2.3} in Part \ref{part1} but for 
  readers' convenience we exhibit it also here. We remark that some of the 
 details below are borrowed from \cite{Ve3} but we believe that our exposition
 and the way of proof are more precise, complete and easier to follow.
 We start with the lower bound. 
 Assume that $S_{0T}(\gam)<\infty$, and so that $\gam$ is
 absolutely continuous, since there is nothing to prove otherwise. Then by 
 (\ref{2.2.9}), $L(\gam_s,\dot{\gam}_s)<\infty$ for Lebesgue almost all 
 $s\in[0,T].$ By (\ref{2.2.1}) and Assumption \ref{ass2.2.1},
 \begin{equation}\label{2.3.29}
 H(x,\be)\leq\tilde K|\be|
 \end{equation}
 for some $\tilde K>0$,
 and so if $L(\gam_s,\dot{\gam}_s)<\infty$ it follows from (\ref{2.2.4}) that
 $|\dot{\gam}|\leq\tilde K$. Suppose that $\cD(L_s)=\{\al:\, L(\gam_s,\al)<
 \infty\}
 \ne\emptyset$ and let ri$\cD(L_s)$ be the interior of $\cD(L_s)$ in its affine
 hull (see \cite{Roc}). Then either ri$\cD(L_s)\ne\emptyset$ or $\cD(L_s)$ (by
 its convexity) consists of one point and recall that $\dot{\gam}_s\in\cD(L_s)$
 for Lebesgue almost all $s\in[0,T]$.  By (\ref{2.2.6}) and (\ref{2.3.29}),
\begin{equation}\label{2.3.30}
0=H(\gam_s,0)=\inf_{\al\in\bbR^d}L(\gam_s,\al).
\end{equation}
This together with the nonnegativity and lower semi-continuity of 
$L(\gam_s,\cdot)$ yield that there exists $\hat\al_s$ such that 
$L(\gam_s,\hat\al_s)=0$ and by a version of the measurable selection (of the
implicit function) theorem (see \cite{CV}, Theorem III.38), $\hat\al_s$ can be
chosen to depend measurably in $s\in[0,T]$. Of course, if ri$\cD(L_s)=\emptyset$
then $\cD(L_s)$ contains only $\hat\al_s$ and in this case $\hat\al_s=
\dot{\gam}_s$ for Lebesgue almost all $s\in[0,T]$. Taking $\al_s=\hat\al_s$ 
and $\be_s=0$ we obtain
\begin{equation}\label{2.3.31}
L(\gam_s,\al_s)=\langle\be_s,\al_s\rangle-H(\gam_s,\be_s).
\end{equation}

Observe that $\ell(s,\al)=L(\gam_s,\al)$ is measurable as a function of $s$
and $\al$ since it is obtained via (\ref{2.2.4}) as a supremum in one argument 
of a family of continuous functions, and so this supremum can be taken there
over a countable dense set of $\be$'s. Hence, the set $A=\{ (s,\al):\,
s\in[0,T],\,\al\in\cD(L_s)\}=\ell^{-1}[0,\infty)$ is measurable, and so the set
$B=A\setminus\{(s,\dot{\gam_s}),\, s\in[0,T]\}$ is measurable, as well. Its
projection $V=\{ s\in[0,T]:\, (s,\al)\in B\,\mbox{for some}\,\al\in\bbR^d\}$
on the first component of the product space is also measurable and $V$ is the
set of $s\in[0,T]$ such that $\cD(L_s)$ contains more than one point. Employing
Theorem III.22 from \cite{CV} we select $\bar\al_s\in\bbR^d$ measurably in
$s\in V$ and such that $(s,\bar\al_s)\in B$. By convexity and lower 
semicontinuity of $L(\gam_s,\cdot)$ it follows from Corollary 7.5.1 in
\cite{Roc} that 
\begin{equation}\label{2.3.32}
L(\gam_s,\dot{\gam}_s)=\lim_{p\uparrow\infty}L(\gam_s,\al^{(p)}_s)\,\,
\,\mbox{where}\,\,\,\al^{(p)}_s=(1-p^{-1})\dot{\gam}_s+p^{-1}\bar\al_s.
\end{equation}
For each $\del>0$ set 
\[
n_\del(s)=\min\{ n\in\bbN:\, |L(\gam_s,\dot{\gam}_s)-
L(\gam_s,\al^{(n)}_s)|+|\dot{\gam}_s-\al^{(n)}_s|<\del\}.
\]
 Then, clearly, $n_\del(s)$ is a measurable function of $s$, and so 
 $\al_s=\al_s^{(\del)}=\al^{(n_\del(s))}_s$ and $L(\gam_s,\al_s)$ are
  measurable in $s$, as well. By Theorems 23.4 and 23.5 from \cite{Roc}
  for each $\al_s=\al_s^{(\del)}$ there exists $\be_s=\be_s^{(\del)}\in\bbR^d$
  such that (\ref{2.3.31}) holds true. Given $\del',\la>0$ take $\del=\min(\del',
  \la/3)$ and for $s\in[0,T]\setminus V$ set $\al_s=\hat\al_s$. Then
\begin{equation}\label{2.3.33}
\int_0^T\big\vert L(\gam_s,\dot{\gam}_s)-L(\gam_s,\al_s)\big\vert ds<\la/3\,
\mbox{and}\,\int_0^T|\dot{\gam}_s-\al_s|ds<\del'.
\end{equation}

For each $b>0$ set $\al_s^b=\al_s$ if the corresponding $\be_s$ in (\ref{2.3.31})
satisfies $|\be_s|\leq b$ and $\al_s^b=\hat\al_s$, otherwise. Note, that 
(\ref{2.3.31}) remains true with $\al_s^b$ in place of $\al_s$ with $\be_s=0$ if 
$\al_s^b=\hat\al_s$. As observed above $|\al|\leq K$ whenever $L(z,\al)<\infty$,
and so $|\hat\al_s|\leq K$ for Lebesgue almost all $s\in[0,T]$. We recall also
that $|\dot{\gam}_s-\al_s|<\del$ and $\dot{\gam}_s\leq K$ for Lebesgue
almost all $s\in[0,T]$. Since $S_{0T}(\gam)<\infty$, $|L(\gam_s,\dot{\gam}_s)-
L(\gam_s,\al_s)|<\del$, and $L(\gam_s,\al_s^b)\uparrow L(\gam_s,\al_s)$ as
$b\uparrow\infty$ for Lebesgue almost all $s\in[0,T]$, we conclude from 
(\ref{2.3.33}) and the above observations that for $b$ large enough
\begin{equation}\label{2.3.34}
\int_0^T\big\vert L(\gam_s,\al_s)-L(\gam_s,\al_s^b)\big\vert ds<\la/3\,\,\,
\mbox{and}\,\,\int_0^T|\al_s-\al_s^b|ds<\del'.
\end{equation}
Next, we apply Lemma \ref{lem2.3.3} to conclude that there exists a sequence
$m_j\to\infty$ such that for each $\Del_j=T/m_j$ and Lebesgue almost all 
$c\in[0,T)$,
\begin{equation}\label{2.3.35}
\int_0^T\big\vert L(\gam_s,\al_s^b)-L(\gam_{q_j(s,c)},\al_{q_j(s,c)}^b)\big
\vert ds<\la/3\,\mbox{and}\,\int_0^T|\al_s^b-\al_{q_j(s,c)}^b|ds<\del'.
\end{equation}
where $q_j(s,c)=[(s+c)\Del^{-1}_j]\Del_j-c$, $[\cdot]$ denotes the integral part
and we assume $L(\gam_s,\al_s^b)=0$ and $\al_s^b=0$ if $s<0$.

Choose $c=c_j\in[\frac 13\Del_j,\frac 23\del_j]$ and set $\hat\gam_s=x+\int_0^s
\al^b_{q_j(u,c)}du$, $\psi_s=\gam_{q_j(s,c)}$ where $\gam_u=\gam_0$ if $u<0$,
$x_0=\tilde x_0=x$, $x_N=\hat\gam_T,$ $\tilde x_N=\gam_T$ and $x_k=\hat\gam_
{k\Del_j-c}$, $\tilde x_k=\gam_{k\Del_j-c}$ for $k=1,...,N-1$ and $\xi_k=
\al^b_{(k-1)\Del_j-c}$ for $k=1,2,...,N$ where $N=\min\{ k:k\Del_j-c>T\}$.
Since $|\dot{\gam}_s|\leq\tilde K$ for Lebesgue almost all $s\in[0,T]$ then
$\bfr_{0T}(\gam,\psi)\leq\tilde K\Del_j$ and, in addition, $\bfr_{0T}(\gam,
\hat\gam)
\leq 3\del'$ by (\ref{2.3.33})--(\ref{2.3.35}). This together with (\ref{2.3.13}) and
(\ref{2.3.14}) yield that for $v=(x,y)$,
\begin{eqnarray}\label{2.3.36}
&\,\,\,\,\,\,\,\bfr_{0T}(Z^\ve_v,\gam)\leq\bfr_{0T}(Z^\ve_v,\psi)+
\tilde K\Del_j\leq(KTe^{KT}+1)
\bfr_{0T}(Z^{\ve,\psi}_v,\psi)+\tilde K\Del_j\\
&\leq(KTe^{KT}+1)\big(3\del'+\tilde K\Del_j+(T+1)\max_{1\leq k\leq N}\big\vert
\Xi^\ve_k(v,\tilde x_{k-1})-\xi_k\big\vert\big)+\tilde K\Del_j \nonumber
\end{eqnarray}
provided $\Del_j\leq 1$ where $Z_v^{\ve,\psi}$ and $\Xi^\ve_k(v,x)$ are the same
as in Lemma \ref{lem2.3.4}, the latter is defined with $t_k=k\Del_j-c$, $k=1,...,N-1$
and $t_N=T$. Choose $\del'$ so small and $m_j$ so large that
\[
(KTe^{KT}+1)\big(3\del'+\tilde K\Del_j+(T+1)\del'\big)+\tilde K\Del_j<\del
\]
then by (\ref{2.3.36}),
\begin{equation}\label{2.3.37}
\big\{ \bfr_{0T}(Z^\ve_{x,y},\gam)<\del\big\}\supset\big\{ 
\max_{1\leq k\leq N}\big\vert\Xi^\ve_k(v,\tilde x_{k-1})-\xi_k\big\vert 
<\del'\big\}.
\end{equation}
By (\ref{2.3.33})--(\ref{2.3.35}),
\begin{equation}\label{2.3.38}
\sum_{k=1}^N(t_k-t_{k-1})L(\tilde x_{k-1},\xi_k)\leq S_{0T}(\gam)+\la
\end{equation}
and by the construction above the conditions of the assertion (ii) of 
Proposition \ref{prop2.3.5} hold true, so choosing $m_j$ sufficiently large we 
derive (\ref{2.2.6}) (with $2\la$ in place of $\la$) from (\ref{2.3.20}), 
(\ref{2.3.37}) and (\ref{2.3.38}) provided $\ve$ is small enough.

Next, we pass to the proof of the upper bound (\ref{2.2.7}). Assume that 
(\ref{2.2.7}) is not true, i.e. there exist $a,\la,\del>0$ and $x\in\cX_T$ such
that for some sequence $\ve_k\to 0$ as $k\to\infty$,
\begin{equation}\label{2.3.39}
P\big\{ \bfr_{0T}\big( Z^{\ve_k}_{x,y},\Psi^a_{0T}(x)\big)\geq 
3\del\big\}>\exp\big(-\frac 1{\ve_k}(a-\la)\big).
\end{equation}
Since $| B(x,y)|\leq K$ by (\ref{2.2.1}) all paths of $Z^{\ve}_{x,y}(t),\, 
t\in[0,T]$ and of $Z^{\ve,\psi}_{v,x}(t),\, t\in[0,T]$ given by (\ref{2.3.11}) 
(the latter for any measurable $\psi$) belong to a compact set 
$\tilde \cK^x\subset C_{0T}$ which consists of curves starting at $x$ and 
satisfying the Lipschitz condition with the constant
$K$. Let $\tilde U^x_{\rho}$ denotes the open $\rho$-neighborhood of the compact set
$\Psi^a_{0T}(x)$ and $\cK^x_\rho=\tilde \cK^x\setminus\tilde U^x_{\rho}$. 
For any small
$\del'>0$ choose a $\del'$-net $\gam_1,...,\gam_n$ in $\cK^x_{2\del}$ where 
$n=n(\del')$. Since
\begin{equation*}
\big\{ \bfr_{0T}\big( Z^{\ve_k}_{x,y},\Psi^a_{0T}(x)\big)\geq 3\del
\big\}\subset\bigcup_{n\geq j\geq 1}\big\{ \bfr_{0T}\big( 
Z^{\ve_k}_{x,y},\gam_j)\leq\del'\big\}
\end{equation*}
then there exists $j$ and a subsequence of $\{\ve_k\}$, for which we use the 
same notation, such that
\begin{equation}\label{2.3.40}
P\big\{ \bfr_{0T}( Z^{\ve_k}_{x,y},\gam_j)\leq\del'\big\}>
n^{-1}\exp\big(-\frac 1{\ve_k}(a-\la)\big).
\end{equation}
Denote such $\gam_j$ by $\gam^{\del'}$, choose a sequence $\del_l\to 0$ and set
$\gam^{(l)}=\gam^{\del_l}$. Since $\cK^x_{2\del}$ is compact there exists a 
subsequence $\gam^{(l_j)}$ converging in $C_{0T}$ to $\hat\gam\in \cK^x_{2\del}$ 
which together with (\ref{2.3.40}) yield
\begin{equation}\label{2.3.41}
\limsup_{\ve\to 0}\ve\ln P\big\{\bfr_{0T}(Z^{\ve}_{x,y},\hat\gam)\leq
\del'\big\}>-a+\la
\end{equation}
for all $\del'>0$.

We claim that (\ref{2.3.41}) contradicts (\ref{2.3.12}) and the assertion (i) of
Proposition \ref{prop2.3.5}. Indeed, set
\[
S^\psi_{b,0T}(\gam)=\int_0^TL_b(\psi(s),\dot{\gam}(s))ds\,\,\mbox{and}\,\,
S_{b,0T}(\gam)=S^\gam_{b,0T}(\gam).
\]
By the monotone convergence theorem
\begin{equation}\label{2.3.42}
S^\psi_{b,0T}(\gam)\uparrow S^\psi_{0T}(\gam)\,\,\mbox{and}\,\,
S_{b,0T}(\gam)\uparrow S_{0T}(\gam)\,\,\mbox{as}\,\, b\uparrow\infty.
\end{equation}
Similarly to our remark in Section \ref{sec2.2}
it follows from the results of Section 9.1 of \cite{IT} that the functionals
$S^\psi_{b,0T}(\gam), S^\psi_{0T}(\gam)$ and $S_{b,0T}(\gam),S_{0T}(\gam)$ are
lower semicontinuous in $\psi$ and $\gam$ (see also Section 7.5 in \cite{FW}).
This together with (\ref{2.3.42}) enable us to apply Lemma \ref{lem2.3.2} in
 order to conclude that
\begin{equation}\label{2.3.43}
\lim_{b\to\infty}S_{b,0T}(\cK_\del^x)=S_{0T}(\cK_\del^x)=\inf_{\gam\in 
\cK^x_\del}S_{0T}(\gam)>a
\end{equation}
where $S_{b,0T}(\cK_\del^x)=\inf_{\gam\in \cK^x_\del}S_{b,0T}(\gam)$. The last
inequality in (\ref{2.3.43}) follows from the lower semicontinuity of $S_{0T}$.
Thus we can and do choose $b>0$ such that
\begin{equation}\label{2.3.44}
S_{b,0T}(\cK_\del^x)>a-\la/8.
\end{equation}

By the lower semicontinuity of $S_{b,0T}^\psi(\gam)$ in $\psi$ there exists a
function $\del(\gam)>0$ on $\cK^x_\del$ such that for each $\gam\in \cK^x_\del$,
\begin{equation}\label{2.3.45}
S_{b,0T}^\psi(\gam)>a-\la/4\,\,\mbox{provided}\,\,\bfr_{0T}(\gam,\psi)<\del_\la
(\gam).
\end{equation}
Next, we restrict the set of functions $\psi$ to make it compact. Namely, we 
allow from now on only functions $\psi$ for which there exists $\gam\in 
\cK^x_\del$ such that either $\psi\equiv\gam$ or $\psi(t)=\gam(kT/m)$ for
 $t\in[kT/m,(k+1)T/m)$, $k=0,1,...,m-1$ and $\psi(T)=\gam(T)$ where $m$ is
 a positive integer. It is easy to see that the set of such functions $\psi$ 
 is compact with respect to the uniform convergence topology in $C_{0T}$ and 
 it follows that $\del_\la(\gam)$ in (\ref{2.3.45}) constructed with such 
 $\psi$ in mind is lower semicontinuous in $\gam$. Hence
\begin{equation}\label{2.3.46}
\del_\la=\inf_{\gam\in \cK^x_\del}\del_\la(\gam)>0.
\end{equation}

Now take $\hat\gam$ satisfying (\ref{2.3.41}) and for any integer $m\geq 1$ set 
$\Del=\Del_m=T/m$, $x_k=x_k^{(m)}=\hat\gam(k\Del)$, $k=0,1,...,m$ and 
$\xi_k=\xi_k^{(m)}=\Del^{-1}\big(\hat\gam(k\Del)-\hat\gam((k-1)\Del)\big),$ 
$k=1,...,m$. Define a piecewise linear 
$\chi_m$ and a piecewise constant $\psi_m$ by
\begin{equation}\label{2.3.47}
\chi_m(t)=x_k+\xi_k\Del\,\,\mbox{and}\,\,\psi_k(t)=x_k\,\,\mbox{for}\,\, 
t\in[k\Del,(k+1)\Del)
\end{equation}
and $k=0,1,...,m-1$ with $\chi_m(T)=\psi_m(T)=\hat\gam(T)$. Since $\hat\gam$ is
 Lipschitz continuous with the constant $\tilde K$ then
\begin{equation}\label{2.3.48}
\bfr_{0T}(\chi_m,\psi_m)\leq\tilde K\Del\,\,\mbox{and}\,\,\bfr_{0T}
(\hat\gam,\psi_m)\leq\tilde K\Del.
\end{equation}
If $m$ is large enough and $\ve>0$ is sufficiently small then
\begin{equation}\label{2.3.49}
\Del<\tilde K^{-1}\min(\del/2,\del_\la)\,\,\mbox{and}\,\,\eta_{b,T}(\ve,\Del)
<\la/8
\end{equation}
where $\eta_{b,T}(\ve,\Del)$ is the same as in (\ref{2.3.16}). Since 
$\hat\gam\in\cK^x_{2\del}$
it follows from (\ref{2.3.48}) and (\ref{2.3.49}) that $\chi_m\in \cK^x_{\del}$ and
 by (\ref{2.3.45})
and the first inequality in (\ref{2.3.49}) we obtain that
\begin{equation}\label{2.3.50}
S^{\psi_m}_{b,0T}(\chi_m)=\Del\sum_{k=0}^{m-1}L_b(x_k,\xi_k)>a-\frac \la{4}.
\end{equation}
Hence, by (\ref{2.3.16}) and the second inequality in (\ref{2.3.49}) for all $\ve$ small enough,
\begin{equation}\label{2.3.51}
 P\big\{ \max_{1\leq k\leq m}\big\vert\Xi^\ve_k((x,y),
 x_{k-1})-\xi_k|<\rho\big\}\leq e^{-\frac 1\ve(a-\la/2)}
 \end{equation}
provided $C_T(b)\rho<\la/8$ (taking into account that $x_0=x$). By (\ref{2.3.12}) 
and the definition of vectors $\xi_k$ for any $v\in\cW$,
\begin{eqnarray}\label{2.3.52}
&\big\vert\Xi^\ve_k(v,x_{k-1})-\xi_k\big\vert\leq\big\vert\Xi^\ve_k(v,x_{k-1})
-\Del^{-1}
\big(Z_v^\ve(k\Del)-Z_v^\ve((k-1)\Del)\big)\big\vert\\
&+2\Del^{-1}\bfr_{0T}(Z^\ve_v,\hat\gam)\leq(K+2\Del^{-1})\bfr_{0T}(Z^\ve_v,
\hat\gam)+\frac 12\tilde K^2\Del.\nonumber
\end{eqnarray}
Therefore,
\begin{eqnarray}\label{2.3.53}
&\big\{ \bfr_{0T}(Z^\ve_{x,y},\hat\gam)\leq\del'\big\}\\
&\subset\big\{\max_{1\leq k\leq m}\big\vert\Xi^\ve_k((x,y), x_{k-1})-\xi_k|
\leq (\tilde K+2\Del^{-1})\del'+\frac 12\tilde K^2\Del\big\}.\nonumber
\end{eqnarray}
Choosing, first, $m$ large enough so that $\Del$ satisfies (\ref{2.3.49}) with
 all sufficiently
small $\ve$ and also that $8C_T(b)\tilde K^2\Del<\la$, and then choosing 
$\del'$ so small that $16C_T(b)(\tilde K+2\Del^{-1})\del'<\del$, we conclude 
that (\ref{2.3.51}) together with (\ref{2.3.53})
contradicts (\ref{2.3.41}), and so the upper bound (\ref{2.2.7}) holds true,
 completing the proof of Theorem \ref{thm2.2.2}.\qed
 
 \begin{remark}\label{rem2.3.6} In view of examples from \cite{BK1} in the 
 fully coupled setup we should not expect convergence (\ref{2.2.22}) in the
  averaging principle with probability one in spite of exponentially
 fast convergence in probability (\ref{2.2.22+}) provided by the upper large
 deviations bound (\ref{2.2.11}). Still, when derivatives of $X^\ve$ and $Y^\ve$
 in $\ve$ grow not too fast we can derive convergence with probability one 
 from (\ref{2.2.22+}). Indeed, consider, for instance, the following example
 \begin{eqnarray}\label{2.3.54}
 &X^\ve(t)=x+\ve\int_0^tB(X^\ve(s),Y^\ve(s))ds\,\,\mbox{and}\\
 &Y^\ve(t)=y+cw_t+\int_0^tb(X^\ve(s))ds\,\,(\mbox{mod}\, 1)\nonumber
 \end{eqnarray}
 where $c\ne 0$ is a constant, $w_t$ is the standard one dimensional 
 Brownian motion, $B(x,y)$ satisfies (\ref{2.2.1}) and it is 
 $1$-periodic in $y$ and $b$ has a bounded derivative in $x$. Set 
 \[
 x^\ve(t)=\frac {dX^\ve(t)}{d\ve}\,\,\,\mbox{and}\,\,\, 
 y^\ve(t)=\frac {dY^\ve(t)}{d\ve}.
 \]
 Then
 \begin{eqnarray*}
&\frac {d}{dt}\left(\begin{matrix}x^\ve(t)\\ y^\ve(t)\end{matrix}\right)
=\left(\begin{matrix}
\ve\frac {\partial B\big(X^\ve(t),Y^\ve(t)\big)}{\partial x}&
\ve\frac {\partial B\big(X^\ve(t),Y^\ve(t)\big)}{\partial y}\\
\frac {\partial b\big(X^\ve(t)\big)}{\partial x}&0\end{matrix}\right)
\left(\begin{matrix}x^\ve(t)\\ y^\ve(t)\end{matrix}\right)\\
&+\left(\begin{matrix}B\big(X^\ve(t),Y^\ve(t)\big)\\ 0\end{matrix}\right).
\end{eqnarray*}
 The solution of this linear equation is easy to estimate which yields
 that for some constant $C>0$,
 \begin{equation}\label{2.3.55}
 \sup_{0\leq t\leq T/\ve}\big\vert \frac {dX^\ve(t)}{d\ve}\big\vert\leq 
 C\exp(C/\sqrt {\ve}).
 \end{equation}
 Let $\mu^x$ be the invariant measure of the diffusion $Y_x(t)=y+w_t+tb(x)$
 (mod 1) (which is unique since the Doeblin condition is satisfied here) and
 assume that
\begin{equation}\label{2.3.56}
\int B(x,y)d\mu^x(y)=0\quad\mbox{for all}\,\, x
\end{equation}
which does not harm the generality since we always can consider
$B(x,y)-\int B(x,y)d\mu^x(y)$ in place of $B(x,y)$. 
Set $\ve_k=\al(\del)/2\ln k$ where $\al(\del)$ is the same as in (\ref{2.2.22+})
written for our specific situation. Then $e^{-\al(\del)/\ve_k}=k^{-2}$ and by 
the Borel--Cantelli lemma we obtain that there exists $k_\del(\om)$, finite with
probability one, so that for all $k\geq k_\del(\om)$,
\begin{equation*}
\max_{0\leq t\leq T/\ve_k}|X^{\ve_k}_{x,y}(t)-x|<\del.
\end{equation*}
By (\ref{2.3.55}) for $\ve_{k+1}<\ve\leq\ve_k$ and $k\geq 2$,
\begin{eqnarray*}
&\max_{0\leq t\leq T/\ve_{k+1}}|X^{\ve_k}_{x,y}(t)-X^\ve_{x,y}(t)|\leq
Ce^{C/\sqrt {\ve_k}}(\ve_k-\ve_{k+1})\\
&\leq C\exp\big(C\sqrt {2(\ln k)/\al(\del)}\big)\ln(1+\frac 1k)(\ln k)^{-2}
\longrightarrow 0\,\,\mbox{as}\,\, k\to\infty.\nonumber
\end{eqnarray*}
It follows that with probability one,
\begin{equation*}
\max_{0\leq n\leq T/\ve}|X^\ve_{x,y}(t)-x|\to 0\,\,\mbox{as}\,\, k\to\infty
\end{equation*}
 which is what we need since in our case  $X^\ve_{x,y}(t)\equiv x$ 
 in view of (\ref{2.3.56}).  
 \end{remark}

\section{Verifying assumptions for random evolutions}\label{sec2.3+}
\setcounter{equation}{0}

In this section we will prove Proposition \ref{prop2.2.4}. Observe that
$H(x,x',\be)$ obtained by (\ref{2.2.18}) is the principal eigenvalue of
the operator $\cL^x+\langle\be,B(x',\cdot)\rangle$ acting on $C^2$ vector
functions $f=(f_1,...,f_N)$ on the manifold $M$ by the formula (see \cite{Ki0}),
\[
\big((\cL^x+\langle\be,B(x',\cdot)\rangle)f\big)_k=\cL^x_kf_k+\langle\be,
B_k(x',\cdot)\rangle f_k
\]
where $x,x'$ and $\be$ are considered as parameters. According to \cite{PW}
this operator satisfies the strong maximum principle. Thus, the first part 
of Proposition \ref{prop2.2.4} follows from the well known results on operators
satisfying the maximum principle (see \cite{DV1}, \cite{DV2} and \cite{Ki0})
and the results on the principle eigenvalue of positive operators (see
\cite{Kra}, \cite{KLS} and \cite{HH}) and of its smooth dependence on
parameters which can be derived from the general perturbation theory of 
linear operators (see \cite{Ka}).

Now we obtain from (\ref{2.2.18}) that for $k=1,...,N$ uniformly in 
$z,x'\in\bar\cX,$ $v\in M$ and $|\be|\leq b$,
\begin{equation}\label{2.3+.1}
\big\vert\frac 1s\log E\exp\langle\be,\int_0^sB_{\nu_{z,v,k}(u)}(x',\hat
Y_{z,v,k}(u))du\rangle-H(z,x',\be)\big\vert\leq\rho_b(s)
\end{equation}
where $\rho_b(s)\to 0$ as $s\to\infty$. Next, we want to compare
\begin{eqnarray*}
&Q_{z,v,k}(s)=E\exp\langle\be,\int_0^sB_{\nu_{z,v,k}(u)}(x',
\hat Y_{z,v,k}(u))du\rangle\,\,\mbox{and}\\
& Q^\ve_{z,v,k}(s)=E\exp\langle\be,\int_0^sB_{\nu^\ve_
{z,v,k}(u)}(x',\hat Y^\ve_{z,v,k}(u))du\rangle.
\end{eqnarray*}
In order to do this we introduce auxiliary random evolutions $W_{z,v,k}(s)=
(\hat W_{z,v,k}(s),\eta_k(s))$ and $W^\ve_{z,v,k}(s)=
(\hat W^\ve_{z,v,k}(s),\eta_k(s))$ governed by the stochastic differential
equations
\begin{equation}\label{2.3+.2}
    d\hat W_{z,v,k}(s)=\sig_{\eta_k(s)}\big(z,\hat W_{z,v,k}(s)\big)dw_s
    +b_{\eta_k(s)}\big(z,\hat W_{z,v,k}(s)\big)ds
\end{equation}
and
 \begin{equation}\label{2.3+.3}
    d\hat W^\ve_{z,v,k}(s)=\sig_{\eta_k(s)}\big(\tilde X^\ve_{z,v,k}(s),
    \hat W^\ve_{z,v,k}(s)\big)dw_s +b_{\eta_k(s)}\big(\tilde X^\ve_{z,v,k}(s),
    \hat W^\ve_{z,v,k}(s)\big)ds,
\end{equation} 
respectively, where $\hat W_{z,v,k}(0)=\hat W^\ve_{z,v,k}(0)=v$,
 \begin{equation}\label{2.3+.4} 
  \frac {d\tilde X^\ve(s)}{dt}=\ve B\big(\tilde X^\ve(s), W^\ve(s)),\,\,\,
  \tilde X^\ve_{z,v,k}(0)=z
 \end{equation} 
 and for $i\ne j$,
\begin{equation}\label{2.3+.5}  
P\big\{\eta_k(s+\Del)=j\big\vert\eta_k(s)=i\big\}=\Del +o(\Del),\,\,\,
\eta_k(0)=k.
 \end{equation} 
 According to \cite{EF} (which relies on Theorem 2 in \S 6, Ch. VII of 
 \cite{GS}) the distributions in the path space of the processes $Y_{z,v,k}$
 and $Y^\ve_{z,v,k}$ are absolutely continuous with respect to the 
 distributions in the path space of the processes $W_{z,v,k}$ and 
 $W^\ve_{z,v,k}$, respectively, with the densities
 \begin{eqnarray}\label{2.3+.6}
 &p_s(\hat W_{z,v,k}(\cdot),\eta)=\prod_{i=0}^{n(s)-1}q_{\eta_k(\zeta_i)
 \eta_k(\zeta_{i+1})}(z,\hat W_{z,v,k}(\zeta_i))\\
 &\times\exp\big(-\sum_{i=0}^{n(s)}\int_{\zeta_i}^{\zeta_{i+1}\wedge s}
 (q_{\eta_k(\zeta_i)}(z,\hat W_{z,v,k}(u))-N+1)du\big)\nonumber
 \end{eqnarray}
 and
 \begin{eqnarray}\label{2.3+.7}
 &\,\,\,\,\,\,p^\ve_s(\tilde X^\ve_{z,v,k}(\cdot),\hat W^\ve_{z,v,k}(\cdot),
 \eta)=\prod_{i=0}^{n(s)-1}q_{\eta_k(\zeta_i)\eta_k(\zeta_{i+1})}
 (\tilde X^\ve_{z,v,k}(\zeta_i),\hat W^\ve_{z,v,k}(\zeta_i))\\
 &\exp\big(-\sum_{i=0}^{n(s)}\int_{\zeta_i}^{\zeta_{i+1}\wedge s}
 (q_{\eta_k(\zeta_i)}(\tilde X^\ve_{z,v,k}(u),\hat W^\ve_{z,v,k}(u))-N+1)du
 \big),\nonumber
 \end{eqnarray}
 respectively, where
 \[
 q_k(z,y)=\sum_{l=1,l\ne k}^Nq_{kl}(z,y),
 \]
 $\zeta_0=0,\,\zeta_{i+1}=\inf\{u>\zeta_i:\,\eta_k(u)\ne\eta_k(\zeta_i)\}$ and
 $n(s)=\max\{ i:\eta_i\leq s\}$.
 
 Thus, we have to compare
 \begin{equation}\label{2.3+.8}
 Q_{z,v,k}(s)=Ep_s(\hat W_{z,v,k}(\cdot),\eta)\exp\langle\be,\int_0^s
 B_{\eta_k(u)}(x',\hat W_{z,v,k}(u))du\rangle
 \end{equation}
 and
  \begin{equation}\label{2.3+.9}
  Q^\ve_{z,v,k}(s)=Ep^\ve_s(\tilde X^\ve_{z,v,k}(\cdot),\hat W_{z,v,k}(\cdot),
  \eta)\exp\langle\be,\int_0^s B_{\eta_k(u)}(x',\hat W^\ve_{z,v,k}(u))du\rangle.
 \end{equation}
 Observe that by (\ref{2.2.1}),
 \begin{equation}\label{2.3+.10}
 \sup_{0\leq u\leq s}|\tilde X^\ve_{z,v,k}(u)-z|\leq\ve sK.
 \end{equation}
 Let $K_q$ be both an upper bound for $|q_{ij}(x,y)|$ and their Lipschitz
 constant then we see from (\ref{2.2.1}) and (\ref{2.3+.6})--(\ref{2.3+.10})
 that
 \begin{eqnarray}\label{2.3+.11}
 &|Q_{z,v,k}(s)-Q^\ve_{z,v,k}(s)|\leq e^{(N+K_q+K|\be|)s}E\bigg( 
 n(s)K^{n(s)}_q\bigg(\ve sKK_q\\
 &+K_q\sup_{0\leq u\leq s}\mbox{dist}(\hat W_{z,v,k}(s),\hat 
 W^\ve_{z,v,k}(s))+2\big\vert\exp\big(\ve sK(K_q+|\be|)\nonumber\\
 &+s(K_q+K|\be|)\sup_{0\leq u\leq s}
 \mbox{dist}(\hat W_{z,v,k}(s),\hat W^\ve_{z,v,k}(s))\big)-1\big\vert\bigg)
 \bigg). \nonumber\end{eqnarray}
 Employing the Witney theorem embed smoothly $M$ as a compact submanifold 
 in an Euclidean space $\bbR^D$ of a sufficiently high dimension $D$ and
 extend the operator $\cL^x+\langle\be,B(x',\cdot)\rangle$ from $M$ to
 $\bbR^d$ so that its coefficients remain $C^2$ and they vanish outside a
  relatively compact set containing $M$ (cf. \cite{Hs}). Now we can view
  (\ref{2.3+.2}) and (\ref{2.3+.3}) as stochastic differential equations in
  $\bbR^d$ keeping the same notations for their coefficients and processes
  there. Then using standard martingale moment estimates for stochastic 
  integrals (see, 
  for instance, \cite{IW}) together with (\ref{2.3+.10}) and the Lipschitz
  continuity of coefficients in (\ref{2.3+.2}) and (\ref{2.3+.3}) we obtain
  \begin{eqnarray*}
  &E\sup_{0\leq u\leq s}|\hat  W_{z,v,k}(u)-\hat  W^\ve_{z,v,k}(u)|^2\\
  &\leq C_1(1+s)\big(\ve^2 s^2K^2+\int_0^sE\sup_{0\leq r\leq u}|\hat  W_{z,v,k}
  (r)-\hat  W^\ve_{z,v,k}(r)|^2du\big)
  \end{eqnarray*}
  for some $C_1>0$ independent of $t,x,v,k$ and $\ve$. Hence, by the
  Gronwall inequality
  \begin{eqnarray}\label{2.3+.12}
  &E\sup_{0\leq u\leq s}\big(\mbox{dist}(\hat W_{z,v,k}(u),
  \hat W^\ve_{z,v,k}(u))\big)^2\\
  &=E\sup_{0\leq u\leq s}|\hat  W_{z,v,k}(u)-\hat  W^\ve_{z,v,k}(u)|^2\leq
  C_1(1+s)\ve^2 s^2K^2e^{C_1(1+s)s}.\nonumber
  \end{eqnarray}
  Observe also that the distribution of $n(s)$ can be written explicitly as
  (see \S 55 in \cite{Gn}),
  \begin{equation}\label{2.3+.13}
  P\{ n(s)=k\}=e^{-(N-1)s}\frac {((N-1)s)^k}{k!}.
  \end{equation}
  
  In order to estimate the last expression in the right hand side of 
  (\ref{2.3+.11}) we note that for any random variable $\xi$,
  \[
  |e^\xi -1|\leq 2|\xi|+(1+e^\xi)\bbI_{|\xi|>1},
  \]
  and so by the Cauchy--Schwarz and the Chebyshev's inequalities
  \begin{equation}\label{2.3+.14}
  E(e^\xi-1)^2\leq 4E\xi^2+2\big(E(1+e^\xi)^4\big)^{1/2}(E\xi^2)^{1/2}.
  \end{equation}
  Now by (\ref{2.3+.11})--(\ref{2.3+.14}) together with the Cauchy--Schwarz
  inequality we obtain that for $k=1,...,N$ uniformly in $x,x'\in\bar\cX$,
   $v\in M$ and $|\be|\leq b$,
  \begin{equation}\label{2.3+.15}
   |Q_{z,v,k}(s)-Q^\ve_{z,v,k}(s)|\leq C_2(1+b+s)(\ve s+\sqrt {\ve s})
   e^{C_2(1+b+s)s}
   \end{equation}
   for another constant $C_2>0$ independent of $z,x',v,\be,t$ and $\ve$.
   
   Choose $s=s(\ve)=(\log(1/\ve))^{1/3}$ and set $l(\ve)=
   [t/(\log(1/\ve))^{1/3}]$. By (\ref{2.2.21}),
   \begin{equation}\label{2.3+.16}
   e^{-Kbs(\ve)}Q^\ve_{x,y,k}(l(\ve)s(\ve))\leq Q^\ve_{x,y,k}(t)\leq
 e^{Kbs(\ve)}Q^\ve_{x,y,k}(l(\ve)s(\ve)).
 \end{equation}
 If $|z-x|\leq K\ve t$ and $s\leq t$ then by (\ref{2.3+.1}), (\ref{2.3+.15})
 and the Lipschitz continuity of $H$ we obtain that 
 \begin{eqnarray}\label{2.3+.17}
 &\,\,\,\,\,\,\,-C_2(1+b+s)(\ve s+\sqrt {\ve s})e^{C_2(1+b+s)s}+
 \exp\big(s(H(x,x',\be)\\
 &-C(b)K\ve t-\rho_b(s))\big)
 \leq Q^\ve_{z,v,k}(s)\leq C_2(1+b+s)(\ve s+\sqrt {\ve s})e^{C_2(1+b+s)s}
 \nonumber\\
 &+\exp\big(s(H(x,x',\be)-C(b)K\ve t-\rho_b(s))\big)\nonumber
 \end{eqnarray}
 for some constant $C(b)>0$ independent of $z,x,x'\in\bar\cX$, $v\in M$, 
 $|\be|\leq b,\, t$ and $\ve$. Observe that by the Markov property,
 \begin{eqnarray}\label{2.3+.18}
 &Q^\ve_{x,y,k}(ls(\ve))=E\exp\langle\be,\int_0^{(l-1)s(\ve)}
 (B_{\nu^\ve_{x,y,k}(u)}(x',\hat Y^\ve_{x,y,k}(u))du\rangle\\
 &\times Q^\ve_{X^\ve_{x,y,k}((l-1)s(\ve)),Y^\ve_{x,y,k}((l-1)s(\ve)),
 \nu^\ve_{x,y,k}((l-1)s(\ve))}\big(s(\ve)\big).\nonumber
 \end{eqnarray}
 Now by (\ref{2.3+.10}) and (\ref{2.3+.17}) applying (\ref{2.3+.18}) for
 $l=l(\ve),l(\ve)-1,...,2$ we obtain that
 \begin{equation}\label{2.3+.19}
 \big\vert\frac 1t\log Q^\ve_{x,y,k}(t)-H(x,x',\be)\big\vert\leq\tilde C(b)
 \big(\ve^{1/3} +\ve t+\rho_b(\min (t,s(\ve))\big)
 \end{equation}
 for some $\tilde C(b)>0$ independent of $x,x'\in\bar\cX$, $y\in M$, $|\be|\leq
 b,\, t$ and $\ve$, which yields (\ref{2.2.3}) completing the proof of 
 Proposition \ref{prop2.2.4}. \qed

 \section{Further properties of $S$-functionals}\label{sec2.4}
\setcounter{equation}{0}

In this section we study essential properties of the functionals $S_{0T}$
which will be needed in the proofs of Theorems \ref{thm2.2.5} and \ref{thm2.2.7}
in the next sections. The following result which follows from \cite{Pi} is a
 basic step in our analysis of functionals $S_{0t}(\gam)$ and our thanks go
 to R. Pinsky who quickly produced on our request \cite{Pi} deriving some
 properties of functionals $I_x(\mu)$ needed here.

\begin{lemma}\label{lem2.4.1} For each $x\in\bar\cX$ and any vector measure 
$\mu=(\mu_1,...,\mu_N)$ on $M$ with $\sum^N_{k=1}\mu_k(M)=1$,
 $I_x(\mu)<\infty$ if and only if each $\mu_k,\, k=1,...,N$ has density
 $g_k=d\mu_k/dm$ with respect to the Riemannian volume $m$ on $M$ such that
 \begin{equation}\label{2.4.1}
 \int_M\big\|\nabla\sqrt {g_k}\|^2dm<\infty,\,\, k=1,...,N
 \end{equation}
 where $\nabla$ is the Riemannian gradient and $\|\cdot\|$ is a corresponding
 norm. Furthermore, there exists $C>0$ such that for any $x\in\bar\cX$ and
 each $\mu$ as above for which (\ref{2.4.1}) holds true,
  \begin{equation}\label{2.4.2}
  C^{-1}\sum_{k=1}^Na_k\int_M\|\nabla\sqrt {g_k}\|^2dm-C\leq
  I_x(\mu)\leq C\sum_{1\leq k\leq N}a_k\int_M\|\nabla\sqrt {g_k}\|^2dm+C
  \end{equation}
 where $a_k=\mu_k(M)$, and if $z\in\bar\cX$ is another point then
 \begin{equation}\label{2.4.3}
 |I_x(\mu)-I_z(\mu)|\leq C|x-z|
 \sum_{k=1}^Na_k\int_M\|\nabla\sqrt {g_k}\|^2dm.
 \end{equation}
\end{lemma}

Next, using Lemma \ref{lem2.4.1} we are able to show that each point where $B$
is complete can be connected with close points by curves with small 
$S$-functionals which, in particular, enables us to obtain important examples
of $S$-compacts.
\begin{lemma}\label{lem2.4.2} (i) There exists $C>0$ and for each $x\in\bar\cX$ 
where the vector field $B$ is complete there exists $r=r(x)>0$ such that if 
$|z_1-x|<r$ and $|z_2-x|<r$ then we can construct $\gam\in C_{0t}$ with 
$t\leq C|z_1-z_2|$ satisfying 
\[
\gam_0=z_1,\,\,\gam_{t}=z_2\,\,\mbox{and}\,\, S_{0t}(\gam)\leq C|z_1-z_2|.
\]
It follows that $R(\tilde z,z)$ and $R(z,\tilde z)$ are locally Lipschitz
continuous in $z$
belonging to the open $r$-neighborhood of $x$ when $\tilde z$ is fixed.

(ii) Let $\cO\subset\cX$ be a compact $\Pi^t$-invariant 
set which either contains a dense in $\cO$ orbit of $\Pi^t$ or $R(x,z)=0$ for
any pair $x,z\in\cO$. Suppose that $B$ is complete at each point of $\cO$. 
Then $\cO$ is an $S$-compact.

(iii) Assume that for any $\eta>0$ there exists $T(\eta)>0$ such that for
each $x\in\cO$ its orbit $\{\Pi^tx,\, t\in[0,T(\eta)]\}$ of length $T(\eta)$
forms an $\eta$-net in $\cO$ and suppose that $B$ is complete at a point
of $\cO$. Then $\cO$ is an $S$-compact.
\end{lemma}
\begin{proof} (i) Fix some $x\in\bar\cX$ and assume that $B$ is complete at $x$.
 Then we can find a simplex $\Del_x$ with vertices in $\Gam_x=\{\bar B_\mu(x):
 I_x(\mu)<\infty\}$ such that $\{\al\Del_x,\,\al\in[0,1]\}$ contains an open
neighborhood of 0 in $\bbR^d$ and 
\[
\Del_x=\{\sum^k_{i=1}\la_i\bar B_{\mu^{(i)}}(x):\,\sum^k_{i=1}\la_i=1,\,
\la_i\geq 0\,\forall i\}
\]
for some $\mu^{(i)}$ with $I_x(\mu^{(i)})<\infty$. By compactness of $\Del_x$ 
it follows that
\[
\mbox{dist}(\Del_x,0)=d_x>0.
\]
 By (\ref{2.2.1}) there exists a small $r(x)>0$ such that if $|z-x|\leq r(x)$ 
 then each simplex 
\[
\Del_z=\{\sum^k_{i=1}\la_i\bar B_{\mu^{(i)}}(z):\,\sum^k_{i=1}\la_i=1,\,
\la_i\geq 0\,\forall i\}
\]
intersects and not at 0 with any ray emanating from $0\in\bbR^d$ or, in other
words, $\{\al\Del_z,\,\al\in[0,1]\}$ contains an open neighborhood of 0 in 
$\bbR^d$ and, moreover,
\[
\mbox{dist}(0,\Del_z)\geq\frac 12d_x.
\]
It follows that for
any $z$ in the $r(x)$-neighborhood of $x$ and any vector $\xi$ there exist
$\la_1,...,\la_k\geq 0$ with $\la_1+\cdots +\la_k=1$ such that 
\[
\sum_{i=1}^k\la_i\bar B_{\mu^{(i)}}(z)=\bar B_{\sum_{1\leq i\leq k}
\la_i\mu^{(i)}}(z)=\xi.
\]
Observe that by (\ref{2.4.2}) and convexity of $I_z$,
\[
I_z(\sum_{1\leq i\leq k}\la_i\mu^{(i)})\leq\max_{1\leq i\leq k}I_z(\mu^{(i)})
\leq\tilde C\big(\max_{1\leq i\leq k}I_z(\mu^{(i)})+1\big)
\]
for some $\tilde C>0$. Hence, any two points
$z_1$ and $z_2$ from the open $r(x)$-neighborhood of $x$ can be connected
by a curve $\gam$ lying on the interval connecting $z_1$ and $z_2$ with
$K\geq |\dot {\gam}_s^{(1)}|\geq\frac 12d_x$, i.e. $\gam_0=z_1,\,\gam_t=z_2$
with some $t\in[K^{-1}|z_1-z_2|,2d^{-1}_x|z_1-z_2|]$ and by (\ref{2.2.9}),
\[
S_{0t}(\gam)\leq 2\tilde Cd_x^{-1}|z_1-z_2|(\max_{1\leq i\leq k}I_x(\mu^{(i)})
+1).
\]
In view of the triangle inequality for $R$ what we have proved yields the
continuity of $R(\tilde z,z)$ and $R(z,\tilde z)$ in $z$ belonging to the
open $r(x)$-neighborhood of $x$ when $\tilde z$ is fixed. Covering $\bar\cX$
by $r(x)-$neighborhoods of points $x\in\bar\cX$ and choosing a finite
subcover we obtain (i) with the same constant $C>0$ for all points in $\bar\cX$.

For the proof of sufficient conditions (ii) and (iii) of $S$-compacthess
see Lemma \ref{lem1.6.2}(ii)--(iii) in Part \ref{part1}.
 \end{proof}

\begin{lemma}\label{lem2.4.3} For any $\eta>0$ and $T>0$ there exists $\zeta>0$
such that if $\gam\in C_{0T},\,\gam\subset\cX$, $S_{0T}(\gam)<\infty$, 
$\gam_0=x_0$, and $|z_0-x_0|<\zeta$ then we can find $\tilde\gam\in C_{0T}$, 
$\tilde\gam\subset\cX$ with $\tilde\gam_0=z_0$ satisfying
\begin{equation}\label{2.4.4}
\bfr_{0T}(\gam,\tilde\gam)<\eta\,\,\mbox{and}\,\,|S_{0T}(\tilde\gam)-
S_{0T}(\gam)|<\eta.
\end{equation}
\end{lemma}
\begin{proof} By (\ref{2.2.9}), (\ref{2.2.19}) and the lower semicontinuity of
 the functionals $I_z(\nu)$ there exist measures $\nu_t\in\cM_{\gam_t},\,
t\in[0,T]$ such that $\dot\gam_t=\bar B_{\nu_t}(\gam_t)$ for Lebesgue almost
all $t\in[0,T]$ and $I_{\gam_t}(\nu_t)=L(\gam_t,\dot\gam_t)$ for Lebesgue 
almost all $t\in[0,T]$. Recall also that $\dot\gam_t$ is measurable in $t$. 
Introduce the (measurable) map
$q:\,[0,T]\times\cP(\bar\cW)\to\bbR\cup\{\infty\}\times\bbR^d$ defined
by $q(t,\nu)=\big( I_{\gam_t}(\nu),\bar B_\nu(\gam_t)\big)$. Recall that
$\dot\gam_t$ is measurable in $t$, and so another map $r:\,[0,T]\to
\bbR\cup\{\infty\}\times\bbR^d$ defined by $r(t)=\big(L(\gam_t,\dot\gam_t\big),
\dot\gam_t\big)$ is also measurable in $t\in[0,T]$. Then $q(t,\nu_t)=r(t)$
 and it follows from the measurable selection in the implicit function theorem 
(see \cite{CV}, Theorem III.38) that measures $\nu_t$ satisfying this 
condition can be chosen to depend measurably on $t\in[0,T]$. Since
$S_{0T}(\gam)<\infty$ and the $I$-functionals are nonnegative then
$I_{\gam_t}(\nu_t)<\infty$ for Lebesgue almost all $t\in[0,T]$ (and,
actually, without loss of generality we can assume that $I_{\gam_t}(\nu_t)$
is finite for all $t\in[0,T]$).

Now let
\[
\tilde\gam_t=z_0+\int_0^t\bar B_{\nu_s}(\tilde\gam_s)ds,\,\, t\in[0,T],
\]
which in view of (\ref{2.2.1}) determines $\tilde\gam\in C_{0T}$. Then by 
(\ref{2.2.1}),
\[
\bfr_{0t}(\gam,\tilde\gam)\leq\zeta+K\int_0^t\bfr_{0s}(\gam,\tilde\gam)ds
\]
and by Gronwall's inequality
\[
\bfr_{0T}(\gam,\tilde\gam)\leq\zeta e^{KT}.
\]
This together with (\ref{2.4.2}) and (\ref{2.4.3}) yields that
\[
\big\vert\int_0^TI_{\gam_t}(\nu_t)dt-\int_0^TI_{\tilde\gam_t}(\nu_t)dt\big\vert
\leq\tilde C\zeta e^{KT}S_{0T}(\gam)
\]
for some $\tilde C>0$ independent of $\zeta$ and $\gam$. Exchanging $\gam$
and $\tilde\gam$, applying the same argument and using the inequality
$S_{0T}(\tilde\gam)\leq\int_0^TI_{\tilde\gam_t}(\nu_t)dt$ we conclude that
\[
\big\vert S_{0T}(\gam)-S_{0T}(\tilde\gam)\big\vert\leq \tilde C\zeta e^{KT}
\max\big(S_{0T}(\gam),S_{0T}(\tilde\gam)\big)\leq \tilde C\zeta e^{KT}
(1+\tilde C\zeta e^{KT})S_{0T}(\gam).
\]
Choosing $\zeta$ small enough we arrive at (\ref{2.4.4}).
\end{proof}

The following result will enable us to control the time which the slow 
 motion can spend away from the $\om$-limit set of the averaged motion.
 \begin{lemma}\label{lem2.4.4} Let $G\subset\cX$ be a compact set not 
    containing entirely any forward semi-orbit of the flow $\Pi^t$.
    Then there exist positive
     constants $a=a_G$ and $T=T_G$ such that for any $x\in G$ and $t\geq 0$,
     \begin{equation*}
     \inf\big\{ S_{0t}(\gam):\,\gam\in C_{0t}\,\,\mbox{and}\,\,\gam_s\in G\,\,
     \mbox{for all}\,\, s\in[0,t]\big\}\geq a[t/T]
     \end{equation*}
     where $[c]$ denotes the integral part of $c$.
     \end{lemma}
     \begin{proof} The result is a simple consequence of lower semicontinuity
     of functionals $S_{0t}$ and the fact that $S_{0T}(\gam)=0$ if and only
     if $\gam$ is a part of an orbit of the flow $\Pi^t$. Further details of
     the argument can be found in Lemma \ref{lem1.6.4} in Part \ref{part1}
     and in Lemma 2.2(a), Chapter 4 of \cite{FW}.
     \end{proof}
     
     For the proof of the following result see Lemma \ref{lem1.6.5} in 
     Part \ref{part1}.
\begin{lemma}\label{lem2.4.5} 
Let $V$ be a connected open set with a piecewise smooth boundary and assume
that (\ref{2.2.24}) holds true. Then the function $R_\partial(x)$ is upper
semicontinuous at any $x_0\in V$ for which $R_\partial(x_0)<\infty$. 
Let $\cO\subset V$ be an $S$-compact. 

(i) Then for each $z\in\bar V$ the function $R(x,z)$ takes on the same
value $R^\cO(z)$ for all $x\in\cO$, and so $R_\partial(x)$ takes on the same
value $R_\partial$ for all $x\in\cO$ and the set $\partial_{\min}(x)=\{ z\in
\partial V:\, R(x,z)=R_\partial\}$ coincides with the same (may be empty)
set $\partial_{\min}$ for all $x\in\cO$. Furthermore, for each $\del>0$ 
there exists $T(\del)>0$ such that for any $x\in\cO$ we can construct 
$\gam^x\in C_{0t_x}$ with $t_x\in(0,T(\del)]$ satisfying
\begin{equation}\label{2.4.5}
\gam^x_0=x,\,\,\gam^x_{t_x}\in\partial V\,\,\mbox{and}\,\, S_{0t_x}(\gam^z)
\leq R_\partial+\del.
 \end{equation}
 
 (ii) Suppose that $R_\partial<\infty$ and dist$(\Pi^tx,\cO)
 \leq d(t)$ for some $x\in V$ and $d(t)\to 0$ as $t\to\infty$. Then
 $R_\partial(x)\leq R_\partial$ and for any $\del>0$ there exist $T_{\del,d}
 >0$ (depending only on $\del$ and the function $d$ but not on $x$) and
 $\hat\gam^x\in C_{0s_x}$ with $s_x\in(0,T_{\del,d}]$ satisfying
 \begin{equation}\label{2.4.6}
\hat\gam^x_0=x,\,\,\hat\gam^x_{s_x}\in\partial V\,\,\mbox{and}\,\, S_{0s_x}
(\hat\gam^x)\leq R_\partial+\del.
 \end{equation}
 In particular, if $R_\partial <\infty$ then $R_\partial(x)<\infty$ and
 if $\cO$ is an $S$-attractor of the flow $\Pi^t$ then $R_\partial(x)<\infty$
 for all $x\in V$.
 
 (iii) Suppose that for any open set $U\supset\cO$ the compact set
 $\bar V\setminus U$ does not contain entirely any forward semi-orbit of
 the flow $\Pi^t$. Then the function $R^\cO(z)$ is lower semicontinuous
 in $z\in\bar V$, $R^\cO(z)\to 0$ as dist$(z,\cO)\to 0$, and $\partial_{\min}$
 is a nonempty compact set.
\end{lemma}

\section{"Very long" time behavior: exits from a domain}\label{sec2.5}
\setcounter{equation}{0}

We start with the following result which will not only yield Theorem 
\ref{thm2.2.5} but also will play an important role in the proof of 
Theorem \ref{thm2.2.7}. 
\begin{proposition}\label{prop2.5.1} Let $V$ be a connected open set with
a piecewise smooth boundary $\partial V$ such that $\bar V=V\cup\partial V
\subset\cX$. Assume that for each $z\in\partial V$ there exist $\iota=
\iota(z)>0$ and a probability measure $\mu$ with $I_z(\mu)<\infty$ so that
\begin{equation}\label{2.5.1}
z+s\bar B_\mu(z)\in\bbR^d\setminus\bar V\,\,\mbox{for all}\,\, s\in(0,\iota],
\end{equation}
i.e. $\bar B_\mu(z)\ne 0$ and it points out into the exterior of $\bar V$.

(i) Suppose that for some $A_1,T>0$ and any $z\in\bar V$ there exists 
$\vf^z\in C_{0T}$ such that for some $t=t(z)\in(0,T]$,
\begin{equation}\label{2.5.2}
\vf_0^z=z,\,\vf^z_t\not\in V\,\,\mbox{and}\,\, S_{0t}(\vf^z)\leq A_1.
\end{equation}
Then for any $x\in V$ uniformly in $y\in\bfM$,
\begin{equation}\label{2.5.3}
\limsup_{\ve\to 0}\ve\log E\tau^\ve_{x,y}(V)\leq A_1
\end{equation}
and for any $\al>0$ there exists $\la(\al)=\la(x,\al)>0$ such that 
uniformly in $y\in\bfM$ for all small $\ve>0$,
\begin{equation}\label{2.5.4}
P\big\{ \tau^\ve_{x,y}(V)\geq e^{(A_1+\al)/\ve}\big\}\leq 
e^{-\la(\al)/\ve}.
\end{equation}

(ii) Assume that there exists an open set $G$ such that $V$ contains its 
closure $\bar G$ and the intersection of $\bar V\setminus G$ with the 
$\om$-limit set of the flow $\Pi^t$ is empty. Let $\Gam$
be a compact subset of $\partial V$ such that
\begin{equation}\label{2.5.5}
\inf_{x\in G,z\in\Gam}R(x,z)\geq A_2
\end{equation}
for some $A_2>0$. Then for some $T>0$ and any $\be>0$ there exists 
$\la(\be)>0$ such that uniformly in $y\in\bfM$ for each $x\in V$ and 
any small $\ve>0$,
\begin{eqnarray}\label{2.5.6}
&P\big\{ Z^\ve_{x,y}(\tau^\ve_{x,y}(V))\in\Gam,\,
\tau^\ve_{x,y}(V)\leq e^{(A_2-\be)/\ve}\big\}\\
&\leq P\big\{ Z^\ve_{x,y}(\tau^\ve_{x,y}(V))\in\Gam,\,
\tau^\ve_{x,y}(V)<T\big\}+e^{-\la(\be)/\ve}.\nonumber
\end{eqnarray}
Suppose that for some $x\in V$,
\begin{equation}\label{2.5.7}
a(x)=\inf_{t\geq 0}\mbox{dist}(\Pi^tx,\partial V)>0.
\end{equation}
Then $R_\partial(x)>0$ and for each $T>0$ there exists $\hat\la(T)=
\hat\la(T,x)>0$ such that uniformly in $y\in\bfM$ for all small $\ve>0$,
\begin{equation}\label{2.5.8}
P\{ \tau^\ve_{x,y}(V)<T\}\leq e^{-\hat\la(T)/\ve}
\end{equation}
and if the set $\Gam$ from (\ref{2.5.5}) coincides with the whole
$\partial V$ then for all $x\in V$ uniformly in $y\in\bfM$, 
\begin{equation}\label{2.5.9}
\liminf_{\ve\to 0}\ve\log E\tau^\ve_{x,y}(V)\geq A_2.
\end{equation}
\end{proposition}
\begin{proof}
    In order to prove (i) we observe, first, that the assumption (\ref{2.5.1})
    above together with Lemma \ref{lem2.4.2}(i) and the compactness of 
    $\partial V$
    considerations enable us to extend any $\vf^z,\, z\in V$ slightly so that 
    it will exit some fixed neighborhood of $V$ with only slight increase in
    its $S$-functional. Hence, from the beginning we assume that for each
    $\be>0$ there exists $\del=\del(\be)>0$ such that for any $z\in V$ we can
    find $T>0$, $\vf^z\in C_{0T}$ and $t=t(z)\in(0,T]$ satisfying
    \[
    \vf_0^z=z,\,\vf^z_t\not\in V_\del\,\,\mbox{and}\,\, S_{0t}(\vf^z)\leq
    A_1+\be
    \] 
    where $V_\del=\{ x:\,$dist$(x,V)\leq\del\}$. Employing the Markov
    property \index{Markov property} we obtain that for any 
    $x\in V,\, y\in\bfM,\, n\geq 1$,
    \begin{eqnarray}\label{2.5.10}
    &P\big\{ \tau^\ve_{x,y}(V)>nT\big\}=P\big\{ Z^\ve_{x,y}(t)\in V,\,
    \forall\,t\in[0,nT]\big\}\\
    &=P\big\{\tau^\ve_{Z^\ve_{x,y}(kT)}(V)>T,\,\,
    \forall\, k=0,1,...,n-1\big\}\leq\big(\sup_{w\in V\times\bfM}
    P\big\{\tau^\ve_{w}(V)>T\big\}\big)^n. \nonumber
    \end{eqnarray}
    From (\ref{2.2.10}) and (\ref{2.5.2}) it follows that
    \begin{equation}\label{2.5.11}
    P\big\{\tau^\ve_w(V)>T\big\}\leq P\big\{\bfr_{0T}(Z^\ve_w,\vf^z)\geq
    \del\,\,\mbox{for any}\,\, z\in V\big\}\leq 1-\exp
    \big(-(A_1+\be+\la)/\ve\big).
    \end{equation}
    By (\ref{2.5.10}) and (\ref{2.5.11}),
     \begin{equation}\label{2.5.12}
    P\big\{\tau_w^\ve(V)>e^{(A_1+\be)/\ve}\big\}
    <e^{-c(\be)/\ve}
    \end{equation}
    and
    \begin{eqnarray}\label{2.5.13}
    &E\tau_w^\ve(V)\leq\sum_{n=0}^\infty (n+1)T
    \big( P\big\{\tau_w^\ve(V)>nT\big\}- P\big\{\tau_w^\ve(V)\\
    &>(n+1)T\big\}\big) 
    =T\sum_{n=0}^\infty P\big\{\tau_w^\ve(V)>nT\big\}
    \leq T\exp\big(\frac 1\ve(A_1+\be+\la)\big)
    \nonumber\end{eqnarray}
    yielding (\ref{2.5.3}) and (\ref{2.5.4}) since $\be$ and $\la$ in 
    (\ref{2.5.13}) can be chosen arbitrarily small as $\ve\to\infty$.
    
     Next, we derive the assertion (ii). Let $t>0$ and $n$ be 
    the integral part of $t/T$ where $T>0$ will be chosen later. Let, 
    again, $w=(x,y)$ with $x\in V$ and $y\in\bfM$. Then
    \begin{eqnarray}\label{2.5.14}
    & P\{ Z^\ve_w(\tau^\ve_w)\in\Gam,\,\tau_w^\ve(V)<t\}\\
     &\leq P\{ Z^\ve_w(\tau^\ve_w(V))\in\Gam,\,
     \tau^\ve_w(V)<(n+1)T\}\nonumber\\
     &=\sum_{k=0}^nP\{ Z^\ve_w(\tau^\ve_w(V))
     \in\Gam,\, kT\leq\tau^\ve_w(V)<(k+1)T\}.\nonumber 
    \end{eqnarray}
    Let $K$ be the intersection of the $\om$-limit set of the flow $\Pi^t$
    with $\bar V$. Then $K$ is a compact set and by our assumption 
    $K\subset G$. Hence, 
    \[
    \del=\frac 13\inf\{ |x-z|:\, x\in K,\,\, z\in\bar V\setminus G\}>0
    \]
    and if we set $U_\eta=\{ z\in V:\,\mbox{dist}(z,K)<\eta\}$ then
    $U_{3\del}\subset G$. Now suppose that 
    $kT\leq\tau^\ve_{x,v}(V)<(k+1)T$ for some $k\geq 1$ and
    $Z^\ve_{x,v}(\tau^\ve_{x,v}(V))\in\Gam$ with $x\in V$ and $v\in\bfM$.
    Then either there is $t_1\in[(k-1)T,kT]$ such that $Z^\ve_{x,v}(t)\in
    \bar V\setminus U_{2\del}$ for all $t\in[t_1,t_1+T]$ or there exist
    $t_2,t_3>0$ such that $(k-1)T\leq t_2<t_3<(k+1)T$ and $Z^\ve_{x,v}(t_2)
    \in U_{2\del}$ while $Z^\ve_{x,v}(t_3)\in\Gam$. Set $\cT_z=
    \{\gam\in C_{0,2T}:\,\gam_0=z$ and either there is $t_1\in[0,T]$ so that
    $\gam_t\in\bar V\setminus U_{2\del}$ for all $t\in[t_1,t_1+T]$ or 
    $\gam_{t_2}\in U_{2\del}$ and $\gam_{t_3}\in\Gam$ for some
    $0\leq t_2<t_3<2T\}$. Then for any $k\geq 1$,
    \begin{eqnarray}\label{2.5.15}
    &\{ Z^\ve_w(\tau^\ve_w(V))\in\Gam,\, kT\leq\tau^\ve_w(V)<(k+1)T\}\\
    &\subset\big\{ Z^\ve_w(\tau^\ve_w(V))\in\Gam,\,
    Z^\ve_{Z^\ve_w((k-1)T)}\in\cT_{Z^\ve_w((k-1)T)}\big\}.\nonumber
    \end{eqnarray}
   
   For each $q>0$
    set $\cT_z^q=\{\gam\in C_{0,2T}:\gam_0=z\,\,\mbox{and}\,\,\bfr_{0,2T}
    (\gam,\cT_z)\leq q\}$ and suppose that for some $\eta>0$ there is
    $d_\eta\geq 0$ so that
    \begin{equation}\label{2.5.16}
    \inf_{z\in V}\inf_{\gam\in\cT_z^{2\eta}}S_{0,2T}(\gam)> d_\eta .
    \end{equation}
    Then $\cT_z^{2\eta}\cap\Psi^{d_\eta}_{0,2T}(z)=\emptyset$, 
    where $\Psi^a_{0,t}(z)$ is the same as in Theorem \ref{thm2.2.2}, and so
    \begin{equation}\label{2.5.17}
    \cT^\eta_z\subset\big\{\gam\in C_{0,2T}:\,\gam_0=z\,\,\mbox{and}\,\,
    \bfr_{0,2T}(\gam,\Psi^{d_\eta}_{0,2T}(z))\geq\eta\big\}.
    \end{equation}
    From (\ref{2.2.10}) and (\ref{2.5.15})--(\ref{2.5.17}) we obtain that for
    any $\be>0$ and all sufficiently small $\ve$,
    \begin{equation}\label{2.5.18}
    P\{ Z^\ve_w(\tau^\ve_w(V))\in\Gam,\, kT\leq\tau^\ve_w(V)<(k+1)T\}\leq
    \hat Ce^{-(d_\eta-\be)/\ve}
     \end{equation}
     for some $\hat C>0$.
     
      Next, we will specify $d_\eta$ in (\ref{2.5.16}) choosing 
     $\eta\leq\frac 12\del$. For each $z\in V$ we can write
    \begin{equation}\label{2.5.19}
    \cT_z^{2\eta}\subset\tilde\cT_z^\eta\cup\hat\cT_z^\eta
    \end{equation}
    where $\tilde\cT_z^\eta=\{\gam\in C_{0,2T}:\,\gam_0=z,\,
    \gam_{t_2}\in U_{3\del}$ and $\gam_{t_3}\in\Gam_{2\eta}$ for some 
    $0\leq t_2<t_3<2T\}$ with $\Gam_{r}=\{z:$ dist$(z,\Gam)\leq r\}$ and 
    $\hat\cT_z=\{\gam\in C_{0,2T}:\,\gam_0=z$ and there is $t_1\in[0,T]$ so
    that $\gam_t\in V_{2\eta}\setminus U_\del$ for all $t\in[t_1,t_1+T]\}$.
    By (\ref{2.5.5}) and the lower semicontinuity of the functional $S_{0,2T}$
    it follows that for any $\zeta>0$ we can choose $\eta>0$ small enough
    so that 
    \begin{equation}\label{2.5.20}
    \inf_{z\in V}\inf_{\gam\in\tilde\cT_z^\eta}S_{0,2T}(\gam)> A_2-\zeta.
    \end{equation}
    
    Since $\bar V\setminus U_\del$ is disjoint with the $\om$-limit set of the
    flow $\Pi^t$ and the latter is closed then if $\eta$ is sufficiently small
    $V_{2\eta}\setminus U_\del$ is also disjoint with this $\om$-limit set
    and, in particular, it does not contain any forward semi-orbit
    of $\Pi^t$. Hence we can apply Lemma \ref{lem2.4.4} which in view of
    (\ref{2.2.9}) implies that there exists $a>0$ such that for all small
     $\eta>0$,
    \begin{equation}\label{2.5.21}
    \inf_{z\in V}\inf_{\gam\in\hat\cT_z}S_{0,2T}(\gam)> aT
    \end{equation}
    which is not less than $A_2$ if we take $T=A_2/a$. Now, (\ref{2.5.20}) and 
    (\ref{2.5.21}) produce (\ref{2.5.16}) with $d=A_2-\zeta$, and so 
    (\ref{2.5.18}) follows with such $d_\eta$.
    This together with (\ref{2.5.14}) yield that for any $\be>0$ we can choose
    sufficiently small $\zeta,\la >0$ and then $\eta>0$ so that for all $\ve$
    small enough
    \begin{eqnarray}\label{2.5.22}
    &P\big\{ Z^\ve_v(\tau^\ve_v)\in\Gam,\,
    \tau^\ve_v(V)\leq e^{(A_2-\be)/\ve}\big\}\\
    &\leq P\{ Z^\ve_v(\tau^\ve_v)\in\Gam,\,\tau_v^\ve(V)<T\}+e^{-\la/2\ve}
    \nonumber\end{eqnarray}
    and (\ref{2.5.6}) follows.
    
    Now assume that (\ref{2.5.7}) holds true for some $x\in V$. 
    Recall, that $S_{0T}(\gam)=0$ implies that $\gam$ is a piece of an
    orbit of the flow $\Pi^t$. Since no $\gam\in C_{0T}$ satisfying
    \begin{equation}\label{2.5.23}
    \gam_0=x\,\,\mbox{and}\,\,\inf_{t\in[0,T]}\,\mbox{dist}
    (\gam_t,\partial V)\leq a(x)/2
    \end{equation}
    can be such piece of an orbit we conclude by the lower semicontinuity of
    $S_{0T}$ that $S_{0T}(\gam)> c(x)$ whenever (\ref{2.5.23}) holds true 
    for some $c(x)>0$ independent of $\gam$ (but depending on $x$).
     Hence, by (\ref{2.2.11}), 
    \begin{eqnarray}\label{2.5.24}
    &P\{ \tau^\ve_v<T\}\leq P\big\{ \bfr_{0T}
    \big(Z^\ve_v, \Psi^{c(x)}_{0T}(x)\big)\\
    &\geq a(x)/2\big\}\leq\exp(-c(x)/2\ve)\nonumber
    \end{eqnarray}
    provided $\ve$ is small enough and (\ref{2.5.8}) follows. 
    Observe also that any $\gam\in C_{0t}$
    with $\gam_0=x\in V$ and $\gam_t\in\partial V$ should contain a piece
    which either belongs to some $\tilde\cT_z^\eta$ or to $\hat\cT_z^\eta$, 
    as above, or to satify (\ref{2.5.23}). By (\ref{2.5.20}), (\ref{2.5.21}), 
    and the above remarks it follows that $S_{0t}(\gam)\geq q(x)$ for such 
    $\gam$ where $q(x)>0$ depends only on $x$, and so $R_\partial(x)\geq q(x)$.
    If $\Gam=\partial V$ then by (\ref{2.5.6}) and (\ref{2.5.8}),
    \begin{eqnarray}\label{2.5.25}
    &E\tau^\ve_{x,y}(V)\geq e^{(A_2-\be)/\ve}\, P\big\{ 
     \tau^\ve_{x,y}(V)\geq e^{(A_2-\be)/\ve}\big\}\\
    &\geq e^{(A_2-\be)/\ve}(1-e^{-\la(\be)/\ve}-e^{-\hat\la(T)/\ve})\nonumber
    \end{eqnarray}
    and, since $\be>0$ is arbitrary, (\ref{2.5.9}) follows completing the
    proof of Proposition \ref{prop2.5.1}.
    \end{proof}
    
    Now we will derive Theorem \ref{thm2.2.5} from Proposition \ref{prop2.5.1}.
   Assume, first, that $R_\partial <\infty$. Then by Lemma \ref{lem2.4.5}, 
   $R_\partial(x)$ is finite in the whole $V$. Moreover,
   since $\cO$ is an $S$-attractor the conditions of Lemma \ref{lem2.4.5}
   are satisfied with some $d(t)\to 0$ as $t\to\infty$ the same for all
   points of $V$ which yields the conditions of Proposition \ref{prop2.5.1}(i)
   with $A_1=R_\partial +\del$ for any $\del>0$. Hence, (\ref{2.5.3}) 
   and (\ref{2.5.4})
   hold true with $A_1=R_\partial$. Since $\cO$ is an $S$-attractor of the
   flow $\Pi^t$ and its basin contains $\bar V$ then the intersection of
   $\bar V\setminus\cO$ with the $\om$-limit set of $\Pi^t$ is empty.
   By the definition of an $S$-attractor for any $\zeta>0$ there exists
   an open set $U_\zeta\supset\cO$ such that $R(x,z)\leq\zeta$ whenever
   $x\in\cO$ and $z\in U_\zeta$. Hence, by the triangle inequality for
   the function $R$ and Lemma \ref{lem2.4.5} for any set 
   $\Gam\subset\partial V$,
   \begin{equation}\label{2.5.26}
   \inf_{z\in U_\zeta,\tilde z\in\Gam}R(z,\tilde z)\geq\inf_{\tilde z\in\Gam}
   R^\cO(\tilde z)-\zeta.
   \end{equation}
   If $\Gam=\partial V$ then by Lemma \ref{lem2.4.5} the right hand side of 
   (\ref{2.5.26}) equals
   $A_2=R_\partial-\zeta$. Assuming that $R_\partial<\infty$ we can apply
   Proposition \ref{prop2.5.1}(ii) with such $A_2$ yielding (\ref{2.5.6}), 
   (\ref{2.5.8}) and since $\zeta>0$ is arbitrary (\ref{2.2.25}) and 
   (\ref{2.2.26})
   follow in this case. If $R_\partial=\infty$ then (\ref{2.2.26}) is trivial 
   and by (\ref{2.5.26}), $R(z,\tilde z)=\infty$ for any $z\in U_\zeta$ and
   $\tilde z\in\partial V$, and so we can apply Proposition \ref{prop2.5.1}(ii)
   with any $A_2$ which sais that the left hand side in (\ref{2.5.9}) equals
   $\infty$, and so (\ref{2.2.25}) holds true in this case, as well.
   
    Next, we establish (\ref{2.2.27}). For small $\del,\be>0$ and large $T>0$
  which will be specified later on set $t_\ve=T+e^{\be/\ve}$
  and define the event 
  \[
  \Xi^\ve_T(n)=\{ \tau^\ve_{Z_v^\ve(t_\ve n+T),Y^\ve_v((t_\ve n+T)/\ve)}
  (U_\del(\cO))\leq e^{\be/\ve}\}.
  \]
  Then
  \begin{eqnarray}\label{2.5.27}
  &\Te^\ve_v((n+1)t_\ve\wedge\tau^\ve_v(V))-
  \Te^\ve_v(nt_\ve\wedge\tau^\ve_v(V))\\
  &\leq T+t_\ve\bbI_V(Z^\ve_v(t_\ve n))
  \big(\bbI_{V\setminus U_{\del/2}(\cO)}(Z_v^\ve(t_\ve n+T))+
  \bbI_{U_{\del/2}(\cO)}(Z_v^\ve(t_\ve n+T))\bbI_{\Xi^\ve_T(n)}\big).
  \nonumber\end{eqnarray}
  If $\del$ is sufficiently small then $V_\del$ is still contained in the
  basin of $\cO$ with respect to the flow $\Pi^t$, and so we can choose
  $T$ (depending only on $\del$) so that
  \[
  \Pi^TV_\del\subset U_{\del/4}(\cO).
  \]
  Then for some $a>0$,
  \[
  \inf\big\{ S_{0T}(\gam):\,\gam\in C_{0T},\,\gam_0\in V_\del,\,\gam_T\not\in
  U_{\del/3}(\cO)\big\}>a,
  \]
  and so if $\gam_0\in V_\del$ and $\gam_T\not\in U_{\del/2}(\cO)$ then 
  dist$(\gam,\Psi^a_{0T}(z))\geq\del/6$ for any $z\in V_\del$. Relying on
  (\ref{2.2.11}) and the Markov property we obtain that for any 
  $v=(z,y)$ with $z\in V$,
  \begin{equation}\label{2.5.28}
  P\big\{ Z_v^\ve(t_\ve n)\in V\,\,\mbox{and}\,\, Z_v^\ve(t_\ve n+T)\in
  V\setminus U_{\del/2}(\cO) \big\}\leq e^{-a/2\ve}
  \end{equation}
  provided $\ve$ is small enough. Next, the same arguments which yield 
  (\ref{2.5.22}) and (\ref{2.5.24}) together with the Markov property
  enable us to conclude that if $\be>0$ is
  small enough then for any $v=(z,y)$ with $z\in V$,
  \begin{equation}\label{2.5.29}
   P\big\{ Z_v^\ve(t_\ve n+T)\in U_{\del/2}(\cO)\,\,\mbox{and}\,\,
   \Xi^\ve_T(n)\big\}\leq e^{-\be/\ve}.
  \end{equation}
  
  Applying (\ref{2.5.27})--(\ref{2.5.29}) we conclude that for sufficiently 
  small $\be$ and any much smaller $\ve$,
  \begin{equation}\label{2.5.30}
  E\big(\Te^\ve_v((n+1)t_\ve\wedge\tau^\ve_v(V))-
  \Te^\ve_v(nt_\ve\wedge\tau^\ve_v(V))\big)\leq t_\ve e^{-\be/\ve}(T+1).
  \end{equation}
  Finally, (\ref{2.2.26}) and (\ref{2.5.30}) together with the Chebyshev inequality
  yield that for $n(\ve)=[e^{(R_\partial+\be/4)/\ve}t_\ve^{-1}]$, each $x\in V$,
  a small $\be>0$ and any much smaller $\ve>0$,
  \begin{eqnarray}\label{2.5.31}
  &P\big\{ \Te^\ve_{x,w}(\tau^\ve_{x,w}(V))\geq e^{-\be/4\ve}
  \tau^\ve_{x,w}(V)\big\}\\
  &\leq P\big\{ \Te^\ve_{x,w}((n(\ve)+1)t_\ve)
  \geq e^{-\be/4\ve}e^{(R_\partial-\be/4)/\ve}\big\}\nonumber\\
  &+P\big\{ 
  \tau_{x,w}^\ve(V)<e^{(R_\partial-\be/4)/\ve}\,\,\mbox{or}\,\,
  \tau_{x,w}^\ve(V)>e^{(R_\partial+\be/4)/\ve}\big\}\nonumber\\
  &\leq\tilde Ce^{-\be/4\ve}\big(1+e^{-(R_\partial+\be/4)/\ve}
  (T+e^{\be/\ve})\big)+e^{-\la(\be/4)/\ve}.\nonumber
  \end{eqnarray}
   Since $R_\partial>0$ and we can choose $\be$ to be arbitrarily small, 
 (\ref{2.5.31}) yields (\ref{2.2.27}).
 
 In order to complete the proof of Theorem \ref{thm2.2.5} we have to derive
  (\ref{2.2.28}). If $\partial_{\min}=\partial V$ then 
    there is nothing to prove, so we assume that $\partial_{\min}$ is a
    proper subset of $\partial V$ and in this case, clearly,
    $R_\partial<\infty$. Since $\Gam=\{ z\in\partial V:\,\mbox{dist}
    (z,\partial_{\min})\geq\del\}$ is compact and disjoint with 
    $\partial_{\min}$ which is also compact then by the lower semicontinuity
    of $R^\cO(z)$ established in Lemma \ref{lem2.4.5}(iii) it follows that
    $R^\cO(z)\geq R_\partial +\be$ for some $\be>0$ and all $z\in\Gam$.
    Then by (\ref{2.5.26}), $R(z,\tilde z)\geq R_\partial+\be/2$ for any
    $z\in U_{\be/2}$ and $\tilde z\in\Gam$. Hence, applying Proposition 
    \ref{prop2.5.1} we obtain that
    \[
  P\big\{ \tau^\ve_{x,y}(V)\geq e^{(R_\partial+\frac 13\be)/\ve}
  \big\}\leq e^{-\la/\ve}
  \]
  and
  \[
  P\big\{ Z^\ve_{x,y}(\tau^\ve_{x,y}(V))\in\Gam,\,
  \tau^\ve_{x,y}(V)\leq
  e^{(R_\partial+\frac 13\be)/\ve}\big\}< e^{-\la/\ve}
  \]
  for some $\la>0$ and all $\ve$ small enough yielding (\ref{2.2.28})
  and completing the proof of Theorem \ref{thm2.2.5}. 
   \qed

\section{Adiabatic transitions between basins of attractors}\label{sec2.6}
\setcounter{equation}{0}

 In this section we will prove Theorem \ref{thm2.2.7} relying, again, on 
    Proposition \ref{prop2.5.1} together with Markov and strong Markov
    property of the Markov process $(X^\ve(t),Y^\ve(t))$.
    In view of (\ref{2.2.31}) and Lemma \ref{lem2.4.2}(i) any curve 
    $\gam\in C_{0t}$
     starting at $\gam_0=x\in V_{j_1}$ and ending at
    $\gam_t=z\in\cap_{1\leq i\leq k}\partial V_{j_i},\, k\leq\ell$ can be
    extended into each $V_{j_i},\, i=1,...,k$ with arbitrarily small increase
    in its $S$-functional. Hence, 
    \begin{equation}\label{2.6.1}
    R^{(i)}_\partial=\min_{j\ne i}R_{ij}
    \end{equation}
    where $R^{(i)}_\partial=\inf\{ R(x,z):\, x\in\cO_i,\, z\in\partial V_i\}$.
    Let $Q$ be an open ball of radius at least $r_0$ centered at the origin
    of $\bbR^d$. By Assumption \ref{ass2.2.6} the slow motion $Z^\ve_{x,y}$
    cannot exit $Q$ provided $x\in Q$ and $y\in\cW$. Furthermore, it is clear
    that $Q$ contains the $\om$-limit set of the averaged flow $\Pi^t$.
    Assumption \ref{ass2.2.6} enables us to deal only with restricted basins
    $V^Q_i=V_i\cap Q$ and though the boundaries $\partial V^Q$ of $V_i^Q$
    may include now parts of the boundary $\partial Q$ of $Q$ it makes no
    difference since $Z^\ve$ cannot reach $\partial Q$ if it starts in $Q$.
    Set $V^{(i)}=Q\setminus\cup_{j\ne i}U_\del(\cO_j)$ where $\del>0$ is
    small enough. We claim that in view of (\ref{2.2.31}) each $V_i$ satisfies 
  conditions of Proposition \ref{prop2.5.1}(i) for any $\be>0$ with $A_1=
  R_\partial^{(i)}+\be$ and some $T=T_\be$ depending on $\be$. Indeed, set 
 \[
 \partial(\eta)=\{ v\in Q:\,\mbox{dist}(v,\cup_{1\leq j\leq\ell}\partial V_j)
 \leq\eta\},\,\eta>0.
 \]
 In view of (\ref{2.2.31}) and Lemma \ref{lem2.4.1} there exists $L>0$ such 
 that 
 if $\eta$ is small enough and $z\in\partial(\eta)$ we can construct a curve
 $\vf^z\in C_{0,L\eta}$ with $S_{0,L\eta}(\vf^z)\leq L\eta$, $\vf^z_0=z,\,
 \vf_t^z\in V_j\setminus\partial(\eta)$ for some $t\in[0,L\eta]$ and 
 $j=1,...,\ell$. Since $V_j$ is the basin of $\cO_j$ there exists 
 $T=T_{\eta,\del}$ such that $\Pi^T\vf^z_t\in U_\del(\cO_j)$ and extending
 $\vf^z$ by the piece of the orbit of $\Pi^t$ we obtain a curve $\tilde\vf^z
 \in C_{0,L\eta+T}$ starting at $z$, entering $U_\del(\cO_j)$ and satisfying
 $S_{0,L\eta+T}(\tilde\vf^z)\leq L\eta$. Hence, for $z\in\partial(\eta)$
 the condition (\ref{2.5.2}) holds true with $V=V^{(i)}$ and $A_1=L\eta$.
 Since the $\om$-limit set of the flow $\Pi^t$ is contained in $Q\cap
 \big(\cup_{1\leq j\leq\ell}(\partial V_j\cup\cO_j)\big)$ it follows from
 Assumption \ref{ass2.2.6} and compactness considerations that there exists
 $\tilde T=\tilde T_{\eta,\delta}$ such that for any $z\in Q\setminus V_i$
 we can find $t_z\in[0,\tilde T]$ with $\Pi^{t_z}z\in\partial(\eta)\cup
 \big(\cup_{j\ne i}U_\del(\cO_j)\big)$. If $\Pi^{t_z}z\in
 \cup_{j\ne i}U_\del(\cO_j)$ then we take $\vf^z_t=\Pi^tz,\, t\in[0,\tilde T]$
 to satisfy (\ref{2.5.2}) for $V=V^{(i)}$ and $A_1=0$. If $\Pi^{t_z}z\in
 \partial(\eta)$ then we extend the curve $\vf^z_t=\Pi^tz,\, t\in[0,t_z]$
 as in the above argument which yields a curve $\tilde\vf^z$ starting at
 $z$, ending in some $U_\del(\cO_j),\, j\ne i$ and having its $S$-functional
 not exceeding $L\eta$. Finally, in the same way as in the proof of 
 Theorem \ref{thm2.2.5} for any $\be>0$ there exists $\hat T=\hat T_{\eta,\del,
 \be}$ such that whenever $z\in V_i(\eta)=V_i\cap Q\setminus\partial(\eta)$
 we can construct $\vf^z\in C_{0\hat T}$ such that (\ref{2.5.2}) holds true 
 with $V=V_i(\eta)$ and $A_1=R^{(i)}_\partial +\be/2$ and, moreover,
 dist$(\vf_t^z,V_j)\leq\eta$ for some $t\leq\hat T$ and $j\ne i$ with $R_{ij}
 =R^{(i)}_\partial$. Then in the same way
 as above we can extend $\vf^z$ to some $\tilde\vf^z\in C_{\hat T+\tilde T}$
 so that $\tilde\vf_t^z\in U_\del(V_j)$ for some $j$ as above,
 $t\leq\hat T+\tilde T$ and $S_{0,\hat T+\tilde T}(\tilde\vf^z)
 \leq R^{(i)}_\partial+\be/2+L\eta$ which gives (\ref{2.5.2}) for all 
 $z\in V=V^{(i)}$ with $A_1=R^{(i)}_\partial +\be$ provided $\eta$ is small 
 enough. Hence, Proposition \ref{prop2.5.1}(i) yields the estimates 
 (\ref{2.5.3}) 
 and (\ref{2.5.4}) for $\tau^\ve_{x,y}(i)$ in place of $\tau^\ve_{x,y}(V)$ with 
   $A_1=R_\partial^{(i)}$. In order to obtain the corresponding bounds in the
   other direction observe that in view of (\ref{2.2.31}),
   \begin{equation}\label{2.6.2}
   R^{(i)}_\partial(\del)=\inf\{ R(x,z):\, x\in\cO_i,\, z\not\in V_i(\eta)\}
   \to R^{(i)}_\partial\,\,\mbox{as}\,\,\del\to 0.
   \end{equation}
    Since $\overline {V_i(\eta)}$ is contained in the basin of $\cO_i$ we
    can apply to $V_i(\eta)$ the same estimates as in Theorem \ref{thm2.2.5}
    which together with (\ref{2.6.2}) and the fact that the exit time of
    $Z^\ve$ from $V_i(\eta)$ is smaller than its exit time from $V_i$ provide
    the remaining bounds yielding (\ref{2.2.32}) and (\ref{2.2.33}).
    
    Next, we derive (\ref{2.2.34}) similarly to (\ref{2.2.27}) but taking into
    account that $\cup_{1\leq j\leq\ell}\partial V_j$ may contain parts of
    the $\om$-limit set of the flow $\Pi^t$ which allows the slow motion
    $Z^\ve$ to stay long time near these boundaries. Still, set
    \[
    \te^\ve_v=\inf\{ t\geq 0:\, Z^\ve_v(t)\in\cup_{1\leq j\leq\ell}U_{\del/3}
    (\cO_j)\}.
    \]
    Using the same arguments as above we conclude that for any $\eta>0$
    there exists $T=T_{\eta,\del}$ such that whenever $z\in Q$ we can 
    construct $\vf^z\in C_{0T}$ with $\vf^z_0=z,\,\vf_T^z\in\cup_{1\leq j
    \leq\ell}U_\del(\cO_j)$ and $S_{0T}(\vf^z)\leq\eta$. This together with
    (\ref{2.5.12}) and Assumption \ref{ass2.2.6} yield that 
    \[
    P\{\te^\ve_v\geq e^{2\eta/\ve}\}\leq e^{-\la(\eta)/\ve}
    \]
    for some $\la(\eta)=\la(x,\eta)>0$ and all small $\ve$. Set 
    \[
    \Gam_1(v)=\big\{  Z^\ve_v(e^{2\eta/\ve})\in Q\setminus
    \cup_{1\leq j\leq\ell}U_{\del/2}(\cO_j)\big\},
    \]
    \[
    \Gam_2(v)=\big\{ \tau^\ve_v
    \big(\cup_{1\leq j\leq\ell}U_\del(\cO_j)\big)\leq e^{\be/\ve}\big\}
    \]
    and $t_\ve=e^{2\eta/\ve}+e^{\be/\ve}$ where $\eta$ is much smaller than 
    $\be$. Then proceeding similarly to the proof of (\ref{2.2.27}) as in
    (\ref{2.5.28})--(\ref{2.5.31}) above we arrive at (\ref{2.2.34}).
     
    Next, we obtain (\ref{2.2.35}) relying on additional assumptions specified
    in the statement of Theorem \ref{thm2.2.7}. Let $V^Q_i$ be the same as 
    above and $\partial^{(i)}_0(x)=\{ z\in\partial V_i^Q:\, 
  R(x,z)=R_\partial^{(i)}\}$. Since $\cO_i$ is an $S$-attractor it follows from
  Lemma \ref{lem2.4.5}(i) that $R(x,z)$ and $\partial_0^{(i)}(x)$ coincide with 
  the same function $R^{\cO_i}(z)$ and the same (in general, may be empty) set 
   $\partial_0^{(i)}$, respectively, for all $x\in\cO_i$. By Lemma 
   \ref{lem2.4.2}(i), our assumption that $B$ is complete on $\partial V_{i}$
   implies that $R^{\cO_{i}}(z)$ is continuous in a neighborhood of
   $\partial V_{i}$, and so $\partial_0^{(i)}$ is a nonempty compact set.
   Since we assume that $\io(i)\ne i$ is the unique index $j$ for which 
   $R_{ij}=R_{i\io(i)}=R_\partial^{(i)}$ then by (\ref{2.2.31}),
    \begin{equation*}
    \min_{j\ne i,\io(i)}\inf_{z\in\partial_0^{(i)}}\mbox{dist}
    (z,\partial V_j)>0.
    \end{equation*}
    Observe that if $\tilde\cO\subset\partial V_{i}$ is an $S$-compact then
    either $\tilde\cO\subset\partial_0^{(i)}$ or $\tilde\cO\cap
    \partial_0^{(i)}=\emptyset$. Denote by $L_\Pi$ the $\om$-limit set of
    the averaged flow $\Pi^t$. Since $L_\Pi\cap\partial V_{i}$ consists
    of a finite number of $S$-compacts it follows that
    \[
    \inf\{ |z-\tilde z|:\, z\in L_\Pi\cap\partial_0^{(i)},\,\tilde z\in
    L_\Pi\setminus\partial_0^{(i)}\}>0.
    \]
    By the continuity of $R^{\cO_{i}}(z)$ in $z\in\partial 
    V_{i}$ there exists $a>0$ such that 
    \[
    \inf\big\{  R^{\cO_{i}}(z):\, z\in\big(\cup_{j\ne i,\io(i)}(\partial 
    V_{i}\cap\partial V_j)\big)\cup\big((L_\Pi\setminus\partial_0^{(i)})
    \cap\partial V_{i}\big)\big\}\geq R_\partial^{(i)}+9a.
    \]
   These considerations enable us to construct a connected open set $G$
   with a piecewise smooth boundary $\partial G$ such that
   \[
   \bar G\subset V_{i}\cup(V_{\io(i)}\setminus\cO_{\io(i)})
   \cup\big((\partial V_{i}\cap
   \partial V_{\io(i)})\setminus (L_\Pi\setminus\partial_0^{(i)})\big)
   \]
   and for $\Gam=\partial G\setminus U_\del(\cO_{\io(i)})$ and some $a(\del)>0$,
   \begin{equation}\label{2.6.3}
   \inf_{z\in\Gam}R^{\cO_{i}}(z)\geq R_\partial^{(i)}+8a
   \end{equation}
   provided $a\leq a(\del)$. The idea of this construction is that if
   $Z^\ve_{x,y}(\tau^\ve_{x,y}(i))\not\in V_{\io(i)}$ then the slow motion 
   should exit $G$ through the part $\Gam$ of its boundary. Somewhat
   similarly to the proof of Proposition \ref{prop2.5.1}(ii) we will show
   that "most likely" this can only occur after the
   time $\exp\big((R_\partial^{(i)}+2a)/\ve\big)$ and, on the
   other hand, we conclude from (\ref{2.2.33}) that except for small 
   probability the exit time $\tau^\ve_{x,y}(i)$ does not exceed
   $\exp\big((R_\partial^{(i)}+a)/\ve\big)$.
   
   Let $U_0$ be a sufficiently small open neighborhood of $\partial_0^{(i)}$ 
   so that, in particular,
   \[
   \sup_{z\in U_0}R^{\cO_{i}}(z)\leq R_\partial^{(i)}+a
   \]
   and set
   \[
   \tau^\ve_{x,y}(G)=\inf\{ t\geq 0:\, Z^\ve_{x,y}(\tau^\ve_{x,y}(G))
   \not\in G\}.
   \]
   Then
   \begin{eqnarray}\label{2.6.4}
   &\big\{ \tau^\ve_v(G)\leq e^{(R_\partial^{(i)}+a)/\ve},\,
   Z^\ve_v(\tau^\ve_v(G))\in\Gam\big\}\subset\bigcup_{0\leq n\leq n(\ve)+1}
   \big(A^{(1)}(n)\cup\\
   &\bigcup_{(n-1)t_\ve\leq k\leq(n+1)t_\ve}\big( A^{(2)}(k)+A^{(3)}(k)
   +\bigcup_{k-2t_\ve\leq m\leq k-2T} A^{(4)}(m)\cap A^{(5)}(k)\big)\big)
   \nonumber\end{eqnarray}
   where $t_\ve=e^{\be/\ve}$ for some small $\be>0$, $n(\ve)=
   \big[e^{(R_\partial^{(i)}+a-\be)/\ve}\big]$, 
   $A^{(1)}(n)=\{ Z^\ve_v(t)\in G\setminus\big( U_\eta(\cO_{i})
   \cup U_\del(\cO_{\io(i)})\big)\,\,\mbox{for all}\,\, t\in[(n-1)t_\ve,nt_\ve]
   \big\}$ for a sufficiently small $\eta>0$, $A^{(2)}(k)=\big\{ 
   \exists t_1,t_2\,\mbox{with}\, k\leq t_1<t_2<k+3T,\, Z^\ve_v(t_1)
   \in U_\eta(\cO_{i}),\, Z^\ve_v(t_2)\in\Gam\big\}$,\,
   $A^{(3)}(k)=\{  Z^\ve_v(t)\in G\setminus\big(U_0\cup
    U_\eta(\cO_{i})\cup U_\del(\cO_{\io(i)})\big)\,\,\mbox{for all}\,\, 
    t\in[k,k+T]\big\}$,\, $A^{(4)}(m)=\big\{ 
   \exists t_1,t_2\,\mbox{with}\, m\leq t_1<t_2<m+T,\, Z^\ve_v(t_1)
   \in U_\eta(\cO_{i}),\, Z^\ve_v(t_2)\in U_0\big\}$,\, and
   $A^{(5)}(k)=\big\{ \exists t_3,t_4\,\mbox{with}\, 
   k\leq t_3<t_4<k+T,\, Z^\ve_v(t_3)\in U_0,\, Z^\ve_v(t_4)\in\Gam\big\}$.
   Observe that $G\setminus\big(U_\eta(\cO_{i})\cup U_\del(\cO_{\io(i)})\big)$
   satisfies conditions of Proposition \ref{prop2.5.1}(i) with arbitrarily 
   small $A_1$, so similarly to (\ref{2.5.12}) we can estimate
   \begin{equation}\label{2.6.5}
   P(A^{(1)}(n))\leq\exp(-\frac 12e^{\be/\ve}).
   \end{equation}
   Similarly to the proof of Proposition \ref{prop2.5.1}(ii) we obtain also 
   that
   \begin{equation}\label{2.6.6}
   \max\big(P(A^{(2)}(k)),P(A^{(3)}(k))\big)\leq 
   e^{-(R_\partial^{(i)}+3a)/\ve}
   \end{equation}
   where we, first, choose $\eta$ small and then $T$ large enough. 
   
   Next, relying on the Markov property \index{Markov property} 
   and the arguments similar
   to the proof of Proposition \ref{prop2.5.1}(ii) we estimate
    \begin{equation}\label{2.6.7}
    P\big(A^{(4)}(m)\cap A^{(5)}(k)\big)\leq 
    e^{-(R^{(i)}_\partial+3a)/\ve}
    \end{equation}
    provided $m\leq k-2T$ and $\ve$ is small enough. Summing in $m,\,k$ and 
    $n$ we obtain from (\ref{2.6.4})--(\ref{2.6.7}) that
    for a small $\be$ and all sufficiently small $\ve$,
     \begin{equation}\label{2.6.8}
     P\big\{ \tau^\ve_v(G)\leq e^{(R^{(i)}_\partial+a)/\ve},\,
     Z^\ve_v(\tau^\ve_v(G))\in\Gam\big\}\leq e^{-a/\ve}.
     \end{equation}
     Employing Proposition \ref{prop2.5.1}(i) we derive that
       \[
        P\big\{ \tau^\ve_v(G)>e^{(R^{(i)}_\partial+a)/\ve}\big\}
        \leq e^{-\la /\ve}
     \]
     for some $\la>0$ and all $\ve$ small enough which together with 
     (\ref{2.6.8}) yield (\ref{2.2.35}). 
     
     In order to complete the proof of Theorem \ref{thm2.2.7} it remains to
     derive (\ref{2.2.36}) and (\ref{2.2.37}). Both statements hold true for 
     $n=1$ in view of (\ref{2.2.33}) and (\ref{2.2.35}) and we proceed by
     induction.
    Set
      \[
      H(n,\al)=\big\{ \Sigma^\ve_i(k,-\al)\leq\tau_v(i,k)\leq
      \Sigma^\ve_i(k,\al)\,\,\,\forall k\leq n\big\}
      \]
      and
      \[
      G(n)=\big\{ Z^\ve_v(\tau_v(i,k))\in V_{\io_k(i)}
      \,\,\,\forall k\leq n\big\}.
      \]
      As the induction hypotesis we assume that for any $\al>0$ there exist
      $\la(\al)>0$ and $\la>0$ such that for all small $\ve$,
      \begin{equation}\label{2.6.9}
      P\big(H(n,\al)\big)\geq 1-ne^{-\la(\al)/\ve}\,\,\mbox{and}\,\,
      m\big(G(n)\big)\geq 1-ne^{-\la/\ve}.
      \end{equation}
      By (\ref{2.2.35}) and the strong Markov property 
      \index{strong Markov property}
      \begin{eqnarray}\label{2.6.10}
      &P\big(G(n)\setminus G(n+1)\big)=P\big(\big\{ Z^\ve_v(\tau_v(i,n+1))
      \not\in V_{\iota_{n+1}(i)}\big\}\cap G(n)\big)\\
      &E\bbI_{Z^\ve_v(\tau_v(i,n))\in\partial U_\del(\cO_{\iota_n(i)})}
      P\big\{ Z^\ve_{Z^\ve_v(\tau_v(i,n))}(\tau_{Z^\ve_v(\tau_v(i,n))}(\iota_n
      (i))\not\in V_{\iota_{n+1}(i)}\big\}\leq e^{-\la/\ve}\nonumber
      \end{eqnarray}
      which implies (\ref{2.2.37}). Similarly, by (\ref{2.2.33}) and the strong
      Markov property
      \begin{eqnarray}\label{2.6.11}
      &P\big((H(n,\al)\setminus H(n+1,\al))\cap G(n)\big)\leq
      E\bbI_{Z^\ve_v(\tau_v(i,n))\in\partial U_\del(\cO_{\iota_n(i)})}\\
      &\times P\big\{\tau_{Z^\ve_v(\tau_v(i,n))}(\iota_n(i))>
      \Sigma^\ve_i(n,\al)\,
      \mbox{or}\,\tau_{Z^\ve_v(\tau_v(i,n))}(\iota_n(i))<\Sigma^\ve_i(n,-\al)
      \big\}< e^{-\la(\al)/\ve}\nonumber
      \end{eqnarray}
      proving (\ref{2.2.36}) and completing the proof of Theorem \ref{thm2.2.7}.
      \qed
      
      Finally, we prove Theorem \ref{thm2.2.8} employing the arguments similar
      to \S 2 and \S 3 in Ch. 6 of \cite{FW}. Namely, in order to obtain the
      upper bound in (\ref{2.2.38}) observe that for any $h>0$ there are
      $\rho_0,\del_0>0$ such that if $\rho<\rho_0,\,\del<\del_0$ and a
      curve $\gam\in C_{0t}$ satisfies $\gam_0\in\partial U_\del(\cO_i)$
      and dist$(\gam_t,\partial U_\del(\cO_j))<\rho$ then $S_{0t}(\gam)\geq
      R_{ij}-h$. Using Lemma \ref{lem2.4.4} and the upper bound of large
      deviations (\ref{2.2.11}) we can choose $t=T_1$ such that for all small
       $\ve$ and any $v\in\cup_{1\leq j\leq\ell}\partial U_{2\del}(\cO_j)$,
       \begin{equation}\label{2.6.12}
       P\big\{\sig_v^{\ve,\del}(1)>T_1\big\}\leq e^{-\frac {R_{ij}}\ve}.
       \end{equation}
       Any path of $Z^\ve$ starting at a point of $\partial U_{2\del}(\cO_i)$ and
       reaching $\partial U_\del(\cO_j)$ at time $\sig_v^{\ve,\del}(1)$
       either spends the time $T_1$ without touching the set $\cup_{1\leq
       k\leq\ell}\partial U_\del(\cO_k)$ or arrives at $\Gam_j$ during the time 
       $T_1$. In the latter case $\bfr_{0T_1}\big(Z^\ve,\Psi_{0T_1}^{R_{ij}-h}
       (x)\big)\geq\rho$ and by (\ref{2.2.11}) and (\ref{2.6.12}) for any 
       $v=(x,y)$ with
       $x\in\partial U_{2\del}(\cO_i)$, all $\ve$ small enough and $j\ne i$,
       \begin{eqnarray}\label{2.6.13}
       &P\big\{Z^\ve_v(\sig_v^{\ve,\del}(1))\in\partial U_\del(\cO_j)\big\}
       \leq P\big\{\sig_v^{\ve,\del} (1)>T_1\big\}\\
       &+P\big\{\bfr_{0T_1}(Z^\ve,\Psi_{0T_1}^{R_{ij}-h}(x))\geq\rho\big\}
       \leq\exp\big(-(R_{ij}-h-\tilde\be)/\ve\big)\nonumber
       \end{eqnarray}
       for some $\tilde\be>0$ independent of $\ve$. Any path of $Z^\ve_v$
       starting at $x\in\partial U_\del(\cO_i)$ and reaching $\partial U_\del
       (\cO_j)$ at the time $\sig_v^{\ve,\del}(1)$ must first hit at time
       $\hat\sig_v^{\ve,\del}(1)$ the set $\partial U_{2\del}(\cO_i)$, and so
       (\ref{2.6.13}) together with the Markov property yields the upper bound
       in (\ref{2.2.38}).
       
       In order to derive the lower bound in (\ref{2.2.38}) observe that using
       the definition of $S$-attractors and Lemma \ref{lem2.4.3} (similarly to
       the proof of Lemma \ref{lem1.6.5}(ii) in Part \ref{part1} and see also
       \S 2 in Ch. 6 of \cite{FW}) we conclude 
       that for any $h>0$ there exists $\del_0$ such that if $\del<\del_0$ then
       for any $v=(x,y)\in\Gam_i$ there exists a curve $\gam\in C_{0t}$ such
       that $\gam_0=x,\gam_s\in U_{2\del}(\cO_i)$ for $s\in[0,s_1]$, 
       $\gam_s\not\in
       \cup_{1\leq k\leq\ell,k\ne j}\partial U_\del(\cO_k)$ for $s>s_1$, 
       $\gam_t\in U_{\del/2}(\cO_j)$ and, finally, $S_{0t}(\gam)\leq R_{ij}+h$.
       Then by (\ref{2.2.10}) for all small $\ve>0$,
       \begin{equation}\label{2.6.14}
       P\big\{Z^\ve_v(\sig_v^{\ve,\del}(1))\in\partial U_\del(\cO_j)\big\}
       \geq P\big\{\bfr_{0t}(Z^\ve,
       \gam)<\del/2\big\}\geq\exp\big(-(R_{ij}+h+\tilde\be)/\ve\big)
       \end{equation}
       for some $\tilde\be>0$ independent of $\ve$ which together with
       (\ref{2.6.13}) yields (\ref{2.2.38}).
       
       Now, (\ref{2.2.39}) follows from (\ref{2.2.38}) and the estimates for
       invariant measures of Markov chains from \S 3, Ch. 6 in \cite{FW}. \qed

       \section{Averaging in difference equations}\label{sec2.7}
\setcounter{equation}{0}

Theorem \ref{thm2.2.10} follows by a slight modification (essentially, by
simplification) of the proof of Theorems \ref{thm2.2.2}, in particular, the
standard Gronwall inequality required in the proof of Lemma \ref{lem2.3.4} 
should be replaced by its discrete time version from \cite{El}. We have also 
to check that (\ref{2.2.41}) holds true here which is easier to do than in the
 continuous time case. Indeed,
\begin{eqnarray*}
&Q_k^\ve(x',x,y)=E\exp\langle\be,\sum_{j=1}^kB(x',Y^\ve_{x,y}(j))\rangle
=\int_\bfM\cdots\int_\bfM dm(y_1)p^x(y,y_1)\\
&\times\exp\langle\be,B(x',y_1)\rangle
 dm(y_2)p^{x_1}(y_1,y_2)\exp\langle\be,B(x',y_2)\rangle\\
 &\times\cdots\times 
 dm(y_k)p^{x_{k-1}}(y_{k-1},y_k)\exp\langle\be,B(x',y_k)\rangle
\end{eqnarray*}
where $x_{k+1}=x_k+\ve B(x_k,y_k)$, $k=0,1,...,k-1$, $x_0=x,y_0=y$. By 
(\ref{2.2.1}),
\[
|x_j-x|\leq Kk\ve\,\,\mbox{for}\,\, j=1,...,k.
\]
Since all $y_j\in\bfM$ which is compact and $|x_j-x|\leq KT$, i.e. all $x_j$
stay also in a compact set, we obtain from our assumptions on transition
densities that for all $y,\tilde y\in\bfM$ and $j=0,1,...,k-1$,
\begin{equation*}
1-CKk\ve\leq\frac {p^{x_j}(y,\tilde y)}{p^x(y,\tilde y)}\leq 1+CKk\ve
\end{equation*}
for some $C>0$ independent of $\ve,k,y,\tilde y$ and $x$ staying in a compact
set. Hence,
\begin{equation*}
(1-CKk\ve)^k\leq\frac {Q^\ve_k(x',x,y)}{Q_k(x',x,y)}\leq (1+CKk\ve)^k
\end{equation*}
where $Q_k(x',x,y)=Q_k^0(x',x,y)$ is obtained from $Q_k^\ve(x',x,y)$ by
replacing $Y^\ve_{x,y}$ in the latter by $Y^0_{x,y}=Y_{x,y}$. It follows
from standard facts on principal eigenvalues of positive operators 
(see, for instance, \cite{KLS} and \cite{HH}) that uniformly in $y\in\bfM$ 
and $x,x'\in\bar\cX$ the limit
\[
\lim_{k\to\infty}\frac 1k\log Q_k((x',x,y)=H(x,x',\be),
\]
exists and it satisfies the conditions of Assumption \ref{ass2.2.9},
and so taking the logarithm in the ineguality above and dividing by $k$
we arrive at (\ref{2.2.41}).

Theorem \ref{thm2.2.12} also follows by a slight modification of proofs of
Theorems \ref{thm2.2.5} and \ref{thm2.2.7}, only we have to derive a result
which replaces Lemma \ref{lem2.4.1} providing required properties of 
$I$-functionals given by (\ref{2.2.44}). Since, without loss of generality,
we can assume that $C^{-1}\leq p^x(y,v)\leq C$ for some $C>0$ and by 
(\ref{2.2.44}),
\[
I_x(\mu)=\sup_{u>0}\int_{\bfM}\log\frac {u(y)}{\int_\bfM p^x(y,v)u(v)dm(v)}
d\mu(y),
\]
where the supremum is taken over positive continuous functions $u$. Then
\[
\sup_{u>0}\int_{\bfM}\log\frac {C^{-1}u(y)}{\int_\bfM u(v)dm(v)}d\mu(y)\leq
I_x(\mu)\leq\sup_{u>0}\int_{\bfM}\log\frac {Cu(y)}{\int_\bfM u(v)dm(v)}d\mu(y).
\]
It is easy to see from here that $I_x(\mu)<\infty$ if and only if
$d\mu(y)=g(y)dm(y)$ and the density $g$ is bounded. Hence, in this case,
\begin{eqnarray*}
&I_x(\mu)\leq\sup_{u>0}\int_\bfM\sup g\log\frac {Cu(y)}{\int_\bfM u(v)dm(v)}
dm(y)\leq\sup g\big(\log C\\
&+\sup_{u>0}(\int_\bfM\log u(y)dm(y)-\log\int_\bfM
u(y)dm(y))\big)\leq\sup g\log C.
\end{eqnarray*}
Since
\[
|p^z(y,v)-p^x(y,v)|\leq D|x-z|\leq CD|x-z|\min(p^z(y,v),p^x(y,v))
\]
for some $D>0$, we obtain
\begin{eqnarray*}
&|I_x(\mu)-I_z(\mu)|\leq\sup_{u>0}\int_\bfM\big\vert\log\big(\frac
{\int_\bfM p^z(y,v)u(v)dm(v)}{\int_\bfM p^x(y,v)u(v)dm(v)}\big)\big\vert 
d\mu(y)\\
&\leq\log(1+CD|x-z|).
\end{eqnarray*}
Two last inequalities provide all properties of  $I$-fuctionals which are
 needed in order to replace Lemma \ref{2.4.1} and to proceed with arguments 
 of Sections \ref{sec2.4}--\ref{sec2.6} in the discrete time case.

Theorem \ref{thm2.2.10}
provides, in particular, an approximation of the slow motion by the averaged 
one in probability but, in general, we do not have convergence in (\ref{2.1.4})
also with probability one (see \cite{BK1}). Sometimes, we can derive
this almost sure convergence from the upper large deviations bound estimating
the derivative in $\ve$ of the slow motion as in the following example.
Let $B(x,y)$ be a bounded 1-periodic in $y$ function on $\bbR^1\times
\bbR^1$ with bounded derivatives and let $\xi_1,\xi_2,...$ be a sequence
of independent identically distributed (i.i.d.) random variables. Define
recursively
\begin{eqnarray}\label{2.7.1}
&X^\ve_v(n+1)=X^\ve_v(n)+\ve B\big(X^\ve_v(n),Y^\ve_v(n)\big),\\
&Y^\ve_v(n+1)=Y^\ve_v(n)+X^\ve_v(n)+\xi_{n+1}\nonumber
\end{eqnarray}
where $v=(z,w)$ and $X^\ve_v(0)=z$, $Y^\ve_v(0)=w$. Then
\begin{eqnarray}\label{2.7.2}
&\frac {dX^\ve_v(n+1)}{d\ve}=\frac {dX^\ve_v(n)}{d\ve}+
B\big(X^\ve_v(n),Y^\ve_v(n)\big)\\
&+\ve\frac {\partial B\big(X^\ve_v(n),Y^\ve_v(n)\big)}{\partial x}
\frac {dX^\ve_v(n)}{d\ve}+\ve\frac {\partial B\big(X^\ve_v(n),Y^\ve_v(n)\big)}
{\partial y}\frac {dY^\ve_v(n)}{d\ve},\nonumber\\
&\frac {dY^\ve_v(n+1)}{d\ve}=\frac {dY^\ve_v(n)}{d\ve}+
\frac {dX^\ve_v(n)}{d\ve}.\nonumber
\end{eqnarray}
Set
\[
A^\ve_v(n)=\left(\begin{matrix}
\frac {\partial B\big(X^\ve_v(n),Y^\ve_v(n)\big)}{\partial x}&
\frac {\partial B\big(X^\ve_v(n),Y^\ve_v(n)\big)}{\partial y}\\
0&0\\
\end{matrix}\right)\,\,\mbox{and}\,\, q^\ve_v(n)=\left(\begin{matrix}
B\big(X^\ve_v(n),Y^\ve_v(n)\big)\\0\end{matrix}\right)
\]
which are sequences of bounded matrices and vectors. Taking into account
the equalities
\[
\frac {dX^\ve_v(0)}{d\ve}=\frac {dY^\ve_v(0)}{d\ve}=0\,\,\mbox{and}\,\,
\left(\begin{matrix}1&0\\1&1\end{matrix}\right)^k=
\left(\begin{matrix}1&0\\k&1\end{matrix}\right)
\]
we obtain from (\ref{2.7.2}) by induction (with the agreement $\prod^{n-1}_{j=n}
=1$) that
\begin{eqnarray*}
&\left(\begin{matrix}\frac {dX^\ve_v(n)}{d\ve}\\\frac {dY^\ve_v(n)}{d\ve}
\end{matrix}\right)=\sum_{k=0}^{n-1}\prod_{j=k+1}^{n-1}\big(\left(
\begin{matrix}1&0\\1&1\end{matrix}\right)+\ve A^\ve_v(j)\big)q^\ve_v(k)\\
&=\sum_{k=0}^{n-1}\big(\sum_{l=0}^{n-k-2}\left(\begin{matrix}1&0\\l&1
\end{matrix}\right)\\
&\times\ve^{n-k-l-2}\sum_{k+1\leq j_1<...<j_{n-k-l-2}\leq n-1}
\prod_{i=1}^{n-k-l-2}A^\ve_v(j_i)\big)q^\ve_v(k).
\end{eqnarray*}
Since $A^\ve_v(j)$ and $q^\ve_v(k)$ are bounded we obtain that
\begin{equation}\label{2.7.3}
\big\vert\frac {dX^\ve_v(n)}{d\ve}\big\vert\leq Cn\sum_{k=0}^{n-1}
(1+C\ve)^{n-k-2}\leq n\ve^{-1}(1+C\ve)^n
\end{equation}
for some $C>0$ independent of $n$ and $\ve$. 

Since $B(x,y)$ is 1-periodic in $y$ we can replace the second equality
in (\ref{2.7.1}) by
\begin{equation}\label{2.7.4}
Y^\ve_v(n+1)=Y^\ve_v(n)+X^\ve_v(n)+\xi_{n+1}\,\,(\mbox{mod 1}),
\end{equation}
i.e. we consider now $Y^\ve_v(n)$ evolving on the interval $[0,1]$ with
0 and 1 identified which makes it the circle of radius $1/2\pi$. Suppose
that the distribution of $\xi_1$ has a $C^1$ density $p(y)$ with respect 
to the Lebesgue measure which is positive on $[0,1]$. Now we have the 
family of Markov chains $Y_{x,y}(n),\, n\geq 0$ with transition probabilities
\begin{equation}\label{2.7.5}
P_x(y,\Gam)=P\{ Y_{x,y}(1)\in\Gam\}=P\{ x+y+\xi_1\, (\mbox{mod 1})\in\Gam\}
=\int_\Gam p(z-x-y)dz.
\end{equation}
Thus we are in the framework of our main model satisfying Assumption 
\ref{ass2.2.9}, and so the assertion of Theorem \ref{thm2.2.10} holds true.
Let $\mu^x$ be the invariant measure of the Markov chain $Y_x$ (which is
unique since the Doeblin condition is satisfied here) and assume that
\begin{equation}\label{2.7.6}
\int B(x,y)d\mu^x(y)=0\quad\mbox{for all}\,\, x
\end{equation}
which is, essentially, not a restriction since we always can consider
$B(x,y)-\int B(x,y)d\mu^x(y)$ in place of $B(x,y)$. This means that
$\bar X^\ve_x(n)\equiv x$ and we derive from Theorem \ref{thm2.2.10} that
for any $\del>0$ there exists $\al(\del)>0$ such that for all small $\ve$,
\begin{equation}\label{2.7.7}
P\{\max_{0\leq n\leq T/\ve}|X^\ve_{x,y}(n)-x|\geq\del\}\leq e^{-\al(\del)/\ve}.
\end{equation}
Set $\ve_k=\al(\del)/2\ln k$ then $e^{-\al(\del)/\ve_k}=k^{-2}$ and by the
Borel--Cantelli lemma we obtain that there exists $k_\del(\om)$ finite with
probability one so that for all $k\geq k_\del(\om)$,
\begin{equation}\label{2.7.8}
\max_{0\leq n\leq T/\ve_k}|X^{\ve_k}_{x,y}(n)-x|<\del.
\end{equation}
By (\ref{2.7.3}) for $\ve_{k+1}<\ve\leq\ve_k$ and $k\geq 2$,
\begin{eqnarray}\label{2.7.9}
&\,\,\,\,\,\,\,\,\,\max_{0\leq n\leq T/\ve_{k+1}}|X^{\ve_k}_{x,y}(n)-
X^\ve_{x,y}(n)|\leq
T\ve_{k+1}^{-2}(1+C\ve_k)^{T/\ve_{k+1}}(\ve_k-\ve_{k+1})\\
&\leq 2Te^{2CT}(\al(\del))^{-1}\ln(1+\frac 1k)\longrightarrow 0\,\,\mbox{as}
\,\, k\to\infty.\nonumber
\end{eqnarray}
It follows that with probability one,
\begin{equation}\label{2.7.10}
\max_{0\leq n\leq T/\ve}|X^\ve_{x,y}(n)-x|\to 0\,\,\mbox{as}\,\, k\to\infty.
\end{equation}
The conditions above can be relaxed a bit but this method will not already
work if, for instance, the second equality in (\ref{2.7.1}) is replaced by
\[
Y^\ve_v(n+1)=2Y^\ve_v(n)+X^\ve_v(n)+\del\xi_{n+1}
\]
since in this case the derivative $\frac {dX^\ve_v(n)}{d\ve}$ may grow
exponentially in $n$ and, indeed, we show in \cite{BK1} that for the latter
example there is no convergence with probability one in (\ref{2.7.10})
provided $\del>0$ is small enough.

Next, we exhibit two examples of computations which demonstrate adiabatic
transitions between attractors of the averaged system via the statistics of
proportions of time the slow motion spends in basins of different attractors.
The fast motions $Y^\ve_v$ in both examples are given by the second equation
in (\ref{2.7.1}) where $\xi_1,\xi_2,...$ are i.i.d. random variables 
with the uniform distribution on $[0,1]$. The slow motion $X^\ve_v$ is given
 by the first equation in (\ref{2.7.1}) where in the first example 
 \[
 B(x,y)=B_1(x,y)=x(x^2-4)(1-x^2)+50\sin 2\pi y
 \]
 and in the second example
 \[
 B(x,y)=B_2(x,y)=x(x^2-4)(1-x)(1.5+x)+50\sin 2\pi y.
 \]
 The Markov chains $Y_x$ preserve here the Lebesgue measure on $[0,1]$
 which is the unique invariant measure for them, and so the averaged
 equation (\ref{2.1.6}) for $\bar Z(t)=\bar X^\ve(t/\ve)$ has the right
 hand side $\bar B(x)=\bar B_1(x)=x(x^2-4)(1-x^2)$ in the first case
 and, $\bar B(x)=\bar B_2(x)=x(x^2-4)(1-x)(1.5+x)$ in the second case.
 The one dimensional vector field $\brB(x)$ has three attracting fixed
 points $\cO_1=2,\cO_2=0,\cO_3=-2$ and two repelling fixed points 1 and $-1$,
 while $\bar B_2(x)$ has the same attracting fixed points but one repelling
 fixed point moves now from $-1$ to $-3/2$ making the basin of $-2$ smaller
 which makes it easier for the slow motion to escape from there. It is easy
  to see that $B_1$ and $B_2$ are complete at the fixed points of the averaged 
  system, and so Theorem \ref{thm2.2.12} is applicable in this situation.
 According to the corresponding part of Theorem \ref{thm2.2.12} the transitions
 between $\cO_1,\cO_2,$ and $\cO_3$ are determined by  $R_{ij},\, i,j=1,2,3$ 
 which are obtained via the functionals $S_{0t}(\gam)$ given by (\ref{2.2.9}) 
 but even here these functionals are not easy to compute. The functionals
  $S_{0t}(\gam)$ yield non classical variational problems and the effective
 ways of their computation remain for further research . 
 
 In the first example
 we plot above the histogram with $10^4$ intervals of a single orbit of 
 the slow motion $X^\ve_{x,y}(n),\, n=0,1,2,...,10^8$ with $\ve=10^{-3}$ 
  and the initial values $x=0,\, y=0$.
  The histogram shows that most of the points of the orbit stay near the
  attractors $\cO_1,\,\cO_2$ and $\cO_3$
  and $X^\ve_{x,y}(n)$ hops between basins of attraction of these points. 
  The form of the histogram indicates the equality $R_{21}=R_{23}$, which
  follows also by the symmetry considerations, but in this case Theorem 
  \ref{thm2.2.12} cannot specify whether the slow motion mostly exits from
  the basin of $\cO_2$ to the basin of $\cO_1$ or to the basin of $\cO_3$.

  \begin{figure} 
\centerline{\epsfig{file=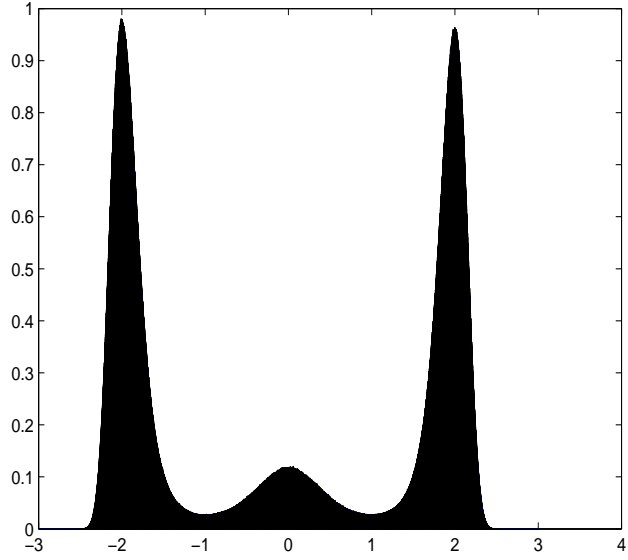, width=10cm, height=8.5cm}}
\caption{Symmetrical basins case}
\label{fig2.8.1:fig}
\end{figure}

In the second example the basin of attraction of $-2$ becomes smaller while 
 the left interval of the basin of attraction of $0$ becomes larger. The latter
leads to the inequality $R_{23}>R_{21}$ which according to Theorem 
\ref{thm2.2.12} makes it more difficult for the slow motion to exit to 
the left from the basin of $\cO_2$ than to the right. In the histogram below
(which has again $10^4$ intervals)
we plot $X^\ve_{x,y}(n),\, n=0,1,2,...,10^8$ with $\ve=10^{-3}$ 
and the initial values $x=-2,\, y=0$. In compliance with Theorem \ref{thm2.2.12} 
the histogram demonstrates that the slow
motion leaves the basin of $\cO_3$ and after arriving at the basin of $\cO_2$
it exits mostly to the basin of $\cO_1$, and so the slow motion hops mostly
between basins of $\cO_1$ and $\cO_2$ staying most of the time in small
neighborhoods of these points. Still, a complete rigorous explanation of these
histograms even for our simple examples requires nontrivial additional
arguments. It is interesting to observe that these 
histograms have the same form as in Section \ref{sec1.9} of Part \ref{part1}
where randomness is generated by the expanding (chaotic) map $y\,\to\, 3y$ 
instead of adding uniformly distributed random variables as we do it here.

\begin{figure}
\centerline{\epsfig{file=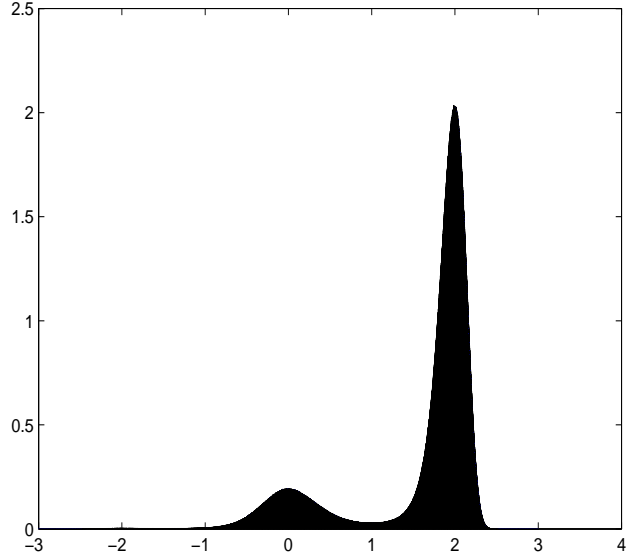, width=10cm, height=8.5cm}}
\caption{Asymmetrical basins case}
\label{fig2.8.2:fig}
\end{figure}

\section{Extensions: stochastic resonance}\label{sec2.8}\setcounter{equation}{0}

The scheme for the stochastic resonance type phenomenon described below is a
slight modification of the model 
suggested by M.Freidlin (cf. \cite{Fre+}) and it can be demonstrated in the 
setup of three scale systems
 \begin{eqnarray}\label{2.8.1}   
&   \frac {dV^{\ve,\del}(t)}{dt}=\del\ve A(V^{\ve,\del}(t),X^{\ve,\del}(t),
Y^{\ve,\del}(t))\nonumber\\
&\frac {dX^{\ve,\del}(t)}{dt}=\ve B(V^{\ve,\del}(t),X^{\ve,\del}(t),
Y^{\ve,\del}(t))\\
&dY^{\ve,\del}(t)=\sig(V^{\ve,\del}(t),X^{\ve,\del}(t),Y^{\ve,\del}(t))dw_t+
b(V^{\ve,\del}(t),X^{\ve,\del}(t),Y^{\ve,\del}(t))dt,\nonumber
\end{eqnarray}
$V^{\ve,\del}=V^{\ve,\del}_{v,x,y}$, $X^{\ve,\del}=X^{\ve,\del}_{v,x,y}$, 
$Y^{\ve,\del}=Y^{\ve,\del}_{v,x,y}$ with initial conditions $V^{\ve,\del}(0)=v$,
$X^{\ve,\del}(0)=x$ and $Y^{\ve,\del}(0)=y$ and the last equation in 
(\ref{2.8.1})
is a stochastic differential equation coupled with first two ordinary 
differential equations though together they should be considered as a system
of stochastic differential equations (with a degeneration in the first two). 
We assume that $V^{\ve,\del}\in\bbR^l$, $X^{\ve,\del}\in\bbR^d$ while
 $Y^{\ve,\del}$ evolves on a compact $n_\bfM$-dimensional $C^2$ Riemannian 
 manifold $\bfM$ and the coefficients $A$, $B$, $b$ are bounded 
 smooth vector fields on $\bbR^l$, $\bbR^d$ and $\bfM$, respectively,
 depending on other variables as parameters. We suppose also that $a=\sig\sig^*$
 is a uniformly positive definite smooth matrix field on $\bfM$. In the same
 way as in Section \ref{sec2.2} we can generalize the setup taking $Y^{\ve,\del}$
  to be random evolutions but in order to simplify the notations we restrict 
  ourselves to fast motions $Y^{\ve,\del}$ being diffusions. The solution of
 (\ref{2.8.1}) determines a Markov diffusion process which the triple
 $(V^{\ve,\del},X^{\ve,\del},Y^{\ve,\del})$. Taking $\ve=\del=0$ we arrive at
  the (unperturbed) process $(v,x,Y^{0,0}_{v,x,y})$ where $Y=Y_{v,x,y}=
  Y^{0,0}_{v,x,y}$ solves the unperturbed stochastic differential equation
 \begin{equation}\label{2.8.2}
 dY(t)=\sig(v,x,Y(t))dw_t+b(v,x,Y(t))dt.
 \end{equation} 
 It is natural to view the diffusion $Y(t)$ as describing an idealized
  physical system where parameters $v=(v_1,...,v_l)$ and $x=(x_1,...,x_d)$ are
  assumed to be constants of motion while the perturbed process 
  $(V^{\ve,\del},X^{\ve,\del},Y^{\ve,\del})$ is regarded as describing a real 
  system where evolution of these parameters is also taken into consideration 
  but unlike the averaging setup (\ref{2.1.1}) we have now two sets of 
  parameters moving with very different speeds.
  
  Let $\mu_{v,x}$ be the unique invariant measure of the diffusion $Y(t)=
  Y_{v,x,y}(t)$. Set
  \begin{equation}\label{2.8.3}
\bar B_v(x)=\bar B(v,x)=\int B(v,x,y)d\mu_{v,x}(y)
\end{equation}
and let $\bar X^{(v)}$ be the solution of the averaged equation
\begin{equation}\label{2.8.4} 
\frac {d\brX^{(v)}(t)}{dt}=\brB_v(\brX^{(v)}(t)).
\end{equation}
First, we apply averaging and large deviations estimates in averaging from
the previous sections to two last equations in (\ref{2.8.1}) freezing the 
slowest
variable $v$ (i.e. taking for a moment $\del=0$). Namely, set 
$\hat X^\ve(t)=X^{\ve,0}_{v,x,y}(t/\ve)$ and $\hat Y^\ve(t)=Y^{\ve,0}_{v,x,y}
(t/\ve)$ 
so that
\begin{eqnarray}\label{2.8.5}
  &\frac {d\hat X^\ve(t)}{dt}= \ve B(v,\hat X^\ve(t),\hat Y^\ve(t))\\
  &\frac {d\hat Y^\ve(t)}{dt}=\sig(v,\hat X^\ve(t),Y^\ve(t))dw_t+
  b(v,\hat X^\ve(t),\hat Y^\ve(t))dt.
  \nonumber\end{eqnarray}
Suppose for simplicity that $l=d=1$ (i.e. both $V^{\ve,\del}$ and
 $X^{\ve,\del}$ are one dimensional) and that the solution $\brX^{(v)}(t)$
 of (\ref{2.8.4}) has the limit set consisting of two attracting points
 $\cO_1$ and $\cO_2$, which for simplicity we assume to be independent
 of $v$, and a repelling fixed point $\cO_0^v$ depending on $v$ and separating
 their basins. As an example of $\brB$ we may have in mind
  $\brB_v(x)=(x-v)(1-x^2), -1<v<1$. Let $S^v_{0T}(\gam)$ be the 
 large deviations rate functional for the system of last two equations in
 (\ref{2.8.1}) defined in (\ref{2.2.9}) and set for $i,j=1,2$,
 \begin{equation}\label{2.8.6}
 R_{ij}(v)=\inf\{ S^v_{0T}(\gam):\,\gam\in C_{0T},\,\gam_0=\cO_i,\,\gam_T
 =\cO_j,\, T\geq 0\}
 \end{equation}
 (cf. with $R_{ij}$ in Theorem \ref{thm2.2.7}).
 
  Set
 \begin{equation}\label{2.8.7}
 \bar A_i(v)=\int A(v,\cO_i,y)d\mu_{v,\cO_i}(y)
 \end{equation}
 and assume that for all $v$,
  \begin{equation}\label{2.8.8}
  \bar A_1(v)<0\quad\mbox{and}\quad\bar A_2(v)>0
  \end{equation}
  which means in view of the averaging principle (see Theorem \ref{thm2.2.2}
  and the following it discussion) that $V^{\ve,\del}_{v,x,y}(t)$ decreases 
  (increases) with high probability while $X^{\ve,\del}_{v,x,y}(t)$ stays close to $\cO_1$ 
  (to $\cO_2$).
  
  The following statement suggests a "nearly" periodic behavior of the
  slowest motion.
  \begin{conjecture}\label{conj2.8.1} Suppose that there exist strictly 
  increasing and decreasing functions $v_{-}(r)$ and $v_{+}(r)$, respectively,
   so that
  \[
  R_{12}(v_{-}(r))=R_{21}(v_{+}(r))=r
  \]
  and $v_{-}(\la)=v_{+}(\la)=v^*$ for some $\la>0$ while
   $v_{-}(r)<v^*<v_{+}(r)$
  for $r<\la$. Assume that $\del\to 0$ and $\ve\to 0$ in such a way that 
  \begin{equation}\label{2.8.9}
  \lim_{\ve,\del\to 0}\ve\ln(\del\ve)=-\rho>-\la.
  \end{equation}
  Then for any $v,x$ there exists $t_0>0$ so that the slowest motion
  $\tilde V^{\ve,\del}_{v,x,y}(t+t_0)=V^{\ve,\del}_{v,x,y}
  \big((t+t_0)/\del\ve\big),\, t\geq 0$ converges in distribution
  (as $\ve,\del\to 0$ so that (\ref{2.8.9}) holds true) to a periodic
  function $\psi(t)$, $\psi(t+T)=\psi(t)$ with 
  \[
  T=T(\rho)=\int_{v_{-}(\rho)}^{v_{+}(\rho)}\frac {dv}{|\bar A_1(v)|}+
  \int_{v_{-}(\rho)}^{v_{+}(\rho)}\frac {dv}{|\bar A_2(v)|}.
  \]
  \end{conjecture}
  
  The argument supporting this conjecture goes as follows.
   Since $V^{\ve,\del}$ moves much slower than $X^{\ve,\del}$ we can freeze 
   the former and in place of (\ref{2.8.1}) we can study, first, (\ref{2.8.5}).
Applying the arguments of Theorem \ref{thm2.2.7} to the pair $\hat X,\hat Y$ 
from (\ref{2.8.5}) we conclude from (\ref{2.2.34}) that the intermediate
   motion $X^{\ve,\del}$ most of the time stays very close to either
   $\cO_1$ or $\cO_2$ before it exits from the corresponding basin, and so
   in view of an appropriate averaging principle (which follows, for
   instance, from Theorem \ref{2.2.2}) on bounded time intervals the slowest 
   motion $V^{\ve,\del}$ mostly stays close to the corresponding 
   averaged motion determined by the vector fields $\bar A_1$ and $\bar A_2$
   given by (\ref{2.8.7}). When $X^{\ve,\del}$ is close to $\cO_1$ 
    the slowest motion $V^{\ve,\del}$ decreases until $v=v_{-}(\rho)$ 
    where $R_{12}(v)=\rho$. In view of (\ref{2.2.33}) and the scaling 
    (\ref{2.8.9}) between $\ve$ and $\del$, a moment later $R_{12}(v)$ becomes 
    less than $\rho$ and $X^{\ve,\rho}$ jumps immediately close to 
    $\cO_2$. There $\bar A_2(v)>0$, and so $V^{\ve,\del}$ starts to 
    grow until it reaches $v=v_{+}(\rho)$ where $R_{21}(v)=\rho$. A moment
    later $R_{21}(v)$ becomes  smaller than $\rho$ and in view of 
    (\ref{2.2.33}) the intermediate motion $X^{\ve,\del}$ jumps immediately 
    close
    to $\cO_1$. This leads to a nearly periodic behavior of $V^{\ve,\del}$.
    In order to make these arguments precise we have to deal here with
    an additional difficulty in comparison with the two scale setup considered
    in previous sections since now the large deviations $S$-functionals from 
    Theorem \ref{thm2.2.2} and the $R$-functions describing adiabatic 
    fluctuations
    and transitions of Theorems \ref{thm2.2.5} and \ref{thm2.2.7} depend on 
    another
    very slowly changing parameter. Still, we can use the technique of 
    Sections \ref{sec2.5} and \ref{sec2.6} above applied on time intervals 
    where changes in the $v$-variable can be neglected should work here
     but the details of this approach have not been worked out yet.  
    
    On the other hand, when the fast motion $Y^{\ve,\del}$ does not depend
     on the slow motions, i.e. when the coefficients $\sig$ and
     $b$ in (\ref{2.8.1}) depend 
     only on the coordinate $y$ (but not on $v$ and $x$), then the above 
     arguments can be made precise without much effort. Indeed, we can obtain
     estimates for transition times $\tau^\ve(1)$ and $\tau^\ve(2)$ of 
     $X^{\ve,\del}(t/\ve)$ between neighborhoods of $\cO_1$ and $\cO_2$ as
     in Theorem \ref{thm2.2.7} applying the latter to $\hat X^\ve$ and
     $\hat Y^\ve$ from (\ref{2.8.5}) with freezed $v$-variable. This is
     possible since the method of Proposition \ref{prop2.5.1} requires us
     to make large deviations estimates, essentially, only for probabilities
      $P\{ kT\leq\tau(i)<(k+1)T\}$, i.e. on bounded time intervals,
      and then combine them with the Markov property arguments. During
      such times the slowest motion $V^{\ve,\del}$ can move only a distance of
      order $\del T$. Thus freezing $v$ and using the Gronwall inequality
      for the equation of $X^\ve$ in order to estimate the resulting error
      we see that the latter is small enough for our purposes. Observe, that
      it would be much more difficult to justify freezing $v$ in the
      coefficients $\sig$ and $b$ of $Y^\ve$, if we allow the latter to depend 
    on $v$, since a strightforward application of the Gronwall inequality there
      would yield an error estimate of an exponential in $1/\ve$ order which 
      is comparable with $1/\del$. Still, it may be possible to take care 
      about the general case using methods of Section \ref{sec2.3} 
       since we produce large deviations estimates there by gluing 
      large deviations estimates on smaller time intervals where the 
      $x$-variable (and so, of course, $v$-variable) can be freezed. Next, set
      \[
      V_{v,y}^{\ve,\del,i}(t)=v+\del\ve\int_0^tA(V_v^{\del,\ve,i}(s),\cO_i,
      Y(s))ds
      \]
      where now $Y$ does not depend on $\ve$ and $\del$. Then by (\ref{2.8.1})
      together with the Gronwall inequality we obtain that
      \[
      |V_{v,x,y}^{\ve,\del}(t)-V_{v,y}^{\ve,\del,i}(t)|\leq K\del\ve
      e^{\del\ve Kt}\int_0^t|X^{\ve,\del}_{v,x,y}(s)-\cO_i|ds
      \]
      where $K$ is the Lipschitz constant of $A$. If $x$ belongs to the 
      basin $\cO_i$ then according to Theorem \ref{thm2.2.7} $\hat X^\ve$,
      and so also $X^{\ve,\del}$, stays most of the time near $\cO_i$ up
      to its exit from the basin of the latter which yields according to
      the above inequality that $V^{\ve,\del}$ stays close to $V^{\ve,\del,i}$
      during this time. But now we can employ the averaging principle for
      the pair $V^{\ve,\del,i}(t),Y(t)$ which sais that $V^{\ve,\del,i}(t)$
      stays close on the time intervals of order $1/\del\ve$ to the averaged
      motion $\bar V^{\ve,\del,i}_v(t)$ defined by
      \[
      \bar V^{\ve,\del,i}_v(t)=v+\int_0^t\bar A_i(\bar V^{\ve,\del,i}_v(s))ds
      \]
      and in view of (\ref{2.8.9}), $\bar V^{\ve,\del,1}_v(t)$ decreases while
      $\bar V^{\ve,\del,2}_v(t)$ increases which leads to the behavior
      described in Conjecture \ref{conj2.8.1}.

    A similar conjecture can be made under the corresponding conditions for
    the discrete time case determined by a three scale difference system of
    equations of the form
    \begin{eqnarray}\label{2.8.10}
  &V^{\ve,\del}(n+1)-V^{\ve,\del}(n)=\ve\del A(V^{\ve,\del}(n),X^{\ve,\del}(n),
  Y^{\ve,\del}(n)),\,\,\,V^{\ve,\del}(0)=v,\nonumber\\
 &\,\,\,\,\,\,\,\, X^{\ve,\del}(n+1)-X^{\ve,\del}(n)=\ve B(V^{\ve,\del}(n),
 X^{\ve,\del}(n),Y^{\ve,\del}(n)),\,\,\, X^{\ve,\del}(0)=x\nonumber
 \end{eqnarray}
 where $A$ and $B$ are smooth vector functions and 
 $Y^{\ve,\del}(n))=Y^{\ve,\del}_{v,x,y}(n),\, Y^{\ve,\del}_{v,x,y}(0)=y$
are coupled with $V^{\ve,\del}(n)$ and $X^{\ve,\del}(n)$ perturbations of a
parametrized by $v$ and $x$ appropriate family of Markov chains having
smooth transition densities similar to those considered in Theorem 
\ref{thm2.2.12}.

\section{Young measures approach to averaging}\label{sec2.9}
\setcounter{equation}{0}

In this section we derive the averaging principle and discuss corresponding
large deviations in the sense of convergence of 
Young measures \index{Young measures} adapted to our probabilistic setup. For
 more detailed information about Young measures we refer the reader to 
 \cite{AG} and references there.
 
 Let $\mu$ belongs to the space $\cP(\bbR^d\times\bfM)$ of probability measures
on $\bbR^d\times\bfM$. We consider a random Young measure $\zeta^\ve$ from 
$([0,T]\times\bbR^d\times\bfM,\ell_T\times\mu)$ to $\cP(\bbR^d\times\bfM)$ 
which we define by
\[
\zeta^\ve(t,x,y)=\del_{X^\ve_{x,y}(t/\ve),Y^\ve_{x,y}(t/\ve)}
\]
where $\ell_T$ is the Lebesgue measure on $[0,T]$, $\del_w$ is the unit mass
at $w$, and $X^\ve,\, Y^\ve$ are the same as in (\ref{2.1.1}).

Suppose that $\mu\in\cP(\bbR^d\times\bfM)$ has a disintegration
\begin{equation}\label{2.9.1}
d\mu(x,y)=d\mu_x(y)d\la(x),\,\,\,\la\in\cP(\bbR^d)
\end{equation}
such that for each Lipschitz continuous function $g$ on $\bfM$ and any
$x,z\in\bbR^d$,
\begin{equation}\label{2.9.2}
|\int gd\mu_x-\int gd\mu_z|\leq K_{L(g)}|x-z|
\end{equation}
for some $K_L>0$ depending only on $L$ where $L(g)$ is both a Lipschitz
constant of $g$ and it also bounds $|g|$. Set
\begin{equation}\label{2.9.3}
\bar B(x)=\int B(x,y)d\mu_x(y)
\end{equation}
and assume that (\ref{2.2.1}) holds true which together with (\ref{2.9.2}) yields
that $\bar B$ is bounded and Lipschitz continuous, and so there exists a 
unique solution $\bar X^\ve(t)=\bar X^\ve_x(t)$ of (\ref{2.1.3}). For any bounded
continuous function $g$ on $\bbR^d\times\bfM$ define
\[
\cE^g_\ve(t,\del)=\big\{ (x,y)\in\bbR^d\times\bfM:\, E\big\vert\frac 1t\int_0^t
g(x,Y^\ve_{x,y}(u))du-\bar g(x)\big\vert >\del\big\}
\]
where $\bar g(x)=\int g(x,y)d\mu_x(y)$.

In the spirit of \cite{AG}) we say that the Young measures $\zeta^\ve$ 
converge as $\ve\to 0$ to a Young measure $\zeta^0$ defined by
\[
\zeta^0(t,x,y)=\del_{\bar Z_x(t)}\times\mu_{\bar Z_x(t)}\in
\cP(\bbR^d\times\bfM), 
\]
 $\bar Z_x(t)=\bar X^\ve_x(t/\ve)$,
if for any bounded continuous function $f(t,x,y)$ on 
$([0,T]\times\bbR^d\times\bfM$,
\[
E\big\vert\int_0^T\big(f(s,X^\ve_{x,y}(s/\ve),Y^\ve_{x,y}(s/\ve))-
\bar f(s,\bar Z_x(s))\big)ds\big\vert\to 0\,\,\mbox{as}\,\,\ve\to 0.
\]
The following result provides a verifiable (in some interesting cases) criterion
for even stronger convergence.
\begin{theorem}\label{thm2.9.1} Let $\mu\in\cP(\bbR^d\times\bfM)$ has the
 disintegration (\ref{2.9.1}) satisfying (\ref{2.9.2}). Then
 \begin{equation}\label{2.9.4}
 \lim_{\ve\to 0}\int_{\bbR^d}\int_\bfM E\sup_{0\leq t\leq T}\big\vert\int_0^t
 \big(f(s,X^\ve_{x,y}(s/\ve),Y^\ve_{x,y}(s/\ve))-\bar f(s,\bar Z_x(s))\big)ds
 \big\vert d\mu(x,y)=0
 \end{equation}
for any bounded continuous function $f=f(t,x,y)$ on $[0,T]\times\bbR^d\times
 \bfM$ where $\bar f(t,x)=\int f(t,x,y)d\mu_x(y)$ if and only if for each
  $N\in\bbN$ and any finite collection $g_1,...,g_N$ of bounded Lipschitz
  continuous functions on $\bbR^d\times\bfM$ there exists an integer valued 
function $n=n(\ve)\to\infty$ as $\ve\to 0$ such that for any $\del>0$ and
$l=1,...,N$,
    \begin{equation}\label{2.9.5}
    \lim_{\ve\to 0}\max_{0\leq j< n(\ve)}\int_{\bbR^d}\int_\bfM P\big\{ 
    \big(X^\ve_{x,y}(jt(\ve)),Y^\ve_{x,y}(jt(\ve))\big)\in\cE^{g_l}_\ve
    (t(\ve),\del)\}d\mu(x,y)=0,
    \end{equation}
    where $t(\ve)=\frac T{\ve n(\ve)}$. 
 \end{theorem}
 \begin{proof} First, we prove that (\ref{2.9.4}) implies (\ref{2.9.5}). Let
 $g_1,...,g_N$ be bounded Lipschitz continuous functions on $\bbR^d\times
 \bfM$ and set
 \begin{equation}\label{2.9.6}
 \rho_{x,y}^{\ve,l}(t)=\ve\int_0^t\big(g_l(X^\ve_{x,y}(s),Y^\ve_{x,y}(s))
 -\bar g_l(\bar X^\ve_x(s))\big)ds.
 \end{equation}
 If
 \[
 \rho_{x,y}^{\ve,l}=\sup_{0\leq t\leq T/\ve}|\rho_{x,y}^{\ve,l}(t)|
 \]
 then by (\ref{2.9.4}) for each $l=1,...,N$,
 \begin{equation}\label{2.9.7}
 \rho^\ve_l=\int_{\bbR^d}\int_\bfM E\rho_{x,y}^{\ve,l}d\mu(x,y)\to 0\,\,\,
 \mbox{as}\,\,\ve\to 0.
 \end{equation}
 Choose an integer valued function $n(\ve)\to\infty$ as $\ve\to 0$ so that
  \begin{equation}\label{2.9.8}
  n(\ve)\max_{1\leq l\leq N}\rho^\ve_l\to 0\,\,\,\mbox{as}\,\,\ve\to 0
 \end{equation}
 and let $t(\ve)=T/\ve n(\ve)$. Set $x_k^\ve=X^\ve_{x,y}(kt(\ve))$, $y_k^\ve=
 Y^\ve_{x,y}(kt(\ve))$ and $\bar x_k^\ve=\bar X^\ve_{x}(kt(\ve))$, $k=0,1,...$.
 Then by (\ref{2.9.6}),
 \begin{equation}\label{2.9.9}
 \rho_{x,y}^{\ve,l}((j+1)t(\ve))-\rho_{x,y}^{\ve,l}(jt(\ve))=
 \ve\int_0^{t(\ve)}\big(g_l(X^\ve_{x^\ve_j,y^\ve_j}(u),Y^\ve_{x^\ve_j,y^\ve_j}
 (u))-\bar g_l(\bar X^\ve_{\bar x_j^\ve}(u))\big)du
 \end{equation}
 where $X^\ve_{x^\ve_j,y^\ve_j}(u)=X^\ve_{x,y}(jt(\ve)+u)$ and
 $Y^\ve_{x^\ve_j,y^\ve_j}(u)=Y^\ve_{x,y}(jt(\ve)+u)$. By (\ref{2.2.1}),
 \begin{eqnarray}\label{2.9.10}
 &\ve\big\vert\int_0^{t(\ve)}\big(g_l(X^\ve_{x^\ve_j,y^\ve_j}(u),
 Y^\ve_{x^\ve_j,y^\ve_j}(u))-g_l(x^\ve_j,Y^\ve_{x^\ve_j,y^\ve_j}(u))\big)du
 \big\vert\\
 &\leq\ve L_l\int_0^{t(\ve)}|X^\ve_{x^\ve_j,y^\ve_j}(u)-x_j^\ve|du\leq
 L_lK(\ve n(\ve))^2\nonumber
 \end{eqnarray}
 where $L_l$ is the Lipschitz constant of $g_l$. Similarly, by (\ref{2.2.1})
 and (\ref{2.9.2}),
 \begin{equation}\label{2.9.11}
 \ve\big\vert\int_0^{t(\ve)}\big(\bar g_l(\bar X^\ve_{\bar x_j^\ve}(u))-
 \bar g_l(\bar x_j^\ve)\big)du\big\vert\leq (L_l+K_{L_l})K(\ve t(\ve))^2
 \end{equation}
 and
 \begin{equation}\label{2.9.12}
 |\bar g_l(\bar x^\ve_j)-\bar g_l(x^\ve_j)|\leq (L_l+K_{L_l})|\bar x_j^\ve-
 x^\ve_j|\leq (L_l+K_{L_l})\rho^\ve_{x,y}.
 \end{equation}
 It follows from (\ref{2.9.9})--(\ref{2.9.12}) that
 \begin{eqnarray}\label{2.9.13}
 &\big\vert\frac 1{t(\ve)}\int_0^{t(\ve)}g_l(x^\ve_j,Y^\ve_{x^\ve_j,y^\ve_j}(u))
 du-\bar g_l(x^\ve_j)\big\vert\\
 &\leq TK(2L_l+K_{L_l})/n(\ve)+(L_l+K_{L_l}+2T^{-1}
 n(\ve))\rho^\ve_{x,y}.\nonumber
 \end{eqnarray}
 Given $\del >0$ choose $\ve_\del>0$ such that for all $\ve\leq\ve_\del$ and
 $l=1,...,N$,
 \[
 TK(2L_l+K_{L_l})/n(\ve)\leq\del/2.
 \]
  Then by (9.13) and the Markov property,
 \begin{eqnarray*}
 &\big\{ (x,y)\in\bbR^d\times\bfM:\,\big( X^\ve_{x,y}(jt(\ve)),
 Y^\ve_{x,y}(jt(\ve))\big)\in\cE_\ve^{g_l}(t(\ve),\del)\subset A_\ve(\del)\\
 &=\big\{ (x,y)\in\bbR^d\times\bfM:\, (L_l+K_{L_l}+2T^{-1}n(\ve))
 E_{x^\ve_j,y^\ve_j}\rho^{\ve,l}_{x,y} >\del/2\big\}
 \end{eqnarray*}
 where 
 \[
 E_{x^\ve_j,y^\ve_j}=E\big(\cdot\big\vert X^\ve_{x,y}(jt(\ve)),
 Y^\ve_{x,y}(jt(\ve))\big)
 \]
 is the conditional expectation. By Chebyshev's inequality
  \begin{equation}\label{2.9.14}
 E\mu(A_\ve(\del))\leq\frac 2{\del}(L_l+K_{L_l}+2T^{-1}n(\ve))E\rho_l^\ve.
 \end{equation}
 By (\ref{2.9.8}) the right hand side of (\ref{2.9.14}) tends to 0 as $\ve\to 0$
 yielding (\ref{2.9.5}).
 
  Next, we derive (\ref{2.9.4}) from (\ref{2.9.5}).
 Since $f$ in (\ref{2.9.4}) is a bounded function and $\la$ is
 a probability measure it is easy to see that it suffices to prove 
 (\ref{2.9.4}) when the integration in $x$ there is restricted to
 compact subsets of $\bbR^d$. But if $x$ belongs to
  a compact set $G\subset\bbR^d$ in view of (\ref{2.1.1}) and (\ref{2.2.1})
  the slow motion $X^\ve_{x,y}(s)$, as well as the averaged one
  $\bar X^\ve_x(s)$, stays during the time $T/\ve$ in a $KT-$neighborhood 
  $G_{KT}$ of $G$. But on $[0,T]\times G_{KT}\times\bfM$ we can approximate
  $f$ uniformly by Lipschitz continuous functions. Thus, in place of 
  (\ref{2.9.4}) it suffices to show that for any compact set $G\subset\bbR^d$ and a bounded Lipschitz
 continuous function $f$ on $[0,T]\times G_{KT}\times\bfM$ with a Lipschitz
 constant $L=L(f)$ in all variables,
 \begin{eqnarray}\label{2.9.15}
 &\lim_{\ve\to 0}\ve\int_G\int_\bfM E\sup_{0\leq t\leq T/\ve}\big\vert\int_0^t
 \big(f(\ve s,X^\ve_{x,y}(s),Y^\ve_{x,y}(s))\\
 &-\bar f(\ve s,\bar X^\ve_x(s))\big)ds\big\vert d\mu_x(y)d\la(x)=0.\nonumber
 \end{eqnarray}
 
  By (\ref{2.2.1}), (\ref{2.9.2}) and (\ref{2.9.3}),
 \begin{eqnarray}\label{2.9.16}
 &\ve\big\vert\int_0^t\big(f(\ve s,X^\ve_{x,y}(s),Y^\ve_{x,y}(s))-
 \bar f(\ve s,\bar X^\ve_x(s))\big)ds\big\vert\\
 &\leq\ve\big\vert\int_0^t\big(f(\ve s,X^\ve_{x,y}(s),Y^\ve_{x,y}(s))-
 \bar f(\ve s,X^\ve_x(s))\big)ds\big\vert\nonumber\\
 &+(L+K_L)T\sup_{0\leq s\leq T/\ve}|X^\ve_x(s)-\bar X^\ve_x(s)|.\nonumber
 \end{eqnarray}
 By (\ref{2.1.1}), (\ref{2.2.1}), (\ref{2.9.2}) and (\ref{2.9.3}),
 \begin{eqnarray*}
 &|X^\ve_{x,y}(t)-\bar X^\ve_x(t)|=\ve\big\vert\int_0^t\big(B(X^\ve_{x,y}(s),
 Y^\ve_{x,y}(s))-\bar B(\bar X^\ve_x(s))\big)ds\big\vert\\
 &\leq\ve\big\vert\int_0^t\big(B(X^\ve_{x,y}(s),Y^\ve_{x,y}(s))-
 \bar B(X^\ve_{x,y}(s))\big)ds\big\vert\\
 &+(K+K_K)\ve\int_0^t|X^\ve_{x,y}(s)-\bar X^\ve_x(s)|ds.
 \end{eqnarray*}
 This together with the Gronwall inequality gives
 \begin{eqnarray}\label{2.9.17}
 &\sup_{0\leq t\leq T/\ve}|X^\ve_{x,y}(t)-\bar X^\ve_x(t)|\\
 &\leq e^{(K+K_K)T}\ve\sup_{0\leq t\leq T/\ve}\big\vert\int_0^t
 \big(B(X^\ve_{x,y}(s),Y^\ve_{x,y}(s))-\bar B(X^\ve_{x,y}(s))\big)ds\big\vert.
 \nonumber\end{eqnarray}
 Now we see that the integral term in the right hand side of (\ref{2.9.17})
 is a particular case of the integral term in the right hand side of 
 (\ref{2.9.16}) with $f=B$, and so it suffices to estimate only the latter.
 
 Set, again, $x_k^\ve=X^\ve_{x,y}(kt(\ve))$, $y_k^\ve=
 Y^\ve_{x,y}(kt(\ve))$, $\bar x_k^\ve=\bar X^\ve_{x}(kt(\ve))$, $k=0,1,...$
 and fix a large $N\in\bbN$. Let $l=[\ve jt(\ve)N/T]=[jN/n(\ve)]$ then by 
 (\ref{2.2.1}), (\ref{2.9.1}) and (\ref{2.9.2}),
 \begin{eqnarray}\label{2.9.18}
 &\quad\,\,\,\,\,\,\,\,\,\,\,\,\,\,\ve\big\vert\int_0^{t(\ve)}\big(f(\ve jt(\ve)
 +\ve u,X^\ve_{x^\ve_j,y^\ve_j}(u),Y^\ve_{x^\ve_j,y^\ve_j}(u))\\
 &-f(lT/N,x^\ve_j,Y^\ve_{x^\ve_j,y^\ve_j}(u))\big)ds\big\vert
 \leq LT^2\big((Nn(\ve))^{-1}+(n(\ve))^{-2}\big)\nonumber\\
 &+L\ve\int_0^{t(\ve)}|X^\ve_{x^\ve_j,y^\ve_j}(u)-x^\ve_j|du
 \leq LT^2/Nn(\ve)+LT^2(1+K))(n(\ve))^{-2}\nonumber
 \end{eqnarray}
 and
 \begin{eqnarray}\label{2.9.19}
 &\ve\big\vert\int_0^{t(\ve)}\big(\bar f(\ve jt(\ve)+\ve u,
 X^\ve_{x^\ve_j,y^\ve_j}(u))-\bar f(lT/N,x_j^\ve))du\big\vert\\
 &\leq LT^2/Nn(\ve)+T^2(L+LK+KK_L)(n(\ve))^{-2}.\nonumber
 \end{eqnarray}
  Now using (\ref{2.9.18}), (\ref{2.9.19}) together with the Markov property
   and assuming that $|f|\leq\hat L_f$ for some constant $\hat L_f>0$ we obtain
 \begin{eqnarray}\label{2.9.20}
 &\,\,\,\,\,\,\,\,\,\ve E\sup_{0\leq t\leq T/\ve}\big\vert\int_0^t\big(f(\ve s,
 X^\ve_{x,y}(s),
 Y^\ve_{x,y}(s))-\bar f(\ve s, X^\ve_{x,y}(s))\big)ds\big\vert\leq\\
 & 2\hat L_f\ve t(\ve)+\ve E\sum_{j=0}^{n(\ve)-1}
 \big\vert\int_{jt(\ve)}^{(j+1)t(\ve)}\big(f(\ve s,X^\ve_{x,y}(s),
 Y^\ve_{x,y}(s))-\bar f(\ve s,\bar X^\ve_x(s))\big)ds\big\vert\nonumber\\
 &\leq 2\hat L_f\ve t(\ve)+\ve E\sum_{j=0}^{n(\ve)-1}
 \big\vert\int_0^{t(\ve)}\big(f(\ve jt(\ve)+\ve s,X^\ve_{x^\ve_j,y^\ve_j}(s),
 Y^\ve_{x^\ve_j,y^\ve_j}(s))\nonumber\\
 &-\bar f(\ve jt(\ve)+\ve s, X^\ve_{x^\ve_j,y^\ve_j}(s))\big)ds\big\vert
 \nonumber\\
 &\leq 2LT^2/N+2\big(\hat L_fT+T^2(L+LK+KK_L)\big)/n(\ve)+\ve t(\ve)\nonumber\\
 &\times\sum_{l=0}^{N-1}\sum_{ln(\ve)/N\leq j<(l+1)n(\ve)/N,j\leq n(\ve)}
 E\big\vert\frac 1{t(\ve)}\int_0^{t(\ve)}f(lT/N,x^\ve_j,Y^\ve_{x^\ve_j,y^\ve_j}
 (s))ds\nonumber\\
 &-\bar f(lT/N,x_j^\ve)\big\vert\leq 2LT^2/N+2\big(\hat L_fT+T^2(L+LK+KK_L)
 \big)/n(\ve)\nonumber\\
 &+\ve t(\ve)n(\ve)\del
+2\hat L_f\ve t(\ve)\sum_{l=0}^{N-1}\sum_{ln(\ve)/N\leq j<(l+1)n(\ve)/N,j
 \leq n(\ve)}\nonumber\\
 &P\big\{\big(X^\ve_{x,y}(jt(\ve)),Y^\ve_{x,y}(jt(\ve))\big)\in
 \cE^{f_l}_\ve(t(\ve),\del)\big\}\nonumber
 \end{eqnarray}
 where $f_l(z,v)=f(lT/N,z,v).$ Integrating against $\mu$ both parts of 
 (\ref{2.9.20}) over $G\times\bfM$ we obtain
 \begin{eqnarray}\label{2.9.21}
 &\,\,\,\,\,\,\,\,\,\,\,\,\,\,\,\ve\int_G\int_\bfM E\sup_{0\leq t\leq T/\ve}
 \big\vert\int_0^t\big(f(\ve s,X^\ve_{x,y}(s),Y^\ve_{x,y}(s))\\
 &-\bar f(\ve s,X^\ve_{x,y}(s))\big)ds\big\vert d\mu(x,y)
 \leq 2\big(\hat L_fT+T^2(L+LK+KK_L)\big)/n(\ve)\nonumber\\
 &+2LT^2/N+T\del+2\hat L_f\max_{0\leq l\leq N-1}\eta_l(\ve,\del)\nonumber
 \end{eqnarray}
 where 
 \[
 \eta_l(\ve,\del)=\max_{0\leq j\leq n(\ve)-1}\int_G\int_\bfM P\big\{
 \big(X^\ve_{x,y}(jt(\ve)),Y^\ve_{x,y}(jt(\ve))\big)\in
 \cE^{f_l}_\ve(t(\ve),\del)\big\}d\mu(x,y).
 \]
 By the assumption there exists an integer valued function $n(\ve)\to\infty$
 as $\ve\to 0$ such that (\ref{2.9.5}) holds true for all $g=f_0,f_1,...,f_{N-1}$
 and then $\max_{0\leq l\leq N-1}\eta_l(\ve,\del)\to 0$ as $\ve\to 0$. Hence,
 letting first $\ve\to 0$, then $\del\to 0$ and, finally, $N\to 0$ we obtain 
 (\ref{2.9.15}) in view of (\ref{2.9.16}) and (\ref{2.9.17}), completing the proof 
 of Theorem \ref{thm2.9.1}.
 \end{proof}
 
 Observe that (\ref{2.9.4}) holding true for all bounded continuous functions
is, in principle, stronger than the averaging principle in the form 
\begin{equation}\label{2.9.22}
\lim_{\ve\to 0}\ve\int_G\int_\bfM E\sup_{0\leq t\leq T/\ve}
 \big\vert X^\ve_{x,y}(t)
 -\bar X^\ve_x(t)\big\vert d\mu(x,y)\to 0\,\,\mbox{as}\,\,\ve\to 0
 \end{equation}
 since (\ref{2.9.22}) is equivalent to (\ref{2.9.15}) with $f=B$.
In fact, if we require (\ref{2.9.5}) only for one function $g=B$ then
it will be equivalent to (\ref{2.9.22}) which follows in the same way as
the proof of Theorem \ref{thm2.9.1} above. Still, the main
interesting classes of systems, we are aware of, for which (\ref{2.9.4}) holds
true are the same for which (\ref{2.9.22}) is satisfied though it is easy to
construct examples of (somewhat degenerate) right hand sides $B$ in
(\ref{2.1.1}) for which (\ref{2.9.22}) holds true but (\ref{2.9.4}) fails (since
in the latter we require convergence for all functions $f$ and in the 
former only for $f=B$).

It follows from \cite{Ki3+} that the assumptions of Theorem \ref{thm2.9.1} hold
true when the unperturbed fast motions $Y_{x,y}(t)$ are diffusion processes
on $\bfM$ so that $\mu_x$ is an invariant measure of $Y_x$ on $\bfM$ ergodic
for $\la$-almost all $x$, where $\la$ is the normalized 
    Lebesgue measure on a large compact in $\bbR^d$, and $\mu_x(U)=
    \int_Uq(x,y)dm(y)$ with $q(x,y)>0$ differentiable in $x$ and $y$.
    This can be extended to random evolutions considered in previous sections.
    
    Observe that under assumptions of Theorem \ref{thm2.2.2} we can
     obtain also large deviations bounds in the form (\ref{2.2.10}) and
     (\ref{2.2.11}) for
     \[
     \tilde Z^\ve_{x,y}(t)=\int_0^tf(s,X^\ve_{x,y}(s/\ve),Y^\ve_{x,y}(s/\ve))
     ds
     \]
     with the functional
     \begin{eqnarray*}
     &\tilde S_{0T}(\tilde\gam)=\inf\big\{ S_{0T}(\gam):\, S_{0T}(\gam)=
     \int_0^TI_{\gam_t}(\nu_t)dt,\\
     &\dot {\gam}_t=\bar B_{\nu_t}(\gam_t),\,\tilde\gam_t=\int_0^t\bar 
     f_{\nu_s}(s,\gam_s)ds\,\,\forall t\in[0,T]\big\},\,\,\bar f_\nu(s,x)=
     \int f(s,x,y)d\nu(y),
     \end{eqnarray*}
     where $f$ is a bounded Lipschitz continuous vector function. The proof 
     can be carried out quite similarly to the proof of Theorem \ref{thm2.2.2}.
     Analogous results can be obtained in the discrete time setup of difference
     equations (\ref{2.1.7}).

\backmatter

\printindex
\end{document}